\documentclass[11pt]{amsart}

\usepackage{geometry}
\usepackage{soul}

\usepackage{todonotes}
\usepackage{amsfonts, adjustbox, amssymb, amscd}
\numberwithin{equation}{section}

\usepackage{bm}
\usepackage{verbatim}
\usepackage{mathrsfs}
\usepackage{graphicx}
\usepackage{tikz-cd}
\usepackage{subcaption}
\usepackage{listings}
\usepackage{subfiles}
\usepackage[toc,page]{appendix}
\usepackage{mathtools}
\usepackage{comment}
\usepackage{enumerate}
\usepackage{enumitem}
\usepackage[all]{xy}

\usepackage{graphicx}
\graphicspath{{images/}}

\usepackage{appendix}
\usepackage{hyperref}
\hypersetup{
    colorlinks=true,
    citecolor=red,
    linkcolor=blue,
    filecolor=magenta,      
    urlcolor=red,
}
\newcommand{\bb}{\bm{b}}
\newcommand{\Mm}{{\bf{M}}}
\newcommand{\Nn}{{\bf{N}}}

\newcommand{\Pp}{{\bf{P}}}
\newcommand{\NN}{{\bf{N}}}
\newcommand{\Dd}{{\bf{D}}}

\newcommand{\Qq}{\mathbb{Q}}

\newcommand{\Rr}{\mathbb{R}}

\newcommand{\Center}{\operatorname{center}}

\newcommand{\Exc}{\operatorname{Exc}}

\newcommand{\rk}{\operatorname{rank}}

\newcommand{\Nklt}{\operatorname{Nklt}}

\newcommand{\Bs}{{\operatorname{Bs}}}
\newcommand{\Src}{{\operatorname{Src}}}
\newcommand{\Spr}{{\operatorname{Spr}}}
\newcommand{\num}{{\operatorname{num}}}
\newcommand{\tang}{{\operatorname{tang}}}

\newcommand{\lct}{\operatorname{lct}}

\newcommand{\Supp}{\operatorname{Supp}}

\newcommand{\Nlc}{\operatorname{Nlc}}

\newcommand{\mult}{\operatorname{mult}}

\newcommand{\cont}{\operatorname{cont}}

\newcommand{\lf}{\lfloor}
\newcommand{\rf}{\rfloor}

\newcommand{\Aa}{{\bf{A}}}

\newcommand{\Bb}{{\bf{B}}}
\newcommand{\Ff}{\mathcal{F}}

\newcommand{\Oo}{\mathcal{O}}
\newcommand{\Ii}{\Gamma}

\newcommand{\Sing}{\mathrm{Sing}}


\newcounter{parentnumber}

\newtheorem{thm}{Theorem}[subsection]
\newtheorem{conj}[thm]{Conjecture}
\newtheorem{cor}[thm]{Corollary}
\newtheorem{lem}[thm]{Lemma}
\newtheorem{prop}[thm]{Proposition}

\newtheorem{claim}[thm]{Claim}

\newtheorem{alphthm}{Theorem}

\theoremstyle{definition}
\newtheorem{defn}[thm]{Definition}

\theoremstyle{definition}
\newtheorem{rem}[thm]{Remark}

\newtheorem{deflem}[thm]{Definition-Lemma}
\newtheorem{defthm}[thm]{Definition-Theorem}

\newtheorem{sce}[thm]{Scenario}
\newtheorem{nota}[thm]{Notation}

\newtheorem{cons}[thm]{Construction}

\theoremstyle{definition}

\begin{document}

\title{Minimal model program for algebraically integrable foliations and generalized pairs}
\author{Guodu Chen, Jingjun Han, Jihao Liu, and Lingyao Xie}

\subjclass[2020]{14E30, 37F75}
\keywords{Generalized pairs, Algebraically integrable foliations, Minimal model program}
\date{\today}

\begin{abstract}
Using techniques from the theory of foliations, we establish the cone theorem and the contraction theorem for lc generalized pairs in full generality, and meanwhile develop the minimal model program for $\mathbb Q$-factorial foliated dlt algebraically integrable foliations. As an application, we obtain the canonical bundle formula for generalized pairs completely, together with several further consequences, including answering a question of Cascini and Spicer.
\end{abstract}

\address{School of Mathematical Sciences, Shanghai Jiao Tong University, 800 Dongchuan RD Shanghai, Minhang District Shanghai 200240, China}
\email{chenguodu@sjtu.edu.cn}

\address{Shanghai Center for Mathematical Sciences \& School of Mathematical Sciences, Fudan University, Shanghai, 200438, China}
\email{hanjingjun@fudan.edu.cn}

\address{Department of Mathematics, Peking University, No. 5 Yiheyuan Road, Haidian District, Peking 100871, China}
\email{liujihao@math.pku.edu.cn}

\address{Department of Mathematics, University of California, San Diego, 9500 Gilman Drive \# 0112, La Jolla, CA 92093-0112, USA}
\email{l6xie@ucsd.edu}

\maketitle

\pagestyle{myheadings}\markboth{\hfill G. Chen, J. Han, J. Liu, and L. Xie \hfill}{\hfill MMP for algebraically integrable foliations and generalized pairs\hfill}

\tableofcontents

\section{Introduction}\label{sec:Introduction}

We work over the field of complex numbers $\mathbb C$.

\subsection{Main theorems} Generalized pairs (g-pairs) and algebraically integrable foliations are two structures that play important roles in modern birational geometry, particularly in the minimal model program (MMP). 

\medskip

\noindent{\emph{Generalized pairs.}} 
G-pairs arise naturally in the study of higher dimensional birational geometry, especially when applying dimension induction arguments like the canonical bundle formula. Birkar and Zhang introduced them in \cite{BZ16} while studying the effective Iitaka fibration conjecture. Though they might appear technical, the MMP for g-pairs has proven to be a powerful tool in other topics in algebraic geometry in recent years:

\begin{itemize}
    \item Borisov-Alexeev-Borisov conjecture \cite{Bir19,Bir21a}. 
    \item The termination of flips for pseudo-effective pairs in dimension four \cite{HL22,CT23}. 
    \item Hacon-Han's conjecture on the connectedness principle of non-klt loci \cite{FS23,Bir20}. 
    \item Boundedness of polarized Calabi-Yau varieties \cite{Bir23a}.
    \item The minimal model program for K\"ahler 3-folds \cite{DH23}.
\end{itemize}

Log canonical (lc) singularities are the largest class for which the traditional MMP makes sense. As pointed out by \cite{KM98}, it is very hard to work with them as they are rather complicated from the cohomological point of view. Our first main result is the following which answers questions in \cite[6.1]{Bir21} and \cite[3.1, 3.3]{HL22}.

\begin{alphthm}\label{thm: main mmp gpair}
Let $(X,B,\Mm)/U$ be an lc generalized pair. Then:
\begin{enumerate}
    \item (Theorem \ref{thm: cone theorem nonnqc gpair}) The cone theorem and contraction theorem hold for $(X,B,\Mm)/U$. %In particular, we can run a $(K_X+B+\Mm_X)$-MMP.
    \item (Theorem \ref{thm: semi-ampleness intro}) If $K_X+B+A+\Mm_X$ is nef$/U$ for some ample$/U$ $\Rr$-divisor $A$, then $K_X+B+A+\Mm_X$ is semi-ample$/U$.
\end{enumerate}
\end{alphthm}

We note here that \cite{HL21a,Xie22,CLX23,LX23a,LX23b} proved Theorem \ref{thm: main mmp gpair} when $(X,B,\Mm)/U$ satisfies a technical condition ``NQC'' (see Definition \ref{defn: b divisors} for details) which was introduced in \cite{HL22}.

Hacon suggested to us that Theorem \ref{thm: main mmp gpair} might have essential implications on K\"ahler varieties in higher dimensions as the associated K\"ahler class on a K\"ahler variety can be considered as a nef $\Rr$-class and is suitable for the nef part of a g-pair (cf. \cite{DH23,DHY23}).

\medskip 

In our proof of Theorem \ref{thm: main mmp gpair}, we built a bridge between the theory of g-pairs and the theory of algebraic integrable foliations. We introduced the theory of \emph{generalized foliated quadruples}, allowing us to develop frameworks for both theories, although they seem to have different origins and flavors.

\medskip

\noindent{\emph{Algebraically integrable foliations.}} The MMP for (algebraically integrable) foliations is a generalization of the classical MMP. Instead of examining the structures associated with the canonical divisor of the ambient variety $K_X$, the foliations theory concentrates on the structures connected to the foliated canonical divisor $K_{\Ff}$. Unlike the usual pairs, Bertini type theorems, as well as the abundance conjecture fail for algebraically integrable foliations.

The MMP for foliations has been established for surfaces (\cite{McQ08, Bru15}) and threefolds (\cite{CS20, Spi20, CS21, SS22}). However, for higher-dimensional foliations, the MMP is still widely open.

In this article, we focus on the MMP for algebraically integrable foliations, where the general leaves are algebraic varieties. In other words, they are induced by dominant rational maps. The theory of algebraically integrable foliations holds an important place in number theory and birational geometry, e.g. in Miyaoka's proof of the non-vanishing conjecture in dimension three \cite{Miy87}, in the study of varieties admitting a non-trivial fibration with rationally connected fibers \cite{BM16}, and in the Grothendieck-Katz conjecture and the Ekedahl-Shepherd-Barron-Taylor conjecture \cite{Bos01}.

Our second result is as follows: we establish the MMP for algebraically integrable foliations with F-dlt singularities (see Definition \ref{defn: fdlt}) in arbitrary dimensions, except the termination part. 

\begin{alphthm}\label{thm: main mmp foliation}
Let $(X,\Ff,B)/U$ be a $\Qq$-factorial F-dlt foliated triple such that $\Ff$ is algebraically integrable. Let $A$ be an ample$/U$ $\Rr$-divisor on $X$. Then:
\begin{enumerate}
    \item (\cite[Theorem 3.9]{ACSS21}+Theorem \ref{thm: fdlt is acss}) The cone theorem, contraction theorem, and the existence of flips hold for $(X,\Ff,B)/U$. In particular, we can run a $(K_\Ff+B)$-MMP$/U$.
    \item (Theorem \ref{thm: bpf fdlt}) If $K_{\Ff}+B+A$ is nef$/U$, then $K_{\Ff}+B+A$ is semi-ample$/U$\footnote{\cite[Theorem 1.2]{CD23} claimed a proof of some special cases of Theorem \ref{thm: main mmp foliation}(2) and other results that are similar to some results in this paper. However, the current proofs in \cite{CD23} seem to be incomplete, mainly because of the failure of \cite[Lemma 2.4]{CD23} and some gaps in the proof of \cite[Theorem 3.5]{CD23}. In this paper, we will avoid using any results in \cite{CD23}.}.
    \item (Theorem \ref{thm: +a gmm fdlt}) If $B\geq A\geq 0$, then $(X,\Ff,B)/U$ has a good minimal model$/U$ or a Mori fiber space$/U$.
\end{enumerate}
\end{alphthm}

As pointed out by \cite{CS21}, F-dlt foliated triples play the same role as dlt pairs in the classical MMP, making them a natural class of singularities to study in the theory of foliations. When $\Ff=T_X$ and $\lfloor B\rfloor=0$, Theorem \ref{thm: main mmp foliation} becomes the classical result of the existence of good minimal models of varieties of general type and the finite generation of the canonical ring \cite[Theorem 1.2]{BCHM10}. We remark that Theorem \ref{thm: main mmp foliation}(2) does not hold in general without the polarization of the ample$/U$ $\Rr$-divisor $A$ (cf. \cite[Example 5.4]{ACSS21}).

\subsection{Idea of the proofs of Theorems \ref{thm: main mmp gpair} and \ref{thm: main mmp foliation}} The proofs of  Theorems \ref{thm: main mmp gpair} and \ref{thm: main mmp foliation} rely crucially on a larger framework: the theory of \emph{generalized foliated quadruples}. 
\begin{defn}[{cf. \cite[Definition 1.2]{LLM23}}]\label{defn: gfq intro}
A \emph{generalized foliated quadruple} (\emph{gfq} for short) $(X,\Ff,B,\Mm)/U$ consists of a normal quasi-projective variety $X$, a foliation $\Ff$ on $X$, an $\mathbb R$-divisor $B\geq 0$ on $X$, a projective morphism $X\rightarrow U$, and a nef$/U$ $\Rr$-divisor $\Mm_{X'}$ on a high model $X'$ of $X$, such that $K_{\Ff}+B+\Mm_X$ is $\Rr$-Cartier. Here $\Mm_X$ is the image of $\Mm_{X'}$ on $X$.
\end{defn}
The notation $\Mm$ in Definition \ref{defn: gfq intro} is considered as a $\bb$-divisor on $X$. We refer the reader to Definition \ref{defn: b divisors} for the definition of $\bb$-divisors, and to Definition \ref{defn: gfq} for a more detailed definition of generalized foliated quadruples. It is clear that when $\Mm=\bm{0}$ is the trivial $\bb$-divisor, a generalized foliated quadruple is just a foliated triple $(X,\Ff,B)/U$; on the other hand, when $\Ff=T_X$, a generalized foliated quadruple is a generalized pair $(X,B,\Mm)/U$ (\cite[Definition 1.4]{BZ16}, \cite[Definition 2.1]{HL22}). Therefore, generalized foliated quadruples can be considered as a mixture of foliated triples and generalized pairs. 

The concept of generalized foliated quadruples was introduced by the third author, Luo, and Meng in their study of the global ACC for foliated threefolds \cite{LLM23}. The main part of this paper is dedicated to a systematic study of this new structure, to establish foundational theorems on it, and to apply them to prove key results and conjectures in both foliations and generalized pairs. It is worth mentioning that our definition slightly differs from \cite[Definition 1.2]{LLM23}, as the latter requires $\Mm$ to be NQC$/U$ (see Definition \ref{defn: b divisors}), while we only require it to be nef$/U$. This will be crucial for us to use this structure to prove Theorem \ref{thm: main mmp gpair}. We refer the reader to Subsection \ref{sec: reason to consider gfq} for a detailed explanation of why this new structure is vital not only for this paper but also for future studies of foliations and generalized pairs.

Under the framework of generalized foliated quadruples, the proofs of Theorems  Theorems \ref{thm: main mmp gpair} and \ref{thm: main mmp foliation} proceed simultaneously.

The first result to prove is the cone theorem. As a positive beginning, the cone theorem for algebraically integrable foliations is already known, as seen in \cite[Theorem 3.9]{ACSS21}. With some adjustments to the details of the proofs, the same approach used in the proof of \cite[Theorem 3.9]{ACSS21} also works for algebraically integrable generalized foliated quadruples (Theorem \ref{thm: cone theorem gfq}). Since the proof of the cone theorem depends on induction on dimension by adjunction, a crucial part of our proof is the adjunction formulas for generalized foliated quadruples, proven in Section \ref{sec: adjunction}. Indeed, we present a precise adjunction formula in Section \ref{sec: adjunction} as we aim to prove the ACC for lc thresholds (Theorem \ref{thm: acc lct alg int gfq}) and the global ACC (Theorem \ref{thm: global acc alg int gfq}) for algebraically integrable generalized foliated quadruples (Section \ref{sec: acc gfq}). These two results will not be used in the rest of the paper, but are expected to be useful for future research. 

Now, the cone theorem for algebraically integrable generalized foliated quadruples implies the cone theorem for generalized pairs by letting $\Ff=T_X$. It is worth mentioning that our approach to prove the cone theorem is different from the one in the proof of \cite{HL21a}. More precisely, the proof in \cite{HL21a} depends on the subadjunction formula for NQC generalized pairs, but we do not have these formulas for non-NQC generalized pairs yet (although we will prove them in the later part of the paper).

With the cone theorem established, we can move on to prove the rest of Theorem \ref{thm: main mmp foliation}(1). We only need to show that each step of a $(K_{\Ff}+B)$-MMP is also a $(K_X+\Delta)$-MMP for some lc pair $(X,\Delta)$. To do this, we need to show that $(X,\Ff,B)$ satisfies a property called ``ACSS" (see Definition \ref{defn: ACSS f-triple}) and that this property is preserved under each step of the MMP. The latter is established in Lemma \ref{lem: ACSS mmp can run}, which is based on the cone theorem. To prove the former, i.e., that F-dlt implies ACSS, we only need to show the termination of MMP with scaling for $\Qq$-factorial F-dlt algebraically integrable foliations with very exceptional foliated log canonical divisors. This primarily relies on the general negativity lemma \cite[Lemma 3.3]{Bir12} and the cone theorem, and is proven in Theorem \ref{thm: mmp very exceptional alg int fol}. We note that the same approach to the proof also works for $\Qq$-factorial generalized foliated quadruples with F-dlt singularities.

Our next goal is to establish Theorem \ref{thm: main mmp gpair}(2), the base-point-freeness theorem for generalized pairs. An approach highlighted in \cite{CLX23} suggests that the base-point-freeness theorem for generalized pairs should only depend on the subadjunction formula of generalized pairs (Theorem \ref{thm: subadj intro}), which in turn relies only on the canonical bundle formula of generalized pairs (Theorem \ref{thm: cbf gpair nonnqc}). Sections \ref{sec: du bois} and \ref{sec: vanishing gpair} confirm this fact. Consequently, we only need to establish the canonical bundle formula for generalized pairs. More specifically, we only need to prove that the moduli part of the canonical bundle formula for generalized pairs is nef.

A key observation is that when the lc-trivial fibration structure satisfies certain stability conditions (e.g., BP stable), the moduli part of the canonical bundle formula corresponds to the log canonical divisor of the induced foliation (Proposition \ref{prop: weak cbf gfq}), and it is nef when the log canonical divisor of the associated pair is nef (Proposition \ref{prop: bp stable nef}). This stability condition is closely related to the singularity of the induced foliation (Proposition \ref{prop: bp semistable foliation lc}, Theorem \ref{thm: lc+weak acc=bpstable}), so it can be attained by taking a foliated log resolution, leaving us to run an MMP to achieve the nef condition. To run the MMP, it suffices to consider the existence of log minimal models for generalized foliated quadruples with numerical dimension zero (Propositions \ref{prop: weak ss num 0 mmp}, \ref{prop: projective num 0 mmp}). The existence of such log minimal models depends only on the cone theorem and the Nakayama-Zariski decomposition, for which we already have the respective results.

%Having settled the base-point-freeness theorem (Theorem \ref{thm: main mmp gpair}(2)), the contraction theorem in Theorem \ref{thm: main mmp gpair}(1) becomes an immediate corollary. Specifically, for a $(K_X+B+\Mm_X)$-flipping contraction $X\rightarrow Z$ and an $\Rr$-divisor $D$ on $X$, $D\equiv_Z0$ if and only if $D\sim_{\mathbb R,Z}0$. Given this, and with a specific perturbation in the coefficients (see Lemma \ref{lem: flip reduce special gpair to pair}), we establish the existence of flips when $X$ is $\Qq$-factorial. This concludes the proof of Theorem \ref{thm: main mmp gpair}.

Finally, we turn to proving Theorem \ref{thm: main mmp foliation}(2-3). Although the Bertini-type theorem fails for foliations, by employing the structure of generalized foliated quadruples, we can, roughly speaking, reduce Theorem \ref{thm: main mmp foliation}(3) to Theorem \ref{thm: main mmp foliation}(2) (Lemma \ref{lem: +a keep under mmp}, Theorem \ref{thm: gmm polarized gfq}). The proof of Theorem \ref{thm: main mmp foliation}(2) is segmented into three steps: First, we utilize the already proven contraction theorem in Theorem \ref{thm: main mmp foliation}(2) to construct a contraction $X\rightarrow T$, where the general fibers of $X\rightarrow Z$ are tangent to $\Ff$. Next, we apply the canonical bundle formula for generalized foliated quadruples to derive a generalized pair structure polarized with an ample divisor on $T$. This canonical bundle formula (Definition-Theorem \ref{defthm: cbf lctrivial morphism}) can be derived using the canonical bundle formula for lc-trivial fibrations of generalized pairs. Lastly, we apply the cone theorem for generalized foliated quadruples to show that the generalized foliated log canonical divisor on $T$ is ample, therefore the foliated log canonical divisor on $X$ is semi-ample, completing the proof of Theorem \ref{thm: main mmp foliation}(2). This concludes the proof of Theorem \ref{thm: main mmp foliation}.

\medskip

\noindent\textbf{Structure of the Paper} In Section \ref{sec: statement of main results}, we list the main results of this paper and explain the importance of the structure of generalized foliated quadruples. The rest of the paper is divided into four parts, each of which we will introduce below. For the convenience of the reader, we have prepared the following flowchart (Table \ref{tbl: flowchart}) to illustrate the streamlined process involved in the proofs of our main theorems.

\begin{table}[ht]
\caption{Structure of the paper}\label{tbl: flowchart}
\begin{adjustbox}{width=1\textwidth,right}
%\begin{center}
$\xymatrix{
 &  *+[F]\txt{Preliminary results on foliations\\ and generalized pairs\\ (Sections \ref{sec: preliminaries}, \ref{sec: basic property gpair}, \ref{sec: stability gpair}, \ref{sec: acss gfq})} & \\
& *+[F]\txt{Precise adjunction formulas\\ for gfqs (Section \ref{sec: adjunction})}\ar@{->}[d]\ar@{->}[r] & *+[F]\txt{ACC and global ACC\\ for gfqs (Section \ref{sec: acc gfq})}\\
 & *+[F]\txt{Cone theorem for gfqs\\ (Section \ref{sec: cone})}\ar@{->}[d]\ar@{->}[r]\ar@{->}[dr]& *+[F]\txt{Existence of Mori fiber space\\ for gfqs (Subsection \ref{subsec: eomfs})}\ar@{->}[u]\\
*+[F]\txt{Very exceptional MMP\\ (Subsection \ref{subsec: very exceptional})}\ar@/_8pc/[ddddd] & *+[F]\txt{MMP with scaling\\ (Subsection \ref{subsec: eomfs})}\ar@{->}[d]\ar@{->}[l] & *+[F]\txt{Theorem \ref{thm: main mmp gpair}}
 \\
 & *+[F]\txt{$\kappa_{\sigma}=0$ MMP\\ (Subsection \ref{sub: num 0 mmp})}\ar@{->}[d]  &\\
 *+[F]\txt{Stability of gfqs\\ (Subsection \ref{subsec: stability gfq})}\ar@{->}[r] &  *+[F]\txt{Canonical bundle formula\\ for generalized pairs \\(Subsection \ref{subsec: cbf gpair})}\ar@{->}[d]\ar@{->}[dl] & \\
  *+[F]\txt{Canonical bundle formula for gfqs \\ (Subsection \ref{subsec: cbf gfq}, Section \ref{sec: subadj})}\ar@/_2pc/[ddr] & *+[F]\txt{Subadjunction\\ for generalized pairs \\(Section \ref{sec: subadj})}\ar@{->}[d] & \\
 &  *+[F]\txt{Base-point-freeness\\ and Contraction theorem\\ for generalized pairs\\ (Sections \ref{sec: du bois}, \ref{sec: vanishing gpair})}\ar@/_2pc/[uuuur]\ar@{->}[d] &  
\\*+[F]\txt{$\Qq$-factorial F-dlt\\ $\Rightarrow$ACSS\\ (Theorem \ref{thm: fdlt is acss})}\ar@{->}[r] &  *+[F]\txt{Theorem \ref{thm: main mmp foliation}}&  \\
&  &  \\
}$ 
\end{adjustbox}
\end{table}

\smallskip

\noindent\textbf{Part \ref{part:prelim}}. Preliminaries: Sections \ref{sec: preliminaries}, \ref{sec: basic property gpair}, and \ref{sec: stability gpair}. This part contains preliminary results and definitions that will be utilized throughout the remainder of the paper, with few foliation structures involved. In Section \ref{sec: preliminaries}, we introduce some basic definitions, including the concept of generalized foliated quadruples and their singularities. In Section \ref{sec: basic property gpair}, we establish basic properties of generalized pairs. Section \ref{sec: stability gpair} is parallel to \cite[Section 2]{ACSS21}, studying the stability of generalized pairs and introducing the concept of Property $(*)$ for generalized pairs.

\smallskip

\noindent\textbf{Part \ref{part:cone}}. Cone Theorem and the minimal model program for algebraically integrable foliations: Sections \ref{sec: adjunction}, \ref{sec: acss gfq}, \ref{sec: cone}, \ref{sec: mmp gfq}, and \ref{sec: acc gfq}. This part focuses on the cone theorem for algebraically integrable generalized foliated quadruples and its applications. In Section \ref{sec: adjunction}, we prove a precise adjunction formula for algebraically integrable generalized foliated quadruples. Section \ref{sec: acss gfq} defines ACSS generalized foliated quadruples and studies its fundamental behaviors. Section \ref{sec: cone} proves the cone theorem for algebraically integrable generalized foliated quadruples. In Section \ref{sec: mmp gfq}, we prove most results of this paper on the minimal model program for generalized foliated quadruples with $\Qq$-factorial F-dlt singularities, excluding the existence of good minimal models. Section \ref{sec: acc gfq} verifies the ACC, the global ACC, and the existence of uniform rational polytopes for algebraically integrable generalized foliated quadruples.

\smallskip

\noindent\textbf{Part \ref{part:cbf}}. Canonical bundle formula and minimal model program for generalized pairs: Sections \ref{sec: cbf}, \ref{sec: subadj}, \ref{sec: du bois}, and \ref{sec: vanishing gpair}. This part presents the canonical bundle formula and applies it to establish the minimal model program for generalized pairs. In Section \ref{sec: cbf}, we state and prove the canonical bundle formula for lc-trivial fibrations of generalized foliated quadruples and, in particular, generalized pairs. Section \ref{sec: subadj} studies lc-trivial morphisms of generalized foliated quadruples and proves the subadjunction formula for generalized pairs. Section \ref{sec: du bois} shows that lc generalized pairs have Du Bois singularities. Section \ref{sec: vanishing gpair} proves the Kodaira vanishing theorem, the Kawamata-Viehweg vanishing theorem, the base-point-freeness theorem, and the contraction theorem for generalized pairs. 

\smallskip

\noindent\textbf{Part \ref{part:gmm}}. Good minimal model, applications, and proofs of the main theorems: Sections \ref{sec: gmm fdlt} and \ref{sec: proof of the main theorems}. In Section \ref{sec: gmm fdlt}, we prove the existence of good minimal models for $\Qq$-factorial F-dlt generalized foliated quadruples polarized with an ample divisor. Similar arguments imply the base-point-freeness theorem for such quadruples, leading to a special case of the Prokhorov-Shokurov $\bb$-semi-ampleness conjecture. Lastly, in Section \ref{sec: proof of the main theorems}, we prove all the main theorems of the paper that are not covered in the previous sections.

\medskip

\noindent\textbf{Postscript}.
We remark that after this paper was posted on Arxiv, subsequent works (further development on the MMP of foliations \cite{LMX24} and its adjoint structures \cite{CHL+24,CHL+25}) have emerged, which build on the framework presented here. 

\medskip

\noindent\textbf{Acknowledgement}. The authors would like to thank Caucher Birkar, Paolo Cascini, Priyankur Chaudhuri, Omprokash Das, Christopher D. Hacon, Chen Jiang, Junpeng Jiao, Jie Liu, Yuchen Liu, Roktim Mascharak, Fanjun Meng, Wenhao Ou, Vyacheslav V. Shokurov, Chenyang Xu, and Qingyuan Xue for fruitful discussions. 
The work is supported by the National Key R\&D Program of China \#2024YFA1014400. The first author and the second author are supported by the National Key R\&D Program of China \#2025YFA1018100. The second author is supported by the National Key R\&D Program of China \#2023YFA1010600. The first author was sponsored by the NSFC (\#12401055) and the NSF of Shanghai (\#24ZR1430000). The second author was supported by NSFC for Excellent Young Scientists (\#12322102), and he is a member of LMNS, Fudan University. The fourth author had been partially supported by NSF research grants no. DMS-1801851 and DMS-1952522, as well as a grant from the Simons Foundation (Award Number: 256202).

\section{Statement of main results}\label{sec: statement of main results}

In this section, we provide the statements of the main results of this paper.

\subsection{Minimal model program for generalized pairs}\label{subsec: mmp gpair} In the past few years, there has been significant advancement in the minimal model program for NQC generalized pairs. The cone theorem, as well as the $\mathbb Q$-factorial cases of the contraction theorem and the proof of the existence of flips, were established in \cite[Theorem 1.1]{HL21a}. Later, the existence of flips for (potentially non-$\mathbb Q$-factorial) NQC generalized pairs was verified in \cite{LX23a}, while the contraction theorem for these pairs was proven in \cite{Xie22}. Additionally, \cite{CLX23} confirmed the Kodaira and the Kawamata-Viehweg vanishing theorems for NQC generalized pairs, offering an alternative proof for the contraction theorem. These developments form the foundation of the minimal model program for NQC generalized pairs, with numerous corollaries and applications already provided in \cite{LT22,TX23}.

The structure of NQC generalized pairs has naturally arisen in the study of the canonical bundle formulas, making them a fundamental structure in the study in the minimal model program. For a considerable amount of time, it has been presumed that the realm of NQC generalized pairs would be the most extensive category necessary to establish in the minimal model program. This is because the structure of NQC generalized pairs is maintained under the canonical bundle formula, adjunction formula, and each stage of the minimal model program, thereby eliminating the need to consider the minimal model program for non-NQC generalized pairs or other larger categories.

However, recent studies on the minimal model program for K\"ahler varieties \cite{DH23,DHY23} have emphasized the critical role the structure of non-NQC generalized pairs plays in the minimal model program for K\"ahler varieties. In the case of K\"ahler varieties, the selection of ample divisors is restricted, preventing many procedures, such as the minimal model program with scaling and general hyperplane section cuttings. Nevertheless, the associated K\"ahler class $\omega$ on a K\"ahler variety serves as a substitute for ample divisors. Although $\omega$ cannot be categorized as an $\Rr$-divisor, it can be considered as an nef $\Rr$-class and is suitable for the nef part of a generalized pair. As explained in \cite{DHY23}, it is now possible to formally define ``running MMP with scaling of the nef $\Rr$-$(1,1)$-class $\overline{\omega}$". Given that $\omega$ is only an $\Rr$-class and NQC cannot be assured, the study of the structure of non-NQC generalized pairs immediately becomes vital for the minimal model program on K\"ahler vareties.

Although little was known about the minimal model program for non-NQC generalized pairs, we have been able to establish the cone theorem and the contraction theorem for non-NQC generalized pairs, thanks to the cone theorem and the canonical bundle formula for generalized foliated quadruples. %With additional effort, we also prove the existence of flips for $\Qq$-factorial non-NQC generalized pairs. These results collectively lay the groundwork for the minimal model program for $\Qq$-factorial generalized pairs. The detailed theorems are as follows:

\begin{thm}[Cone and contraction theorems]\label{thm: cone theorem nonnqc gpair}
Let $(X,B,\Mm)/U$ be a generalized pair and $\pi: X\rightarrow U$ the associated morphism. Let $\{R_j\}_{j\in\Lambda}$ be the set of $(K_{X}+B+\Mm_X)$-negative extremal rays in $\overline{NE}(X/U)$ that are rational. Then:
\begin{enumerate}
    \item $$\overline{NE}(X/U)=\overline{NE}(X/U)_{K_{X}+B+\Mm_X\geq 0}+\overline{NE}(X/U)_{\Nlc(X,B,\Mm)}+\sum_{j\in\Lambda} R_j.$$
    In particular, any $(K_{X}+B+\Mm_X)$-negative extremal ray in $\overline{NE}(X/U)$ is rational.
    \item Each $R_j$ is spanned by a rational curve $C_j$ such that $\pi(C_j)=\{pt\}$ and 
    $$0<-(K_{X}+B+\Mm_X)\cdot C_j\leq 2\dim X.$$
    \item For any ample$/U$ $\Rr$-divisor $A$ on $X$,
    $$\Lambda_A:=\{j\in\Lambda\mid R_j\subset\overline{NE}(X/U)_{K_{X}+B+A+\Mm_X<0}\}$$
    is a finite set. In particular, $\{R_j\}_{j\in\Lambda}$ is countable, and is a discrete subset in $\overline{NE}(X/U)_{K_{X}+B+\Mm_X<0}$. Moreover, we may write
    $$\overline{NE}(X/U)=\overline{NE}(X/U)_{K_{X}+B+A+\Mm_X\geq 0}+\overline{NE}(X/U)_{\Nlc(X,B,\Mm)}+\sum_{j\in\Lambda_A}R_j.$$
    \item Let $F$ be a $(K_X+B+\Mm_X)$-negative extremal face in $\overline{NE}(X/U)$ that relatively ample at infinity (cf. Definition \ref{defn: basics of cone theorem}) with respect to $(X,B,\Mm)$. Then $F$ is a rational extremal face, and there exists a contraction$/U$ $\cont_F: X\rightarrow Z$ of $F$ satisfying the following.
\begin{enumerate}
    \item For any integral curve $C$ on $X$ such that the image of $C$ in $U$ is a closed point, $\cont_F(C)$ is a point if and only if $[C]\in F$.
    \item $\mathcal{O}_Y=(\cont_F)_*\mathcal{O}_X$. In other words, $\cont_F$ is a contraction.
    \item For any Cartier divisor $D$ on $Y$ such that $D\cdot C=0$ for any curve $C$ contracted by $\cont_F$, there exists a Cartier divisor $D_Y$ on $Y$ such that $D=\cont_F^*D_Y$.
\end{enumerate}
\end{enumerate}
\end{thm}
When $\Mm$ is NQC$/U$ and $(X,B,\Mm)$ is lc, Theorem \ref{thm: cone theorem nonnqc gpair}(1-3) was proven in \cite[Theorem 1.1]{HL21a} and Theorem \ref{thm: cone theorem nonnqc gpair}(4) was proven in \cite[Theorem 1.5]{Xie22} (see also \cite[Theorem 1.7]{CLX23}).

There are several other important results on the structure of generalized lc pairs. The first two are the Kodaira vanishing theorem and the Kawamata-Viehweg vanishing theorem:

\begin{thm}[Kodaira vanishing theorem for lc generalized pairs]\label{thm: kod vanishing gpair intro}
Let $(X,B,\Mm)$ be a projective lc generalized pair, and let $D$ be a Cartier divisor on $X$ such that $D-(K_X+B+\Mm_X)$ is ample. Then $H^i(X,\mathcal{O}_X(D))=0$ for any positive integer $i$.
\end{thm}

\begin{thm}[Relative Kawamata-Viehweg vanishing for lc generalized pairs]\label{thm: kv vanishing gpair intro}
Let $(X,B,\Mm)/U$ be an lc generalized pair associated with morphism $f: X\rightarrow U$, and let $D$ be a Cartier divisor on $X$ such that $D-(K_X+B+\Mm_X)$ is nef$/U$ and log big$/U$ with respect to $(X,B,\Mm)$. Then $R^if_*\mathcal{O}_X(D)=0$ for any positive integer $i$.
\end{thm}

When $\Mm$ is NQC$/U$, Theorem \ref{thm: kod vanishing gpair intro} was proven in \cite[Theorem 1.3]{CLX23} while Theorem \ref{thm: kv vanishing gpair intro} was proven in \cite[Theorem 1.4]{CLX23}.

\smallskip

Next, we have the base-point-freeness theorem and the semi-ampleness theorem for lc generalized pairs:

\begin{thm}[Base-point-freeness theorem]\label{thm:base-point-freeness intro}
Let $(X,B,\Mm)/U$ be an lc generalized pair and $D$ a nef$/U$ Cartier divisor on $X$, such that $aD-(K_X+B+\Mm_X)$ is ample$/U$ for some positive real number $a$. Then $\mathcal{O}_X(mD)$ is globally generated over $U$ for any integer $m\gg 0$.
\end{thm}

\begin{thm}[Semi-ampleness theorem]\label{thm: semi-ampleness intro}
Let $(X,B,\Mm)/U$ be an lc generalized pair and $D$ a nef$/U$ $\mathbb R$-Cartier $\Rr$-divisor on $X$, such that $D-(K_X+B+\Mm_X)$ is ample$/U$. Then $D$ is semi-ample$/U$.
\end{thm}

When $\Mm$ is NQC$/U$, Theorem \ref{thm:base-point-freeness intro} was proven in \cite[Theorem 1.4]{Xie22}, \cite[Theorem 1.5]{CLX23} while Theorem \ref{thm: semi-ampleness intro} was proven in \cite[Theorems 1.2]{Xie22}, \cite[Theorem 1.6]{CLX23}.

We also prove the canonical bundle formula and the subadjunction formula for generalized pairs. As the canonical bundle formula's statement is very technical and is a special case of Theorem \ref{thm: cbf gfq} below (by setting $\Ff=T_X$), we will omit it here and only state the subadjunction formula.

\begin{thm}[Subadjunction formula]\label{thm: subadj intro}
    Let $(X,B,\Mm)/U$ be an lc generalized pair an $V$ an lc center of $(X,B,\Mm)$ such that $\dim V\geq 1$. Let $W$ be the normalization of $V$. Then there exists an lc generalized pair $(W,B_W,\Mm^W)/U$ such that
    $$K_W+B_W+\Mm^W_W\sim_{\mathbb R}(K_X+B+\Mm_X)|_W.$$
    Moreover, the image of any lc center of $(W,B_W,\Mm^W)$ in $X$ is an lc center of $(X,B,\Mm)$.
\end{thm}
The main part of Theorem \ref{thm: subadj intro} was proven in \cite[Theorem 5.1]{HL21b} when $\Mm$ is NQC$/U$.

Finally, we can show that lc generalized pairs are Du Bois:

\begin{thm}\label{thm: glc sings are Du Bois}
Let $(X,B,\Mm)$ be an lc generlaized pair. Then any union of lc centers of $(X,B,\Mm)$ is Du Bois. In particular, $X$ is Du Bois.
\end{thm}
Theorem \ref{thm: glc sings are Du Bois} was proven in \cite[Theorem 1.6]{LX23b} when $\Mm$ is NQC$/X$.

\subsection{Minimal model program for algebraically integrable foliations}\label{subsec: mmp aif} 

The theory of foliations holds a significant place in birational geometry. Most notably, it has played a critical role in Miyaoka's proof of several key cases of the abundance conjecture in dimension three \cite{Miy87}. In recent developments, foliations have been used by Bogomolov and McQuillan to analyze projective varieties which admit a non-trivial fibration with rationally connected fibers \cite{BM16}. Furthermore, foliation theory has strong connections with other areas of algebraic geometry, such as the algebraic geometry in characteristic $p>0$ and number theory as highlighted by the Grothendieck-Katz conjecture and the Ekedahl-Shepherd-Barron-Taylor conjecture. Its importance is also highlighted in hyperbolicity theory, where it was essential in McQuillan's proof of a specific case of the Green-Griffiths-Lang conjecture \cite{McQ98}.

In recent years, it has been discovered that many structures and results in classical birational geometry can be extended to foliations, especially, within the context of the minimal model program. Instead of examining the structures associated with the canonical divisor of the ambient variety $K_X$, the foliations theory concentrates on the structures connected to the foliated canonical divisor $K_{\Ff}$. This approach offers greater flexibility in practice. Specifically, when $\Ff = T_X$, we find that $K_{\Ff}=K_X$, bringing us back to the classical setting.

The foundational work for the minimal model program for foliations has been established for foliated surfaces (cf. \cite{McQ08,Bru15}) and foliated threefolds (cf. \cite{CS20,Spi20,CS21,SS22}). Moreover, several classic questions from the minimal model programs, such as the ascending chain condition (ACC) conjecture for minimal log discrepancies, the ACC conjecture for lc thresholds, the global ACC, and the index theorems, have been adapted to foliations and verified in dimensions 2 and/or 3, as indicated in \cite{Che22,Che23,LLM23,LMX23a,LMX23b}.

Given these developments, it is natural to ask whether the minimal model program for foliations could extend to higher dimensions. Unfortunately, this seems to be a challenging question, with limited information available, even in dimension 4. However, from the perspective of the minimal model program, it seems sufficient to focus on a subset of foliations that have an additional structure: algebraically integrable foliations.

Algebraically integrable foliations are foliations where the general leaves are algebraic varieties; in other words, they are induced by dominant rational maps. These foliations naturally come into play when a fibration structure is established. Notably, Miyaoka's study of the abundance conjecture in dimension $3$ primarily utilized algebraically integrable foliations \cite{Miy87}, as opposed to arbitrary ones. This approach has been reflected in recent research into the abundance conjecture for Kähler threefolds \cite{DO23a,DO23b} and threefolds over fields of characteristic $p>3$ with numerical dimension $2$ \cite{Xu23}. In these studies, the algebraic integrability of foliations is guaranteed; indeed, all the foliations addressed in these papers are induced by MRC fibrations, making them automatically algebraically integrable. Given this, algebraically integrable foliations are expected to be crucial in future research of the minimal model program, particularly in questions related to the abundance conjecture.

The first objective of this paper is to develop the minimal model program for algebraically integrable foliations of arbitrary rank with ``mild" singularities in arbitrary dimensions. Here ``mild" singularity is usually referred to as ``F-dlt" (see Definition \ref{defn: fdlt}). As explained in \cite{CS21,SS22}, F-dlt foliated triples play the same role as dlt pairs in the classical MMP and is a natural class of singularities to study. Moreover, any terminal foliated singularity is F-dlt.

Recall that a foliated triple $(X,\Ff,B)/U$ consists of a normal quasi-projective variety $X$ associated with a projective surjective morphism $X\rightarrow U$, a foliation $\Ff$ on $X$, and an $\Rr$-divisor $B\geq 0$ on $X$, such that $K_{\Ff}+B$ is $\Rr$-Cartier. The first result of this paper shows that we can run a $(K_{\Ff}+B)$-MMP$/U$ provided that it is $\Qq$-factorial F-dlt:

\begin{thm}[Minimal model program]\label{thm: mmp fdlt}
    Let $(X,\Ff,B)/U$ be a $\Qq$-factorial foliated triple. Assume that $\Ff$ is algebraically integrable and $(X,\Ff,B)$ is F-dlt. Then we may run a $(K_{\Ff}+B)$-MMP$/U$.
\end{thm}

We remark that when $\dim X\leq 3$, Theorem \ref{thm: mmp fdlt} is known when $\rk\Ff=2$ (\cite[Corollary 2.3]{SS22}; \cite[Theorem 1.1]{CS21} when $U=\{pt\}$) and when $\rk\Ff=1$ and $U=\{pt\}$ (\cite[Theorems 1.1, 2.36, Section 6]{CS21}), even without the algebraically integrable condition. When assuming the termination of klt flips in dimension $\leq\rk\Ff$, \cite[Theorem 1.1]{CS23a} proves Theorem \ref{thm: mmp fdlt} even without the F-dlt condition, but requires that $(X,B)$ is klt. In particular, when $(X,B)$ is klt and $\dim X\leq 4$, Theorem \ref{thm: mmp fdlt} can be deduced from \cite[Theorem 1.1]{CS23a}.

Proceeding further, we demonstrate the termination of MMP with scaling as well as the existence of good minimal models for algebraically integrable foliations polarized with an ample divisor. The polarization of the ample divisor is a natural condition to add, as can be seen in \cite[Theorem 1.2]{CS21} and \cite[Theorem 1.3]{CS20}. It is worth noting that, even within the framework of the classical MMP, the existence of good minimal models in higher dimensions is only known when polarized with an ample divisor (\cite[Theorem C]{BCHM10}, \cite[Theorem 1.5]{HH20}),  while the general case remains an open conjecture.

\begin{thm}[Good minimal model]\label{thm: +a gmm fdlt}
    Let $(X,\Ff,B)/U$ be a $\Qq$-factorial foliated triple. Assume that $\Ff$ is algebraically integrable, $B\geq A\geq 0$ for some ample$/U$ $\Rr$-divisor $A$, and $(X,\Ff,B)$ is F-dlt. Then we may run a $(K_{\Ff}+B)$-MMP$/U$ with scaling of an ample$/U$ $\Rr$-divisor $H$, and any such MMP terminates 
    \begin{enumerate}
        \item with a Mori fiber space of $(X,\Ff,B)/U$ if $K_{\Ff}+B$ is not pseudo-effective$/U$, and
        \item with a good minimal model of $(X,\Ff,B)/U$ if $K_{\Ff}+B$ is pseudo-effective$/U$.
    \end{enumerate}
\end{thm}

We also have the following result on the abundance of algebraically integrable foliations polarized with an ample divisor.

\begin{thm}[Abundance]\label{thm: +a abundance fdlt}
  Let $(X,\Ff,B)/U$ be a $\Qq$-factorial foliated triple and $A$ an ample$/U$ $\Rr$-divisor on $X$.  Assume that $\Ff$ is algebraically integrable and $(X,\Ff,B)$ is F-dlt. Then 
  $$\kappa_{\sigma}(X/U,K_{\Ff}+B+A)=\kappa_{\iota}(X/U,K_{\Ff}+B+A).$$
\end{thm}

It is important to note that Theorem \ref{thm: +a abundance fdlt} is not a direct consequence of Theorem \ref{thm: +a gmm fdlt}. This is because Bertini type theorems fail for foliations, and it is possible that $(X,\Ff,B+H)$ is not lc for any $H\in |A/U|_{\mathbb R}$ (see \cite[Example 3.4]{DLM23}).

We also prove a base-point-freeness theorem for algebraically integrable foliations.

\begin{thm}[Base-point-freeness]\label{thm: bpf fdlt}
  Let $(X,\Ff,B)/U$ be a $\Qq$-factorial foliated triple. Assume that $\Ff$ is algebraically integrable and $(X,\Ff,B)$ is F-dlt. Let $A$ be an ample$/U$ $\Rr$-divisor on $X$ such that $K_{\Ff}+B+A$ is nef$/U$. Then:
  \begin{enumerate}
      \item $K_{\Ff}+B+A$ is semi-ample$/U$.
      \item Suppose that there exists a positive integer $m$ such that $m(K_{\Ff}+B+A)$ is Cartier. Then
      $$\mathcal{O}_X(mn(K_{\Ff}+B+A))$$
      is globally generated$/U$ for any integer $n\gg 0$.
    \end{enumerate}
\end{thm}
In the literature, the semi-ampleness of $K_\Ff+B+A$ is known when $(X,\Ff,B+A)$ is $\mathbb Q$-factorial F-dlt, $U=\{pt\}$, and $\dim X\leq 3$, even without the algebraically integrable condition (see \cite[Theorem 1.3]{CS20}, \cite[Theorem 1.3]{CS21}). However, there was no result on the base-point-freeness theorem of foliations in dimensions $\geq 3$. It is worth mentioning that the base-point-freeness theorem Theorem \ref{thm: bpf fdlt}(2) is crucial for us to prove a special case of the Prokhorov-Shokurov's $\bb$-semi-ampleness conjecture later in this paper (Theorem \ref{thm: ps intro}).

An important application of Theorem \ref{thm: +a gmm fdlt} is the existence of Mori fiber spaces for foliated triples, even with, at worst, lc singularities. We note that in this paper, Mori fiber spaces and log minimal models are in the sense of Birkar-Shokurov; that is, we allow the extraction of lc centers. See Definition \ref{defn: models I} for details.

\begin{thm}\label{thm: eomfs}
    Let $(X,\Ff,B)/U$ be an lc foliated triple. Assume that $\Ff$ is algebraically integrable and $K_{\Ff}+B$ is not pseudo-effective$/U$. Then $(X,\Ff,B)/U$ has a Mori fiber space.
\end{thm}

Another interesting type of foliations is the class of foliations with numerical dimension zero. For example, based on \cite[Theorem 1.4]{CS20} and \cite[Theorem 1.7]{CS21}, \cite[Theorem 1.9]{LLM23} has shown the existence of good minimal models for numerical dimension zero foliations in dimension $\leq 3$. In this paper, we obtain the existence of good minimal models for algebraically integrable foliations with numerical dimension zero:
\begin{thm}\label{thm: gmm ai num0}
    Let $(X,\Ff,B)$ be a projective lc foliated triple. Assume that $\Ff$ is algebraically integrable and $\kappa_{\sigma}(K_{\Ff}+B)=0$. Then:
   \begin{enumerate}
       \item $(X,\Ff,B)$ has a good minimal model.
       \item $\kappa_{\iota}(K_{\Ff}+B)=0$.
       \item If $(X,\Ff,B)$ is $\Qq$-factorial dlt, then we may run a $(K_{\Ff}+B)$-MMP with scaling of an ample $\Rr$-divisor, and any such MMP terminates with a good minimal model of $(X,\Ff,B)$.
   \end{enumerate}
\end{thm}

Siu \cite{Siu10} has used Eckl’s construction of numerically trivial foliations \cite{Eck04} to sketch a plan to solve the abundance conjecture. One step of Siu's approach, \cite[(4.1)]{Siu10}, focuses on the abundance conjecture for smooth projective varieties associated with an ``algebraically integrable numerically trivial foliation". Though the concept of ``numerically trivial foliation" in \cite{Siu10}, which was defined analytically in \cite{Eck04}, seems to differ from the concept of ``foliations whose canonical divisor has numerical dimension zero", these two types of foliations are closely connected. Hence, studying the abundance properties of numerical dimension zero algebraically integrable foliations (potentially with singularities that are worse than lc) on smooth projective varieties becomes intriguing, as it may have implications for the abundance conjecture. With this in mind, we prove the following theorem in this paper:

\begin{thm}\label{thm: abundance num 0 no restriction to f}
    Let $(X,\Ff,B)$ be a projective algebraically integrable f-triple such that $\kappa_{\sigma}(K_{\Ff}+B)=0$. Assume that $K_X+B$ is pseudo-effective and $(X,B)$ is lc. Then $\kappa_{\iota}(K_{\Ff}+B)=0$.
\end{thm}

Finally, we recall the following conjecture of Cascini and Spicer:

\begin{conj}[{\cite[Conjecture 4.2]{CS23a}}]\label{conj: cs23 4.2(1)}
Let $(X,\Ff,B)$ be a $\Qq$-factorial projective foliated triple, such that $\Ff$ is algebraically integrable, $B$ is a $\Qq$-divisor, $(X,B)$ is klt, and one of the following cases hold:
\begin{enumerate}
    \item $(X,\Ff,B)$ is F-dlt.
    \item $(X,\Ff,B)$ is canonical. 
\end{enumerate}
 Then there exists a morphism $f: X\rightarrow Y$ which induces $\Ff$. 
\end{conj}
In this paper, we provide a positive answer to Conjecture \ref{conj: cs23 4.2(1)}(1) with weaker assumptions and stronger results:
\begin{thm}\label{thm: cs23 4.2(1)}
    Let $(X,\Ff,B)$ be a $\Qq$-factorial foliated triple such that $\Ff$ is algebraically integrable and $(X,\Ff,B)$ is F-dlt. Then:
    \begin{enumerate}
    \item $(X,B)$ is qdlt (cf. Definition \ref{defn: qdlt}). In particular, if $\lfloor B\rfloor=0$, then $(X,B)$ is klt.
    \item There exists a morphism $f: X\rightarrow Y$ to a smooth variety which induces $\Ff$.
    \end{enumerate}
\end{thm}

We also prove a weaker form of Conjecture \ref{conj: cs23 4.2(1)}(2) without assuming that $(X,B)$ is klt.
\begin{thm}\label{thm: canonical almost holomorphic}
    Let $(X,\Ff,B)$ be a $\Qq$-factorial canonical foliated triple such that $\Ff$ is algebraically integrable. Then $\Ff$ is induced by an almost holomorphic map.
\end{thm}
A very recent result \cite[Theorem 1.4]{CS23b} shows that the algebraic part of a foliation $\Ff$ on a projective variety $X$ is induced by an almost holomorphic map, provided that $X$ is $\Qq$-factorial klt and $\Ff$ is canonical. In particular, \cite[Theorem 1.4]{CS23b} implies Theorem \ref{thm: canonical almost holomorphic} when $X$ is projective klt and $B=0$.

\subsection{Generalized foliated quadruples}

As explained above, to establish the minimal model program for algebraically integrable foliations and generalized pairs, we need to broaden the category of objects for our study and consider the structure of \emph{generalized foliated quadruples} $(X,\Ff,B,\Mm)$, as defined in Definition \ref{defn: gfq intro}.

Most of the main theorems of this paper on algebraically integrable foliations can also be extended to the category of algebraically integrable generalized foliated quadruples. Two results related to this structure that are particularly worth mentioning are the cone theorem and the canonical bundle formula. These two results will be essential in other main theorems of the paper, the statements of most of which do not rely on the language of generalized foliated quadruples.

\subsubsection{Cone theorem} We establish the cone theorem for algebraically integrable generalized foliated quadruples in full generality. 

\begin{thm}[Cone theorem for algebraically integrable generalized foliated quadruples]\label{thm: cone theorem gfq}
Let $(X,\Ff,B,\Mm)/U$ be a generalized foliated quadruple and $\pi: X\rightarrow U$ the associated morphism. Let $\{R_j\}_{j\in\Lambda}$ be the set of $(K_{\Ff}+B+\Mm_X)$-negative extremal rays in $\overline{NE}(X/U)$ that are rational. Assume that $\Ff$ is algebraically integrable. Then:
\begin{enumerate}
    \item $$\overline{NE}(X/U)=\overline{NE}(X/U)_{K_{\Ff}+B+\Mm_X\geq 0}+\overline{NE}(X/U)_{\Nlc(X,\Ff,B,\Mm)}+\sum_{j\in\Lambda} R_j.$$
    Here $\Nlc(X,\Ff,B,\Mm)$ is the non-lc locus of $(X,\Ff,B,\Mm)$ (cf. Definition \ref{defn: gfq singularity}). In particular, any $(K_{\Ff}+B+\Mm_X)$-negative extremal ray in $\overline{NE}(X/U)$ is rational.
    \item Each $R_j$ is spanned by a rational curve $C_j$ such that $\pi(C_j)=\{pt\}$, $C_j$ is tangent to $\Ff$, and 
    $$0<-(K_{\Ff}+B+\Mm_X)\cdot C_j\leq 2\dim X.$$
    \item For any ample$/U$ $\Rr$-divisor $A$ on $X$,
    $$\Lambda_A:=\{j\in\Lambda\mid R_j\subset\overline{NE}(X/U)_{K_{\Ff}+B+A+\Mm_X<0}\}$$
    is a finite set. In particular, $\{R_j\}_{j\in\Lambda}$ is countable, and is a discrete subset in $\overline{NE}(X/U)_{K_{\Ff}+B+\Mm_X<0}$. Moreover, we may write
    $$\overline{NE}(X/U)=\overline{NE}(X/U)_{K_{\Ff}+B+A+\Mm_X\geq 0}+\overline{NE}(X/U)_{\Nlc(X,\Ff,B,\Mm)}+\sum_{j\in\Lambda_A}R_j.$$
    \item Let $F$ be a $(K_X+B+\Mm_X)$-negative extremal face in $\overline{NE}(X/U)$ that relatively ample at infinity (cf. Definition \ref{defn: basics of cone theorem}) with respect to $(X,\Ff,B,\Mm)$. Then $F$ is a rational extremal face.
\end{enumerate}
\end{thm}
 When $\Ff=T_X$ and $\Mm=\bm{0}$, Theorem \ref{thm: cone theorem gfq} follows from \cite[Theorem 5.10]{Amb03} and \cite[Theorem 4.5.2]{Fuj17}. However, whenever either $\Ff\not=T_X$ or $\Mm\not=\bm{0}$, Theorem \ref{thm: cone theorem gfq} becomes new. More precisely, there are two cases that worth to mention:
 \begin{enumerate}
 \item When $\Ff=T_X$, we get the cone theorem for generalized pairs, Theorem \ref{thm: cone theorem nonnqc gpair}, which is new.
 \item When $U=\{pt\}$ and $\Mm=\bm{0}$, (2) and a large part of (1) (the part without considering he rationality of $R_j$) were proven in \cite[Theorem 3.9]{ACSS21}, but the rest parts are missing. Therefore, we cannot directly use \cite[Theorem 3.9]{ACSS21} to prove Theorem \ref{thm: mmp fdlt} and Theorem \ref{thm: cone theorem gfq} becomes necessary.
 \end{enumerate}
We would like to note that the contraction theorem, the existence of flips, and the base-point-freeness theorem are still valid for generalized foliated quadruples that possess nice singularities (e.g. F-dlt). However, since these theorems do not hold the same level of importance in proving our other main theorems as the cone theorem does, we choose to omit them here.

\subsubsection{Canonical bundle formula} The canonical bundle formula for foliated triples, as established in \cite[Theorem 1.3]{LLM23}, plays a crucial role in proving the global ACC for foliated threefolds. However, due to technical challenges, the work presented in \cite{LLM23} could not prove the canonical bundle formula for generalized foliated quadruples $(X,\Ff,B,\Mm)$ unless the nef part $\Mm$ is $\bb$-semi-ample. In this study, we overcome these technical difficulties with innovative approaches, successfully proving the canonical bundle formula for generalized foliated quadruples in a more comprehensive manner.

\begin{thm}\label{thm: cbf gfq}
Let $(X,\Ff,B,\Mm)/U$ be a sub-generalized foliated quadruple and $f: X\rightarrow Z$ a contraction$/U$, such that $f: (X,\Ff,B,\Mm)\rightarrow Z$ is an lc-trivial morphism (see Definition \ref{defn: lc trivial morphism}). Let $B_Z$ and $\Mm^Z$ be the discriminant part and the base moduli part of $f: (X,\Ff,B,\Mm)\rightarrow Z$ (see Definition-Theorem \ref{defthm: cbf lctrivial morphism}). Then $\Mm^Z$ is $\bb$-nef$/U$ and
$$K_{\Ff}+B+\Mm_X\sim_{\mathbb R}f^*(K_{\Ff_Z}+B_Z+\Mm^Z_Z).$$
Moreover, we have the following properties:
\begin{enumerate}
\item $B_Z$ is uniquely determined and $\Mm^Z$ is uniquely determined up to $\Rr$-linear equivalence.
\item If $B\geq 0$, then $B_Z\geq 0$.
\item If $(X,\Ff,B,\Mm)$ is sub-lc, then $(Z,\Ff_Z,B_Z,\Mm^Z)$ is sub-lc.
\item If $(X,\Ff,B,\Mm)$ is lc, then $(Z,\Ff_Z,B_Z,\Mm^Z)$ is lc.
\item If $(X,\Ff,B,\Mm)$ is sub-lc or $f$ has connected fibers, then any lc center of $(Z,\Ff_Z,B_Z,\Mm^Z)$ is the image of an lc center of $(X,\Ff,B,\Mm)$ on $Z$.
\item If $f$ has connected fibers, then the image of any lc center of $(X,\Ff,B,\Mm)$ on $Z$ is an lc center of $(Z,\Ff_Z,B_Z,\Mm^Z)$.
\item If $f$ has connected fibers, then for any prime divisor $D$ on $X$, 
$$\mult_DB_Z=\epsilon(D)-\sup\{t\mid (X,\Ff,B+tf^*D,\Mm)\text{ is lc over the generic point of }D\}$$
where $\epsilon(D)=0$ if $D$ is $\Ff_Z$-invariant, and $\epsilon(D)=1$ otherwise (see Definition \ref{defn: special divisors on foliations}). 
\item The $\Rr$-linear equivalence class of $\Mm^Z$ only depends on $(X,B,\Mm)$ over the generic point of $Z$.
\item If $\Mm$ is NQC$/U$, then $\Mm^Z$ is NQC$/U$.
\end{enumerate}
\end{thm}

We want to emphasize that Theorem \ref{thm: cbf gfq} is applicable to any foliation, not just those that are algebraically integrable. Consequently, we anticipate that Theorem \ref{thm: cbf gfq} will play a significant role in future studies, encompassing both algebraically integrable foliations and those that are not necessarily algebraically integrable.

Next, we revisit the history of partial results that have contributed to the main part of Theorem \ref{thm: cbf gfq}, i.e. the nefness$/U$ of $\Mm^Z$. 
\begin{enumerate}
    \item When $\Ff=T_X$ and $\Mm=\bm{0}$, the main part of Theorem \ref{thm: cbf gfq} is\cite[Theorem 1.2]{JLX22} (or \cite[Theorem 2.23]{JLX22} combined with \cite[Lemma 1.1]{FG12}). For other related references, see \cite{Kaw98,Amb05,Kol07,Flo14,FG14}.
    \item When $\Ff=T_X$ and $\Mm$ is NQC$/U$, previously we only knew the cases where either $B\geq 0$ at the generic point of $Z$ or $\Mm$ is $\bb$-semi-ample$/Z$ (\cite[Theorem 2.23]{JLX22}+\cite[Theorem 4.5]{HL21b}). We direct the reader to \cite{Fil19,Fil20,FS23} for other related references. It is worth noting that no results were known when $\Mm$ is not NQC$/U$.
    \item When $\Ff \neq T_X$, we only knew the cases where $f$ is a contraction and $\Mm$ is NQC$/U$ and $\bb$-semi-ample$/Z$ (\cite[Proposition 6.14]{LLM23}).
\end{enumerate}
In this paper, we not only prove Theorem \ref{thm: cbf gfq} in full generality but also clarify why \cite{Kol07} was able to address the horizontal negative coefficients, while \cite{Fil19,Fil20,JLX22,FS23} cannot deal with this issue. Further details on this are provided in Remark \ref{rem: lc trivial fibration definition}. Following this discussion, we refine the definition of lc-trivial fibrations and lc-trivial morphisms, which are elaborated in Definition \ref{defn: lc trivial fibration gfq}.

Finally, we note that the proof of Theorem \ref{thm: cbf gfq} does not depend on the mixed Hodge structure, as opposed to what is done in \cite{Kol07}. Remark \ref{rem: lc trivial fibration definition} also explains why the mixed Hodge structure cannot be applied to our case. Instead, our approach is based on the structure of algebraically integrable foliations, a method similar to the one used in \cite[Proof of Theorem 1.3]{ACSS21}. However, the canonical bundle formula in that reference is not complete as the ``BP stable" condition is required. Moreover, \cite[Theorem 1.3]{ACSS21} also has the additional requirement that $B\geq 0$ over the generic point of $Z$, a condition we aim to avoid.

\subsection{Singularities of algebraically integrable generalized foliated quadruples}

The cone theorem and the canonical bundle formula are primarily concerned with understanding the global behavior of algebraically integrable generalized foliated quadruples. However, it is equally important to examine the local behavior, particularly the singularities of these structures. In this paper, we will concentrate on two key aspects that are tied to the singularity structure of generalized foliated quadruples: the precise adjunction formula and the ACC for lc thresholds.

\subsubsection{Adjunction formulas}
\cite[Proposition 3.2]{ACSS21} proves the adjunction formula for algebraically integrable foliations provided that the ambient variety is $\Qq$-factorial and that the foliation is induced by a contraction. In this paper, we remove these two technical conditions and prove the adjunction formula for algebraically integrable foliations in full generality:

\begin{thm}\label{thm: fol adj intro}
    Let $(X,\Ff,B)$ be an foliated triple such that $\Ff$ is algebraically integrable. Let $S$ be a prime divisor on $X$, such that $\mult_SB=0$ if $S$ if $\Ff$-invariant and $\mult_SB=1$ otherwise.  Let $S^\nu$ be the normalization of $S$ and $\Ff_S$ the restricted foliation (see Definition \ref{defn: restricted foliation}) of $\Ff$ on $S^\nu$. Then
    $$K_{\Ff_S}+B_S=(K_{\Ff}+B)|_{S^\nu}$$
    for some $\Rr$-divisor $B_S\geq 0$. Moreover, if $(X,\Ff,B)$ is lc, then $(S^\nu,\Ff_S,B_S)$ is lc.
\end{thm}
We remark that \cite[Theorem 3.16]{CS23b} proves the adjunction formula to non-invariant divisors for any foliation with a minor requirement that the boundary has $\Qq$-coefficients. In particular, when $B$ has rational coefficients and $\mult_SB=1$, Theorem \ref{thm: fol adj intro} is implied by \cite[Theorem 3.16]{CS23b}.

Theorem \ref{thm: fol adj intro} can be extended to the category of algebraically integrable generalized foliated quadruples:

\begin{thm}\label{thm: lc adjunction foliation nonnqc}
  Let $(X,\Ff,B,\Mm)/U$ be a generalized foliated quadruple such that $\Ff$ is algebraically integrable. Let $S$ be a prime divisor on $X$, such that $\mult_SB=0$ if $S$ if $\Ff$-invariant, and $\mult_SB=1$ otherwise. Let $S^\nu$ be the normalization of $S$, $\Mm^S:=\Mm|_{S^\nu}$ (see Definition \ref{defn: restriction b divisor}), and $\Ff_S$ the restricted foliation of $\Ff$ on $S^\nu$. Then
  $$K_{\Ff_S}+B_S+\Mm^S_{S^\nu}:=(K_{\Ff}+B+\Mm_X)|_{S^\nu}$$
  for some $\Rr$-divisor $B_S\geq 0$. Moreover, if $(X,\Ff,B,\Mm)$ is lc, then $(S^\nu,\Ff_S,B_S,\Mm^S)$ is lc.
\end{thm}

In \cite[Theorem 1.6]{DLM23}, a precise adjunction formula was introduced for algebraically integrable foliated triples, playing a crucial role in proving the ACC for lc thresholds and the global ACC for algebraically integrable foliated triples. Building upon this concept, we formulate and establish a precise adjunction formula for algebraically integrable generalized foliated quadruples in this paper. We then apply it to prove the ACC for lc thresholds and the global ACC for algebraically integrable generalized foliated quadruples. We have the following theorem:

\begin{thm}\label{thm: dcc adjunction is dcc}
Let $\Ii\subset [0,+\infty)$ be a set of real numbers. Let $(X,\Ff,B,\Mm)/U$ be an lc generalized foliated quadruple, $S$ a prime divisor on $X$ with normalization $S^\nu$, such that $\mult_SB=0$ if $S$ is $\Ff$-invariant and $\mult_SB=1$ otherwise. Assume that the coefficients of $B$ belong to $\Ii$ and $\Mm$ is a $\Ii$-linear combination of nef$/U$ $\bb$-Cartier $\bb$-divisors. Let
 $$K_{\Ff_S}+B_S+\Mm^S_{S^\nu}:=(K_{\Ff}+B+\Mm_X)|_{S^\nu}$$
where $\Ff_S$ is the restricted foliation of $\Ff$ on $S^\nu$ and $\Mm^S=\Mm|_{S^\nu}$.

Then the coefficients of $B_S$ belong to $D(\Ii)$ (see Definition \ref{defn: derived set}). In particular, if $\Ii$ is a DCC set, then the coefficients of $B_S$ belong to a DCC set.
\end{thm}

We offer a more detailed version of Theorem \ref{thm: dcc adjunction is dcc} in Theorem \ref{thm: precise adj gfq}. Given its highly technical nature, we omit it from the introduction. 

\subsubsection{ACC and the global ACC} Theorem \ref{thm: dcc adjunction is dcc} leads to the proof of the ACC and the global ACC for algebraically integrable generalized foliated quadruples. For foliated triples, these results were proven in \cite[Theorem 1.1]{DLM23} and \cite[Theorem 1.2]{DLM23} respectively.

\begin{thm}[ACC for lc thresholds for algebraically integrable generalized foliated quadruples]\label{thm: acc lct alg int gfq}
Let $r$ be a positive integer and $\Ii\subset [0,+\infty)$ a DCC set. Then there exists an ACC set $\Ii'$ depending only on $r$ and $\Ii$ satisfying the following. Let $(X,\Ff,B,\Mm)/X$ be an lc generalized foliated quadruple, such that 
\begin{enumerate}
    \item $\Ff$ is algebraically integrable of rank $r$, 
    \item the coefficients of $B$ belong to $\Ii$, and
    \item  $\Mm$ is a $\Ii$-linear combination of nef$/X$ $\bb$-Cartier $\bb$-divisors.
\end{enumerate}
Then the lc threshold
$$\lct(X,\Ff,B,\Mm;D,\Nn):=\sup\{t\mid t\geq 0, (X,\Ff,B+tD,\Mm+t\Nn)\text{ is lc}\}$$
is contained in $\Ii'$.
\end{thm}

\begin{thm}[Global ACC for algebraically integrable generalized foliated quadruples]\label{thm: global acc alg int gfq}
Let $r$ be a positive integer and $\Ii\subset [0,1]$ a DCC set. Then there exists a finite set $\Ii_0\subset\Ii$ depending only on $r$ and $\Ii$ satisfying the following. Let $(X,\Ff,B,\Mm)$ be a projective lc generalized foliated quadruple such that 
\begin{enumerate}
\item $\Ff$ is algebraically integrable of rank $r$, 
\item the coefficients of $B$ belong to $\Ii$,
\item  $\Mm=\sum\gamma_j\Mm_j$, where each $\gamma_j\in\Ii$ and each $\Mm_j$ is a nef $\bb$-Cartier $\bb$-divisor,
\item $\Mm_j\not\equiv\bm{0}$ if $\gamma_j\not=0$, and
\item $K_{\Ff}+B+\Mm_X\equiv 0$.
\end{enumerate}
Then the coefficients of $B$ belong to $\Ii_0$, and $\gamma_j\in\Ii_0$ for each $j$.
\end{thm}
The proof of Theorem \ref{thm: global acc alg int gfq} is harder than the proof of the global ACC for algebraically integrable foliated triples \cite[Theorem 1.2]{DLM23}. This is because our proof heavily relies on the existence of Mori fiber spaces (Theorem \ref{thm: eomfs}), unlike the proof in \cite[Theorem 1.2]{DLM23}.

As a straightforward corollary of Theorem \ref{thm: global acc alg int gfq}, we obtain the global ACC for rank one generalized foliated quadruples:

\begin{cor}[Global ACC for rank one generalized foliated quadruples]\label{cor: global acc rank 1 gfq}
Let $\Ii\subset [0,1]$ be a DCC set. Then there exists a finite set $\Ii_0\subset\Ii$ depending only on $\Ii$ satisfying the following. Let $(X,\Ff,B,\Mm)$ be an lc generalized foliated quadruple such that
\begin{enumerate}
\item $\rk\Ff=1$,
\item the coefficients of $B$ belong to $\Ii$, 
\item  $\Mm=\sum\gamma_j\Mm_j$, where each $\gamma_j\in\Ii$ and each $\Mm_j$ is a nef $\bb$-Cartier $\bb$-divisor,
\item $\Mm_j\not\equiv\bm{0}$ if $\gamma_j\not=0$, and
\item $K_{\Ff}+B+\Mm_X\equiv 0$.
\end{enumerate}
Then the coefficients of $B$ belong to $\Ii_0$, and $\gamma_j\in\Ii_0$ for each $j$.
\end{cor}

When $\Mm=\bm{0}$, Corollary \ref{cor: global acc rank 1 gfq} was proven in \cite[Corollary 1.3]{DLM23}.

\subsubsection{Uniform rational polytopes} As a direct consequence of Theorems \ref{thm: acc lct alg int gfq} and \ref{thm: global acc alg int gfq}, we establish the existence of uniform lc rational polytopes for algebraically integrable generalized foliated quadruples. Despite their complex nature, these polytopes are powerful tools in birational geometry. They are notably used in several applications on the ACC conjecture for minimal log discrepancies and the boundedness of complements. These polytopes are essential for the formal definition of KSBA moduli spaces \cite[6.27.3, Theorem 11.49]{Kol23}, and play a crucial role in proving the global ACC for foliated threefolds \cite{LMX23b}. In this paper, we prove the existence of uniform lc rational polytopes for algebraically integrable generalized foliated quadruples:

\begin{thm}\label{thm: uniform rational polytope gfq}
Let $r$ be a positive integer, $v_1^0,\dots,v_m^0,u_1^0,\dots,u_n^0$ positive real numbers, $\bm{v}_0:=(v_1^0,\dots,v_m^0)$, and $\bm{u}_0:=(u_1^0,\dots,u_n^0)$. Then there exists an open set $U\ni (\bm{v}_0,\bm{u}_0)$ of the rational envelope of $(\bm{v}_0,\bm{u}_0)$ in $\mathbb R^{m+n}$ depending only on $r$ and $\bm{v}_0$, $\bm{u}_0$ satisfying the following. Let $$\left(X,\Ff,B=\sum_{j=1}^mv_j^0B_j,\Mm=\sum_{k=1}^nu_k^0\Mm_k\right)\Bigg{/}X$$ be an lc generalized foliated quadruple, such that $\Ff$ is algebraically integrable, $\rk\Ff=r$, $B_j\geq 0$ are distinct Weil divisors, and $\Mm_k$ are nef$/X$ $\bb$-Cartier $\bb$-divisors. Then $$\left(X,\Ff,B=\sum_{j=1}^mv_jB_j,\sum_{k=1}^nu_k\Mm_k\right)$$ is lc for any $(v_1,\dots,v_m,u_1,\dots,u_n)\in U$.
\end{thm}

\subsection{Miscellaneous results on the minimal model program and foliations} We also prove several other interesting theorems that can be useful for further applications.

\subsubsection{Analogues of dlt models} We establish the existence of $(*)$ models and ACSS models for algebraically integrable generalized foliated quadruples (see Definitions \ref{defn: ACSS f-triple} and \ref{defn: acss model}). As detailed in \cite{ACSS21,CS23a,DLM23}, these models play the same role as dlt models in the classic MMP. Moreover, $\Qq$-factorial dlt algebraically integrable generalized foliated quadruples always satisfy the property ``ACSS" and the property $(*)$ (see Theorem \ref{thm: fdlt is acss}).

\begin{thm}[Existence of ACSS model]\label{thm:  ACSS model}
    Let $(X,\Ff,B,\Mm)/U$ be an lc generalized foliated quadruple. Assume that $\Ff$ is algebraically integrable. Then $(X,\Ff,B,\Mm)/U$ has an ACSS model which is also a $(*)$ model.
    
    In particular, there exists a birational morphism $f: Y\rightarrow X$ and a contraction $\pi: Y\rightarrow Z$ satisfying the following. Let $\Ff_Y:=f^{-1}\Ff$ and
    $$K_{\Ff_Y}+B_Y+\Mm_Y=f^*(K_{\Ff}+B+\Mm_X),$$
    then
    \begin{enumerate}
        \item  $(Y,B_Y,\Mm)$ is $\Qq$-factorial qdlt, 
        \item $\pi$ is equi-dimensional and  $\Ff_Y$ is induced by $\pi$, and
        \item any prime $f$-exceptional divisor is an lc place of $(X,\Ff,B,\Mm)$.
    \end{enumerate}
\end{thm}

\subsubsection{Minimal model program for very exceptional divisors} When running the relative minimal model program, especially the birational minimal model program, we often encounter the minimal model program for very exceptional divisors \cite[Theorem 1.8]{Bir12}. In this paper, we establish the minimal model program for algebraically integrable generalized foliated quadruples whose generalized foliated log canonical divisor is very exceptional. 

\begin{thm}\label{thm: gfq mmp very exceptional intro}
    Let $(X,\Ff,B,\Mm)/U$ be a $\Qq$-factorial F-dlt generalized foliated quadruple and $E\geq 0$ and $\Rr$-divisor on $X$, such that $E$ is very exceptional$/U$ and
    $$K_{\Ff}+B+\Mm_X\sim_{\mathbb R,U}E.$$
    Then we may run a $(K_{\Ff}+B+\Mm_X)$-MMP$/U$ with scaling of an ample$/U$ $\Rr$-divisor, and any such MMP terminates with a good minimal model $(X',\Ff',B',\Mm)/U$ of $(X,\Ff,B,\Mm)/U$, such that $K_{\Ff'}+B'+\Mm_{X'}\sim_{\mathbb R,U}0$. 
\end{thm}
Theorem \ref{thm: gfq mmp very exceptional intro} is vital for proving that $\Qq$-factorial dlt implies ACSS (Theorem \ref{thm: fdlt is acss}). This is essential for proving Theorem \ref{thm: main mmp foliation}.

\subsubsection{A special case of Prokhorov-Shokurov's base-point-freeness conjecture}
Prokhorov-Shokurov's base-point-freeness conjecture \cite[Conjecture 7.13]{PS09} is a major conjecture in birational geometry. It has been verified when the relative dimension of the fibration is $1$ (\cite[Theorem 8.1]{PS09}) or $2$ (\cite[Theorem 1.4]{ABBDILW23}). However, for the relative dimension of the fibration $\geq 3$, the conjecture is still largely open. Through the application of foliation theory, we prove a special case of the Prokhorov-Shokurov base-point-freeness conjecture:
\begin{thm}\label{thm: ps intro}
    Let $d$ and $m$ be two positive integers. Then there exists a positive integer $I$ depending only on $d$ and $I$ satisfying the following.
    
    Let $(X,B)$ be a projective klt pair, and $\pi: X\rightarrow Z$ a contraction to a smooth variety $Z$. Let $B^h$ and $B^v$ 
 be the horizontal$/Z$ part of $B$ and the vertical$/Z$ part of $B$ respectively, and let $\Mm$ be the moduli part of $\pi: (X,B)\rightarrow Z$. Assume that the following conditions hold.
    \begin{enumerate}
        \item (Semi-stability) $(X,B)$ is BP semi-stable$/Z$.
        \item (Klt-trivial) $K_X+B\sim_{\mathbb Q,Z}0$.
        \item (Coefficient control) $mB$ is a Weil divisor.
        \item (Fano type) There exists an ample $\Rr$-divisor $H$ such that $B^h\geq H\geq 0$.
        \item (Snc condition) There exists a reduced divisor $\Sigma_Z$ on $Z$, such that $B^v=\pi^{-1}(\Sigma_Z)$, $(Z,\Sigma_Z)$ is log smooth, and for any reduced divisor $H\geq 0$ such that $(Z,\Sigma+H)$ is log smooth, $(X,B+\pi^*H)$ is lc.
    \end{enumerate}
   Then $\Mm$ descends to $X$, $I\Mm_X$ is Cartier, and $nI\Mm_X$ is base-point-free for any integer $n\gg 0$.
\end{thm}
While the conditions of Theorem \ref{thm: ps intro} are quite restrictive, mainly because condition (4) is not preserved under birational transformations, there are no requirements on the dimension of the varieties or the relative dimension of $\pi$. This makes the theorem potentially useful for future applications.

\subsection{Why should we care about generalized foliated quadruples?}\label{sec: reason to consider gfq}

Before we move on to the main part of the paper, we would like to briefly explain why we need to consider the structure of generalized foliated quadruples and why it is essential for us to prove some main theorems of the paper, even if these theorems' statements do not explicitly mention this structure. To clarify this, we present five scenarios where generalized foliated quadruples play a crucial role. Four out of these five scenarios are unavoidable in the proofs of this paper.

\begin{sce}[Canonical bundle formula of foliations]
\cite[Theorem 1.2]{LLM23} established the canonical bundle formula for foliations. More precisely, given a projective lc foliated triple $(X,\Ff,B)$ and a contraction $f: X\rightarrow Z$ such that the general fibers of $f$ are tangent to $\Ff$ and $K_{\Ff}+B\sim_{\mathbb R,Z}0$, we have
$$K_{\Ff}+B\sim_{\mathbb R}f^*(K_{\Ff_Z}+B_Z+\Mm_Z)$$
where $(Z,\Ff_Z,B_Z,\Mm)$ is a projective lc generalized foliated quadruple. Therefore, if we want to study the behavior of $(X,\Ff,B)$, then it is necessary to study the structure of $(Z,\Ff_Z,B_Z,\Mm)$. When $\Ff=T_X$ and $(X,B)$ is klt, the classical approach is to find an $\Rr$-divisor
$$0\leq\Delta_Z\sim_{\mathbb R}B_Z+\Mm_Z$$
such that $(Z,\Delta_Z)$ is klt \cite[Theorem 0.2]{Amb05} and use the structure of $(Z,\Delta_Z)$ instead of $(Z,B_Z,\Mm)$. This approach is essential for the proof of the finite generation of the canonical ring (\cite[Proof of Corollary 1.1.2]{BCHM10}, \cite{FM00}). 

However, when $\Ff\not=T_X$, it is generally not possible for us to combine $B_Z$ and $\Mm_Z$ and get an lc triple structure $(Z,\Ff_Z,\Delta_Z\sim_{\mathbb R}B_Z+\Mm_Z)$. This is because of the following two reasons:
\begin{enumerate}
    \item ``Klt" is almost an empty condition for foliations when $\Ff\not=T_X$. Actually, unless the foliation is purely transcendental, there are always (infitely many) lc centers of $(X,\Ff,B)$ as long as $\Ff\not=T_X$. This will cause trouble when we try to perturb the coefficients of $B_Z+\Mm_Z$ and get $\Delta_Z$. In fact, even when $\Ff=T_X$, if $(X,B)$ is not klt, then we do not know whether there exists such $\Delta_Z$ so that $(Z,\Delta_Z)$ is lc, and we usually need the $\bb$-semi-ampleness of $\Mm$ to show this fact. The $\bb$-semi-ampleness of $\Mm$, on the other hand, is the Prokhorov-Shokurov conjecture \cite[Conjecture 7.13]{PS09} as mentioned above, which is widely open when $\dim X-\dim Z\geq 3$. Indeed, the unconfirmed status of the Prokhorov-Shokurov conjecture is one key reason why Birkar-Zhang introduced the concept of generalized pairs in \cite{BZ16}.
    \item Even if $\Mm$ is semi-ample, the existence of such $\Delta_Z$ is also unknown as Bertini type theorems fail for foliations in general, even for surfaces (cf. \cite[Example 3.4]{DLM23}). In other words, it is possible that $(Z,\Ff_Z,B_Z+G_Z)$ is not lc for any $G_Z\in |\Mm_Z|_{\mathbb R}$. 
\end{enumerate}
Therefore, in many situations, we must analyze the structure of generalized foliated quadruples rather than foliated triples. Furthermore, due to (2), the concept of generalized foliated quadruples becomes essential for studying foliations, even in lower dimensions. This is a key reason why \cite{LLM23,LMX23b} rely on the theory of generalized foliated quadruples to establish the global ACC for foliated triples in dimension $3$.
\end{sce}

\begin{sce}[MMP with scaling]
We recall how we run the minimal model program for with scaling for usual pairs. For simplicity, we only consider the projective case. Given a projective lc pair $(X,B)$ and an  $\Rr$-divisor $A\geq 0$ on $X$ such that $K_X+B+A$ is nef, we consider the scaling numbers
$$\lambda:=\inf\{t\mid t\geq 0, K_X+B+tA\text{ is nef}\}.$$
If $t=0$ then we are done. Otherwise, we contract a $(K_X+B)$-negative extremal ray $R$ such that $(K_X+B+tA)\cdot R=0$, and let $f: (X,B)\dashrightarrow (X',B')$ be a corresponding divisorial contraction, flip, or Mori fiber space associated to the contraction of $R$. We may replace $(X,B)$ with $(X',B')$ and $A$ with $A'$ and continue this process.

Although we only need $K_X+B+A$ to be nef to run the MMP, in practice, we usually also need the additional condition that $(X,B+A)$ is lc. This is helpful in many situations: since we do not know the termination of the MMP, it is likely for us to consider pairs $(X,B+\mu A)$ where $\mu$ is related to the scaling numbers $\lambda$. In this case, we usually need $(X,B+A)$ to be lc in order to guarantee that $(X,B+\mu A)$ is lc. For this reason, for the very first step of the MMP, we usually require $A$ to be a general ample $\Rr$-divisor, or a general base-point-free big and nef $\Rr$-divisor when $(X,B)$ is klt.

Now we consider the minimal model program for foliations. We definitely want to consider the minimal model program with scaling of ample divisors as well. However, as we have explained above, Bertini type theorems fail for foliations in general, even for surfaces (cf. \cite[Example 3.4]{DLM23}). Therefore, it is possible that for any ample $\Rr$-divisor $A\geq 0$ on $X$, $(X,\Ff,B+A)$ is not lc. Now the minimal model program of $(X,\Ff,B)$ with scaling of $A$ becomes weird: we can still run the minimal model program, but it will become difficult to study the intermediate outputs $(X',B'+\lambda A')$ with $\lambda>0$ after each step of the MMP, where $\lambda$ is the scaling number. This causes a lot of inconvenience for the minimal model program of foliations.

The structure of generalized foliated quadruples, however, can easily resolve this issue: if we identify $(X,\Ff,B)$ with the generalized foliated quadruple $(X,\Ff,B,\overline{0})$, then instead of running an MMP with scaling of an ample $\Rr$-divisor $A$, we can let $\Aa:=\overline{A}$ be the nef $\bb$-divisor associated to $A$. Now may run an MMP with scaling of $(0,\Aa)$. That is, although we still consider
$$\lambda:=\inf\{t\mid t\geq 0, K_{\Ff}+B+tA\text{ is nef}\},$$
the output of the first step of the MMP $\phi: X\dashrightarrow X'$ becomes $$(X',\Ff'=\phi_*\Ff,B'=\phi_*B,\lambda\Aa)$$ which is still an lc generalized foliated quadruple. Therefore, by using the theory of generalized foliated quadruples, we can bypass the failure of Bertini type theorems of foliations straightforwardly.
\end{sce}

\begin{sce}[Minimal model program on K\"ahler varieties]\label{sce: kahler}
It is well-known that foliations, especially algebraically integrable foliations, has a tight connection with the minimal model program for K\"ahler varieties. As we have mentioned above, Das and Ou essentially use the structure of algebraically integrable foliations to prove the abundance conjecture for K\"ahler manifolds in dimension $3$ \cite{DO23a,DO23b}. There is no doubt that foliations are expected to be useful in the study of K\"ahler minimal model program in the future.

On the other hand, generalized pairs is also known to have a tight connection with the minimal model program for K\"ahler varieties \cite{DH23,DHY23}. The key reason is due to the minimal model program with scaling of K\"ahler classes. K\"ahler classes cannot be considered as divisors, but by considering K\"ahler classes as $\bb$-nef classes and move it to the nef part of the generalized pair, we can formally define the minimal model program with scaling of K\"ahler classes.

In summary, the study of K\"ahler varieties seems to be a natural place for foliations and generalized pairs to get mixed together. We therefore can expect the structure of generalized foliated quadruples, particularly the algebraically integrable ones, to play a crucial role in the study of K\"ahler varieties in the future. 
\end{sce}

\begin{sce}[Cone theorem and semi-ampleness theorem]
In \cite[Theorem 3.9]{ACSS21}, a version of the cone theorem for algebraically integrable foliated triples $(X,\Ff,B)$ is proved. When $(X,\Ff,B)$ is lc, the main part of the cone theorem, i.e. the formula
$$\overline{NE}(X)=\overline{NE}(X)_{K_{\Ff}+B\geq 0}+\sum R_j$$
was proved in  \cite[Theorem 3.9]{ACSS21}. However, \cite[Theorem 3.9]{ACSS21} did not prove the countableness of $R_j$ nor the finiteness of $R_j$ when polarizing $(X,\Ff,B)$ with an ample divisor $A$. That is, the formula
$$\overline{NE}(X)=\overline{NE}(X)_{K_{\Ff}+B+A\geq 0}+\sum_{\text{finite}} R_j$$
is missing. One key reason for this seems to be the issue that $(X,\Ff,B+A)$ may no longer be lc, and this is, again, due to the failure of the Bertini type theorems for foliations. However, if we consider $(X,\Ff,B,\bar A)$ instead of $(X,\Ff,B+A)$, then $(X,\Ff,B,\bar A)$ becomes an lc generalized pair and we have immediately have more flexibility.

Similar issues appear when we consider the semi-ampleness theorem for foliations. For usual pairs, the semi-ampleness theorem is usually formulated in the following way:
$$(X,B) \text{ lc}, A\text{ ample}, K_X+B+A\text{ nef}\Rightarrow K_X+B+A\text{ semi-ample}.$$
However, for foliations, the semi-ampleness theorem is usually formulated in the following way: (under suitable conditions)
$$(X,\Ff,B+A) \text{ lc}, B\geq 0, A\geq 0\text{ ample}, K_{\Ff}+B+A\text{ nef}\Rightarrow K_{\Ff}+B+A\text{ semi-ample}.$$
This is again due to the failure of Bertini-type theorems. Nevertheless, with the new concept of generalized foliated quadruples, these arguments can now be strengthened back to:  (under suitable conditions)
$$(X,\Ff,B) \text{ lc}, A\text{ ample}, K_{\Ff}+B+A\text{ nef}\Rightarrow K_{\Ff}+B+A\text{ semi-ample}.$$
\end{sce}

\begin{sce}[Canonical bundle formula for generalized pairs]
The final scenario where the structure of generalized foliated quadruples plays a vital role is in obtaining the canonical bundle formula for generalized pairs. To establish this formula for lc-trivial fibrations in cases involving either non-NQC generalized pairs or NQC generalized pairs with potentially negative coefficients, we cannot rely on the structure of mixed Hodge structure (see Remark \ref{rem: lc trivial fibration definition}). Filipazzi's approach \cite{Fil19,Fil20} is also unsuitable due to its requirements for $\Qq$-coefficients and its inability to handle negative coefficients. Therefore, the only viable approach to achieve such a canonical bundle formula is by utilizing the theory of foliations as in  \cite{ACSS21}. Now, since we are dealing with generalized pairs in this context, the introduction of generalized foliated quadruples becomes necessary. For more details, we refer the reader to the proof of Theorem \ref{thm: cbf gpair nonnqc}.
\end{sce}

\part{Preliminaries}\label{part:prelim}

\section{Basic definitions}\label{sec: preliminaries}

 Throughout the paper, we work primarily with normal quasi-projective varieties to ensure consistency with the references. However, most results should also hold for normal varieties that are not necessarily quasi-projective. Similarly, most results in our paper should hold for any algebraically closed field of characteristic zero. We will adopt the standard notations and definitions in \cite{KM98, BCHM10} and use them freely. For foliations, we generally follow the notations and definitions in \cite{CS20,ACSS21,CS21}, but there may be minor differences. For generalized pairs, we follow the notations and definitions in \cite{HL21a}.

\subsection{Special notations}

\begin{nota}
    In this paper, $\mathbb N$ stands for the set of non-negative integers and $\mathbb N^+$ stands for the set of positive integers. 
\end{nota}

\begin{nota}
In this paper, the notation ``$/$" is always considered a shorthand for ``over". For example, ``$/Z$" means ``over $Z$".
\end{nota}

\begin{nota}\label{nota: general r divisor}
    Let $X\rightarrow U$ be a projective morphism from a normal variety to a variety, and let $A$ be a semi-ample$/U$ $\Rr$-divisor on $X$. An $\Rr$-divisor $H$ on $X$ is said to be \emph{general} in $|A/U|_{\mathbb R}$ if there exist base-point-free$/U$ divisors $A_1,\dots,A_n$ and real numbers $r_1,\dots,r_n\in (0,1)$ such that $A=\sum_{i=1}^nr_iA_i$ and $H=\sum_{i=1}^nr_iH_i$, where $H_i\in |A_i/U|$ are general elements. A \emph{general ample$/U$ $\Rr$-divisor} on $X$ is an ample$/U$ $\Rr$-divisor $D$ on $X$ such that $D$ is general in $|D/U|_{\mathbb R}$.
\end{nota}

\begin{nota}
    Let $\Ii$ be a set of real numbers, $X$ a normal variety, and $B$ an $\Rr$-divisor on $X$. We write $B\in\Ii$ if the coefficients of $B$ belong to $\Ii$.
\end{nota}

\begin{nota}
    Let $\pi: X\rightarrow U$ be a projective morphism between varieties and $D$ an $\Rr$-divisor on $X$. We denote by $\kappa_{\sigma}(X/Z,D)$ (resp. $\kappa_{\iota}(X/Z,D)$, $\kappa(X/Z,D)$) the relative numerical dimension (resp. relative invariant Iitaka dimension, relative Iitaka dimension) of $D$ over $Z$. When $Z=\{pt\}$, we may drop $X/Z$ and use the notation $\kappa_{\sigma}(D)$ (resp. $\kappa_{\iota}(D)$, $\kappa(D)$) instead. We refer the reader to \cite[Section 2]{HH20} for the formal definitions and basic properties of $\kappa_{\sigma}(X/Z,D)$, $\kappa(X/Z,D)$, and $\kappa(X/Z,D)$.
\end{nota}

\begin{defn}[Log big]\label{defn: log big}
Let $(X,B,\Mm)/U$ be a g-pair and $D$ an $\Rr$-Cartier $\Rr$-divisor on $X$. We say that $D$ is \emph{log big$/U$ with respect to $(X,B,\Mm)$} if $D|_V$ is big$/U$ for any lc center $V$ of $(X,B,\Mm)$. In particular, $D$ is big$/U$.
\end{defn}

\begin{defn}[{cf. \cite[Definition 5.3]{Amb03}, \cite[Definition 6.7.2]{Fuj11}}]\label{defn: basics of cone theorem}
Let $(X,\Delta)$ be a (not necessarily lc) pair and $\pi: X\rightarrow U$ a projective morphism. Let $F$ be an extremal face of $\overline{NE}(X/U)$.
\begin{enumerate}
    \item A \emph{supporting function} of $F$ is a $\pi$-nef $\Rr$-divisor $H$ such that $F=\overline{NE}(X/U)\cap H^{\bot}$. If $H$ is a $\Qq$-divisor, we say that $H$ is a \emph{rational supporting function}. Since $F$ is an extremal face of $\overline{NE}(X/U)$, $F$ always has a supporting function.
    \item We say that $F$ is \emph{rational} if $F$ has a rational supporting function.
    \item For any $\Rr$-Cartier $\Rr$-divisor $D$ on $X$, we say that $F$ is $D$-\emph{negative} if $$F\cap\overline{NE}(X/U)_{D\geq 0}=\{0\}.$$
    \item We say that $F$ is \emph{relatively ample at infinity with respect to} $(X,\Delta)$  if $$F\cap\overline{NE}(X/U)_{\Nlc(X,\Delta)}=\{0\}.$$ Equivalently, $H|_{\Nlc(X,\Delta)}$ is $\pi|_{\Nlc(X,\Delta)}$-ample for any supporting function $H$ of $F$.
    \item We say that $F$ is \emph{contractible at infinity with respect to} $(X,\Delta)$ if $F$ has a rational supporting function $H$ and $H|_{\Nlc(X,\Delta)}$ is $\pi|_{\Nlc(X,\Delta)}$-semi-ample.
\end{enumerate}
\end{defn}

\begin{deflem}
Let $K$ be a convex cone containing no lines. A ray $R$ of $K$ is called \emph{exposed} if there is a hyperplane meeting $K$ exactly along $R$. In particular, any exposed ray of $K$ is extremal in $K$. If $K$ does not contain any lines, then $K$ is the closure of the subcone of $K$ spanned by exposed rays \cite[Corollary 18.7.1]{Roc70}.

Let $\pi: X\rightarrow U$ be a projective morphism from a normal quasi-projective variety to a variety. By definition, an extremal ray in $\overline{NE}(X/U)$ is exposed if and only if it has a supporting function (that is not necessarily rational). Moreover, for any subcone $V$ of $\overline{NE}(X/U)$, we have
$$\overline{NE}(X/U)=\overline{V+\sum R_i}$$
where $R_i$ are exposed rays that are not contained in $V$.
\end{deflem}

\subsection{Sets}

\begin{defn}\label{defn: DCC and ACC}
Let $\Ii\subset\Rr$ be a set. We say that $\Ii$ satisfies the \emph{descending chain condition} (DCC) if any decreasing sequence in $\Ii$ stabilizes, and $\Ii$ satisfies the \emph{ascending chain condition} (ACC) if any increasing sequence in $\Ii$ stabilizes. 
\end{defn}

\begin{defn}\label{defn: derived set}
    Let $\Ii\subset [0,+\infty)$ be a set. We define 
    $$\Ii_+:=\{0\}\cup\left\{\sum_{i=1}^l\gamma_i\bigg| \gamma_i\in\Ii,l\in\mathbb N^+\right\}\text{ and }D(\Ii):=\left\{\frac{m-1+\gamma}{m}\bigg|\gamma\in\Ii_+,m\in\mathbb N^+\right\}.$$
\end{defn}

\subsection{Foliations}

\begin{defn}[Foliations, {cf. \cite[Section 2.1]{CS21}}]\label{defn: foliation}
Let $X$ be a normal variety. A \emph{foliation} on $X$ is a coherent sheaf $\Ff\subset T_X$ such that
\begin{enumerate}
    \item $\Ff$ is saturated in $T_X$, i.e., $T_X/\Ff$ is torsion-free, and
    \item $\Ff$ is closed under the Lie bracket.
\end{enumerate}
The \emph{rank} of the foliation $\Ff$ is the rank of $\Ff$ as a sheaf and is denoted by $\rk\Ff$. The \emph{corank} of $\Ff$ is $\dim X-\rk\Ff$. The \emph{canonical divisor} of $\Ff$ is a divisor $K_\Ff$ such that $\mathcal{O}_X(-K_{\mathcal{F}})\cong\mathrm{det}(\Ff)$. We define $N_{\Ff}:=(T_X/\Ff)^{\vee\vee}$ and $N_{\Ff}^*:=N_{\Ff}^{\vee}$.

If $\Ff=0$, then we say that $\Ff$ is a \emph{foliation by points}.
\end{defn}

\begin{defn}[Singular locus]
     Let $X$ be a normal variety and $\Ff$ a rank $r$ foliation on $X$. We can associate to $\Ff$ a morphism $$\phi: \Omega_X^{[r]}\to \mathcal{O}_X(K_{\Ff})$$ defined by taking the double dual of the $r$-wedge product of the map $\Omega^1_X\to \Ff^*$, induced by the inclusion $\Ff\to T_X$. This yields a map $$\phi': (\Omega_X^{[r]}\otimes\mathcal{O}_X(-K_{\Ff}))^{\vee\vee}\to \mathcal{O}_X$$ and we define the singular locus, denoted by $\Sing \Ff$, to be the cosupport of the image of $\phi'$.
\end{defn}

\begin{defn}[Pullbacks and pushforwards, {cf. \cite[3.1]{ACSS21}}]\label{defn: pullback}
Let $X$ be a normal variety, $\Ff$ a foliation on $X$, $f: Y\dashrightarrow X$ a dominant map, and $g: X\dashrightarrow X'$ a birational map. We denote by $f^{-1}\Ff$ the \emph{pullback} of $\Ff$ on $Y$ as constructed in \cite[3.2]{Dru21}. We also say that $f^{-1}\Ff$ is the \emph{induced foliation} of $\Ff$ on $Y$. If $\Ff=0$, then we say $f^{-1}\Ff$ is induced by $f$. In this case, we say $f^{-1}\Ff$ is \emph{algebraically integrable}.

We define the \emph{pushforward} of $\Ff$ on $X'$ as $(g^{-1})^{-1}\Ff$ and denote it by $g_*\Ff$.
\end{defn}

\begin{defn}[Invariant subvarieties, {cf. \cite[3.1]{ACSS21}}]\label{defn: f-invariant}
Let $X$ be a normal variety, $\Ff$ a foliation on $X$, and $S\subset X$ a subvariety. We say that $S$ is \emph{$\Ff$-invariant} if and only if for any open subset $U\subset X$ and any section $\partial\in H^0(U,\Ff)$, we have $$\partial(\mathcal{I}_{S\cap U})\subset \mathcal{I}_{S\cap U}$$ 
where $\mathcal{I}_{S\cap U}$ is the ideal sheaf of $S\cap U$. Note that if $\Ff$ is the foliation induced by a dominant map $f:X\dashrightarrow Z$, then a divisor $D$ is $\Ff$-invariant if and only if $D$ is vertical with respect to $f$.
\end{defn}

\begin{defn}[Special divisors on foliations, cf. {\cite[Definition 2.2]{CS21}}]\label{defn: special divisors on foliations}
Let $X$ be a normal variety and $\Ff$ a foliation on $X$. For any prime divisor $C$ on $X$, we define $\epsilon_{\Ff}(C):=1$ if $C$ is not $\Ff$-invariant, and  $\epsilon_{\Ff}(C):=0$ if $C$ is $\Ff$-invariant. If $\Ff$ is clear from the context, we may use $\epsilon(C)$ instead of $\epsilon_{\Ff}(C)$. For any $\Rr$-divisor $D$ on $X$, we define $$D^{\Ff}:=\sum_{C\mid C\text{ is a component of }D}\epsilon_{\Ff}(C)C.$$
Let $E$ be a prime divisor over $X$ and $f: Y\rightarrow X$ a projective birational morphism such that $E$ is on $Y$. We define $\epsilon_{\Ff}(E):=\epsilon_{f^{-1}\Ff}(E)$. It is clear that $\epsilon_{\Ff}(E)$ is independent of the choice of $f$.
\end{defn}

\subsection{Polarized foliations}

\begin{defn}[$\bb$-divisors]\label{defn: b divisors} Let $X$ be a normal quasi-projective variety. We call $Y$ a \emph{birational model} over $X$ if there exists a projective birational morphism $Y\to X$. 

Let $X\dashrightarrow X'$ be a birational map. For any valuation $\nu$ over $X$, we define $\nu_{X'}$ to be the center of $\nu$ on $X'$. A \emph{$\bb$-divisor} $\Dd$ on $X$ is a formal sum $\Dd=\sum_{\nu} r_{\nu}\nu$ where $\nu$ are valuations over $X$ and $r_{\nu}\in\mathbb R$, such that $\nu_X$ is not a divisor except for finitely many $\nu$. If in addition $r_{\nu}\in\Qq$ for every $\nu$, then $\Dd$ is called a \emph{$\Qq$-$\bb$-divisor}. The \emph{trace} of $\Dd$ on $X'$ is the $\Rr$-divisor
$$\Dd_{X'}:=\sum_{\nu_{X'}\text{ is a divisor}}r_\nu\nu_{X'}.$$
If $\Dd_{X'}$ is $\Rr$-Cartier and $\Dd_{Y}$ is the pullback of $\Dd_{X'}$ on $Y$ for any birational model $Y$ over $X'$, we say that $\Dd$ \emph{descends} to $X'$ and $\Dd$ is the \emph{closure} of $\Dd_{X'}$, and write $\Dd=\overline{\Dd_{X'}}$. 

Let $X\rightarrow U$ be a projective morphism and assume that $\Dd$ is a $\bb$-divisor on $X$ such that $\Dd$ descends to some birational model $Y$ over $X$. If $\Dd_Y$ is nef$/U$ (resp. base-point-free$/U$, semi-ample$/U$), then we say that $\Dd$ is \emph{nef}$/U$ (resp. \emph{base-point-free}$/U$, \emph{semi-ample}$/U$). If $\Dd_Y$ is a Cartier divisor, then we say that $\Dd$ is \emph{$\bb$-Cartier}. If $\Dd_Y$ is a $\Qq$-Cartier $\Qq$-divisor, then we say that $\Dd$ is \emph{$\Qq$-$\bb$-Cartier}. If $\Dd$ can be written as an $\Rr_{\geq 0}$-linear combination of nef$/U$ $\bb$-Cartier $\bb$-divisors, then we say that $\Dd$ is \emph{NQC}$/U$.

Let $X\rightarrow U$ be a projective morphism and assume that $\Dd$ and $\Dd'$ are two $\bb$-divisors over $X$. We write $\Dd\sim_{\mathbb R,U}\Dd'$ (resp. $\Dd\sim_{\mathbb Q,U}\Dd',\Dd\equiv_{\mathbb Q,U}\Dd'$) if for any birational model $Y$ of $X$, $\Dd_Y\sim_{\mathbb R,U}\Dd'_Y$ (resp. $\Dd_Y\sim_{\mathbb Q,U}\Dd'_Y,\Dd_Y\equiv_{\mathbb Q,U}\Dd_Y'$). 

We let $\bm{0}$ be the $\bb$-divisor $\bar{0}$.
\end{defn}

\begin{defn}\label{defn: restriction b divisor}
We will use two types of restrictions of $\bb$-divisors in this paper. Let $X$ be a normal variety and $\Dd$ a $\bb$-divisor on $X$.
\begin{enumerate}
    \item Let $V$ be a non-empty subset of $X$.  We define the \emph{restricted $\bb$-divisor} of $\Dd$ on $V$, which is denoted by $\Dd|_{V}$, in the following way. 

    For any birational morphism $\pi: W\to V$, there exists a birational morphism $\pi': Y\rightarrow X$ such that $W\subset Y$ and $\pi'|_W=\pi$. We let $(\Dd|_V)_{W}=(\Dd_Y)|_W$. It is easy to see that this definition is independent of the choice of $Y$ and defines a $\bb$-divisor.
    \item Suppose that $\Dd$ descends to a birational model of $X$. Let $S$ be a prime divisor on $X$ and $\nu: S^\nu\rightarrow S$ the normalization of $S$. The \emph{restricted $\bb$-divisor} of $\Dd$ on $S^\nu$, which is denoted by $\Dd|_{S^\nu}$, is defined in the following way. 

     Let $f: Y\rightarrow X$ be a log resolution of $(X,S)$ such that $\Dd$ descends to $Y$. Let $S_Y:=f^{-1}_*S$. Then there exists an induced birational morphism $f_S: S_Y\rightarrow S^\nu$ such that $\nu\circ f_S=f|_{S_Y}$.
     We define $$\Dd|_{S^\nu}:=\overline{\Dd_Y|_{S_Y}}.$$
     It is clear that $\Dd|_{S^\nu}$ is well-defined and is independent of the choice of $Y$.
\end{enumerate}
\end{defn}

\begin{defn}[Generalized foliated quadruples]\label{defn: gfq}
A \emph{generalized foliated sub-quadruple} (\emph{sub-gfq} for short) $(X,\Ff,B,\Mm)/U$ consists of a normal quasi-projective variety $X$, a foliation $\Ff$ on $X$, an $\Rr$-divisor $B$ on $X$, a projective morphism $X\rightarrow U$, and a nef$/U$ $\bb$-divisor $\Mm$ over $X$, such that $K_{\Ff}+B+\Mm_X$ is $\mathbb R$-Cartier. If $\Mm$ is NQC$/U$, then we say that $(X,\Ff,B,\Mm)/U$ is \emph{NQC}. If $B\geq 0$, then we say that $(X,\Ff,B,\Mm)/U$ is a \emph{generalized foliated quadruple} (\emph{gfq} for short). If $U=\{pt\}$, we usually drop $U$ and say that $(X,\Ff,B,\Mm)$ is \emph{projective}. 

Let $(X,\Ff,B,\Mm)/U$ be a (sub-)gfq. If $\Mm=\bm{0}$, then we may denote $(X,\Ff,B,\Mm)/U$ by $(X,\Ff,B)/U$ or $(X,\Ff,B)$, and say that $(X,\Ff,B)$ is a \emph{foliated (sub-)triple} (\emph{f-(sub-)triple} for short). If $\Ff=T_X$, then we may denote $(X,\Ff,B,\Mm)/U$ by $(X,B,\Mm)/U$, and say that $(X,B,\Mm)/U$ is a \emph{generalized (sub-)pair} (\emph{g-(sub-)pair} for short). If $\Mm=\bm{0}$ and $\Ff=T_X$, then we may denote $(X,\Ff,B,\Mm)/U$ by $(X,B)/U$ or $(X,B)$, and say that $(X,B)$ is a \emph{(sub-)pair}. 

A (sub-)gfq (resp. f-(sub-)triple, f-(sub-)pair, g-(sub-)pair, (sub-)pair) $(X,\Ff,B,\Mm)/U$ (resp. $(X,\Ff,B)/U$,$(X,B,\Mm)/U$, $(X,B)/U$) is called a \emph{$\mathbb Q$-(sub-)gfq} (resp. \emph{$\mathbb Q$-f-(sub-)triple, $\mathbb Q$-g-(sub-)pair, $\mathbb Q$-(sub-)pair}) if $B$ is a $\mathbb Q$-divisor and $\Mm$ is a $\mathbb Q$-$\bb$-divisor.
\end{defn}

\begin{nota}
In the previous definition, if $U$ is not important, we may also drop $U$. This usually happens when we emphasize the structures of $(X,\Ff,B,\Mm)$ which are independent of the choice of $U$, such as the singularities of $(X,\Ff,B,\Mm)$. In addition, if $B=0$, we may drop $B$.
\end{nota}

\begin{defn}[Singularities of gfqs]\label{defn: gfq singularity}
Let $(X,\Ff,B,\Mm)$ be a (sub-)gfq. For any prime divisor $E$ over $X$, let $f: Y\rightarrow X$ be a birational morphism such that $E$ is on $Y$, and suppose that
$$K_{\Ff_Y}+B_Y+\Mm_Y:=f^*(K_\Ff+B+\Mm_X)$$
where $\Ff_Y:=f^{-1}\Ff$. We define $a(E,\Ff,B,\Mm):=-\mult_EB_Y$ to be the \emph{discrepancy} of $E$ with respect to $(X,\Ff,B,\Mm)$. It is clear that $a(E,\Ff,B,\Mm)$ is independent of the choice of $Y$. If $\Mm=\bm{0}$, we let $a(E,\Ff,B):=a(E,\Ff,B,\Mm)$. If $\Ff=T_X$, we let $a(E,X,B,\Mm):=a(E,\Ff,B,\Mm)$. If $\Mm=\bm{0}$ and $\Ff=T_X$, we let $a(E,X,B):=a(E,\Ff,B,\Mm)$.

We say that $(X,\Ff,B,\Mm)$ is \emph{(sub-)lc} (resp. \emph{(sub-)klt}) if $a(E,\Ff,B,\Mm)\geq -\epsilon_{\Ff}(E)$ (resp. $>-\epsilon_{\Ff}(E)$) for any prime divisor $E$ over $X$. We say that $(X,\Ff,B,\Mm)$ is \emph{(sub-)canonical} (resp. \emph{(sub-)terminal}) if $a(E,\Ff,B,\Mm)\geq 0$ (resp. $>0$) for any prime divisor $E$ that is exceptional over $X$. An \emph{lc place} of $(X,\Ff,B,\Mm)$ is a prime divisor $E$ over $X$ such that $a(E,\Ff,B,\Mm)=-\epsilon_{\Ff}(E)$. An \emph{lc center} of $(X,\Ff,B,\Mm)$ is a subvariety $W$ of $X$ such that either $W$ is the center of an lc place of $(X,\Ff,B,\Mm)$ on $X$, or $W=X$. A \emph{non-trivial lc center} of $(X,\Ff,B,\Mm)$ is an lc center of $(X,\Ff,B,\Mm)$ that is not $X$. A \emph{non-lc place} of $(X,\Ff,B,\Mm)$ is a prime divisor $E$ over $X$ such that $a(E,\Ff,B,\Mm)<-\epsilon_{\Ff}(E)$. A \emph{non-lc center} of $(X,\Ff,B,\Mm)$ is the center of a non-lc place of $(X,\Ff,B,\Mm)$ on $X$. The union of all non-lc centers of $(X,\Ff,B,\Mm)$ is called the \emph{non-lc locus} of $(X,\Ff,B,\Mm)$ and is denoted by $\Nlc(X,\Ff,B,\Mm)$. The union of all non-lc centers and non-trivial lc centers of $(X,\Ff,B,\Mm)$ is called the \emph{non-klt locus} of $(X,\Ff,B,\Mm)$ and is denoted by $\Nklt(X,\Ff,B,\Mm)$.
\end{defn}

\begin{defn}\label{defn: lct}
  Let $(X,\Ff,B,\Mm)$ be a sub-gfq, $D\geq 0$ an $\Rr$-divisor on $X$ and $\Nn$ a nef$/X$ $\bb$-divisor, such that $D+\Nn_X$ is $\Rr$-Cartier. The \emph{lc threshold} (\emph{lct} for short) of $(D,\Nn)$ with respect to $(X,\Ff,B,\Mm)$ is defined as
  $$\lct(X,\Ff,B,\Mm;D,\Nn):=\sup\{t\mid (X,\Ff,B+tD,\Mm+t\Nn)\text{ is sub-lc}\}.$$
  If $\Nn=0$, then we may drop $\Nn$ and denote $\lct(X,\Ff,B,\Mm;D,\Nn)$ by $\lct(X,\Ff,B,\Mm;D)$.
\end{defn}

\begin{defn}[Models, I]\label{defn: models I}
Let $(X,\Ff,B,\Mm)/U$ be an lc gfq, $\phi: X\dashrightarrow X'$ a birational map over $U$, $E:=\Exc(\phi^{-1})$ the reduced $\phi^{-1}$-exceptional divisor, $\Ff':=\phi_*\Ff$, and $B':=\phi_*B+E^{\Ff'}$.
\begin{enumerate}
    \item $(X',\Ff',B',\Mm)/U$ is called a \emph{log birational model} of $(X,\Ff,B,\Mm)/U$. 
    \item $(X',\Ff',B',\Mm)/U$ is called a \emph{weak lc model} of $(X,\Ff,B,\Mm)/U$ if 
\begin{enumerate}
\item $(X',\Ff',B',\Mm)/U$ is a log birational model of $(X,\Ff,B,\Mm)/U$, 
    \item $K_{\Ff'}+B'+\Mm_{X'}$ is nef$/U$, and
    \item for any prime divisor $D$ on $X$ which is exceptional over $X'$, $$a(D,\Ff,B,\Mm)\leq a(D,\Ff',B',\Mm).$$
\end{enumerate}
\item $(X',\Ff',B',\Mm)/U$ is called a \emph{semi-good minimal model} of $(X,\Ff,B,\Mm)/U$ if
\begin{enumerate}
        \item $(X',\Ff',B,\Mm)/U$ is a weak lc model of $(X,\Ff,B,\Mm)/U$, and
        \item $K_{\Ff'}+B'+\Mm_{X'}$ is semi-ample$/U$.
\end{enumerate}
\item  Suppose that there exists a contraction$/U$ $X'\rightarrow Z$. $(X',\Ff',B',\Mm)\rightarrow Z$ is called a \emph{Mori fiber space} of $(X,\Ff,B,\Mm)/U$ if
\begin{enumerate}
    \item  $(X',\Ff',B',\Mm)/U$ is a log birational model of $(X,\Ff,B,\Mm)/U$,
    \item $X'$ is $\Qq$-factorial,
    \item $X'\rightarrow Z$ is a $(K_{\Ff'}+B'+\Mm_{X'})$-Mori fiber space$/U$,
    \item for any prime divisor $D$ on $X$ which is exceptional over $X'$, $$a(D,\Ff,B,\Mm)<a(D,\Ff',B',\Mm).$$
\end{enumerate}
\end{enumerate}
We shall not define ``good minimal models" until Definition \ref{defn: models ii}.
\end{defn}

\begin{nota}
    Let $(X_0,\Ff_0,B_0,\Mm)/U$ be a gfq. When we say that the following
\begin{center}$\xymatrix{
(X_0,\Ff_0,B_0,\Mm)\ar@{-->}[r]^{f_0} & (X_1,\Ff_1,B_1,\Mm)\ar@{-->}[r]^{\ \ \ \ \ \ \ \ f_1} & \dots\ar@{-->}[r] & (X_n,\Ff_n,B_n,\Mm)\ar@{-->}[r]^{\ \ \ \ \ \ \ \ \ f_n} & \dots 
}$
\end{center}
is a (possibly infinite) sequence of steps of a $(K_{\Ff_0}+B_0+\Mm_{X_0})$-MMP$/U$, we mean the following: for any $i$, $f_i: X_{i}\dashrightarrow X_{i+1}$ is a step of a $(K_{\Ff_i}+B_i+\Mm_{X_i})$-MMP$/U$ that is not a Mori fiber space, $\Ff_{i+1}:=(f_i)_*\Ff_i$, and $B_{i+1}:=(f_i)_*B_i$.
\end{nota}

\begin{cons}[MMP with scaling]\label{cons: mmp with scaling}
    Let $(X,\Ff,B,\Mm)/U$ be an lc gfq. Let $D\geq 0$ be an $\Rr$-divisor on $X$ and $\Nn$ a nef$/U$ $\bb$-divisor on $X$ such that $D+\Nn_X$ is $\Rr$-Cartier and $K_{\Ff}+B+\Mm_X+t(D+\Nn_X)$ is nef$/U$ for some positive real number $t$. A step of a $(K_{\Ff}+B+\Mm_X)$-MMP$/U$ with scaling of $(D,\Nn)$ is defined as follows. Let
$$\lambda:=\inf\{s\geq 0\mid K_{\Ff}+B+sD+\Mm_X+s\Nn_X\text{ is nef}/U\}.$$
Assume that the following conditions hold:
\begin{itemize}
  \item There exists an extremal ray $R$ in $\overline{NE}(X/U)$ such that $(K_{\Ff}+B+\lambda D+\Mm_X+\lambda\Nn_X)\cdot C=0$ and $(D+\Nn_X)\cdot C>0$. In particular, $R$ is a $(K_{\Ff}+B+\Mm_X)$-negative extremal ray.
  \item The contraction associated to $R$ exists, and if it is a small contraction, the corresponding $(K_{\Ff}+B+\Mm_X)$-flip exists.
\end{itemize}
Then for any such $R$, we call the divisorial contraction or the Mori fiber space associated to $R$, or the $(K_{\Ff}+B+\Mm_X)$-flip associated to $R$, a step of a $(K_{\Ff}+B+\Mm_X)$-MMP$/U$ with scaling of $(D,\Nn)$. 

A sequence of steps of a $(K_{\Ff}+B+\Mm_X)$-MMP$/U$ with scaling of $(D,\Nn)$ is a sequence of steps of a $(K_{\Ff}+B+\Mm_X)$-MMP$/U$
\begin{center}$\xymatrix{
(X_0,\Ff_0,B_0,\Mm)\ar@{-->}[r]^{f_0} & (X_1,\Ff_1,B_1,\Mm)\ar@{-->}[r]^{\ \ \ \ \ \ \ \ f_1} & \dots\ar@{-->}[r] & (X_n,\Ff_n,B_n,\Mm)\ar@{-->}[r]^{\ \ \ \ \ \ \ \ \ f_n} & \dots 
}$
\end{center}
where $(X_0,\Ff_0,B_0,\Mm)=(X,\Ff,B,\Mm)$, and each $f_i$ is a step of a $(K_{\Ff_i}+B_i+\Mm_{X_i})$-MMP$/U$ with scaling of $(D_i,\Nn)$, where $D_i$ is the image of $D$ on $X_i$. 
$$\lambda_{i}:=\inf\{s\geq 0\mid K_{\Ff_{i}}+B_{i}+sD_{i}+\Mm_{X_{i}}+s\Nn_{X_{i}}\text{ is nef}/U\}$$
are called the \emph{scaling numbers} (of this MMP with scaling of $(D,\Nn)$), which are well-defined.

If $\Nn=\bm{0}$, a (sequence of) step(s) of a $(K_{\Ff}+B+\Mm_X)$-MMP$/U$ with scaling of $(D,\Nn)$ is called a (sequence of) step(s) of a $(K_{\Ff}+B+\Mm_X)$-MMP$/U$ with scaling of $D$. 
\end{cons}
We remark that Construction \ref{cons: mmp with scaling} does not require the condition that $(X,\Ff,B+D,\Mm+\Nn)$ is lc.

\begin{defn}
   Let $(X,\Ff,B,\Mm)$ and $(X',\Ff',B',\Mm')$ be two sub-gfqs.   
   We say that $(X,\Ff,B,\Mm)$ and $(X',\Ff',B',\Mm')$ are \emph{crepant} to each other if $\Mm=\Mm'$, and there exist two birational morphisms $p: W\rightarrow X$ and $q: W\rightarrow X'$ and a foliation $\Ff_W$ on $W$ such that $\Ff_W=p^{-1}\Ff=q^{-1}\Ff'$, $\Mm=\Mm'$, and
   $$p^*(K_{\Ff}+B+\Mm_X)=q^*(K_{\Ff'}+B'+\Mm'_{X'}).$$
\end{defn}

\section{Basic properties of generalized pairs}\label{sec: basic property gpair}

In this section, we present several results concerning the structure of generalized pairs. Although most of these results, or their analogous forms, have already been established in the existing literature, it is somewhat surprising to note that many fundamental results for non-NQC generalized pairs remain unaddressed, despite the significant progress in the field of generalized pairs in recent years. Additionally, there are relatively few references available on this subject. For clarity, for the reader's convenience, and to provide a solid reference for future work, we will present detailed proofs for all the results in this section.

\subsection{Dlt modification}

First, we recall the definition of dlt for generalized pairs.

\begin{defn}[Dlt, {\cite[Definition 2.3]{HL22}}]\label{defn: dlt}
Let $(X,B,\Mm)/U$ be an lc g-pair. We say that $(X,B,\Mm)$ is \emph{dlt} if there exists an open subset $V\subset X$ satisfying the following.
\begin{enumerate}
    \item $(V,B|_V)$ is log smooth. In particular, $B|_V$ is an snc Weil $\Qq$-divisor.
    \item $V$ contains the generic point of any lc center of $(X,B,\Mm)$.
    \item The generic point of any lc center of $(X,B,\Mm)$ is the generic point of an lc center of $(V,B|_V)$.
\end{enumerate}
If $(X,B,\Mm)$ is dlt and $\lfloor B\rfloor$ is normal, then we say that $(X,B,\Mm)$ is plt.
\end{defn}

The following lemma indicates that the definition of dlt in \cite[Definition 2.3]{HL22} is the same as the definition of dlt in \cite{Bir20,FS23}. 
\begin{lem}\label{lem: equi def dlt 1}
Let $(X,B,\Mm)/U$ be an lc g-pair. Then the following two conditions are equivalent
\begin{enumerate}
    \item $(X,B,\Mm)$ is dlt.
    \item For any lc center of $(X,B,\Mm)$ with generic point $\eta$, over a neighborhood of $\eta$, $(V,B|_V)$ is log smooth and $\Mm$ descends to $X$. 
\end{enumerate} 
\end{lem}
\begin{proof}
Since being dlt, the property in (2), and being log smooth are local properties, we may work over a neighborhood of a generic point $\eta$ of an lc center of $(X,B,\Mm)$. (2)$\Rightarrow$ (1) immediately becomes obvious, so we only need to prove (1)$\Rightarrow$ (2). 

By Definition \ref{defn: dlt}, there exists a neighborhood $V$ of $\eta$ such that $(V,B|_V)$ is log smooth and $\eta$ is an lc center of $(V,B|_V)$. Since $(V,B|_V)$ is log smooth, $\Mm_X|_V$ is $\Rr$-Cartier. We let $\Mm^V:=\Mm|_V$ be the restricted $\bb$-divisor of $\Mm$ on $V$, then $\Mm^V$ is nef$/X$ and $\Mm^V_V=\Mm_X|_V$. Suppose that $h: V'\rightarrow V$ is a resolution of $V$ such that $\Mm^V$ descends to $V'$ and there exists a prime divisor $E$ on $V'$ such that $\Center_{V}E=\bar\eta$ and $E$ is an lc place of $(V,B|_V)$. By the negativity lemma,
    $$\Mm^V_{V'}=h^*\Mm^V_V-F$$
    for some $F\geq 0$, such that either $F=0$ over $\bar\eta$ or $\Supp F=\Supp h^{-1}(\bar\eta)$. Since $(X,B,\Mm)$ is lc, $(V,B|_V,\Mm^V)$ is lc. Thus $F=0$ over $\bar\eta$. Possibly shrinking $V$, we may assume that $\Mm$ descends to $V$. The lemma follows.
\end{proof}

Lemma \ref{lem: equi def dlt 1} implies the following result

\begin{deflem}[Dlt modification, {\cite[Theorem 2.9]{FS23}}]\label{deflem: dlt model}
    Let $(X,B,\Mm)/U$ be a g-pair. Then there exists a birational morphism $f: Y\rightarrow X$ satisfying the following. Let $E_1,\dots,E_n$ be the prime $f$-exceptional divisors and $B_Y:=f^{-1}_*(B\wedge\Supp B)+\sum_{i=1}^nE_i$, then
    \begin{enumerate}
        \item $(Y,B_Y,\Mm)$ is $\Qq$-factorial dlt, and
        \item $a(E_i,X,B,\Mm)\leq 0$ for any $i$.
    \end{enumerate}
    In particular, if $(X,B,\Mm)$ is lc, then $a(E_i,X,B,\Mm)=0$ for any $i$, and $$K_Y+B_Y+\Mm_Y=f^*(K_X+B+\Mm_X).$$
    
    For any such $f$, we call $f$ a \emph{dlt modification} of $(X,B,\Mm)$, and say that $(Y,B_Y,\Mm)$ is a \emph{dlt model} of $(X,B,\Mm)$. 
\end{deflem}

We conjecture that dlt has another equivalent definition

\begin{conj}\label{conj: dlt has good resolution}
   Let $(X,B,\Mm)/U$ be an lc g-pair. Then $(X,B,\Mm)$ is dlt if and only if there exists a log resolution  $f: Y\rightarrow X$ of $(X,\Supp B)$ and an open subset $V\subset X$, such that $\Mm$ descends to $Y$, $f$ is an isomorphism over $V$, and $V$ contains the generic point of any lc center of $(X,B,\Mm)$. 
\end{conj}
When $(X,B,\Mm)/U$ is NQC, Conjecture \ref{conj: dlt has good resolution} was proven in \cite[Theorem 6.1]{Has22a}.

\subsection{Perturbation and MMP}

\begin{lem}\label{lem: gklt+ample terminate}
    Let $(X,B+A,\Mm)/U$ be a $\Qq$-factorial lc g-pair such that $X$ is klt, $A\geq 0$ is ample$/U$, and $B\geq 0$. Then any $(K_X+B+A+\Mm_X)$-MMP$/U$ with scaling of an ample$/U$ $\Rr$-divisor terminates with either a semi-good minimal model of $(X,B+A,\Mm)/U$ or a Mori fiber space of  $(X,B+A,\Mm)/U$.
\end{lem}
\begin{proof}
   By \cite[Lemma 3.4]{HL22}, there exists $0\leq\Delta\sim_{\mathbb R,U}B+A+\Mm_X$ such that $(X,\Delta)$ is klt. By \cite[Corollary 1.4.2]{BCHM10}, any $(K_X+\Delta)$-MMP$/U$ with scaling of an ample$/U$ $\Rr$-divisor terminates with either a Mori fiber space of $(X,\Delta)/U$ or a semi-good minimal model of $(X,\Delta)/U$. The lemma follows. 
\end{proof}

\begin{lem}\label{lem: scaling number go to 0}
Let $(X,B,\Mm)/U$ be a $\Qq$-factorial lc g-pair such that $X$ is klt, and $A\geq 0$ an ample$/U$ $\Rr$-divisor on $X$. Then we may run a $(K_X+B+\Mm_X)$-MMP$/U$ with scaling of $A$. Moreover, let $$(X,B,\Mm):=(X_1,B_1,\Mm)\dashrightarrow (X_2,B_2,\Mm)\dashrightarrow\dots\dashrightarrow (X_i,B_i,\Mm)\dashrightarrow\dots$$
be any $(K_X+B+\Mm_X)$-MMP$/U$ with scaling of $A$, and let $\lambda_i$ be the $i$-th scaling number of this MMP for each $i$, i.e.,
$$\lambda_i:=\inf\{t\mid t\geq 0, K_{X_i}+B_i+tA_i+\Mm_{X_i}\text{ is nef/}U\},$$
where $A_i$ is the strict transform of $A$ on $X_i$ for each $i$. Then $\lambda_i\geq\lambda_{i+1}$ for each $i$, and one of the following holds
\begin{enumerate}
    \item This MMP terminates after finitely many steps.
    \item $\lim_{i\rightarrow +\infty}\lambda_i=0$.
\end{enumerate}
\end{lem}
\begin{proof}
Possibly rescaling $A$, we may assume that $K_X+B+A+\Mm_X$ is nef$/U$. We first prove that we may run this MMP by induction on $i$. Let $\lambda_0:=1$ and suppose that there is already a sequence of steps of a $(K_X+B+\Mm_X)$-MMP$/U$ with scaling of $A$
$$(X,B,\Mm):=(X_1,B_1,\Mm)\dashrightarrow (X_2,B_2,\Mm)\dashrightarrow\dots\dashrightarrow (X_k,B_k,\Mm)$$
for some $k\ge1$, such that $\lambda_i\geq\lambda_{i+1}$ for any $i\leq k-2$. If $K_{X_k}+B_k+\Mm_{X_k}$ is nef$/U$, then we are done, so we may assume that $K_{X_k}+B_k+\Mm_{X_k}$ is not nef$/U$. Since nef$/U$ is a closed condition, $\lambda_k>0$. Since $K_X+B+\lambda_0A+\Mm_X$ is nef$/U$ and $K_{X_{k-1}}+B_{k-1}+\lambda_{k-1}A_{k-1}+\Mm_{X_{k-1}}$ is nef$/U$ when $k\geq 2$, $K_{X_{k}}+B_{k}+\lambda_{k-1}A_{k}+\Mm_{X_{k}}$ is nef, hence $\lambda_{k-1}\geq\lambda_{k}$.

By \cite[Lemma 3.4]{HL22}, there exists a klt pair $(X,\Delta)$ such that 
$$K_X+\Delta\sim_{\mathbb R,U}K_X+B+\Mm_X+\frac{\lambda_k}{2}A.$$ 
Possibly replacing $A$, we may assume that $(X,\Delta+(1-\frac{\lambda_k}{2})A)$ is lc. Then we have an induced sequence of steps of a $(K_X+\Delta)$-MMP$/U$ with scaling of $(1-\frac{\lambda_k}{2})A$
$$(X,\Delta):=(X_1,\Delta_1)\dashrightarrow (X_2,\Delta_2)\dashrightarrow\dots\dashrightarrow (X_k,\Delta_k),$$
such that $K_{X_k}+\Delta_k$ is not nef. Let $(X_k,\Delta_k)\dashrightarrow (X_{k+1},\Delta_{k+1})$ be the next step of the $(K_X+\Delta)$-MMP$/U$ with scaling of $(1-\frac{\lambda_k}{2})A$. Then the induced birational map $X_k\dashrightarrow X_{k+1}$ is a step of a $(K_X+B+\Mm_X)$-MMP$/U$ with scaling of $A$. 

It remains to prove that if this MMP does not terminate, then $\lim_{i\rightarrow+\infty}\lambda_i=0$. Suppose that $\lambda:=\lim_{i\rightarrow+\infty}\lambda_i>0$. Then
\begin{align*}
\left(X,B+\frac{\lambda}{2}A,\Mm\right):=&\left(X_1,B_1+\frac{\lambda}{2}A_1,\Mm\right)\dashrightarrow\left(X_2,B_2+\frac{\lambda}{2}A_2,\Mm\right)\dashrightarrow\\
\dots\dashrightarrow&\left(X_i,B_i+\frac{\lambda}{2}A_i,\Mm\right)\dashrightarrow\dots
\end{align*}
is an infinite sequence of steps of a $(K_X+B+\frac{\lambda}{2}A+\Mm_X)$-MMP$/U$, which contradicts Lemma \ref{lem: gklt+ample terminate}.
\end{proof}

The following result seems to be missed in known literature.

\begin{prop}\label{prop: qfact nqc any scaling terminate}
    Let $(X,B+A,\Mm)/U$ be a $\Qq$-factorial NQC lc g-pair such that $A\geq 0$ is ample$/U$ and $B\geq 0$. Then any $(K_X+B+A+\Mm_X)$-MMP$/U$ with scaling of an ample$/U$ $\Rr$-divisor terminates with either a semi-good minimal model of $(X,B+A,\Mm)/U$ or a Mori fiber space of $(X,B+A,\Mm)/U$.
\end{prop}
\begin{proof}
By \cite[Lemma A.5]{LT22}, possibly replacing $A$ with a general element in $|A/U|_{\mathbb R}$, there exists an lc pair $(X,\Delta+\frac{1}{2}A)$ such that $0\leq \Delta\sim_{\mathbb R}B+\frac{1}{2}A+\Mm_X$. By \cite[Theorem 1.5]{HH20} and \cite[Theorem 1.9]{Bir12}, any $(K_X+\Delta+\frac{1}{2}A)$-MMP$/U$ with scaling of an ample$/U$ $\Rr$-divisor terminates. Thus any $(K_X+B+A+\Mm_X)$-MMP$/U$ with scaling of an ample$/U$ $\Rr$-divisor terminates. The remaining part of the proposition follows from \cite[Theorem 1.3]{LX23a} and \cite[Lemma 3.5(1)]{HL21a}.
\end{proof}

\begin{lem}\label{lem: movable num 0 is 0}
Let $X$ be a normal projective variety and $D$ a movable $\Rr$-Cartier $\Rr$-divisor on $X$ such that $\kappa_{\sigma}(D)=0$. Then $D\equiv 0$.
\end{lem}
\begin{proof}
We let $f: Y\rightarrow X$ be a resolution of $X$. Let $P_Y:=P(Y,f^*D)$ and $N_Y:=N(Y,f^*D)$ be the positive and negative parts of the Nakayama-Zariski decomposition of $f^*D$ respectively, and let $P:=P(X,D)$ and $N:=N(X,D)$ be the positive and negative parts of the Nakayama-Zariski decomposition of $D$ respectively. Since $D$ is movable, by \cite[Lemma 3.7(3)]{LX23a}, $N=0$. By \cite[Lemma 3.4(3)]{LX23a}, $f_*N_Y=N$, so $N_Y$ is exceptional$/X$. Since $\kappa_\sigma(f^*D)=\kappa_\sigma(D)=0$, by \cite[V 2.7 Proposition(8)]{Nak04}, $P_Y\equiv 0$. Thus $D=f_*D_Y=f_*(P_Y+N_Y)\equiv 0.$
\end{proof}

\begin{prop}
Let $(X,B,\Mm)$ be a projective $\Qq$-factorial lc g-pair such that $\kappa_{\sigma}(K_X+B+\Mm_X)=0$ and $X$ is klt. Let $A$ be an ample $\Rr$-divisor. Then we may run a $(K_X+B+\Mm_X)$-MMP with scaling of $A$, and any such MMP terminates with a model $(X',B',\Mm)$ of $(X,B,\Mm)$ such that $K_{X'}+B'+\Mm_{X'}\equiv 0$. Moreover, if $\kappa_{\iota}(K_X+B+\Mm_X)=0$, then $K_{X'}+B'+\Mm_{X'}\sim_{\mathbb R}0$.
\end{prop}
\begin{proof}
By Lemma \ref{lem: scaling number go to 0}, we may run a $(K_X+B+\Mm_X)$-MMP with scaling of $A$. Let
$$(X,B,\Mm):=(X_1,B_1,\Mm)\dashrightarrow (X_2,B_2,\Mm)\dashrightarrow\dots\dashrightarrow (X_i,B_i,\Mm)\dashrightarrow\dots$$
be any such MMP with scaling numbers $\lambda_i\ge\lambda_{i+1}$. If this MMP does not terminate, then $\lim_{i\rightarrow+\infty}\lambda_i=0$ by Lemma \ref{lem: scaling number go to 0}. There exists a positive integer $m$ such that $X_i\dashrightarrow X_{i+1}$ is a flip for any $i\geq m$. We may denote by $\phi_i: X_m\dashrightarrow X_i$ the induced birational contraction and $A_i$ the strict transform of $A$ on $X_i$ for any $i> m.$ Since $K_{X_i}+B_i+\lambda_iA_i+\Mm_{X_i}$ is nef for each $i$,
$$K_{X_m}+B_m+\Mm_{X_m}=\lim_{i\rightarrow+\infty}(\phi_i^{-1})_*(K_{X_i}+B_i+\lambda_iA_i+\Mm_{X_i})$$
is movable. Moreover, since $\kappa_{\sigma}(K_X+B+\Mm_X)=0$, $\kappa_{\sigma}(K_{X_m}+B_m+\Mm_{X_m})=0$. By Lemma \ref{lem: movable num 0 is 0}, $K_{X_m}+B_m+\Mm_{X_m}\equiv 0$, a contradiction. Thus this MMP terminates with a model $(X',B',\Mm)$ such that $K_{X'}+B'+\Mm_{X'}$ is nef and $\kappa_{\sigma}(K_{X'}+B'+\Mm_{X'})=0$. By Lemma \ref{lem: movable num 0 is 0} again, $K_{X'}+B'+\Mm_{X'}\equiv 0$. Moreover, if $\kappa_{\iota}(K_X+B+\Mm_X)=0$, then $\kappa_{\iota}(K_{X'}+B'+\Mm_{X'})=0$, hence $K_{X'}+B'+\Mm_{X'}\sim_{\mathbb R}0$.
\end{proof}

\subsection{Lc centers of generalized pairs}

 We will discuss the structure of lc centers of lc g-pairs in this section. For NQC generalized pairs, the structure of their lc centers is well-studied in \cite{LX23a} based on the connectedness principle established in \cite{Bir20,FS23} and the canonical bundle formula \cite{Fil20,HL21b,JLX22,FS23}, but little was known for the non-NQC case.

\begin{defn}[Adjunction for generalized pairs to divisors, cf. {\cite[Definition 4.7]{BZ16}}]\label{defn: adj to lc places}
Let $(X,B,\Mm)/U$ be a g-(sub-)pair and $S$ a component of $B^{=1}$. Let $S^\nu$ be the normalization of $S$. The g-(sub-)pair $(S^\nu,B_S,\Mm^S)/U$ induced by the adjunction
$$K_{S^\nu}+B_S+\Mm^S_S:=(K_X+B+\Mm_X)|_S$$
is given in the following way. Let $f: Y\rightarrow X$ be a log resolution of $(X,\Supp B)$ such that $\Mm$ descends to $Y$, $S_Y$ the strict transform of $S$ on $Y$, and 
$$K_Y+B_Y+\Mm_Y:=f^*(K_X+B+\Mm_X).$$ 
We define $\Mm^S:=\Mm|_{S^\nu}$ and $B_{S_Y}:=(B_Y-S_Y)|_{S_Y}$. We let $f|_{S_Y}: S_Y\rightarrow S^\nu$ be the induced birational morphism and define $B_S:=(f|_{S_Y})_*B_{S_Y}$.
\end{defn}

\begin{lem}[cf. {\cite[Lemma 3.18(2)]{LX23b}}]\label{lem: inversion of adjunction gdlt}
Let $(X,B,\Mm)/U$ be a dlt g-pair, $S$ a component of $\lfloor B\rfloor$, and $(S,B_S,\Mm^S)/U$ the g-pair induced by the adjunction 
$$K_S+B_S+\Mm^S_S:=(K_X+B+\Mm_X)|_S.$$
Then $(S,B_S,\Mm^S)$ is dlt. Moreover
\begin{enumerate}
    \item Any lc center of $(S,B_S,\Mm^S)$ is an lc center of $(X,B,\Mm)$.
    \item Any lc center of $(X,B,\Mm)$ that is contained in $S$ is an lc center of $(S,B_S,\Mm^S)$.
\end{enumerate}
\end{lem}
\begin{proof}
By \cite[Lemma 2.9]{HL22}, $(S,B_S,\Mm^S)$ is dlt.

(1) Let $f: \tilde X\rightarrow X$ be a log resolution of $(X,\Supp B)$ such that $\Mm$ descends to \(\tilde X\). Let $K_{\tilde X}+\tilde B+\Mm_{\tilde X}:=f^*(K_X+B+\Mm_X)$ and let $\tilde S$ be the strict transform of $S$ on \(\tilde X\), then
$f|_{\tilde S}$ is a log resolution of $(S,\Supp B_S)$ such that $\Mm^S$ descends to \(\tilde S\). We have
$$f|_{\tilde S}^*(K_S+B_S+\Mm^S_S)=K_{\tilde S}+B_{\tilde S}+\Mm^S_{\tilde S}:=(K_{\tilde X}+\tilde B+\Mm_{\tilde X})|_{\tilde S}.$$
Let $W_S$ be an lc center of $(S,B_S,\Mm^S)$. Then $W_S$ is the image of an lc center $W_{\tilde S}$ of $(\tilde S,B_{\tilde S},\Mm^S)$ in $S$. Since $(\tilde X,\tilde B)$ is log smooth and $\Mm$ descends to \(\tilde X\), $W_{\tilde S}$ is also an lc center of $(\tilde X,\tilde B,\Mm)$ which is contained in \(\tilde S\), so $W:=f(W_{\tilde S})$ is an lc center of $(X,B,\Mm)$ which is contained in $S$. It is clear that $W_S=W$ under the natural inclusion $S\rightarrow X$. This implies (1).

(2) Let $W$ be an lc center of $(X,B,\Mm)$ that is contained in $S$. Since $(X,B,\Mm)$ is dlt, by Lemma \ref{lem: equi def dlt 1}, possibly shrinking $X$ to a neighborhood of the generic point of $W$, we may assume that $(X,B)$ is log smooth and $\Mm$ descends to $X$. Thus $W$ is an lc center of $(X,B)$, $K_S+B_S=(K_X+B)|_S$, and $\Mm^S$ descends to $S$. Since $(X,B)$ is log smooth, $W$ is an lc center of $(S,B_S)$, hence an lc center of $(S,B_S,\Mm^S)$. This implies (2).
\end{proof}

\begin{deflem}
    Let $(X,B,\Mm)/U$ be a dlt g-pair and $V$ an lc center of $(X,B,\Mm)$ such that $\dim V\geq 1$. Then we may construct a dlt g-pair $(V,B_V,\Mm^V)/U$ on $V$ inductively in the following way. If $V=X$ then we let $(V,B_V,\Mm^V):=(X,B,\Mm)$. Otherwise, let $S$ be a codimension $1$ lc center of $(X,B,\Mm)$ such that $V\subset S$. By \cite[Lemma 2.9]{HL22}, there exists a dlt g-pair $(S,B_{S},\Mm^S)$ induced by adjunction
    $$K_{S}+B_{S}+\Mm^S_{S}=(K_X+B+\Mm_X)|_{S}.$$
    By Lemma \ref{lem: inversion of adjunction gdlt}, $V$ is an lc center of $(S,B_S,\Mm^S)$. By repeating this process and applying induction on dimension, we get a dlt g-pair $(V,B_V,\Mm^V)/U$ on $V$. $(V,B_V,\Mm^V)/U$ is called the dlt g-pair \emph{induced by repeatedly applying adjunction to codimension $1$ lc centers}
    $$K_V+B_V+\Mm^V_V:=(K_X+B+\Mm_X)|_V.$$
\end{deflem}

\begin{defn}
An \emph{lc crepant log structure} is of the form $f: (X,B,\Mm)\rightarrow Z$, where
\begin{enumerate}
    \item $(X,B,\Mm)/Z$ is an lc g-pair,
    \item $K_X+B+\Mm_X\sim_{\Rr,Z}0$, and
    \item $f$ is a contraction. In particular, $f_*\Oo_X=\Oo_Z$.
\end{enumerate}
In addition, if
\begin{enumerate}
    \item[(4)] $(X,B,\Mm)$ is dlt, 
\end{enumerate}
then we say that $f: (X,B,\Mm)\rightarrow Z$ is a \emph{dlt crepant log structure}. 

For any irreducible subvariety $W\subset Z$, we say that $W$ is an \emph{lc center} of an lc crepant log structure $f: (X,B,\Mm)\rightarrow Z$, if there exists an lc center $W_X$ of $(X,B,\Mm)$ such that $W=f(W_X)$. For any (not necessarily closed) point $z\in Z$, we say that $\bar z$ is an \emph{lc center} of $f: (X,B,\Mm)\rightarrow Z$ if $\bar z$ is an lc center of $f: (X,B,\Mm)\rightarrow Z$.
\end{defn}

\begin{rem} 
In Section \ref{sec: cbf} below we will introduce the concept of lc-trivial fibrations. We will see that an lc crepant log structure is an lc-trivial fibration $f: (X,B,\Mm)\rightarrow Z$ such that $B\geq 0$ (see Definition \ref{defn: lc trivial fibration gfq} below). We will also see that an lc center of an lc crepant log structure $f: (X,B,\Mm)\rightarrow Z$ is indeed an lc center of the induced g-pair $(Z,B_Z,\Mm^Z)$ via the canonical bundle formula (see Theorem \ref{thm: cbf gpair nonnqc} below).
\end{rem}

\begin{defn}[Standard $\mathbb P^1$-link, cf. {\cite[Definition 2.21]{FS23}}]\label{defn: standard p1 link}
Let $X\rightarrow U$ be a projective morphism from a normal quasi-projective variety to a variety. A \emph{standard $\mathbb P^1$-link}$/U$ $f: (X,B,\Mm)\rightarrow T$ is an lc g-pair $(X,B,\Mm)/U$ with a projective morphism $f: X\to T$ over $U$ satisfying the following properties
\begin{enumerate}
\item $K_X+B+\Mm_X\sim_{\mathbb R,T}0$,
\item there exists a birational morphism $X'\rightarrow X$ such that $\Mm_{X'}\sim_{\mathbb R,T}0$,
\item $\lfloor B\rfloor=D_1+D_2$, where $D_1,D_2$ are prime divisors and $f|_{D_i}: D_i\rightarrow T$ are isomorphisms,
\item $(X,B,\Mm)$ is plt, and
\item every reduced fiber of $f$ is isomorphic to $\mathbb P^1$.
\end{enumerate}
We call $D_1$ and $D_2$ the \emph{horizontal sections} of $(X,B,\Mm)/T$.
\end{defn}

\begin{defn}[$\mathbb P^1$-link, cf. {\cite[Definition 2.23]{FS23}}]\label{defn: p1 link}
Let $(X,B,\Mm)/U$ be a dlt g-pair associated with a projective morphism $f: X\rightarrow U$, such that $K_X+B+\Mm_X\sim_{\mathbb R,Z}0$. Let $Z_1$, $Z_2$ be two lc centers of $(X,B,\Mm)$. We say that $Z_1$ and $Z_2$ are \emph{directly $\mathbb P^1$-linked$/U$} if there exists an lc center $W$ of $(X,B,\Mm)$ satisfying the following
\begin{enumerate}
        \item $Z_i\subset W$ for each $i$.
        \item $f(W)=f(Z_1)=f(Z_2)$.
        \item Let $(W,B_W,\Mm^W)/U$ be the dlt g-pair induced by repeatedly applying adjunction to codimension $1$ lc centers
$$K_W+B_W+\Mm^W_W:=(K_X+B+\Mm_X)|_W.$$
        Then there exists a standard $\mathbb P^1$-link$/U$ $h: (W',B_{W'},\Mm^W)\rightarrow T$ such that $(W',B_{W'},\Mm^W)$ is crepant to $(W,B_W,\Mm^W)$, and $Z_1|_{W'}$ and $Z_2|_{W'}$ are the horizontal sections of $(W',B_{W'},\Mm^W)/T$.
\end{enumerate}
We say that $Z_1$ and $Z_2$ are \emph{$\mathbb P^1$-linked$/U$} if either $Z_1=Z_2$, or there exists an integer $n\geq 2$ and lc centers $Z_1',\dots,Z_n'$ of $(X,B,\Mm)$, such that $Z_1'=Z_1,Z_n'=Z_2$, and $Z'_i$ and $Z'_{i+1}$ are directly $\mathbb P^1$-linked$/U$ for any $1\leq i\leq n-1$.
\end{defn}

The following theorem is important when characterizing the structure of lc centers of g-pairs. We emphasize that, in the following theorem, we do not require $(X,B,\Mm)$ to be NQC.

\begin{thm}[{\cite[Theorem 3.5]{Bir20}; \cite[Theorem 1.4]{FS23} for the $\Qq$-coefficient case}]\label{thm: P1 link for gdlt crepant log structure} 
Let $(X,B,\Mm)/U$ be a dlt g-pair associated with a projective morphism $f: X\rightarrow U$, such that $K_X+B+\Mm_X\sim_{\mathbb R,U}0$. Let $s\in U$ be a (not necessarily closed) point such that $f^{-1}(s)$ is connected (as a $k(s)$-scheme). Let
$$\mathcal{S}:=\{V\mid V\text{ is an lc center of }(X,B,\Mm), s\in f(V)\}$$
and $Z,W\in\mathcal{S}$ be two elements such that $Z$ is minimal in $\mathcal{S}$ with respect to the inclusion. Then there exists $Z_W\in\mathcal{S}$ such that $Z_W\subset W$, and $Z$ and $Z_W$ are $\mathbb P^1$-linked$/U$. In particular, any minimal elements in $\mathcal{S}$ with respect to inclusion are $\mathbb P^1$-linked$/U$ to each other.
\end{thm}
\begin{proof}
\noindent\textbf{Step 1}. In this step, we show that the theorem holds over an \'etale neighborhood $(s'\in U')\rightarrow (s\in U)$ such that $k(s)=k(s')$. We use induction on $\dim X$ and on $\dim U$.

If $f^{-1}(s)\cap\lfloor B\rfloor$ is disconnected, then by \cite[Theorem 3.5]{Bir20}, after an \'etale base change, there are exactly two non-trivial lc centers of $(X,B,\Mm)$ intersecting $f^{-1}(s)$, and they are $\mathbb P^1$-linked with each other. We are done in this case.

If $f^{-1}(s)\cap\lfloor B\rfloor$ is connected, then we let $D_1,\dots,D_r$ be the irreducible components of $\lfloor B\rfloor$.  By passing to an \'etale neighborhood of $s\in S$ without changing $k(s)$, we may assume that each $D_i$ has connected fiber over $s$, and every lc center of $(X,B,\Mm)$ intersects $f^{-1}(s)$ (cf. \cite[Claim 4.38.1]{Kol13}). Possibly reordering indices, we may assume that $Z\subset D_1$, $W\subset D_r$, and 
$$f^{-1}(s)\cap D_i\cap D_{i+1}\not=\emptyset$$
for any $1\leq i\leq r-1$. Let $(D_i,B_{D_i},\Mm^{D_i})$ be the g-pair induced by adjunction
$$K_{D_i}+B_{D_i}+\Mm^{D_i}_{D_i}:=(K_X+B+\Mm_X)|_{D_i}$$
for each $i$.
\begin{claim}\label{claim: p1link induction}
Let $Z_1:=Z$. For any $2\leq i\leq r$, there exists an lc center $Z_i\subset D_{i-1}\cap D_{i}$ in $\mathcal{S}$ such that
\begin{enumerate}
    \item $Z_i$ and $Z_{i-1}$ are $\mathbb P^1$-linked$/U$ with each other,
    \item $Z_i$ is minimal in $\mathcal{S}$, and
    \item $Z_i$ is an lc center of $(D_i,B_{D_i},\Mm^{D_i})$ and $(D_{i-1},B_{D_{i-1}},\Mm^{D_{i-1}})$.
\end{enumerate}
\end{claim}
\begin{proof}
Suppose we have already constructed $Z_{i-1}$. By Lemma \ref{lem: inversion of adjunction gdlt}, $Z_{i-1}$ and $D_{i-1}\cap D_i$ are lc centers of $(D_{i-1},B_{D_{i-1}},\Mm^{D_{i-1}})$, and $Z_{i-1}$ is minimal among all lc centers of $(D_{i-1},B_{D_{i-1}},\Mm^{D_{i-1}})$ which dominate $s$. By induction hypothesis of $\dim X$ and $\dim U$, there exists an lc center $Z_i\subset D_{i-1}\cap D_i$ that is minimal among all lc centers of $(D_{i-1},B_{D_{i-1}},\Mm^{D_{i-1}})$ which dominate $s$, and $Z_i$ and $Z_{i-1}$ are $\mathbb P^1$-linked with each other. By Lemma \ref{lem: inversion of adjunction gdlt}, $Z_i$ is an lc center of $(X,B,\Mm)$, is minimal in $\mathcal{S}$, and is an lc center of $(D_i,B_{D_i},\Mm^{D_i})$. The claim follows by induction on $i$.
\end{proof}
\noindent\textit{Proof of Theorem \ref{thm: P1 link for gdlt crepant log structure} continued}. By Claim \ref{claim: p1link induction} applied to $i=r$, the theorem holds over an \'etale neighborhood $(s'\in U')\rightarrow (s\in U)$ such that $k(s)=k(s')$ under the induction hypothesis of \(\dim X\) and \(\dim U\).

\medskip

\noindent\textbf{Step 2}. We show that the \'etale base change was not necessary and conclude the proof of the theorem. Let 
$$X\xrightarrow{\tilde f}\tilde U\rightarrow U$$
be the Stein factorization of $f$. Since $f^{-1}(s)$ is connected, there exists a unique pre-image $\tilde s\in\tilde U$ of $s$. Let $Z_i$ be the minimal elements of $\mathcal{S}$. Since lc centers commute with \'etale base change, we see that there is a unique irreducible subvariety 
$$\tilde s\in \tilde V\subset U$$
such that $\tilde V=\tilde f(Z_i)$ for each $i$. 

Let $\tilde v$ be the generic point of \(\tilde V\). By \textbf{Step 1} and induction hypothesis, the theorem holds after an \'etale base change 
$$\tilde\pi: (\tilde v'\in \tilde U')\rightarrow (\tilde v\in\tilde U).$$
Since $\tilde f$ has connected fibers, \(\tilde\pi\) induces an isomorphism of the fibers
$$\tilde\pi: (\tilde f')^{-1}(\tilde v')\cong\tilde f^{-1}(\tilde v).$$
Thus $Z_i$ canonically lift to $Z_i'\cong Z_i$ and the $\mathbb P^1$-links$/U$ between the $Z_i'$ descend
to $\mathbb P^1$-links$/U$ between the $Z_i$.
\end{proof}

\begin{lem}\label{lem: intersection of lc center gpair}
Let $f: (X,B,\Mm)\rightarrow Z$ be an lc crepant log structure and $z\in Z$ a (not necessarily closed) point. Let $$\mathcal{S}_z:=\{V\mid V\text{ is an lc center of }f: (X,B,\Mm)\rightarrow Z, z\in V\}.$$
Then
\begin{enumerate}
    \item There exists a unique element $W\in\mathcal{S}_z$ that is minimal with respect to inclusion. 
    \item $W$ is unibranch at $z$, i.e., the completion $\widehat{W}_z$ is irreducible.
    \item Any intersection of lc centers of $f: (X,B,\Mm)\rightarrow Z$ is a union of lc centers.
\end{enumerate}
\end{lem}
\begin{proof}
By Definition-Lemma \ref{deflem: dlt model}, possibly replacing $(X,B,\Mm)$ with a dlt model, we may assume that $(X,B,\Mm)$ is dlt. Since $f$ is a contraction, $f^{-1}(z)$ is connected. For any element $W\in\mathcal{S}_z$ that is minimal with respect to inclusion, there exists an lc center $Z_W$ of $(X,B,\Mm)$ that is minimal among all lc centers whose image on $Z$ is equal to $W$ with respect to inclusion. By Theorem \ref{thm: P1 link for gdlt crepant log structure}, all such $Z_W$ are $\mathbb P^1$-linked$/Z$ to each other, hence their images on $Z$ are the same. This proves (1). (2) follows from (1) by considering every \'etale neighborhood of $z$. 

For any lc centers $W_1,W_2$ on $Z$, let $z\in W_1\cap W_2$ be any point. By (1), there exists a unique element $W_z$ of \(\mathcal{S}_z\). Then
$$z\in W_z\subset W_1\cap W_2,$$
so 
$$W_1\cap W_2=\cup_{z\in W_1\cap W_2}z\subset \cup_{z\in W_1\cap W_2}W_z\subset W_1\cap W_2.$$
Therefore, 
$$W_1\cap W_2=\cup_{z\in W_1\cap W_2}W_z$$
is a union of lc centers. We get (3).
\end{proof}

\begin{lem}\label{lem: gdlt crepant log structure is compatible under subadjunction}
Let $f: (X,B,\Mm)\to Z$ be a dlt crepant log structure and $Y\subset X$ an lc center. Let
$$
f|_Y: Y\xrightarrow{f_Y}Z_Y\xrightarrow{\pi} Z
$$
be the Stein factorization of $f|_Y$, and $(Y,B_Y,\Mm^Y)/Z$ the dlt g-pair induced \emph{by repeatedly applying adjunction to codimension $1$ lc centers}
$$K_Y+B_Y+\Mm_Y^Y:=(K_X+B+\Mm_X)|_Y.$$
Then
\begin{enumerate}
\item $f_Y: (Y,B_Y,\Mm^Y)\rightarrow Z_Y$ is a dlt crepant log structure.
\item For any lc center $W_Y\subset Z_Y$ of $f_Y: (Y,B_Y,\Mm^Y)\rightarrow Z_Y$, $\pi(W_Y)$ is an lc center of $f: (X,B,\Mm)\rightarrow Z$.
\item For any lc center $W\subset Z$ of $f: (X,B,\Mm)\rightarrow Z$, every irreducible component of $\pi^{-1}(W)$ is an lc center of $f_Y: (Y,B_Y,\Mm^Y)\rightarrow Z_Y$.
\end{enumerate}
\end{lem}

\begin{proof}
(1) We only need to show that $(Y,B_Y,\Mm^Y)$ is dlt, which follows from \cite[Lemma 2.9]{HL22}. 

(2) There exists an lc center $V_Y$ of $(Y,B_Y,\Mm^Y)$ such that $f_Y(V_Y)=W_Y$. By Lemma \ref{lem: inversion of adjunction gdlt}, $V_Y$ is also an lc center of $(X,B,\Mm)$. Thus $\pi(W_Y)=f(V_Y)$ is an lc center of $f: (X,B,\Mm)\rightarrow Z$.

(3) Let $z$ be the generic point of $W$. Since the question is \'etale local, possibly replacing $Z$ by an \'etale neighborhood of $z$ and replacing $Y$ with its irreducible components, we may assume that $f^{-1}(z)\cap Y$ is connected, and we only need to show that there exists an lc center $V_Y$ of $f_Y: (Y,B_Y,\Mm^Y)\rightarrow Z_Y$ such that $f_Y(V_Y)$ is an irreducible component of \(\pi^{-1}(W)\).

Let $V_X$ be a minimal lc center of $(X,B,\Mm)$ which dominates $W$, i.e., $V_X$ is minimal in
$$\{V\mid V\text{ is an lc center of }(X,B,\Mm), V\text{ dominates }W\}$$
with respect to inclusion. Then $f(V_X)=W$. By Theorem \ref{thm: P1 link for gdlt crepant log structure}, there exists an lc center $V_Y\subset Y$ of $(X,B,\Mm)$ that is $\mathbb P^1$-linked$/Z$ to $V_X$. By Lemma \ref{lem: inversion of adjunction gdlt}, $V_Y$ is also an lc center of $(Y,B_Y,\Mm^Y)$. Thus $f_Y(V_Y)\subset Z_Y$ is an lc center of $f_Y: (Y,B_Y,\Mm^Y)\rightarrow Z_Y$. Moreover, since $V_Y$ is $\mathbb P^1$-linked$/Z$ to $V_X$, $(f|_Y)(V_Y)=f(V_X)=W$. Thus $f_Y(V_Y)$ is an irreducible component of \(\pi^{-1}(W)\) and we are done.
\end{proof}

The following lemma is implicitly stated in \cite[Definition 4.7]{BZ16}.
\begin{lem}[{\cite[Definition 4.7]{BZ16}}]\label{lem: lc adjunction keep lc}
Let $(X,B,\Mm)/U$ be a g-pair and $S$ a component of $\lfloor B\rfloor$. Let $S^\nu$ be the normalization of $S$, and let $(S^\nu,B_S,\Mm^S)/U$ be the g-pair induced by the adjunction
$$K_{S^\nu}+B_S+\Mm^S_S:=(K_X+B+\Mm_X)|_S.$$
Assume that $(X,B,\Mm)$ is lc near $S$. Then $(S^\nu,B_S,\Mm^S)$ is lc.
\end{lem}

\begin{lem}\label{lem: lc adjunction keep lc center}
Let $(X,B,\Mm)/U$ be an lc g-pair, $S$ a component of $\lfloor B\rfloor$, $\nu: S^\nu\rightarrow S$ the normalization of $S$, and $(S^\nu,B_S,\Mm^S)/U$ the g-pair induced by the adjunction 
$$K_{S^\nu}+B_S+\Mm^S_{S^\nu}:=(K_X+B+\Mm_X)|_{S^\nu}.$$
Let $V$ be an lc center of $(X,B,\Mm)$. Let $W$ be an irreducible subvariety of $S^\nu$ such that $\nu(W)$ is a component of $V\cap S$. Then $W$ is an lc center of $(S^\nu,B_S,\Mm^S)$.
\end{lem}
\begin{proof}
By Lemma \ref{lem: intersection of lc center gpair}, $V\cap S$ is a union of lc centers of $(X,B,\Mm)$. In particular, \(\nu(W)\) is an lc center of $(X,B,\Mm)$. Possibly replacing $V$ with \(\nu(W)\), we may assume that $V\subset S$ and \(\nu(W)=S\). We let $f: Y\rightarrow X$ be a dlt modification of $(X,B,\Mm)$,
 $$K_Y+B_Y+\Mm_Y:=f^*(K_X+B+\Mm_X),$$
 and let $S_Y:=f^{-1}_*S$.  We let
 $$K_{S_Y}+B_{S_Y}+\Mm^S_{S_Y}:=(K_Y+B_Y+\Mm_Y)|_{S_Y}$$
 and let $f_S: S_Y\rightarrow S^\nu$ be the induced morphism such that $\nu\circ f_S=f|_{S_Y}$. Since $f: (Y,B_Y,\Mm)\rightarrow X$ is a dlt crepant log structure and $V$ is an lc center of $f: (Y,B_Y,\Mm)\rightarrow X$, by Lemma \ref{lem: gdlt crepant log structure is compatible under subadjunction}(1)(3), any irreducible component of \(\nu^{-1}(V)\) is an lc center of $f_S: (S_Y,B_{S_Y},\Mm^S)\rightarrow S$. In particular, $W$ is an lc center of $f_S: (S_Y,B_{S_Y},\Mm^S)\rightarrow S$. Since
 $$K_{S_Y}+B_Y+\Mm^S_{S_Y}=f_S^*(K_{S^\nu}+B_S+\Mm^S_{S^\nu}),$$
 $W$ is an lc center of $(S^\nu,B_S,\Mm^S)$.
\end{proof}

\begin{lem}\label{lem: same lct if adjunction is same}
Let $(X,B,\Mm)/U$ be an lc g-pair and $D_1\geq 0,D_2\geq 0$ two $\Rr$-Cartier $\Rr$-divisors on $X$. Let $S_1,\dots,S_n$ be distinct prime divisors on $X$ such that $B\geq\sum_{i=1}^n S_i$ and $V:=\cap_{i=1}^nS_i$. Assume that for any component $W$ of $V$ with normalization $W^\nu$, $W$ is not contained in $\Supp D_1\cup\Supp D_2$ and
$D_1|_{W^\nu}=D_2|_{W^\nu}.$
Let
$$t_i:=\sup\{t\geq 0\mid (X,B+tD_i,\Mm)\text{ is lc near }V\}$$
for $i=1,2$. Then $t_1=t_2$.
\end{lem}
\begin{proof}
We let $V_k:=\cap_{i=1}^kS_i$ for each $1\le k\le n-1$. Let $W$ be a component of $V$ with normalization $W^\nu$.
Then for each $1\le k\le n-1$, there exists an irreducible component $W_k$ of $V_k$ such that
$$W:=W_n\subset W_{n-1}\subset\dots\subset W_1=S_1.$$
We let $W_k^\nu$ be the normalization of $W_k$ for each $k$ and let $(W_0,B_0,\Mm^{(0)}):=(X,B,\Mm)$. By Lemmas \ref{lem: lc adjunction keep lc} and \ref{lem: lc adjunction keep lc center}, we may repeatedly apply adjunction and get lc g-pairs $(W_k^\nu,B_k,\Mm^{(k)})$, such that for each $1\leq k\leq n$
\begin{itemize}
    \item $W_k$ is an lc center of $(W_{k-1}^\nu,B_{k-1},\Mm^{(k-1)})$, and
    \item $K_{W_k^\nu}+B_k+\Mm^{(k)}_{W_k^\nu}=\left(K_{W_{k-1}^\nu}+B_{k-1}+\Mm^{(k-1)}_{W_{k-1}^\nu}\right)|_{W_k^\nu}=(K_X+B+\Mm_X)|_{W_k^\nu}.$
\end{itemize}
Let
$$t_{i,W}:=\sup\{t\geq 0\mid (X,B+tD_i,\Mm)\text{ is lc near }W\}$$
for $i=1,2$. Then there exists an lc center $T_i\subset\Supp D_i$ such that $T_i\cap W\not=\emptyset$. By Lemma \ref{lem: intersection of lc center gpair}, we may assume that $T_i\subset W$ for each $i$. By Lemma \ref{lem: lc adjunction keep lc}, $(W^\nu,B_n+t_{i,W}D_i|_{W^\nu},\Mm^{(n)})$ is lc for each $i$. By Lemma \ref{lem: lc adjunction keep lc center}, there exists an lc center $T_{W,i}$ of $(W^\nu,B_n+t_{i,W}D_i|_{W^\nu},\Mm^{(n)})$ such that $T_{W,i}\subset\Supp D_i|_{W^\nu}$. Since $D_1|_{W^\nu}=D_2|_{W^\nu}$,
$$t_{1,W}=\lct\left(W^\nu,B_n,\Mm^{(n)};D_1|_{W^\nu}\right)=\lct\left(W^\nu,B_n,\Mm^{(n)};D_2|_{W^\nu}\right)=t_{2,W}.$$
Therefore, $t_1=t_2$.
\end{proof}

\subsection{Inversion of adjunction} In this subsection, we prove the following inversion of adjunction for NQC generalized pairs

\begin{thm}\label{thm: inversion of adjunction}
Let $(X,B,\Mm)$ be an NQC g-pair and $S$ a component of $B^{=1}$. Let $S^\nu$ be the normalization of $S$, and let $(S^\nu,B_S,\Mm^S)/U$ be the g-pair induced by the adjunction
$$K_{S^\nu}+B_S+\Mm^S_S:=(K_X+B+\Mm_X)|_S.$$
Suppose that $(S^\nu,B_S,\Mm^S)$ is lc. Then $(X,B,\Mm)$ is lc near $S$.
\end{thm}
\begin{proof}
First we prove the case when $(X,B,\Mm)$ is a $\Qq$-g-pair. 

By Definition-Lemma \ref{deflem: dlt model}, there exists a birational morphism $f: Y\rightarrow X$ satisfying the following. Let $E$ be the reduced $f$-exceptional divisor and $B_Y:=f^{-1}_*(B\wedge\Supp B)+E$, then
\begin{enumerate}
    \item $(Y,B_Y,\Mm)$ is \(\Qq\)-factorial dlt,
    \item $a(F,X,B,\Mm)\leq 0$ for any prime $f$-exceptional divisor $F$.
\end{enumerate}
We let
$$K_Y+\bar B_Y+\Mm_Y:=f^*(K_X+B+\Mm_X)$$
and let $S_Y$ be the strict transform of $S$ on $Y$. Let $(S_Y,B_{S_Y},\Mm^S)/U$ and $(S_Y,\bar B_{S_Y},\Mm^S)/U$ be the g-pairs induced by adjunction
$$K_{S_Y}+B_{S_Y}+\Mm^S_{S_Y}=(K_Y+B_Y+\Mm_Y)|_{S_Y}$$
and
$$K_{S_Y}+\bar B_{S_Y}+\Mm^S_{S_Y}=(K_Y+\bar B_Y+\Mm_Y)|_{S_Y}$$
respectively. We let $Q:=\bar B_Y-B_Y$. 

Let $A$ be an ample divisor on $Y$ such that $K_Y+B_Y+\Mm_Y+A$ is nef. We may run a $(K_Y+B_Y+\Mm_Y)$-MMP$/X$ with scaling of $A$
$$(Y,B_Y,\Mm):=(X_0,B_0,\Mm)\dashrightarrow (X_1,B_1,\Mm)\dashrightarrow\dots\dashrightarrow (X_n,B_n,\Mm)\dashrightarrow\dots.$$
Let $S_i,A_i,Q_i,\bar B_i$ be the images of $S_Y,A,Q,\bar B_Y$ on $X_i$ for each $i$, $f_i: X_i\rightarrow X$ the induced birational morphism, and 
$$\lambda_i:=\inf\{t\geq 0\mid K_{X_i}+B_i+tA_i+\Mm_{X_i}\text{ is nef}/X\}$$
the scaling numbers. Then $K_{X_i}+\bar B_i+\Mm_{X_i}=f_i^*(K_X+B+\Mm_X)$ and $\bar B_i=B_i+Q_i$ for any $i$. Let
$$K_{S_i}+B_{S_i}+\Mm^S_{S_i}:=(K_{X_i}+B_i+\Mm_{X_i})|_{S_i}$$
and
$$K_{S_i}+\bar B_{S_i}+\Mm^S_{S_i}:=(K_{X_i}+\bar B_i+\Mm_{X_i})|_{S_i}$$
for any $i$. Then $\bar B_{S_i}=B_{S_i}+Q_i|_{S_i}$. Moreover, there exists a birational morphism $g_i: S_i\rightarrow S$ such that
$$K_{S_i}+\bar B_{S_i}+\Mm^S_{S_i}=g_i^*(K_S+B_S+\Mm^S_S).$$
Thus $(S_i,\bar B_{S_i},\Mm^S)$ is lc. Since $(Y,B_Y,\Mm)$ is dlt, $(X_i,B_i,\Mm)$ is dlt. By Lemma \ref{lem: inversion of adjunction gdlt}, $(S_i,B_{S_i},\Mm^S)$ is dlt. By Lemma \ref{lem: inversion of adjunction gdlt}, any lc center of $(X_i,B_i,\Mm)$ is an lc center of $(S_i,B_{S_i},\Mm^S)$. Since all components of $Q_i$ are lc centers of $(X_i,B_i,\Mm)$ and $(S_i,\bar B_{S_i},\Mm^S)$ is lc, $\Supp Q_i$ does not intersect $S_i$ for any $i$.

We pick a non-negative integer $m$ in the following way. If the $(K_Y+B_Y+\Mm_Y)$-MMP$/X$ terminates, then we let $m$ be the index so that $(X_m,B_m,\Mm)/X$ is a log minimal model of $(Y,B_Y,\Mm)/X$ for some non-negative integer $m$. If the $(K_Y+B_Y+\Mm_Y)$-MMP$/X$ does not terminate, then by Lemma \ref{lem: scaling number go to 0}, $\lim_{i\rightarrow+\infty}\lambda_i=0$, so by special termination (cf. \cite[Lemma 2.18]{LX23a}), we may pick a positive integer $m$, such that $S_i\dashrightarrow S_{i+1}$ is an isomorphism in codimension $1$ for any $i\geq m$. We let $I\geq 2$ be any sufficiently divisible positive integer satisfying the following
\begin{itemize}
    \item $IQ$ is a Weil divisor.
    \item 
    $$(f_m)_*\mathcal{O}_{X_m}(A_m-IQ_m)\subset (f_m)_*\mathcal{O}_{X_m}(A_m)$$
    are contained in 
    $$\mathcal{I}_{f_m(\Supp Q)}\cdot (f_m)_*\mathcal{O}_{X_m}(A_m).$$
    \item If $(X_m,B_m,\Mm)/X$ is a log minimal model of $(Y,B_Y,\Mm)/X$ and $m\geq 2$, then $\lambda_{m-1}>\frac{1}{I}$.
    \item If the $(K_Y+B_Y+\Mm_Y)$-MMP$/X$ does not terminate, then \(\lambda_m>\frac{1}{I}\).
\end{itemize}
Since $S_i\dashrightarrow S_{i+1}$ is an isomorphism in codimension $1$ for any $i\geq m$, for any $j\geq m$, we have
$$(f_j)_*\mathcal{O}_{X_j}(A_j-IQ_j)=(f_m)_*\mathcal{O}_{X_m}(A_m-IQ_m).$$
Since $\Supp Q_i$ does not intersect $S_i$ for any $i$, we have an induced homomorphism 
$$(f_i)_*\mathcal{O}_{X_i}(A_i-IQ_i)\rightarrow(f_m|_{S_m})_*\mathcal{O}_{S_m}(A_i)=(f_i|_{S_i})_*\mathcal{O}_{S_i}(A_i)$$
which is not surjective. Therefore,
$$R^1(f_i)_*\mathcal{O}_{X_i}(A_i-IQ_i-S_i)\not=0$$
for any $i\geq m$. 

We let $l:=m$ if $(X_m,B_m,\Mm)/X$ is a log minimal model of $(Y,B_Y,\Mm)/X$, and let $l$ be the unique positive integer such that $\lambda_{l-1}>\frac{1}{I}\geq\lambda_{l}$ if the $(K_Y+B_Y+\Mm_Y)$-MMP$/X$ does not terminate. Then $l\geq m$,
$$X_0\dashrightarrow X_1\dashrightarrow X_{l}$$
is also a sequence of steps of a $(K_{X_0}+B_0+\Mm_{X_0}+\frac{1}{I}A)$-MMP$/X$ with scaling of $A$, and $$K_{X_{l}}+B_{l}+\Mm_{X_{l}}+\frac{1}{I}A_{l}$$ 
is nef$/X$. Since $X_0$ is $\Qq$-factorial klt, we may pick $$0\leq \Delta_0\sim_{\mathbb Q}B_0-S_0+\Mm_{X_0}+\frac{1}{I}A$$ such that $(X_0,\Delta_0)$ is klt and $(X_0,S_0+\Delta_0)$ is plt. We let \(\Delta_l\) be the image of \(\Delta_0\) on \(X_l\), then $(X_l,S_l+\Delta_l)$ is plt, so $(X_l,\Delta_l)$ is klt. Then

$$A_l-IQ_l-S_l\sim_{\mathbb Q,X}K_{X_l}+\Delta_l+(I-1)\left(K_{X_l}+B_l+\Mm_{X_l}+\frac{1}{I}A_l\right),$$
so by the relative Kawamata-Viehweg vanishing \cite[Theorem 1-2-5]{KMM87}, $$R^1(f_{l})_*\mathcal{O}_{X_{l}}(A_{l}-IQ_{l}-S_{l})\not=0,$$
a contradiction. We are done with the case when $(X,B,\Mm)$ is a $\Qq$-g-pair.

\medskip

Now we prove the case when $(X,B,\Mm)$ is not necessarily a $\Qq$-g-pair, hence conclude the proof of the theorem. There exist real numbers $r_1,\dots,r_c$ such that $1,r_1,\dots,r_c$ are linearly independent over $\Qq$, $\bm{r}:=(r_1,\dots,r_c)\in\mathbb R^c$, and $\Qq$-linear functions $s_1,\dots,s_p,t_1,\dots,t_q$, such that
$$B=\sum_{i=1}^ps_i(1,\bm{r})B_i,\Mm=\sum_{i=1}^qt_i(1,\bm{r})\Mm_i,$$
where $B_i\geq 0$ are distinct Weil divisors and \(\Mm_i\) are nef$/X$ $\bb$-Cartier $\bb$-divisors. Let 
$$B(\bm{v}):=\sum_{i=1}^ps_i(1,\bm{v})B_i\text{ and }\Mm(\bm{v}):=\sum_{i=1}^qt_i(1,\bm{v})\Mm_i,$$
for any $\bm{v}\in\mathbb R^c$. 

Since the coefficients of divisors under adjunction are transformed via $\Qq$-linear functions, there are $\Qq$-linear functions $s'_1,\dots,s'_{p'},t'_1,\dots,t'_{q'}$, distinct Weil divisors $B_{S_i}\geq 0$, and nef$/X$ $\bb$-Cartier $\bb$-divisors \(\Mm^S_i\), 
$$B_S(\bm{v}):=\sum_{i=1}^ps_i(1,\bm{v})B_{S,i},\text{ and }\Mm^S(\bm{v}):=\sum_{i=1}^qt_i(1,\bm{v})\Mm^S_i,$$
such that
$$K_{S^\nu}+B_S(\bm{v})+\Mm^S(\bm{v})_{S^\nu}=(K_X+B(\bm{v})+\Mm(\bm{v})_X)|_{S^\nu}$$
for any $\bm{v}\in\mathbb R^c$. Since $$(S^\nu,B_S=B_S(\bm{r}),\Mm^S=\Mm^S(\bm{r}))$$ is lc, 
there exists an open neighborhood $U\ni\bm{r}$ of $\mathbb R^c$ such that
$$(S^\nu,B_S(\bm{v}),\Mm^S(\bm{v}))$$
is lc for any $\bm{v}\in U$. By the $\Qq$-g-pair case,
$$(X,B(\bm{v}),\Mm(\bm{v}))$$
is lc for any $\bm{v}\in U\cap\mathbb Q$. Thus $$(X,B=B(\bm{r}),\Mm=\Mm(\bm{r}))$$
is lc by continuity of log discrepancies.
\end{proof}

\begin{rem}
    We do not need Theorem \ref{thm: inversion of adjunction} in the rest of the paper but we expect it to be useful for future works. We remark that several alternative versions of Theorem \ref{thm: inversion of adjunction} can be found in \cite[Theorems 1.5, 1.6, 6.7]{Fil20} but we cannot apply them directly to prove Theorem \ref{thm: inversion of adjunction} because of the following reasons
\begin{enumerate}
\item All these theorems require that $(X,B,\Mm)$ is a $\Qq$-g-pair.
\item \cite[Theorem 1.5]{Fil20} requires $S$ to be a minimal lc center and $S$ is projective.
\item \cite[Theorem 1.6]{Fil20} requires that $X$ is projective and $(X,B,\Mm)$ is a \(\Qq\)-g-pair. Moreover, the potential g-pair structure constructed on $W^\nu$ \cite[Theorem 1.6]{Fil20} is not known to be identical to the g-pair structure constructed in \cite[Theorem 4.5]{HL21b}.
\item \cite[Theorem 6.7]{Fil20} requires that $X$ is $\Qq$-factorial projective klt.
\end{enumerate}  
\end{rem}

\subsection{Boundedness on the number of components}

\begin{prop}\label{prop: bound number of components}
Let $\gamma_0\le1$ be a positive real number, and $b_1,\dots,b_n\in [\gamma_0,1]$ positive real numbers. Let $(X,B=\sum_{i=1}^nb_iB_i+D,\Mm)/X$ be an lc g-pair and $x\in X$ a point, such that $B_i\geq 0$ is a non-zero $\Qq$-Cartier Weil divisor for each $i$, and $D\geq 0$. Suppose that $\bar x\subset\Supp B_i$ for each $i$. Moreover, assume that one of the following holds
    \begin{enumerate}
        \item $\Mm$ is NQC$/X$.
        \item There exists a klt g-pair $(X,B',\Mm')/X$.
        \item $\gamma_0=1$ and each $B_i$ is Cartier.
    \end{enumerate}
Then 
$$n\leq\frac{\dim X-\dim\bar x}{\gamma_0}.$$
\end{prop}
\begin{proof}
When $\dim X=1$ the proposition is trivial, so we may assume that $\dim X\geq 2$. We may also assume that $n\geq 1$, otherwise there is nothing left to prove. 

Let $B_{n+1},\dots,B_{n+\dim\bar X}$ be general hyperplane sections on $X$ and let $b_i:=1$ when $i\geq n+1$. Possibly replacing $x$ with $\bar x\cap\cap_{i=1}^{\dim X}H_i$ and $B$ with $\sum_{i=1}^{n+\dim\bar x}b_iB_i+D$, we may assume that $x$ is a closed point. 

First we prove the proposition under conditions (1) or (2). Possibly adding general hyperplane sections which pass through $x$, we may assume that $x$ is an lc center of $(X,B,\Mm)$. Let $E$ be an lc place of $(X,B,\Mm)$ such that $\Center_XE=x$.
\begin{claim}\label{claim: extract divisor which is ample}
    There exists a contraction $f: Y\rightarrow X$ of $E$ such that $-E$ is ample$/X$.
\end{claim}
\begin{proof}
    If $\Mm$ is NQC$/X$, then the claim follows from \cite[Theorem 1.7]{LX23b}. Otherwise, the claim follows from \cite[Lemma 2.11]{Bir20}.
\end{proof}
\noindent\textit{Proof of Proposition \ref{prop: bound number of components} continued}. By Claim \ref{claim: extract divisor which is ample}, there exists a contraction $f: Y\rightarrow X$ of $E$ such that $-E$ is ample$/X$. We let $B_{i,Y},D_Y,B_Y$ be the strict transforms of $B_{i},D,B$ on $Y$ respectively. Since $x\in\Supp B_i$ for each $i$, $\mult_EB_i>0$ for each $i$, so $B_{i,Y}$ is ample$/X$ for each $i$. We let $E^\nu$ be the normalization of $E$, $\Mm^E:=\Mm|_{E^\nu}$, and let
$$K_{E^\nu}+B_E+\Mm^E_{E^\nu}=(K_Y+B_Y+E+\Mm_Y)|_{E^\nu}.$$
We let $B_{i,E}:=\Supp(B_{i,Y}|_{E^\nu})$ for each $i$. Then for any component $D_{i,j}$ of $B_{i,E}$, we have
$$\mult_{D_{i,j}}B_E=\frac{n_{i,j}-1+\sum_{k=1}^nb_km_{k,i,j}+\gamma_{i,j}}{n_{i,j}}$$
for some real number $\gamma_{i,j}\geq 0$ and non-negative integers $m_{k,i,j}$, such that $m_{i,i,j}\not=0$. Since $(E,B_E,\Mm^E)$ is lc, $\mult_{D_{i,j}}B_E\leq 1$, so $$B_E\geq\sum_{i=1}^nb_iB_{i,E}.$$
Since $B_{i,Y}$ is ample$/X$, $B_{i,Y}|_{E^\nu}$ is ample, so $B_{i,E}$ is big. The proposition under conditions (1) or (2) follows from \cite[Proposition 5.1]{BZ16}.

Now we prove the proposition under condition (3). Let $S$ be the normalization of an irreducible component of $B_1$ such that $x\in S_1$, and let $(S,B_S,\Mm^S)$ be the g-pair induced by the adjunction
$$K_S+B_S+\Mm^S_S:=(K_X+B+D+\Mm_X)|_S.$$
Since $x\in\Supp B_i$ for each $i$, $B_i|_S\not=0$ for any $i\geq 2$. Since $B_i$ is Cartier and $(S,B_S,\Mm^S)$ is lc, $B_i|_S=\Supp(B_i|_S)$ for any $i\geq 2$, and 
$$B_S\geq\sum_{i=2}^nB_i|_S.$$
Since each $B_i|_S$ is Cartier, by induction on $\dim X$, we have $n\leq\dim X$ and the proposition follows.
\end{proof}

\section{Stability of generalized pairs}\label{sec: stability gpair}

In this section, we discuss the stability properties of g-pairs. We will define the concepts of generically lc, Property $(*)$ BP (semi-)stable, and log stable for g-pairs, and then study the basic properties of g-pairs satisfying these properties. This section is parallel to \cite[Section 2]{ACSS21}.

\subsection{Toroidal generalized pairs}

\begin{defn}[{cf. \cite[Definition 2.1]{ACSS21}}]\label{defn: toroidal g-pairs}
Let $(X,\Sigma_X,\Mm)/U$ be a g-pair. We say that $(X,\Sigma_X,\Mm)$ is \emph{toroidal} if $\Sigma_X$ is a reduced divisor, $\Mm$ descends to $X$, and for any closed point $x\in X$, there exists a toric variety $X_{\sigma}$, a closed point $t\in X_{\sigma}$, and an isomorphism of complete local algebras 
$$\phi_x:\widehat{\mathcal{O}}_{X,x}\cong\widehat{\mathcal{O}}_{X_\sigma,t}$$
such that the ideal of $\Sigma_X$ maps to the invariant ideal of $X_{\sigma}\backslash T_{\sigma}$, where $T_\sigma\subset X_\sigma$ is the maximal torus of $X_{\sigma}$. Any such $(X_\sigma, t)$ will be called a \emph{local model} of $(X,\Sigma_X,\Mm)$ at $x\in X$.

Let $(X,\Sigma_X,\Mm)/U$ and $(Z,\Sigma_Z,\Mm^Z)/U$ be toroidal g-pairs and $f: X\rightarrow Z$ a surjective morphism$/U$. We say that $f: (X,\Sigma,\Mm)\rightarrow (Z,\Sigma,\Mm^Z)$ is \emph{toroidal} if for every closed point $x\in X$, there exist a local model $(X_\sigma,t)$ of $(X,\Sigma_X,\Mm)$ at $x$, a local model $(Z_{\tau},s)$ of $(Z,\Sigma_Z,\Mm^Z)$ at $z:=f(x)$, and a toric morphism $g: X_\sigma\to Z_{\tau}$, so that the diagram of algebras commutes:
\begin{center}$\xymatrix{
    \widehat{\mathcal{O}}_{X,x}\ar@{->}[r]^{\cong}  &     \widehat{\mathcal{O}}_{X_{\sigma},t} \\
     \widehat{\mathcal{O}}_{Z,z}\ar@{->}[r]^{\cong}\ar@{->}[u] & \widehat{\mathcal{O}}_{Z_{\tau},s}\ar@{->}[u]
}$
\end{center}
Here the vertical maps are the algebra homomorphisms induced by $f$ and $g$, respectively. 
\end{defn}

\begin{defthm}[{\cite[Definition-Theorem 6.5]{LLM23}, \cite[Theorem 2.2]{ACSS21}}]\label{defthm: weak ss reduction}
Let $X$ be a normal quasi-projective variety, $X\rightarrow U$ a projective morphism, $X\rightarrow Z$ a contraction, $B$ an $\Rr$-divisor on $X$, $\Mm$ a nef$/U$ $\bb$-divisor on $X$, $D_1,\dots,D_m$ prime divisors over $X$, and $D_{Z,1},\dots,D_{Z,n}$ prime divisors over $Z$. Then there exist a toroidal g-pair $(X',\Sigma_{X'},\Mm)/U$, a log smooth pair $(Z',\Sigma_{Z'})$, and a commutative diagram
 \begin{center}$\xymatrix{
X'\ar@{->}[r]^{h}\ar@{->}[d]_{f'}& X\ar@{->}[d]^{f}\\
Z'\ar@{->}[r]^{h_Z} & Z\\
}$
\end{center}
satisfying the following.
\begin{enumerate}
\item $h$ and $h_Z$ are projective birational morphisms.
\item $f': (X',\Sigma_{X'},\Mm)\rightarrow (Z',\Sigma_{Z'})$ is a toroidal contraction.
\item $\Supp(h^{-1}_*B)\cup\Supp\Exc(h)$ is contained in $\Supp\Sigma_{X'}$.
\item $X'$ has at most toric quotient singularities.
\item $f'$ is equi-dimensional.
\item $\Mm$ descends to $X'$.
\item $X'$ is $\Qq$-factorial klt.
\item The center of each $D_i$ on $X'$ and the center of each $D_{Z,i}$ on $Z'$ are divisors.
\end{enumerate}
We call any such $f': (X',\Sigma_{X'},\Mm)\rightarrow (Z',\Sigma_{Z'})$ (associated with $h$ and $h_Z$) which satisfies (1–7) an \emph{equi-dimensional model} of $f: (X,B,\Mm)\rightarrow Z$. 
\end{defthm}
\begin{proof}
Possibly replacing $X$ and $Z$ with high models, we may assume that $\Mm$ descends to $X$, each $D_i$ is a divisor on $X$, and each $D_{Z,i}$ is a divisor on $Z$. Now the theorem follows from \cite[Theorem 2.2]{ACSS21}, which in turn follows from \cite[Theorem 2.1 and Proposition 4.4]{AK00}. We also refer the reader to \cite[Theorem B.6]{Hu20} for a more detailed explanation.
\end{proof}

\begin{rem}
In Definition-Theorem \ref{defthm: weak ss reduction}, it is important to note that the contraction $X\rightarrow Z$ may not necessarily be over $U$. This kind of phenomenon will appear throughout the rest of the paper.
\end{rem}

\subsection{Discriminant and moduli parts of generalized pairs}

\begin{defn}[Birationally equivalent morphisms, cf. {\cite[Page 4, Paragraph 2]{ACSS21}}]
Let $f: X\rightarrow Z$ and $f': X'\rightarrow Z'$ be surjective morphisms between normal varieties. We say that $f$ and $f'$ are \emph{birationally equivalent} if there exist birational maps $h: X\dashrightarrow X'$ and $h_Z: Z\dashrightarrow Z'$ such that $f'\circ h=h_Z\circ f$.
\end{defn}

\begin{defn}[Generically lc, cf. {\cite[2.2. Discriminant and Moduli Part]{ACSS21}}]\label{defn: generically lc}
Let $(X,B,\Mm)/U$ be a g-sub-pair and $f: X\rightarrow Z$ a contraction. We say that $(X,B,\Mm)$ is \emph{generically (sub-)lc$/Z$} if $(X,B,\Mm)$ is (sub-)lc over the generic point of $Z$. Note that $f$ may not be a contraction$/U$. We remark that we will not use the notation ``GLC" for ``generically lc" as in \cite{ACSS21} since GLC also stands for ``generalized lc" in many references.
\end{defn}

\begin{defn}[Crepant generalized pairs, cf. {\cite[Definition 2.3]{ACSS21}}]\label{defn: crepant g-pairs}
Let $(X,B,\Mm)/U$ and $(X',B',\Mm')/U$ be two g-sub-pairs and $f: X\rightarrow Z$, $f': X'\rightarrow Z'$ two contractions. We say that $(X,B,\Mm)$ and $(X',B',\Mm')$ are \emph{crepant over the generic point of $Z$} if we have the following commutative diagram
 \begin{center}$\xymatrix{
& W\ar@{->}[ld]_{p}\ar@{->}[dr]^{q} &\\
X\ar@{.>}[rr]^{h}\ar@{->}[d]_{f}& & X'\ar@{->}[d]^{f'}\\
Z\ar@{.>}[rr]^{h_Z} & & Z'\\
}$
\end{center}
satisfying the following. Let $$K_W+B_W+\Mm_W:=p^*(K_X+B+\Mm_X)$$ and $$K_W+B'_W+\Mm'_W:=q^*(K_{X'}+B'+\Mm'_X).$$ Then:
\begin{enumerate}
    \item $h$ and $h_Z$ are birational maps. In particular, $f$ and $f'$ are birationally equivalent.
    \item $\Mm$ and $\Mm'$ descend to $W$.
    \item $B_W-B'_W$ and $\Mm_W-\Mm'_W$ are vertical$/Z$.
\end{enumerate}
\end{defn}

\begin{defn}[Discriminant and moduli parts, cf. {\cite[Definition 2.3]{ACSS21}}]
Let $(X,B,\Mm)/U$ be a g-sub-pair and $f: X\rightarrow Z$ a contraction such that $(X,B,\Mm)$ is generically sub-lc$/Z$. In the following, we fix a choice of $K_X$ and a choice of $K_Z$, and suppose that for any birational morphism $g: \bar X\rightarrow X$ and $g_Z: \bar Z\rightarrow Z$, $K_{\bar X}$ and $K_{\bar Z}$ are chosen as the Weil divisors such that $g_*K_{\bar X}=K_X$ and $(g_Z)_*K_{\bar Z}=K_Z$.

Let $f': X'\rightarrow Z'$ be any contraction that is birationally equivalent to $f$ such that the induced birational maps $h: X'\dashrightarrow X$ and $h_Z: Z'\dashrightarrow Z$ are morphisms and $Z'$ is $\Qq$-factorial. We let 
$$K_{X'}+B'+\Mm_{X'}:=h^*(K_X+B+\Mm_X).$$ 
For any prime divisor $D$ on $Z'$, we define
$$b_D(X',B',\Mm;f):=1-\sup\left\{t\mid \left(X',B'+tf'^*D,\Mm\right)\text{ is sub-lc over the generic point of } D\right\}.$$
Since being sub-lc is a property that is preserved under crepant transformations, $b_D(X,B,\Mm;f)$ is independent of the choices of $X'$ and $Z'$ and is also independent of $U$.

Since $(X,B,\Mm)$ is generically sub-lc$/Z$, $(X',B',\Mm)$ is generically sub-lc$/Z$, so we may define
$$B_{Z'}:=\sum_{D\text{ is a prime divisor on }Z'}b_D(X,B,\Mm;f)D$$
and $$N_{X'}:=K_{X'}+B'+\Mm_{X'}-f'^*(K_{Z'}+B_{Z'}).$$
We call $B_{Z'}$ and $N_{X'}$ the \emph{discriminant part} and \emph{trace moduli part} of $f': (X',B',\Mm)\rightarrow Z'$ respectively, and call $B_Z:=(h_Z)_*B_{Z'}$ and $N_X:=h_*B$ the \emph{discriminant part} and \emph{trace moduli part} of $f: (X,B,\Mm)\rightarrow Z$ respectively.

By construction, there exist two $\bb$-divisors $\Bb$ on $Z$ and $\Nn$ on $X$, such that for any contraction $f'': X''\rightarrow Z''$ that is birationally equivalent to $f$ such that the induced birational maps $h': X''\dashrightarrow X'$ and $h_{Z'}: Z''\dashrightarrow Z'$ are morphisms and $Z''$ is $\Qq$-factorial, $\Bb_{Z''}$ is the discriminant part of $f'': (X'',B'',\Mm)\rightarrow Z''$, and $\Nn_{X''}$ is the trace moduli part of $f'': (X'',B'',\Mm)\rightarrow Z''$, where 
$$K_{X''}+B''+\Mm_{X''}:=h'^*(K_{X'}+B'+\Mm_{X'}).$$ 
We call $\Nn$ the \emph{moduli part} of $f: (X,B,\Mm)\rightarrow Z$ and $\Bb$ the \emph{discriminant $\bb$-divisor} of $f: (X,B,\Mm)\rightarrow Z$. By construction, $\Bb$ is uniquely determined and $\Nn$ is uniquely determined for any fixed choices of $K_X$ and $K_Z$.
\end{defn}

\subsection{Boundary property of generalized pairs}

\begin{defn}[BP (semi-)stable, boundary property, cf. {\cite[Definition 2.5]{ACSS21}}]\label{defn: bpstable}
Let $(X,B,\Mm)/U$ be a g-sub-pair and $f: X\rightarrow Z$ a contraction, such that $(X,B,\Mm)$ is generically sub-lc$/Z$. Let $\Bb$ be the discriminant $\bb$-divisor of $f: (X,B,\Mm)\rightarrow Z$.

We say that $f: (X,B,\Mm)\rightarrow Z$ is BP stable (resp. BP semi-stable) if $K_Z+\Bb_Z$ is $\Rr$-Cartier, and for any birational morphism $h_Z: Z'\rightarrow Z$,
$$h_Z^*(K_Z+\Bb_Z)=\text{(resp. }\geq\text{)} K_{Z'}+\Bb_{Z'}.$$
If $f: (X,B,\Mm)\rightarrow Z$ is BP stable (resp. BP semi-stable), then we say that $(X,B,\Mm)$ is BP stable (resp. BP semi-stable) over $Z$.

We say that $f: (X,B,\Mm)\rightarrow Z$ \emph{satisfies the boundary property} if there exists a contraction $f': X'\rightarrow Z'$ that is birationally equivalent to $f$, such that the induced birational maps $h: X'\dashrightarrow X$ and $h_Z: Z'\dashrightarrow Z$ are morphisms, $K_{X'}+B'+\Mm_{X'}:=h^*(K_X+B+\Mm_X)$, and $f': (X',B',\Mm)\rightarrow Z'$ is BP stable.
\end{defn}

\begin{lem}[cf. {\cite[Remark 2.6(1)]{ACSS21}}]\label{lem: bp stable preserved under birational morphism}
Let $(X,B,\Mm)/U$ be a g-sub-pair and $f: X\rightarrow Z$ a contraction, such that $(X,B,\Mm)$ is generically sub-lc$/Z$. Let $h: X'\rightarrow X$ be a birational morphism, $K_{X'}+B'+\Mm_{X'}:=h^*(K_X+B+\Mm_X),$ and $f':=f\circ h$. Let $\Bb$ and $\Bb'$ be the discriminant $\bb$-divisors of $f: (X,B,\Mm)\rightarrow Z$ and $f': (X',B',\Mm)\rightarrow Z$ respectively. Then $\Bb=\Bb'$. In particular, $f: (X,B,\Mm)\rightarrow Z$ is BP stable if and only if $f': (X',B',\Mm)\rightarrow Z$ is BP stable.
\end{lem}
\begin{proof}
Let $f'': X''\rightarrow Z''$ be any contraction that is birationally equivalent to $f$ such that $Z''$ is $\Qq$-factorial and the induced birational map $h': X''\dashrightarrow X'$ is a morphism. Let $$K_{X''}+B''+\Mm_{X''}:=h'^*(K_{X'}+B'+\Mm_{X'}),$$
then $K_{X''}+B''+\Mm_{X''}=(h\circ h')^*(K_X+B+\Mm_X)$. Thus $\Bb_{Z''}=\Bb'_{Z''}$, so $\Bb=\Bb'$. The in particular part follows from the definition of BP stability.
\end{proof}

\begin{lem}[cf. {\cite[Remark 2.6(2)]{ACSS21}}]\label{lem: bp stable implies descend}
Let $(X,B,\Mm)/U$ be a g-sub-pair and $f: X\rightarrow Z$ a contraction, such that $f: (X,B,\Mm)\rightarrow Z$ is BP stable. Let $B_Z$ and $\Nn$ be the discriminant part and the moduli part of $f: (X,B,\Mm)\rightarrow Z$ respectively. Then:
\begin{enumerate}
    \item $\Nn_X=K_X+B+\Mm_X-f^*(K_Z+B_Z).$
    \item $\Nn$ descends to $X$.
\end{enumerate}
\end{lem}
\begin{proof}
For any $f': X'\rightarrow Z'$ that is birationally equivalent to $f$, such that the induced birational maps $h: X'\dashrightarrow X$ and $h_Z: Z'\dashrightarrow Z$ are morphisms, we have $K_{Z'}+B_{Z'}=h_Z^*(K_{Z}+B_{Z})$. Thus
\begin{align*}
    \Nn_{X'}&=K_{X'}+B'+\Mm_{X'}-f'^*(K_{Z'}+B_{Z'})=h^*(K_X+B+\Mm_X-f^*(K_Z+B_Z)),
\end{align*}
where $K_{X'}+B'+\Mm_{X'}:=h^*(K_X+B+\Mm_X)$, and $B_{Z'}$ is the discriminant part of $f': (X',B',\Mm)\rightarrow Z'$. The lemma immediately follows.
\end{proof}

\begin{cor}\label{cor: boundary property imply descend to high model}
Let $(X,B,\Mm)/U$ be a g-sub-pair and $f: X\rightarrow Z$ a contraction, such that $f: (X,B,\Mm)\rightarrow Z$ satisfies the boundary property. Then the moduli part $\Nn$ of $f: (X,B,\Mm)\rightarrow Z$ descends to some birational model $X'$ of $X$.
\end{cor}
\begin{proof}
It immediately follows from Lemma \ref{lem: bp stable implies descend}.
\end{proof}

\subsection{Log stable generalized pairs}

\begin{lem}[cf. {\cite[Lemma 2.9]{ACSS21}}]
Let $(X,B,\Mm)$ be a g-sub-pair and $F$ a prime divisor over $X$. Then there exists an $\Rr$-divisor $\Delta$ on $X$ such that $(X,\Delta,\Mm)$ is a sub-lc g-sub-pair and $F$ is the unique lc place of $(X,\Delta,\Mm)$.
\end{lem}
\begin{proof}
Let $h: W\rightarrow X$ be a log resolution of $(X,
\Supp B)$ such that $F$ is a divisor on $W$ and $\Mm$ descends to $W$, and let
$$K_W+B_W+\Mm_W:=h^*(K_X+B+\Mm_X).$$
Let $C_W:=B_W^{\geq 1}-B_W^{\geq 1}\wedge F$ and let $a:=a(F,X,B,\Mm)$. Let $A_1, A_2\geq 0$ be sufficiently general ample $\mathbb R$-divisors on $Y$ such that 
$$aF-C_W+A_1-A_2\sim_{\mathbb R}0.$$
Then $(W,\Delta_W:=B_W+aF-C_W+A_1-A_2,\Mm)$ is lc and $F$ is the unique lc place of $(W,\Delta_W,\Mm)$. We let $\Delta:=h_*\Delta_W$, then 
$$K_W+\Delta_W+\Mm_W=h^*(K_X+\Delta+\Mm_X)$$
and $(X,\Delta,\Mm)$ has the required properties.
\end{proof}

\begin{defn}[Log stable generalized pairs, cf. {\cite[Definition 2.10]{ACSS21}}]\label{defn: logstable g-pair}
Let $(X,B,\Mm)/U$ be a g-sub-pair and $f: X\rightarrow Z$ a contraction, such that $(X,B,\Mm)$ is generically sub-lc$/Z$. Let $B_Z$ be the discriminant part of $f: (X,B,\Mm)\rightarrow Z$. We say that $(X,B,\Mm)$ is \emph{log stable$/Z$} if for any $\Rr$-Cartier $\Rr$-divisor $H$ on $Z$, $(X,B+f^*H,\Mm)$ is sub-lc if and only if $(Z,B_Z+H)$ is sub-lc.

It is clear that if $(X',B',\Mm)/U$ is a g-sub-pair that is crepant to $(X,B,\Mm)$ with birational morphism $X'\to X$, then $(X,B,\Mm)$ is log stable$/Z$ if and only if $(X',B',\Mm)$ is log stable$/Z$.
\end{defn}

\begin{thm}[cf. {\cite[Theorem 2.11]{ACSS21}}]\label{thm: bp stable is log stable}
Let $(X,B,\Mm)/U$ be a g-sub-pair and $f: X\rightarrow Z$ a contraction such that $(X,B,\Mm)$ is generically sub-lc$/Z$. Then $(X,B,\Mm)$ is BP stable$/Z$ if and only if $(X,B,\Mm)$ is log stable$/Z$.
\end{thm}
\begin{proof}
By Lemma \ref{lem: bp stable preserved under birational morphism}, possibly replacing $(X,B,\Mm)$ by a crepant model, we may assume that $\Mm$ descends to $X$. Then the theorem follows from the definitions and \cite[Theorem 2.11]{ACSS21}.
\end{proof}

\subsection{Property \texorpdfstring{$(*)$}{} generalized pairs}

\begin{lem}[cf. {\cite[Lemma 2.12]{ACSS21}}]\label{lem: dimension of full lc rank}
Let $(X,B,\Mm)/U$ be a g-pair and $f: X\rightarrow Z$ a contraction. Let $d:=\dim X$ and $m:=\dim Z$. Let $z\in Z$ be a closed point, $D_1,\dots,D_m\geq 0$ Cartier divisors on $Z$, such that $z\in\Supp D_i$ for each $i$ and $(X,B+\sum_{i=1}^mf^*D_i,\Mm)$ is lc over $f^{-1}(z)$. Then the dimension of any irreducible component of $f^{-1}(z)$ is $d-m$.
\end{lem}

\begin{proof}
For any irreducible component $G$ of $f^{-1}(z)$, let $H_1,\dots,H_{\dim G}$ be general very ample divisors on $X$, $V:=\cap_{i=1}^{\dim G}H_i$, and $(V,B_V,\Mm^V)/U$ the g-pair induced by the adjunction $$K_V+B_V+\Mm^V_V:=(K_X+B+\Mm_X)|_V.$$ Then $(V,B_V+\sum_{i=1}^mf^*D_i|_V,\Mm^V)$ is lc, $G\cap V$ is a closed point, and $A_i:=f^*D_i|_V$ is Cartier and contains $G\cap V$ for any $i$. By Proposition \ref{prop: bound number of components}, $m\leq\dim V=d-\dim G$. Thus $\dim G\leq d-m$. Therefore, the dimension of any irreducible component of $f^{-1}(z)$ is $\leq d-m$. By \cite[Exercise II 3.22 (a)]{Har77}, the dimension of any irreducible component of $f^{-1}(z)$ is $\geq d-m$. The lemma immediately follows.
\end{proof}

\begin{defn}[Property $(*)$ generalized pairs, cf. {\cite[Definition 2.13]{ACSS21}}]\label{defn: property *}
Let $(X,B,\Mm)/U$ be a g-sub-pair and $f: X\rightarrow Z$ a contraction. We say that $f: (X,B,\Mm)\rightarrow Z$ satisfies \emph{Property $(*)$} if there exists a reduced divisor $\Sigma_Z$ on $Z$ satisfying the following.
\begin{enumerate}
\item $(Z,\Sigma_Z)$ is log smooth. In particular, $Z$ is smooth.
\item The vertical$/Z$ part $B^v$ of $B$ is equal to $f^{-1}(\Sigma_Z)$. In particular, $B^v$ is reduced and $\Sigma_Z$ is the image of $B^v$ on $Z$.
\item For any closed point $z\in Z$ and any reduced divisor $\Sigma\ge \Sigma_Z$ on $Z$ such that $(Z,\Sigma)$ is log smooth near $z$, $(X,B+f^*(\Sigma-\Sigma_Z),\Mm)$ is sub-lc over a neighborhood of $z$.
\end{enumerate}
%If $f$ is clear from the context and $f: (X,B,\Mm)\rightarrow Z$ satisfies Property $(*)$, then we also say that $(X,B,\Mm)/Z$ satisfies \emph{Property $(*)$}. 
By (2), $\Sigma_Z$ is uniquely determined by $f: (X,B,\Mm)\rightarrow Z$. We will temporarily call $\Sigma_Z$ the \emph{base divisor} associated with $f: (X,B,\Mm)\rightarrow Z$. 
\end{defn}

In the following lemma, we will show that $\Sigma_Z$ is actually the discriminant part of $f: (X,B,\Mm)\rightarrow Z$.

\begin{lem}[cf. {\cite[Lemma 2.14]{ACSS21}}]\label{lem: basic property (*) gpair}
Let $(X,B,\Mm)/U$ be a g-sub-pair and $f: X\rightarrow Z$ a contraction such that $f: (X,B,\Mm)\rightarrow Z$ satisfies Property $(*)$. Let $\Sigma_Z$ be the base divisor associated with $f: (X,B,\Mm)\rightarrow Z$. Then:
\begin{enumerate}
\item $(X,B,\Mm)$ is sub-lc.
\item $\Sigma_Z$ is the discriminant part of $f: (X,B,\Mm)\rightarrow Z$.
\item If $B\geq 0$, then $f$ is equi-dimensional over $Z\backslash\Supp\Sigma_Z$.
\end{enumerate}
\end{lem}
\begin{proof}
(1) follows from the definition immediately.
 
(2) Let $B_Z$ be the discriminant part of $f: (X,B,\Mm)\rightarrow Z$. By definition and (1), $\Supp B_Z\ge B_Z\geq\Sigma_Z$. It suffices to show that $B_Z\le\Sigma_Z$. Let $P$ be a prime divisor on $Z$ such that $P\not\subset\Supp\Sigma_Z$, and let $z$ be a general closed point in $P$. Then $(Z,\Sigma_Z+P)$ is log smooth at $z$. Then $(X,B+f^*P,\Mm)$ 
is sub-lc over a neighborhood of $z$ which implies that
$$\sup\{t\mid (X,B+tf^*P,\Mm)\text{ is sub-lc over the generic point of }P\}=1.$$
Thus $P\not\subset\Supp B_Z$ and hence $\Sigma_Z=\Supp\Sigma_Z\geq\Supp B_Z\ge B_Z$. This implies (2).

(3) Let $d:=\dim X$ and $m:=\dim Z$. Let $z\in Z\backslash\Supp\Sigma_Z$ be a closed point, and let $\Sigma_1,\dots,\Sigma_m$ be general hyperplane sections on $Z$ such that $z\in \Sigma_i$ for any $i$. Then $(Z,\Sigma_Z+\sum_{i=1}^m\Sigma_i)$ is log smooth at $z$. Then $(X,B+\sum_{i=1}^m f^*\Sigma_i,\Mm)$ is lc over a neighborhood of $z$. By Lemma \ref{lem: dimension of full lc rank}, the dimension of any irreducible component of $f^{-1}(z)$ is $d-m$. This implies (3).
\end{proof}

\begin{lem}[cf. {\cite[Lemma 2.15]{ACSS21}}]\label{lem: lcc property* pullback} 
Let $(X,B,\Mm)/U$ be a g-pair and $f: X\rightarrow Z$ a contraction such that $f: (X,B,\Mm)\rightarrow Z$ satisfies Property $(*)$. Let $\Sigma_Z$ be the discriminant part of $f: (X,B,\Mm)\rightarrow Z$, and let $\Sigma\geq\Sigma_Z$ be a reduced divisor on $Z$, such that $(Z,\Sigma)$ is log smooth.

Consider $\Sigma$ to be a reduced subscheme of $Z$. Then for any irreducible stratum $V$ of $\Sigma$, any irreducible component of $f^{-1}(V)$ is an lc center of $(X,B+f^{-1}(\Sigma-\Sigma_Z),\Mm)$.
\end{lem}
\begin{proof} 
Let $k:=\dim Z-\dim V$. Since $(Z,\Sigma)$ is log smooth, there exist irreducible components $\Sigma_1,\dots, \Sigma_k$ of $\Sigma$ such that $V=\bigcap_{i=1}^k \Sigma_i$. By Definition \ref{defn: property *}(3), for any $i$ and any general closed point $z\in\Supp\Sigma_i$, $(X,B+f^*(\Sigma-\Sigma_Z),\Mm)$ is sub-lc over a neighborhood of $z$. Thus any irreducible component of $f^{-1}(\Sigma_i)$ is an lc center of $(X,B+f^*(\Sigma-\Sigma_Z),\Mm)$. Therefore, any irreducible component of $f^{-1}(V)$ is an intersection of lc centers of $(X,B+f^*(\Sigma-\Sigma_Z),\Mm)$. The lemma follows from Lemma \ref{lem: intersection of lc center gpair}.
\end{proof}

\begin{prop}[cf. {\cite[Proposition 2.16]{ACSS21}}]\label{prop: weak ss satisfies *}
Let $(X,\Sigma_X,\Mm)/U$ be a toroidal g-pair, $(Z,\Sigma_Z)$ a log smooth pair, and $f: (X,\Sigma_X,\Mm)\rightarrow (Z,\Sigma_Z)$ a toroidal morphism. Let $(X,B,\Mm)/U$ be a sub-lc g-sub-pair such that $\Supp B\subset \Supp\Sigma_X$, and the vertical$/Z$ part of $B$ is equal to $f^{-1}(\Sigma_Z)$. Then $f: (X,B,\Mm)\rightarrow Z$ satisfies Property $(*)$.
\end{prop}
\begin{proof}
Since $\Mm$ descends to $X$, by \cite[Proposition 2.16]{ACSS21}, $f: (X,B)\rightarrow Z$ satisfies Property $(*)$. By Definition \ref{defn: property *}, $f: (X,B,\Mm)\rightarrow Z$ satisfies Property $(*)$.
\end{proof}

The following result indicates that we can always get Property $(*)$ g-pairs up to weak semi-stable reductions.

\begin{prop}[cf. {\cite[Proposition 2.17]{ACSS21}}]\label{prop: weak ss imply *}
Let $(X,B,\Mm)/U$ be a g-sub-pair and $f: X\rightarrow Z$ a contraction, such that $(X,B,\Mm)$ is generically sub-lc$/Z$. Let $f': (X',\Sigma_{X'},\Mm)\rightarrow (Z',\Sigma_{Z'})$ be an equi-dimensional model of $f: (X,B,\Mm)\rightarrow Z$, associated with $h: X'\rightarrow X$ and $h_Z: Z'\rightarrow Z$. Then there exist two $\Rr$-divisors $B'$ and $F'$ on $X'$, such that $F'$ is vertical$/Z'$, $\Supp B'\cup\Supp F'\subset\Sigma_{X'}$, and 
$$K_{X'}+B'+\Mm_{X'}=h^*(K_X+B+\Mm_X)+F'.$$ 
In particular, $(X',B',\Mm)$ and $(X,B,\Mm)$ are crepant over the generic point of $Z$.

Moreover, if $(X,B,\Mm)$ is sub-lc, then $F'\geq 0$ and $f': (X',B',\Mm)\rightarrow Z'$ satisfies Property $(*)$.
\end{prop}
\begin{proof} 
Possibly adding components to $\Sigma_{Z'}$, we may assume that $\Sigma_{Z'}$ coincides with the image of the vertical$/Z'$ part of $\Sigma_{X'}$. We let $G:=f'^{-1}(\Sigma_{Z'})$ and
$$K_{X'}+\tilde B'+\Mm_{X'}:=h^*(K_X+B+\Mm_X),$$
then $G\subset\Supp\Sigma_{X'}$ and $\Supp\tilde B'\subset\Supp\Sigma_{X'}$. We may define an $\Rr$-divisor $B'$ on $X'$ as follows: For any prime divisor $D$ on $X'$,
\begin{itemize}
\item if $D$ is not a component of $\Supp\tilde B'$ nor $G$, then $\mult_DB'=0$,
\item if $D$ is a component of $G$, then $\mult_DB'=1$,
\item if $D$ is a component of $\tilde B'$ but is not a component of $G$, then $\mult_DB'=\mult_D\tilde B'$.
\end{itemize}
Set $F':=B'-\tilde B'$. Since $G$ is the vertical$/Z'$ part of $\Sigma_{X'}$, the vertical$/Z'$ part of $B'$ is equal to $G=f'^{-1}(\Sigma_{Z'})$. By construction $\Supp F'\subset G$ is vertical$/Z'$ and $\Supp B'\cup\Supp F'\subset\Sigma_{X'}$.

Assume that $(X,B,\Mm)$ is sub-lc. As $\mult_DF'=1-\mult_D\tilde B'$ if $D$ is a component of $G$, one can see that $\mult_DF'\geq 0$ and hence $F'\ge0$. By Proposition \ref{prop: weak ss satisfies *}, $f': (X',B',\Mm)\rightarrow Z'$ satisfies Property $(*)$.
\end{proof}

The following proposition shows that Property $(*)$ is stable under the MMP.
\begin{prop}[cf. {\cite[Proposition 2.18]{ACSS21}}]\label{prop: MMP preserves *}
Let $(X,B,\Mm)/U$ be an lc g-pair and $f: X\rightarrow Z$ a contraction, such that $f: (X,B,\Mm)\rightarrow Z$ satisfies Property $(*)$. Let $\phi: (X,B,\Mm)\dashrightarrow (Y,B_Y,\Mm)$ be a sequence of steps of a $(K_X+B+\Mm_X)$-MMP$/Z$ and $f_Y: Y\rightarrow Z$ the induced morphism. Assume that $\phi$ is also a sequence of steps of a $(K_X+B+\Mm_X)$-MMP$/U$. Then:
\begin{enumerate}
    \item $f_Y: (Y,B_Y,\Mm)\rightarrow Z$ satisfies Property $(*)$, and the discriminant part of $f_Y: (Y,B_Y,\Mm)\rightarrow Z$ is equal to the discriminant part of $f: (X,B,\Mm)\rightarrow Z$.
    \item For any closed point $z\in Z$, $\phi^{-1}$ is an isomorphism near the generic point of any irreducible component of $f_Y^{-1}(z)$.
    \item If $f$ is equi-dimensional, then $f_Y$ is equi-dimensional.
\end{enumerate}
\end{prop}
\begin{proof}
Without loss of generality, we may assume that $\phi$ is a step of a $(K_X+B+\Mm_X)$-MMP$/Z$. 

(1) Let $\Sigma_Z$ be the discriminant part of $f: (X,B,\Mm)\rightarrow Z$. Then by assumption $(Z,\Sigma_Z)$ is log smooth. Since the vertical$/Z$ part of $B$ is equal to $f^{-1}(\Sigma_Z)$ and $\phi$ does not extract any divisor, the vertical$/Z$ part of $B_Y$ is equal to $\phi\circ f^{-1}(\Sigma_Z)=f_Y^{-1}(\Sigma_Z)$. For any reduced divisor $\Sigma\geq\Sigma_Z$ on $Z$, $(X,B+f^*(\Sigma-\Sigma_Z),\Mm)/U$ is lc. Since $\phi$ is a step of a $(K_X+B+\Mm_X)$-MMP$/Z$, $\phi$ is also a step of a $(K_X+B+f^*(\Sigma-\Sigma_Z)+\Mm_X)$-MMP$/Z$. Thus 
$$(Y,B_Y+\phi_*f^*(\Sigma-\Sigma_Z)=B_Y+f_Y^*(\Sigma-\Sigma_Z),\Mm)$$
is lc. Therefore, $f_Y: (Y,B_Y,\Mm)\rightarrow Z$ satisfies Property $(*)$. By Lemma \ref{lem: basic property (*) gpair}(2), $\Sigma_Z$ is the discriminant part of $f_Y: (Y,B_Y,\Mm)\rightarrow Z$.

(2) Possibly shrinking $Z$ near $z$, there exists a reduced divisor $\Sigma\geq\Sigma_Z$ on $Z$, such that $(Z,\Sigma)$ is log smooth and $z$ is a stratum of $\Sigma$. By Lemma \ref{lem: lcc property* pullback}, any irreducible component of $f_Y^{-1}(z)$ is an lc center of $(Y,B_Y+f_Y^*(\Sigma-\Sigma_Z),\Mm)$. By Definition \ref{defn: property *}(3), $(X,B+f^*(\Sigma-\Sigma_Z),\Mm)$ is lc. For any irreducible component $G$ of $f^{-1}(z)$, let $D_G$ be an lc place of $(Y,B_Y+f_Y^*(\Sigma-\Sigma_Z),\Mm)$ over the generic point of $G$. Then
$$0\leq a(D_G,X,B+f^*(\Sigma-\Sigma_Z),\Mm)\leq a(D_G,Y,B_Y+f_Y^*(\Sigma-\Sigma_Z),\Mm)=0.$$
Thus
$$a(D_G,X,B+f^*(\Sigma-\Sigma_Z),\Mm)=a(D_G,Y,B_Y+f_Y^*(\Sigma-\Sigma_Z))=0,$$
so $\phi^{-1}$ is an isomorphism near the generic point of $G$.

(3) It immediately follows from (2).
\end{proof}

\begin{lem}\label{lem: snc base change}
Let $X$ be a smooth variety of dimension $d$ and $o\in X$ a closed point. Let $D_1,\dots,D_d$ and $D_1',\dots,D_d'$ be prime divisors on $X$ such that $o\in D_i$, $o\in D_i'$ for each $i$, and $(X,\sum_{i=1}^dD_i)$, $(X,\sum_{i=1}^dD_i')$ are log smooth at $o$. Then there exists an index $j$ such that $(X,\sum_{i=1}^{d-1}D_i+D_j')$ is log smooth at $o$.
\end{lem}
\begin{proof}
There are local coordinates systems $\{x_1,\dots,x_d\}$ and $\{x_1',\dots,x_d'\}$ near $o$, such that $D_i=(x_i=0)$ and $D_i'=(x_i'=0)$ for each $i$. 

Suppose that $(X,\sum_{i=1}^{d-1}D_i+D_j')$ is not log smooth near $o$ for any $j$, then locally analytically near $o$, $(x_j'=0)$ is contained in the subspace spanned by $x_1,\dots,x_{d-1}$ for any $j$. Since locally analytically near $o$, $X$ is spanned by $x_1',\dots,x_d'$, $X$ is contained in its subspace spanned by $x_1,\dots,x_{d-1}$, which is absurd. The lemma follows.
\end{proof}

\begin{prop}[Weak Bertini type theorem for Property $(*)$ modifications, cf. {\cite[Proposition 2.19]{ACSS21}}]\label{prop: weak bertini g-pair property *}
Let $(X,B,\Mm)/U$ be an lc g-pair and $f: X\rightarrow Z$ a contraction, such that $f: (X,B,\Mm)\rightarrow Z$ satisfies Property $(*)$. Let $z\in Z$ be a closed point. Then for any general ample$/U$ $\Rr$-divisor $A$ on $X$, over a neighborhood of $z$,
\begin{enumerate}
    \item $f: (X,B+A,\Mm)\rightarrow Z$ satisfies Property $(*)$, and
    \item the discriminant part of $f: (X,B+A,\Mm)\rightarrow Z$ is equal to the discriminant part of $f: (X,B,\Mm)\rightarrow Z$.
\end{enumerate}
%We remark that the choice of $A$ may depend on the choice of $z$.
\end{prop}

\begin{proof}
Let $\Sigma_Z$ be the discriminant part of $f: (X,B,\Mm)\rightarrow Z$. Then $\Sigma_Z$ is a reduced divisor, $(Z,\Sigma_Z)$ is log smooth, and the vertical$/Z$ part of $B$ is equal to $f^{-1}(\Sigma_Z)$. Since $A$ is general ample$/U$, all components of $A$ are horizontal$/U$. Therefore, the vertical part of $B+A$ is equal to $f^{-1}(\Sigma_Z)$. By Lemma \ref{lem: basic property (*) gpair}, (2) will follow from (1) and we only need to prove (1).

Let $\Sigma_0\geq\Sigma_Z$ be a reduced divisor on $Z$ such that $(Z,\Sigma_0)$ is log smooth near $z$ and $z$ is a stratum of $\Sigma_0$. By Definition \ref{defn: property *}, $(X,B+f^*(\Sigma_0-\Sigma_Z),\Mm)$ is lc over a neighborhood of $z$. Therefore, we may choose $A$ so that $(X,B+A+f^*(\Sigma_0-\Sigma_Z),\Mm)$ is lc over a neighborhood of $z$. By Definition \ref{defn: property *} again, we only need to show that, for any reduced divisor $\Sigma\geq\Sigma_Z$ such that $(Z,\Sigma)$ is log smooth near $z$, $(X,B+A+f^*(\Sigma-\Sigma_Z),\Mm)$ is lc over a neighborhood of $z$.

Possibly adding components to $\Sigma$, we may assume that $z$ is an lc center of $(Z,\Sigma)$. We let $q:=\dim Z$, then $\Sigma$ and $\Sigma_0$ both have $q$ irreducible components near $z$. Let $D^0_1,\dots,D^0_q$ be the components of $\Sigma_0$ and let $D_1,\dots,D_q$ be the components of $\Sigma$. Suppose that $\Sigma$ and $\Sigma_0$ have $k$ different components for any integer $0\leq k\leq q$. Possibly reordering indices, we may assume that $D^0_i\not=D_i$ for any $1\leq i\leq k$, and $D^0_i=D_i$ for any $k+1\leq i\leq q$. If $k=0$ then there is nothing left to prove, so we may assume that $k>0$. By repeatedly applying Lemma \ref{lem: snc base change}, possibly reordering indices, we may assume that
$$\left(Z,\Sigma_j:=\sum_{i=j+1}^qD_i^0+\sum_{i=1}^jD_i\right)$$
is log smooth for any $0\leq j\leq k$, and $\Sigma_k=\Sigma$. Moreover, $\Sigma_j\geq\Sigma_Z$ for each $j$, and $\Sigma_j$ and $\Sigma_{j+1}$ has exactly $1$ different component. By inductively showing that $(X,B+A+f^*(\Sigma_j-\Sigma_Z),\Mm)$ is lc for each $j$, we may assume that $k=1$.

Let $V:=\cap_{i=2}^qD_i^0$. Then we have $D_1|_V=D_1^0|_V=z$. Therefore, for any irreducible component $W$ of $f^{-1}(V)=\cap_{i=2}^qf^{-1}(D_i^0)$ with normalization $W^\nu$, $f^*D_1|_{W^\nu}=f^*D_i^0|_{W^\nu}$. Let $$R:=B+A+f^*\left(\sum_{i=2}^qD_i^0-\Sigma_Z\right).$$ By Lemma \ref{lem: same lct if adjunction is same},
\begin{align*}
    1=&\sup\left\{t\geq 0|\left(X,R+tf^*D_1^0,\Mm\right)\text{ is lc over a neighborhood of }z\right\}\\
=&\sup\left\{t\geq 0|\left(X,R+tf^*D_1^0,\Mm\right)\text{ is lc near }f^{-1}(V)\right\}\\
=&\sup\left\{t\geq 0|\left(X,R+tf^*D_1,\Mm\right)\text{ is lc near }f^{-1}(V)\right\}\\
=&\sup\left\{t\geq 0|\left(X,R+tf^*D_1,\Mm\right)\text{ is lc over a neighborhood of }z\right\}.
\end{align*}
Thus $(X,B+A+f^*(\Sigma-\Sigma_Z),\Mm)$ is lc over a neighborhood of $z$, and the proposition follows.
\end{proof}

\subsection{Maximal moduli}

Our definition of \emph{maximal moduli} has slightly differences with the ones defined in \cite[Proposition-Definition 1]{Sho23} or \cite[Definition 2.20]{ACSS21}.

\begin{defn}[cf. {\cite[Proposition-Definition 1]{Sho23}, \cite[Definition 2.20]{ACSS21}}]\label{defn: maximal moduli}
Let $(X,B,\Mm)/U$ be a g-sub-pair and $f: X\rightarrow Z$ a contraction, such that $(X,B,\Mm)$ is generically sub-lc$/Z$. Let $\Nn$ be the moduli part of $f: (X,B,\Mm)\rightarrow Z$. We say that $f: (X,B,\Mm)\rightarrow Z$ has \emph{maximal moduli} if $\Nn_X$ is $\Rr$-Cartier and the following conditions hold. 

For any g-sub-pair $(X',B',\Mm')/U$ and contraction $f': X'\rightarrow Z'$, such that
\begin{itemize}
\item $(X',B',\Mm')$ is generically sub-lc$/Z'$,
    \item $K_{X'}+B'+\Mm'_{X'}$ is nef$/Z'$, 
    \item $f'$ is birationally equivalent to $f$ with induced morphism $h:X'\to X$, and
    \item $(X,B,\Mm)$ and $(X',B',\Mm')$ are crepant over the generic point of $Z$,
\end{itemize}
then $|h^*\Nn_X-\Nn'_{X'}|\not=\emptyset$, where $\Nn'$ is the moduli part of $f': (X',B',\Mm')\rightarrow Z'$.
\end{defn}

\begin{prop}[{cf. \cite[Proposition 2.21, Remark 2.22]{ACSS21}}]\label{prop: * imply maximal moduli}
Let $(X,B,\Mm)/U$ and $(X',B',\Mm)/U$ be g-sub-pairs and $f: X\rightarrow Z$, $f': X'\rightarrow Z'$ contractions such that
    \begin{enumerate}
    \item  $(X,B,\Mm)$ is generically sub-lc$/Z$, %and  $(X',B',\Mm)$ is generically sub-lc$/Z'$,
    \item $(X,B,\Mm)$ and $(X',B',\Mm)$ are crepant over the generic point of $Z$,
    %\item $\Nn$ is the moduli part of $f: (X,B,\Mm)\rightarrow Z$ and $\Nn'$ is the moduli part of $f': (X',B',\Mm)\rightarrow Z'$.
    \item $f: (X,B,\Mm)\rightarrow Z$ satisfies Property $(*)$ and $\Nn_X$ is $\Rr$-Cartier,
    \item $(X',B',\Mm)$ is BP stable$/Z'$ and $\Nn'_{X'}$ is nef, and
    \item either $f$ is equi-dimensional or $B\geq 0$,
    \end{enumerate}
where $\Nn$ is the moduli part of $f: (X,B,\Mm)\rightarrow Z$ and $\Nn'$ is the moduli part of $f': (X',B',\Mm)\rightarrow Z'$. Then $f: (X,B,\Mm)\rightarrow Z$ has maximal moduli.
\end{prop}
\begin{proof}
Let $\phi_Z: Y\rightarrow Z$ and $\phi_{Z'}: Y\rightarrow Z'$ be a resolution of indeterminacy of the induced birational map $Z\dashrightarrow Z'$. Let $\phi: W\rightarrow X$ and $\phi': W\rightarrow X'$ be any high enough resolution of indeterminacy of the induced birational map $X\dashrightarrow X'$ such that $W\dashrightarrow Y$ is a morphism. Let
$$K_W+B_W'+\Mm_W:=\phi'^*(K_{X'}+B'+\Mm_{X'}).$$ 
Since $(X',B',\Mm)$ is BP stable$/Z'$, $(W,B_W',\Mm)$ is BP stable$/Y$ and $\Nn'$ descends to $X'$ by Proposition \ref{lem: bp stable implies descend}. In particular, $K_W+B_W'+\Mm_W$ is nef$/Z'$. We only need to show that $\phi^*\Nn_X-\Nn'_W\geq 0$. Possibly replacing $(X',B',\Mm)$ and $Z'$ by $(W,B_W',\Mm)$ and $Y$ respectively, we may further assume that both $h:X'\dashrightarrow X$ and $h_Z:Z'\dashrightarrow Z$ are morphisms.

Let $B_Z$ (resp. $B_{Z'}$) be the discriminant part of $f: (X,B,\Mm)\rightarrow Z$ (resp. $f': (X',B',\Mm)\rightarrow Z'$). By Lemma \ref{lem: basic property (*) gpair}(2), $B_Z$ is a reduced divisor, $(Z,B_Z)$ is log smooth, and the vertical$/Z$ part of $B$ is $f^{-1}(B_Z)$. We may write 
$$K_{X'}+\tilde{B}'+\Mm_{X'}=h^*(K_X+B+ \Mm_X)\text{ and }K_{Z'}+\tilde{B}_{Z'}=h_Z^*(K_Z+B_Z)$$
for some $\Rr$-divisors $\tilde{B}'$ and $\tilde{B}_{Z'}$. %Since $(X',B',\Mm)$ is BP stable$/Z'$, $K_{Z'}+B_{Z'}$ is $\Rr$-Cartier and $\Nn'$ descends to $X'$ by Proposition \ref{lem: bp stable implies descend}. 
Then 
\begin{align*}
h^*\Nn_X-\Nn'_{X'}=\tilde{B}'-B'-(f')^*\left(\tilde{B}_{Z'}-B_{Z'}\right)
\end{align*}
which is vertical over $Z$. We claim that $h_*(h^*\Nn_X-\Nn'_{X'})\ge0$. 

Pick any component $D$ of $\Supp h_*(h^*\Nn_X-\Nn'_{X'})$. Let $D'$ be the strict transform of $D$ on $X'$  and $D_Z$ the image of $D$ on $Z$. There are two possibilities.

\medskip

\noindent\textbf{Case 1}. $D_Z$ is contained in $B_Z$.

In this case, since $f: (X,B,\Mm)\rightarrow Z$ satisfies Property $(*)$, $\mult_{D'}\tilde{B}'=\mult_DB=1$ and $(X',\tilde{B}',\Mm)$ is log stable$/Z'$. Since $(Z',\tilde{B}_{Z'})$ is sub-lc and $B_{Z'}$ is the discriminant part of $f': (X',B',\Mm)\rightarrow Z'$, $(X',B'+(f')^*(\tilde{B}_{Z'}-B_{Z'}),\Mm)$ is sub-lc. It follows that
$$\mult_{D'}\left(B'+(f')^*\left(\tilde{B}_{Z'}-B_{Z'}\right)\right)\leq 1=\mult_{D'}\tilde{B}',$$
so $\mult_Dh_*(h^*\Nn_X-\Nn'_{X'})\geq 0.$

 \medskip

\noindent\textbf{Case 2}. $D_Z$ is not contained in $B_Z$. 
    
In this case, by assumption and Lemma \ref{lem: basic property (*) gpair}(3), $f$ is equi-dimensional over the generic point of $D_Z$. Thus $D_Z$ is a prime divisor. Since $f: (X,B,\Mm)\rightarrow Z$ satisfies Property $(*)$, it holds that $1=\mult_D(B+f^*D_Z)$. Hence if we denote by $D_{Z'}$ the strict transform of $D_Z$ on $Z'$, then
$$1=\mult_D(B+f^*D_Z)=\mult_D\left(h_*\left(\tilde{B}'+(f')^*D_{Z'}\right)\right)=\mult_{D'}\left(\tilde{B}'+(f')^*D_{Z'}\right).$$
Let
$$b_{D_{Z'}}(X',B',\Mm;f')=\sup\left\{t\mid \left(X',B'+t(f')^*D_{Z'},\Mm\right)\text{ is sub-lc over the generic point of }D_{Z'}\right\}.$$
Then $\mult_{D_{Z'}}B_{Z'}=1-b_{D_{Z'}}(X',B',\Mm;f')$ and hence 
$$D_{Z'}-B_{Z'}=b_{D_{Z'}}(X',B',\Mm;f')D_{Z'}$$ 
near the generic point of $D_{Z'}$. Moreover, as $\mult_{D_{Z'}}\tilde{B}_{Z'}=\mult_{D_Z}B_Z=0$, we see that 
\begin{align*}
\mult_{D'}\left(h^*\Nn_X-\Nn'_{X'}\right)
=&\mult_{D'}\left(\tilde{B}'-B'-(f')^*\tilde{B}_{Z'}+(f')^*B_{Z'}\right)\\
=&\mult_{D'}\left(\tilde{B}'+(f')^*D_{Z'}-B'-(f')^*(D_{Z'}-B_{Z'})\right)\\
=&1-\mult_{D'}\left(B'+(f')^*(D_{Z'}-B_{Z'})\right)\\
=&1-\mult_{D'}\left(B'+b_{D_{Z'}}(X',B',\Mm;f')\cdot (f')^*D_{Z'}\right)\ge0.
\end{align*}
The claim holds. Then the proposition follows from the negativity lemma immediately. We finish the proof.
\end{proof}

\part{Cone theorem and MMP for algebraically integrable foliations}\label{part:cone}

\section{Precise adjunction formula for algebraically integrable foliations}\label{sec: adjunction}

In this section, we will establish Theorem \ref{thm: precise adj gfq} under the additional assumption that $\Ff$ is induced by a contraction. The complete proof of Theorem \ref{thm: precise adj gfq} will be provided in Section \ref{sec: cone}.

\begin{thm}[Precise adjunction formula for generalized foliated quadruples]\label{thm: precise adj gfq}
Let $(X,\Ff,B,\Mm)/U$ be a gfq. Suppose that 
\begin{itemize}
    \item $\Ff$ is algebraically integrable, and
    \item $B=\epsilon(S)S+\sum_{j=1}^m b_j B_j$, and $\Mm=\sum_{k=1}^n r_k \Mm_k$,
\end{itemize}
where $b_1,\dots,b_m,r_1,\dots,r_n\ge0$, $S,B_1,\dots,B_m$ are distinct prime divisors on $X$, and $\Mm_1,\dots,\Mm_n$ are nef$/U$ $\bb$-Cartier $\bb$-divisors on $X$. Let $S^\nu\rightarrow S$ be the normalization of $S$ and $\Ff_S$ the restricted foliation of $\Ff$ on $S^\nu$.%and let $\Mm_{k}^S:=\Mm_k|_{S^\nu}$ for any $k$.

Then there exist prime divisors $T_1,\dots,T_l,C_1,\dots,C_q$ on $S^\nu$, positive integers $w_1,\dots,w_q$, and non-negative integers $\{w_{i,j}\}_{1\leq i\leq q,1\leq j\leq m}$ and $\{v_{i,k}\}_{1\leq i\leq q, 1\leq k\leq n}$ satisfying the following. For any real numbers $b_1',\dots,b_m'$ and $r_1',\dots,r_n'$ such that 
$$\left(X,\Ff,B':=\epsilon(S)S+\sum_{j=1}^m b_j' B_j,\Mm':=\sum_{k=1}^n r_k' \Mm_{k}\right)$$ 
is a sub-gfq, then
$$\left(K_\Ff+B'+\Mm'\right)|_{S^\nu}=K_{\Ff_S}+B'_{S}+\Mm'^{S}_{S^\nu},$$
where
$$B'_{S}:=\sum_{i=1}^l T_i+\sum_{i=1}^q \frac{w_i-1+\sum_{j=1}^m w_{i,j}b_j'+\sum_{k=1}^n v_{i,k}r_k'}{w_i}C_i\text{ and }\Mm'^{S}:=\Mm'|_{S^\nu}.$$
Moreover, if $(X,\Ff,B',\Mm')$ is lc near $S$, then $(S^\nu,\Ff_S,B'_{S},\Mm'^{S})$ is lc.
\end{thm}

\subsection{Preliminaries for algebraically integrable foliations}

\begin{defn}[Algebraically integrable foliations, {cf. \cite[3.1]{ACSS21}}]\label{defn: algebraically integrable}
Let $X$ be a normal variety and $\Ff$ a foliation on $X$. We say that $\Ff$ is an \emph{algebraically integrable foliation} if there exists a dominant map $f: X\dashrightarrow Y$ such that $\Ff=f^{-1}\Ff_Y$, where $\Ff_Y$ is a foliation by points. In this case, we say that $\Ff$ is \emph{induced by $f$}.
\end{defn}

We will use the following result throughout the paper.

\begin{lem}[{cf. \cite[Lemma 2.7]{DLM23}}]\label{lem:foliation-invariant}
    Let $f:X'\to X$ be a proper birational morphism of normal varieties, $\Ff$ a foliation on $X$, and $\Ff':=f^{-1}\Ff$ the pullback foliation on $X'$. Then $\Ff'$ is algebraically integrable if and only if $\Ff$ is algebraically integrable.
\end{lem}

\begin{defn}[Tangent and transverse, {cf. \cite[Section 3.4]{ACSS21}}]\label{defn: tangent to foliation}
 Let $X$ be a normal variety, $\Ff$ a foliation on $X$, and $V\subset X$ a subvariety. Suppose that $\Ff$ is a foliation induced by a dominant rational map $X\dashrightarrow Z$. We say that $V$ is \emph{tangent} to $\Ff$ if there exists a birational morphism $\mu: X'\rightarrow X$, an equidimensional contraction $f': X'\rightarrow Z$, and a subvariety $V'\subset X'$, such that
    \begin{enumerate}
    \item $\mu^{-1}\Ff$ is induced by $f'$, and
        \item $V'$ is contained in a fiber of $f'$ and $\mu(V')=V$.
    \end{enumerate}
    We say that $V$ is \emph{transverse} to $\Ff$ if $V$ is not tangent to $\Ff$.

For any point $x\in V$, we say that $V$ is \emph{transverse} to $\Ff$ at $x$ if $x\not\in\Sing(X)\cup\Sing(\Ff)\cup\Sing(V)$, and for any analytic neighborhood $U$ of $x$, $T_{V}|_U\rightarrow T_X|_U$ does not factor through $T_{\Ff}|_U$. We say that $V$ is \emph{everywhere transverse} to $\Ff$ if $V$ is transverse to $\Ff$ at $x$ for any $x\in V$ (in particular, $V$ is smooth and $V$ does not intersect $\Sing(X)$ or $\Sing(\Ff)$). We say that $V$ is \emph{generically transverse} to $\Ff$ if $V$ is transverse to $\Ff$ at the generic point $\eta_V$ of $V$.
\end{defn}

\begin{defn}[Tangency of general fibers]
 Let $X$ be a normal variety, $\Ff$ a foliation on $X$, and $\pi: X\dashrightarrow Z$ a dominant map. We say that the general fibers of $\pi$ are \emph{tangent to $\Ff$} if for any general closed point $x$ on a general fiber $F$ of $\pi$, the linear subspace $\Ff_x\subset T_{X,x}$ determined by the inclusion $\Ff\subset T_X$ contains $T_{F,x}$.
\end{defn}

The following results are useful when applying the canonical bundle formula of foliations.

\begin{lem}\label{lem: gen fiber tangent mean induce}
Let $X$ be a normal variety, $\Ff$ a foliation on $X$, and $f: X\rightarrow Z$ a contraction. Suppose that the general fibers of $f$ are tangent to $\Ff$. Then there exists a foliation $\Ff_Z$ on $Z$, such that $\Ff=f^{-1}\Ff_Z$.
\end{lem}
\begin{proof}
Possibly replacing $X$ with a higher model, we may assume that $f$ is a morphism. By Definition-Lemma \ref{defthm: weak ss reduction}, there exists an equidimensional model$/Z$ of $f: X\rightarrow Z$, $f': (X',\Sigma_{X'},\Mm)\rightarrow (Z',\Sigma_{Z'})$, associated with $h: X'\rightarrow X$ and $h_Z: Z'\rightarrow Z$. By \cite[Lemma 6.7]{AD13}, there exists a foliation $\Ff_{Z'}$ on $Z'$ such that $(f')^{-1}\Ff_{Z'}=h^{-1}\Ff$. We may let $\Ff_Z:=(h_Z)_*\Ff_{Z'}$.
\end{proof}

\begin{lem}\label{lem: stein induce same foliation}
Let $\pi: X\rightarrow Z$ be a projective surjective morphism from a normal variety to a variety and $X\xrightarrow{f}Y\xrightarrow{\tau}Z$ the Stein factorization of $\pi$. Let $\Ff$ be the foliation on $X$ induced by $\pi$. Then $\Ff$ is also induced by $f$.
\end{lem}
\begin{proof}
Let $\Ff_Z$ be the foliation by points on $Z$. Then $\Ff_Y:=\tau^{-1}\Ff_Z$ is the foliation by points on $Y$. Since
$$\Ff=(\tau\circ f)^{-1}\Ff_Z=f^{-1}\Ff_Y,$$
 $\Ff$ is induced by $f$.
\end{proof}

\begin{defn}[Restricted foliation]\label{defn: restricted foliation}
Let $X$ be a normal variety, $\Ff$ an algebraically integrable foliation on $X$, $S$ a prime divisor on $X$, and $\nu: S^\nu\rightarrow S$ the normalization of $S$. The \emph{restricted foliation} of $\Ff$ on $S^\nu$ is defined in the following way. 
\begin{enumerate}
  \item If $S$ is $\Ff$-invariant, then we let $U\subset X$ be the largest open subset that does not contain $\Sing(\Ff)\cup\Sing(X)\cup\Sing(S)$ and let $S':=S\cap U$. The natural inclusion of sheaves
  $$\Ff|_{S'}\rightarrow T_X|_{S'}$$
  factors through $T_{S'}$ over $U$, which defines a foliation $\Ff_{S'}$ on $S'$. $\Ff_{S'}$ extends to a foliation $\Ff_S$ on $S^\nu$, and we call $\Ff_S$ the \emph{restricted foliation} of $\Ff$ on $S^\nu$.
  \item If $S$ is not $\Ff$-invariant, then we let $U\subset X$ be the largest open subset that does not contain $\Sing(\Ff)\cup\Sing(X)\cup\Sing(S)$ and where $S$ is transverse to $\Ff$ everywhere in $U$. We let $S':=S\cap U$. The natural inclusion of sheaves $$\Ff|_{S'}\rightarrow T_X|_{S'}$$ induces an inclusion of sheaves $\Ff|_{S'}\cap T_{S'}\rightarrow T_{S'}$. Since $\Ff$ is saturated in $T_X$, $\Ff|_{S'}\cap T_{S'}$ is saturated in $T_{S'}$. Since $\Ff$ is closed under the Lie bracket, $\Ff|_{S'}\cap T_{S'}\subset\Ff$ is closed under the Lie bracket. Thus $\Ff_{S'}:=\Ff|_{S'}\cap T_{S'}$ is a foliation on $S'$. $\Ff_{S'}$ extends to a foliation $\Ff_S$ on $S^\nu$, and we call $\Ff_S$ the \emph{restricted foliation} of $\Ff$ on $S^\nu$.
\end{enumerate}
\end{defn}

Recall the following proposition. %which shows that algebraic integrability is preserved under adjunction:

\begin{prop}[cf. {\cite[Proposition 3.2]{DLM23}}]\label{prop: a.i preserved adjunction}
Let $\Ff$ be an algebraically integrable foliation, $S$ a prime divisor on $X$, and $S^\nu\rightarrow S$ the normalization of $S$. Let $\Ff_S$ be the restriction of the foliation $\Ff$ on $S^\nu$. Then $\Ff_S$ is algebraically integrable and $\rk\Ff_S=\rk\Ff-\epsilon_{\Ff}(S)$.
\end{prop}

\subsection{Foliated log resolution and adjunction formula}

\begin{defn}[{cf. \cite[\S 3.2]{ACSS21}}]\label{defn: foliated log smooth}
Let $(X,\Ff,B,\Mm)/U$ be a sub-gfq such that $\Ff$ is algebraically integrable. We say that $(X,\Ff,B,\Mm)$ is \emph{foliated log smooth} if there exists a contraction $f: X\rightarrow Z$ satisfying the following.
\begin{enumerate}

  \item $X$ has at most quotient toric singularities.
  \item $\Ff$ is induced by $f$.
  \item $(X,\Sigma_X)$ is toroidal for some reduced divisor $\Sigma_X$ such that $\Supp B\subset\Sigma_X$. In particular, $(X,\Supp B)$ is toroidal, and $X$ is $\Qq$-factorial klt.
  \item There exists a log smooth pair $(Z,B_Z)$ such that $$f: (X,\Supp B,\Mm)\rightarrow (Z,\Supp B_Z)$$ is an equidimensional toroidal contraction. %In particular, $(X,\Supp B,\Mm)$ is toroidal, $\Mm$ descends to $X$, and $X$ is $\Qq$-factorial klt.
  \item $\Mm$ descends to $X$.
\end{enumerate}
We may say that the contraction $f: X\rightarrow Z$ is \emph{associated to} $(X,\Ff,B,\Mm)$. It is important to remark that $f$ may not be a contraction$/U$. In particular, $\Mm$ may not be nef$/Z$.
\end{defn}

\begin{lem}[cf. {\cite[Lemma 3.1]{ACSS21}}]\label{lem: foliated log smooth imply lc}
Let $(X,\Ff,B,\Mm)$ be a sub-gfq such that $\Ff$ is algebraically integrable and $(X,\Ff,B,\Mm)$ is foliated log smooth. Then $(X,\Ff,B^\Ff,\Mm)$ is lc.
\end{lem}
\begin{proof}
By \cite[Lemma 3.1]{ACSS21}, $(X,\Ff,B^\Ff)$ is lc. Since $\Mm$ descends to $X$, $(X,\Ff,B^\Ff,\Mm)$ is lc.
\end{proof}

\begin{defn}\label{defn: log resolution}
Let $X\rightarrow U$ be a projective morphism from a normal quasi-projective variety to a variety, $B$ an $\Rr$-divisor on $X$, and $\Ff$ an algebraically integrable foliation on $X$. A \emph{foliated log resolution} of $(X,\Ff,B,\Mm)$ is a projective birational morphism $h: X'\rightarrow X$ such that 
$$(X',\Ff':=h^{-1}\Ff,B':=h^{-1}_*B+\Exc(h),\Mm)$$ 
is foliated log smooth, where $\Exc(h)$ is the reduced $h$-exceptional divisor.
\end{defn}
It is clear that foliated log resolutions exist. Indeed, let $f: X\dashrightarrow Z$ be a dominant map that induces $\Ff$. Let $g: Y\rightarrow X$ be a birational morphism such that $f\circ g$ is a morphism. By Lemma \ref{lem: stein induce same foliation}, we may assume that $f\circ g$ is a contraction. Let $B_Y$ be the reduced divisor supported on \(\Supp(g^{-1}_*B)\cup\Supp\Exc(g)\). By Definition-Theorem \ref{defthm: weak ss reduction}, there exists an equidimensional model$/U$ $f': (X',\Sigma_{X'},\Mm)\rightarrow (Z',\Sigma_{Z'})$ of $f\circ g: (Y,B_Y,\Mm)\rightarrow Z$ associated with $h': X'\rightarrow Y$ and $h_Z: Z'\rightarrow Z$. Let $h:=g\circ h'$, $\Ff':=h^{-1}\Ff$, and $B':=(h')_*^{-1}B_Y+\Exc(h')$.
Then $h$ is a foliated log resolution of $(X,\Ff,B,\Mm)$.

Next we state the adjunction formula for algebraically integrable generalized foliated quadruples. The precise adjunction formula, i.e., the adjunction formula with coefficient control, will be provided later. In particular, we can only prove the sub-lc for the time being.

\begin{thm}\label{thm: not precise adjunction}
    Let $(X,\Ff,B,\Mm)/U$ be an lc gfq such that $\Ff$ is algebraically integrable. Let $S$ be a prime divisor on $X$ such that $\mult_SB=\epsilon_{\Ff}(S)$, and let $\nu: S^\nu\rightarrow S$ be the normalization of $S$, $\Mm^S:=\Mm|_S$, $\Ff_S$ the restricted foliation of $\Ff$ on $S^\nu$, and
    $$K_{\Ff_S}+B_S+\Mm^S_{S^\nu}:=(K_X+B+\Mm_X)|_{S^\nu}.$$
    Then $(S^\nu,\Ff_S,B_S,\Mm^S)$ is sub-lc.
\end{thm}
\begin{proof}
Let $h:X'\to X$ be a foliated log resolution of $(X,\Ff,B+S,\Mm)$ with associated $f':X'\to Z'$. Let $\Ff':=h^{-1}\Ff$ and
$$K_{\Ff'}+{B}'+\Mm_{X'}:=h^*(K_\Ff+B+\Mm_X).$$
By Lemma \ref{lem: foliated log smooth imply lc}, 
$$(X',\Ff':=h^{-1}\Ff,\tilde B':=(B')^{\geq 0},\Mm)$$ is lc. In particular, $(X',\Ff',\tilde B')$ is lc.

We let $S':=h^{-1}_*S$ which is normal. There exists an induced birational morphism $h_S: S'\rightarrow S^\nu$ such that $\nu\circ h_S=h|_{S'}$. Let $\Ff_{S'}$ be the restricted foliation of \(\Ff\) on \(S'\), then \(\Ff_{S'}=h_S^{-1}\Ff_S\). Let
$$K_{\Ff_{S'}}+\tilde B_{S'}:=\left(K_{X'}+\tilde B'\right)|_{S'}\text{ and }K_{\Ff_{S'}}+B_{S'}+\Mm^S_{S'}:=\left(K_{X'}+B'+\Mm_{X'}\right)|_{S'}.$$
By \cite[Proposition 3.2]{ACSS21}, $(S',\Ff_{S'},\tilde B_{S'})$ is lc. Since $\Mm$ descends to $X'$, $\Mm^S$ descends to $S'$, so $(S',\Ff_{S'},\tilde B_{S'},\Mm^S)$ is lc.
Since $\tilde B'\geq B'$, $\tilde B_{S'}\geq B_{S'}$ and thus
$$\left(S',\Ff_{S'},\tilde B_{S'},\Mm^S\right)$$ 
is sub-lc. By construction,
$$K_{\Ff_{S'}}+B_{S'}+\Mm^S_{S'}=h_S^*\left(K_{\Ff_S}+B_S+\Mm^S_{S^\nu}\right),$$
so $(S^\nu,\Ff_S,B_S,\Mm^S)$ is sub-lc. This completes the proof.
\end{proof}

Finally, we recall the following definition of F-dlt.

\begin{defn}[F-dlt]\label{defn: fdlt}
    Let $(X,\Ff,B,\Mm)/U$ be an lc gfq such that $\Ff$ is algebraically integrable. We say that $(X,\Ff,B,\Mm)$ is \emph{F-dlt} if there exists a foliated log resolution $f: Y\rightarrow X$ of $(X,\Ff,B,\Mm)$ such that $a(D,\Ff,B,\Mm)>-\epsilon_{\Ff}(D)$ for any prime $f$-exceptional divisor $D$.
\end{defn}

\subsection{Cutting foliations by general hyperplane sections}

Before we prove the adjunction formulas, we first need to discuss how to cut generalized foliated quadruples with general hyperplane sections. 

\begin{lem}\label{lem: general hyperplane cut keep m}
Let $X$ be a normal quasi-projective variety and $H$ a prime divisor such that $H$ is base-point-free and is a general member of $|H|$. Let $\Mm$ be a $\bb$-divisor on $X$ such that $\Mm$ descends to a birational model $X'$ of $X$ and $\Mm^H:=\Mm|_H$. Then $\Mm^H_H=\Mm_X|_H$.
\end{lem}
\begin{proof}
We may assume that the induced birational map $f: X'\dashrightarrow X$ is a morphism. We let $V:=f(\Supp(\Mm_{X'}-f^{-1}_*\Mm_X))$. Then $\dim X-\dim V\geq 2$. Since $H$ is general, $\dim H-\dim (V\cap H)\geq 2$. Therefore, for any prime divisor $D$ on $H$ and identifying $D$ with its image in $X$, we have that $\Mm$ descends to $X$ near the generic point of $D$. The lemma follows immediately.
\end{proof}

\subsubsection{Cutting by invariant hyperplane sections}

First, we show that we can cut foliations by invariant base-point-free linear systems freely.

\begin{prop}\label{prop: general hyperplane invariant}
Let $(X,\Ff,B,\Mm)/U$ be a sub-gfq and $W$ a proper subvariety of $X$. Suppose that $\Ff$ is induced by a morphism $f: X\rightarrow Z$, $\dim Z>0$, and $W$ is transverse to $\Ff$. Let $H_Z\subset Z$ be a general hyperplane section. Let $H:=f^*H_Z$, $\Mm^H:=\Mm|_H$, and 
$$K_{\Ff_H}+B_H+\Mm_H^H:=(K_{\Ff}+B+\Mm_X)|_H,$$ 
where $\Ff_H$ is the restricted foliation of $\Ff$ on $H$. Then:
\begin{enumerate}
    \item $H$ intersects $W$.
    \item If $(X,\Ff,B,\Mm)$ is (sub-)lc, then $(H,\Ff_H,B_H,\Mm^H)$ is (sub-)lc.
     \item For any component $D$ of $\Supp B$ such that $D$ intersects $H$ and any component $C$ of $D\cap H$, $\mult_CB_H=\mult_DB$.
\end{enumerate}
\end{prop}
\begin{proof}
%We only need to prove the case when $(X,\Ff,B)$ is sub-lc, and the lc case of (2) follows by (3) and the sub-lc case of (2).
By Definition-Theorem \ref{defthm: weak ss reduction}, there exists an equidimensional model$/U$ $f': (X',\Sigma_{X'},\Mm)\rightarrow (Z,\Sigma_Z)$ of $f: (X,B,\Mm)\rightarrow Z$ associated with $h: X'\rightarrow X$ and $h_Z: Z'\rightarrow Z$, such that $h$ is a foliated log resolution of $(X,\Ff,B,\Mm)$ and $\Ff':=h^{-1}\Ff$ is induced by $f'$. We let $H':=h^*H$,
$$K_{\Ff'}+B'+\Mm_{X'}:=h^*(K_{\Ff}+B+\Mm_X),\text{ and }K_{\Ff_{H'}}+B_{H'}+\Mm^H_{H'}:=\left(K_{\Ff'}+B'+\Mm_{X'}\right)|_{H'}.$$ 
Let
$$R(f'):=\sum_{D\mid D\text{ is a prime divisor on }Z'}(f'^*D-f'^{-1}(D))$$
be the ramification divisor of $f'$. Since $H',H_{Z'}$ are general, we have
$$R(f')|_{H'}=R(f')|_{f'^*H_{Z'}}=\sum_{D\mid D\text{ is a prime divisor on }Z'}((f'|_{H'})^*D|_{H_{Z'}}-(f'|_{H'})^{-1}(D|_{H_{Z'}})):=R(f'|_{H'}).$$
Therefore,
$$K_{\Ff'}|_{H'}=(K_{X'/Z'}-R(f'))|_{H'}=K_{H'/H_{Z'}}-R(f'|_{H'})=K_{\Ff_{H'}}.$$
Since $\Mm^H_H=\Mm_X|_H$, we have $B_{H'}=B'|_{H'}$. 

(1) Since $W$ is transverse to $\Ff$, $W':=h^{-1}(W)$ is not tangent to $\Ff'$. Thus $\dim g(W')\geq 1$, so $H_{Z'}:=h_Z^*H_Z$ intersects $g(W')$ and $H_Z$ intersects $h_Z(g(W'))=f(W)$. Hence $H$ intersects $W$.

(3) Since $H$ is general, near the generic point $\eta_C$ of $C$, $h$ is an isomorphism. Since $B_{H'}=B'|_{H'}$, $B|_H=B_H$ near $\eta_C$. We may write $B=\sum b_iB_i$ where $B_i$ are the irreducible components of $B$, then 
$$B_H=B|_H=\sum b_i(B_i\cap H)$$ 
near $\eta_C$. Since $H$ is general, there exists a unique index $i$ such that $B_i\cap H\not=0$ at \(\eta_C\). Then \(B_i\cap H=C\), \(B_i=D\), and hence \(\mult_CB_H=b_i=\mult_DB\).

(2) Assume that $(X,\Ff,B,\Mm)$ is sub-lc. According to Lemma \ref{lem: foliated log smooth imply lc}, $(X',\Ff',\tilde B':=(B')^{\geq 0},\Mm)$ is lc. Let $K_{\Ff_{H'}}+\tilde B_{H'}:=(K_{\Ff'}+\tilde B')|_{H'}$. By \cite[Proposition 3.2]{ACSS21}, $(H',\Ff_{H'},\tilde B_{H'})$ is lc. Since $\Mm$ descends to $X'$, 
$$K_{\Ff_{H'}}+\tilde B_{H'}+\Mm^H_{H'}=(K_{\Ff'}+B+\Mm_X)|_H$$
and $\Mm^H$ descends to $H'$. Thus $(H',\Ff_{H'},\tilde B_{H'},\Mm^H)$ is lc. Since $\tilde B'\geq B'$, $\tilde B_{H'}\geq B_{H'}$, and hence $(H',\Ff_{H'},B_{H'},\Mm^H)$ is sub-lc. Since $H$ is general, one can see that
$$K_{\Ff_{H'}}+B_{H'}+\Mm^H_{H'}=(h|_{H'})^*\left(K_{\Ff_H}+B_H+\Mm^H_H\right),$$
which implies that $(H,\Ff_H,B_H,\Mm^H)$ is sub-lc. If $(X,\Ff,B,\Mm)$ is lc, then $(H,\Ff_H,B_H,\Mm^H)$ is lc by (3).
\end{proof}

\subsubsection{Cutting by non-invariant hyperplanes}

Next we show that, if we only consider the local property of foliations, then we can cut the foliation by non-invariant hyperplane sections. We first prove a lemma.

\begin{lem}\label{lem: toroidal cut general hyperplane still lc}
    Let $f: (X,\Sigma,\Mm)\rightarrow (Z,\Sigma_Z)$ be a toroidal morphism and $z\in Z$ a closed point. Let $\Ff$ be the foliation induced by $f$ and let $B$ be the horizontal$/Z$ part of $\Sigma$. Let $H$ be a general member of a base-point-free linear system on $X$, such that $H$ dominates $Z$. Then $(X,\Ff,B+H,\Mm)$ is lc over a neighborhood of $z$.
\end{lem}
\begin{proof}
By \cite[Lemma 3.6]{DLM23}, $(X,\Ff,B+H)$ is lc over a neighborhood of $z$. Since $\Mm$ descends to $X$, $(X,\Ff,B+H,\Mm)$ is lc over a neighborhood of $z$.
\end{proof}

\begin{prop}\label{prop: general hyperplane non-invariant}
Let $(X,\Ff,B,\Mm)$ be a sub-gfq and $W$ a proper subvariety of $X$. Suppose that $\Ff$ is algebraically integrable, $W$ is tangent to $\Ff$, and $\dim W\geq 1$. Let $H\subset X$ be a general hyperplane section. Let $\Mm^H:=\Mm|_H$ and $$K_{\Ff_H}+B_H+\Mm^H_H:=(K_{\Ff}+B+H+\Mm_X)|_H,$$ where $\Ff_H$ is the restricted foliation of \(\Ff\) on \(H\). Then:
\begin{enumerate}
    \item $H$ intersects $W$.
    \item If $(X,\Ff,B,\Mm)$ is (sub-)lc, then $(H,\Ff_H,B_H,\Mm^H)$ is (sub-)lc near $W|_H$.
    \item For any component $D$ of $\Supp B$ such that $D$ intersects $H$ and any component $C$ of $D\cap H$, $\mult_CB_H=\mult_DB$.
\end{enumerate}
\end{prop}
\begin{proof}
(1) is obvious. By \cite[Proposition 3.6]{Dru21}, $K_{\Ff_H}=(K_{\Ff}+H)|_H$, so $B_H+\Mm^H_H=B|_H+\Mm_X|_H$. We remark that \cite[Proposition 3.6]{Dru21} requires that $\rk\Ff\geq 2$, but the same lines of the proof work for the case when $\rk\Ff=1$ as well. By Lemma \ref{lem: general hyperplane cut keep m}, $B_H=B|_H$ and then (3) follows.

(2) By Definition-Theorem \ref{defthm: weak ss reduction}, there exists an equidimensional model$/U$ $f': (X',\Sigma_{X'},\Mm)\rightarrow (Z,\Sigma_Z)$ of $f: (X,B,\Mm)\rightarrow Z$ associated with $h: X'\rightarrow X$ and $h_Z: Z'\rightarrow Z$, such that $h$ is a foliated log resolution of $(X,\Ff,B,\Mm)$ and $\Ff':=h^{-1}\Ff$ is induced by $f'$. We let $$K_{\Ff'}+B'+\Mm_{X'}:=h^*(K_{\Ff}+B+\Mm_X),$$ $H':=h^*H$, $W':=h^{-1}(W)$, and $\tilde B':=(B')^{\geq 0}$. We let $z$ be the image of $W'$ on $Z'$. Since $(X,\Ff,B,\Mm)$ is lc, by Lemma \ref{lem: foliated log smooth imply lc}, $(X',\Ff',\tilde B',\Mm)$ is lc. Moreover, all components of $\tilde B'$ are horizontal$/Z$. By Lemma \ref{lem: toroidal cut general hyperplane still lc}, $(X',\Ff',\tilde B'+H',\Mm)$ is lc over a neighborhood of $z'$. In particular, $(X',\Ff',\tilde B',\Mm)$ is lc near $W'|_{H'}$.

Let
$$K_{\Ff_{H'}}+\tilde B_{H'}:=(K_{\Ff'}+\tilde B')|_{H'}\text{ and }K_{\Ff_{H'}}+B_{H'}:=(K_{\Ff'}+B')|_{H'}.$$ 
By \cite[Proposition 3.2]{ACSS21}, $(H',\Ff_{H'},\tilde B_{H'})$ is lc near $W'|_{H'}$. Since $\tilde B'\geq B'$, $\tilde B_{H'}\geq B_{H'}$, so $(H',\Ff_{H'},B_{H'})$ is sub-lc near $W'|_{H'}$. Since $\Mm$ descends to $X'$, $\Mm^H$ descends to $H'$, so $(H',\Ff_{H'},B_{H'},\Mm^H)$ is sub-lc near $W'|_{H'}$. Since $H$ is general, $$K_{\Ff_{H'}}+B_{H'}+\Mm^H_{H'}=h|_{H'}^*(K_{\Ff_H}+B_H+\Mm^H_H),$$ 
so $(H,\Ff_H,B_H,\Mm^H)$ is sub-lc near $W|_H$. Finally, if $(X,\Ff,B,\Mm)$ is lc, then by (3), $B_H\geq 0$, so $(H,\Ff_H,B_H,\Mm^H)$ is lc near $W|_H$.
\end{proof}

\subsection{Basic properties of foliated surfaces}
In this subsection, we recall some basic properties of foliated surfaces. Moreover, we introduce the concept of \emph{surface numerical gfqs} and study their basic properties. This is crucial for the proof of adjunction formulas.

\begin{defn}
Let $X$ be a normal surface, $\Ff$ a foliation on $X$, and $x\in X$ a closed point such that $x\not\in\Sing(X)$ and $x\in\Sing(\Ff)$. Let $v$ be a vector field generating $\Ff$ near $x$. By \cite[Page 2, Line 17-18]{Bru15}, $v(x)=0$ and $(Dv)|_x$ has exactly two eigenvalues $\lambda_1$ and $\lambda_2$. 

We say that $x$ is a \emph{reduced singularity} of $\Ff$ if at least one of $\lambda_1$ and $\lambda_2$ is not $0$ (say, $\lambda_2$) and $\frac{\lambda_1}{\lambda_2}\not\in\mathbb Q^+$. We say that $\Ff$ has \emph{at most reduced singularities} if for any closed point $p\in X$, $\Ff$ is either non-singular at $p$ or $p$ is a reduced singularity of $\Ff$. 
\end{defn}

\begin{defn}[Minimal resolution]\label{defn: minimal resolution foliation}
Let $X$ be a normal surface, $\Ff$ a foliation on $X$, $f: Y\rightarrow X$ a projective birational morphism, and $\Ff_Y:=f^{-1}\Ff$.

We say that $f$ is a \emph{resolution} of $\Ff$ if $Y$ is smooth and $\Ff_Y$ has at most reduced singularities. By \cite{Sei68} (we refer to \cite[Pages 908--912]{Can04} for a detailed explanation), a resolution of $\Ff$ always exists.

We say that $f$ is the \emph{minimal resolution} of $\Ff$ if for any resolution $g: W\rightarrow X$ of $\Ff$, $g$ factors through $f$, i.e. there exists a projective birational morphism $h: W\rightarrow Y$ such that $g=f\circ h$. By definition, the minimal resolution of $\Ff$ is unique, and by \cite[Proposition 1.17]{Che23}, the minimal resolution of $\Ff$ exists.
\end{defn}

\begin{defn}
Let $X$ be a normal surface with at most cyclic quotient singularities, $\Ff$ a foliation on $X$, and $C$ a reduced curve on $X$ such that no component of $C$ is $\Ff$-invariant. For any closed point $x\in X$, we define $\tang(\Ff,C,x)$ in the following way. 
\begin{itemize}
\item If $x\notin \Sing(X)$, then we let $v$ be a vector field generating $\Ff$ around $x$, and $f$ a holomorphic function defining $C$ around $x$. We define 
$$\tang(\Ff,C,x):=\dim_{\mathbb{C}}\frac{\mathcal{O}_{X,x}}{\langle f, v(f)\rangle}.$$
\item If $x\in\Sing(X)$, then $x$ is a cyclic quotient singularity of index $r$ for some integer $r\geq 2$. Let $\rho:\tilde X\rightarrow X$ be an index $1$ cover of $X\ni x$, $\tilde x:=\rho^{-1}(x)$, $\widetilde C:=\rho^*C$, and $\tilde\Ff$ the foliation induced by the sheaf $\rho^*\Ff$ near $\tilde x$. Then $\tilde x$ is a smooth point of $\tilde X$, and we define
$$\tang(\Ff,C,x):=\frac{1}{r}\tang(\tilde\Ff,\tilde C,\tilde x).$$
\end{itemize}
We define
$$\tang(\Ff,C):=\sum_{x\in X}\tang(\Ff,C,x)$$
which is well-defined according to \cite[Section 2]{Bru02}.
\end{defn}

\begin{defn}
Let $X$ be a normal surface with at most cyclic quotient singularities, $\Ff$ a foliation on $X$, and $C$ a reduced curve on $X$ such that all components of $C$ are $\Ff$-invariant. For any closed point $x\in X$, we define $Z(\Ff,C,x)$ in the following way.
\begin{itemize}
\item If $x\notin \Sing(X)$, then we let $\omega$ be a $1$-form generating $\Ff$ around $x$, and $f$ a holomorphic function generating $C$ around $x$. Then there are uniquely determined holomorphic functions $g,h$ and a holomorphic $1$-form $\eta$ on $X$ near $x$, such that $g\omega=hdf+f\eta$ and $f,h$ are coprime. We define 
$$Z(\Ff,C,x):= {\rm ord}_x(\frac{h}{g}|_C)$$
the vanishing order of $\frac{h}{g}|_C$ at $x$. By \cite[Chapter 2, Page 15]{Bru15}, $Z(\Ff,C,x)$ is independent of the choice of $\omega$.
\item If $x\in C\cap \Sing(X)$, we define 
$Z(\Ff,C,x):=0.$
\end{itemize}
We define $$Z(\Ff,C):=\sum_{x\in C}Z(\Ff,C,x)$$
which is well-defined according to \cite[Section 2]{Bru02}.
\end{defn}

\begin{defn}[Dual graph]\label{defn: dual graph}
Let $C=\cup_{i=1}^nC_i$ be a collection of irreducible curves contained in the non-singular locus of a normal surface $X$. We define the \emph{dual graph} $\mathcal{D}(C)$ of $C$ as follows.
\begin{enumerate}
    \item The vertices $v_i=v_i(C_i)$ of $\mathcal{D}(C)$ correspond to the curves $C_i$.
    \item For $i\neq j$, the vertices $v_i$ and $v_j$ are connected by $C_i\cdot C_j$ edges.
    \item Each vertex $v_i$ is labeled by $w(C_i):=-C_i^2$. The integer $w(C_i)$ is called the \emph{weight} of $C_i$.
\end{enumerate}
%We sometimes write the name of the curve $C_i$ near the vertex $v_i$.

For any projective birational morphism $f: Y\rightarrow X$ between surfaces, let $E=\cup_{i=1}^nE_i$ be the reduced exceptional divisor. If $E\nsubseteq\Sing Y$, then we define $\mathcal{D}(f):=\mathcal{D}(E)$.
\end{defn}

\begin{defn}\label{defn: num foliated sing}
A \emph{surface numerical sub-gfq} (\emph{surface num-sub-gfq} for short) $(X,\Ff,B,\Mm)/U$ consists of a normal surface $X$, a rank $1$ foliation $\Ff$ on $X$, an $\Rr$-divisor $B$ on $X$, and a nef$/U$ $\bb$-divisor $\Mm$. We say that $(X,\Ff,B,\Mm)$ is a \emph{surface numerical gfq} (\emph{surface num-gfq} for short) if $(X,\Ff,B)$ is a surface num-sub-gfq and $B\geq 0$. A \emph{surface num-gfq germ} $(X\ni x,\Ff,B,\Mm)$ consists of a surface num-gfq $(X,\Ff,B,\Mm)/X$ and a closed point $x\in X$.

Let $(X,\Ff,B,\Mm)$ be a surface num-sub-gfq. Let $f: Y\rightarrow X$ be a resolution of $X$ with prime $f$-exceptional divisors $E_1,\dots,E_n$. Since $\{(E_i\cdot E_j)\}_{n\times n}$ is negative definite, the equation 
$$\begin{pmatrix}
  (E_1\cdot E_1) &\cdots & (E_1\cdot E_n) \\
  \vdots&\ddots & \vdots         \\
  (E_n\cdot E_1) &\cdots & (E_n\cdot E_n) \\
\end{pmatrix} 
\begin{pmatrix}
  a_1 \\
  \vdots \\
  a_n \\
\end{pmatrix} 
 =
 \begin{pmatrix}
-(K_{\Ff_Y}+B_Y+\Mm_Y)\cdot E_1\\
\vdots\\
-(K_{\Ff_Y}+B_Y+\Mm_Y)\cdot E_n\\
\end{pmatrix}$$
has a unique solution $(a_1,\dots,a_n)$, where $\Ff_Y:=f^{-1}\Ff$ and $B_Y:=f^{-1}_*B$. For any prime divisor $E$ on $Y$, we define 
$$a_{\num,f}(E,\Ff,B,\Mm):=-\mult_E\left(B_Y+\sum a_iE_i\right).$$
\end{defn}

\begin{lem}\label{lem: anum same as a}
Let $(X,\Ff,B,\Mm)$ be a sub-gfq such that $\dim X=2$ and $\rk\Ff=1$. Let $f: Y\rightarrow X$ be a resolution of $X$ and $E$ a prime divisor on $Y$. Then $a_{\num,f}(E,\Ff,B,\Mm)=a(E,\Ff,B,\Mm)$.
\end{lem}
\begin{proof}
If $f(E)$ is a prime divisor, then $a_{\num,f}(E,\Ff,B,\Mm)=-\mult_EB=a(E,\Ff,B,\Mm)$. Assume that $E$ is exceptional over $X$. Write
$$K_{\Ff_Y}+\sum a_iE_i+B_Y+\Mm_Y=f^*(K_{\Ff}+B+\Mm_X),$$
where $E_i$ are the prime $f$-exceptional divisors, $\Ff_Y:=f^{-1}\Ff$ and $B_Y:=f^{-1}_*B$. Then
$$a_{\num,f}(E_i,\Ff,B,\Mm)=-a_i=a(E_i,\Ff,B,\Mm)$$ for any $i$. Since $E=E_j$ for some $j$, the lemma holds.
\end{proof}

\begin{lem}\label{lem: anum not depend on resolution}
Let $(X,\Ff,B,\Mm)$ be a surface num-sub-gfq and $f: Y\rightarrow X$, $f': Y'\rightarrow X$ two resolutions of $X$ such that both $\Center_YE$ and $\Center_{Y'}E$ are divisors. Then
$a_{\num,f}(E,\Ff,B,\Mm)=a_{\num,f'}(E,\Ff,B,\Mm).$
\end{lem}
\begin{proof}
%If $\Center_XE$ is a prime divisor, then the lemma is trivial.%$$a_{\num,f}(E,\Ff,B,\Mm)=-\mult_EB=a_{\num,f'}(E,\Ff,B,\Mm),$$ so we may assume that $E$ is exceptional over $X$.
%Let $g: W\rightarrow Y$ and $g': W\rightarrow Y'$ be a common resolution, and $h: W\rightarrow X$ the induced birational morphism. Possibly replacing $f'$ with $h$, 
We may assume that there exists a morphism $g: Y'\rightarrow Y$. Let $E_i$ be the prime $f'$-exceptional divisors, $B_{Y'}:=f'^{-1}_*B-\sum_ia_{\num,f'}(E_i,\Ff,B,\Mm)E_i$, and $B_Y:=g_*B_{Y'}$.
Then $(K_{\Ff_{Y'}}+B_{Y'}+\Mm_{Y'})\cdot E_i=0$ for any $E_i$. Since $Y$ is smooth, $K_{\Ff_Y}+B_Y+\Mm_Y$ is $\Rr$-Cartier. By the negativity lemma, we see that $K_{\Ff_{Y'}}+B_{Y'}+\Mm_{Y'}=g^*(K_{\Ff_Y}+B_Y+\Mm_Y)$. Thus $(K_{\Ff_Y}+B_Y+\Mm_Y)\cdot g_*E_i=0$ for any $E_i$, so 
$$a_{\num,f}(E_i,\Ff,B,\Mm)=-\mult_{g_*E_i}B_Y=\mult_{E_i}B_{Y'}=a_{\num,f'}(E_i,\Ff,B,\Mm)$$
for any $E_i$ such that $g_*E_i\not=0$. In particular, $a_{\num,f}(E,\Ff,B,\Mm)=a_{\num,f'}(E,\Ff,B,\Mm)$.
\end{proof}

\begin{defn}
Let $(X,\Ff,B)$ be a surface num-sub-gfq. We define $a(E,\Ff,B,\Mm):=a_{\num,f}(E,\Ff,B,\Mm)$ for an arbitrary resolution $f: Y\rightarrow X$ of $X$ such that $E$ is a divisor on $Y$. Lemmas \ref{lem: anum same as a} and \ref{lem: anum not depend on resolution} guarantee that there is no abuse of notations.

Let $(X,\Ff,B,\Mm)$ be a surface num-gfq. We say that $(X,\Ff,B,\Mm)$ is \emph{num-lc} if $a(E,\Ff,B,\Mm)\geq-\epsilon_{\Ff}(E)$ for any prime divisor $E$ over $X$.
\end{defn}

\begin{lem}\label{lem: add component worse sing}
Let $(X,\Ff,B,\Mm)$ be a surface num-gfq and $x\in X$ a closed point. Then for any prime divisor $E$ over $X\ni x$, $a(E,\Ff,B,\Mm)\leq a(E,\Ff,B)$, and $a(E,\Ff,B,\Mm)=a(E,\Ff,B)$ if and only if $x\not\in\Supp B$ and $\Mm$ descends to $X$ over a neighborhood of $x$. In particular, if $(X,\Ff,B,\Mm)$ is num-lc, then $(X,\Ff,B)$ is num-lc.
\end{lem}
\begin{proof} 
It follows from \cite[Lemma 3.41]{KM98}.
\end{proof}

\subsection{Adjunction formula for surface generalized foliated quadruples}\label{subsec: adj gfq surface}

\begin{lem}\label{lem:tersurf}
Suppose that $(X,\Ff,B,\Mm)$ is a gfq such that $\dim X=2$ and $\rk\Ff=1$. Let $C$ be an $\Ff$-invariant curve with normalization $C^\nu$, and $x\in C$ a closed point such that $\Ff$ is terminal near $x$. Let $\Mm^C:=\Mm|_{C^\nu}$. Then near $x$, $X$ is $\Qq$-factorial klt and $C$ is non-singular, and $\mult_x(\Mm_X|_{C^\nu}-\Mm^C_{C^\nu})\ge0$. 
\end{lem}
\begin{proof}
By \cite[Theorem 3.19]{LMX23a}, $X$ is $\Qq$-factorial klt and $C$ is non-singular near $x$. Possibly shrinking $X$ near $x$, we may assume that $X$ is $\Qq$-factorial klt, and $C$ is non-singular. Let $f: Y\rightarrow X$ be a birational morphism such that $\Mm$ descends to $Y$. By the negativity lemma, $f^*\Mm_X-\Mm_Y\geq 0$. This implies that $\mu=\mult_x((h|_{C_Y})_*(h^*\Mm_X-\Mm_Y)|_{C_Y})\geq 0$, where $C_Y:=f^{-1}_*C$.
\end{proof}

\begin{lem}\label{lem: surface pia terminal1}
Notation and assumptions as in Lemma \ref{lem:tersurf}. Suppose that $B=\sum_{j=1}^m b_j B_j$ and $B_j$ are the irreducible components of $B$. Then there exist non-negative integers $w_1,\dots,w_m$ satisfying the following.

For any real numbers $b_1',\dots,b_m'$, the vanishing order
$${\rm ord}_x\left(\left(K_{\Ff}+\sum b_j'B_j+\Mm_X\right)|_{C^\nu}-\Mm^C_{C^\nu}\right)=\frac{I-1+\sum_{j=1}^m w_j b_j'}{I}+\mu,$$
where $\Mm^C:=\Mm|_{C^\nu}$, $I$ is the order of the local fundamental group $\pi_1(X\ni x)$, and $\mu:=\mult_x(\Mm_X|_{C^\nu}-\Mm^C_{C^\nu})$.

Moreover, if $(X,\Ff,\sum b_j'B_j,\Mm)$ is lc, then $0\leq \frac{I-1+\sum_{j=1}^m w_j b_j'}{I}+\mu\leq 1.$
    %We note that $\Mm_X|_{C^\nu}$ in (3), $$\left(K_{\Ff}+\sum_{j=1}^mb_j'B_j+\Mm_X\right)\Biggm|_{C^\nu}$$ in (4), and $$\left(K_{\Ff}+\sum_{j=1}^mb_j'B_j+\sum_{k=1}^nr_k'\Mm_{k,X}\right)\Biggm|_{C^\nu}$$ in (6) may not be well-defined, but they are at least well-defined near $x$ so there is no confusion for the statements of the lemma.
\end{lem}
\begin{proof}
By Lemma \ref{lem:tersurf}, possibly shrinking $X$ near $x$, we may assume that $X$ is $\Qq$-factorial klt, $\Ff$ is terminal, and $C$ is non-singular.% In particular, we will identify $C$ with $C^\nu$ in the following arguments.
Then the lemma follows from \cite[Theorem 3.2]{LMX23b}.

Moreover, if $(X,\Ff,\sum_{j=1}^m b_j' B_j,\Mm)$ is lc, then $\mu\geq 0$ by Lemma \ref{lem:tersurf}. Therefore $\frac{I-1+\sum_{j=1}^m w_j b_j'}{I}+\mu\geq 0$. The inequality $\frac{I-1+\sum_{j=1}^m w_j b_j'}{I}+\mu\le1$ follows from \ref{thm: not precise adjunction}.
\end{proof}

\begin{lem}\label{lem: surface pia terminal}
Notation and assumptions as in Lemma \ref{lem:tersurf}. Suppose that 
\begin{itemize}
  \item $B=\sum_{j=1}^m b_j B_j$, where $B_j$ are the irreducible components of $B$, and 
  \item $\Mm=\sum_{k=1}^n r_k \Mm_k$, where $\Mm_k$ are nef$/X$ $\bb$-Cartier $\bb$-divisors. 
\end{itemize} 
Then there exist non-negative integers $w_1,\dots,w_m$, and $v_1,\dots,v_n$ satisfying the following. 

For any real numbers $b_1',\dots,b_m',$ and $r_1',\dots,r_n'$, the vanishing order
$${\rm ord}_x\left(\left(K_{\Ff}+\sum_j b_j' B_j+\sum_k r_k' \Mm_{i,X}\right)|_{C^\nu}-\sum_{k=1}^n r_k' \Mm^C_{k,C^\nu}\right)=\frac{I-1+\sum_{j=1}^m w_j b_j'+\sum_{k=1}^n v_k r_k'}{I}.$$
where $\Mm^C_k:=\Mm_k|_{C^\nu}$ for each $k$ and $I$ is the order of the local fundamental group $\pi_1(X\ni x)$.

Moreover, if $(X,\Ff,\sum_{j=1}^m b_j' B_j,\sum_{k=1}^n r_k' \Mm_k)$ is lc, then
$0\leq\frac{I-1+\sum_{j=1}^m w_j b_j'+\sum_{k=1}^n v_k r_k'}{I}\leq 1.$
\end{lem}

\begin{proof}
By Lemma \ref{lem:tersurf}, possibly shrinking $X$ near $x$, we may assume that $X$ is $\Qq$-factorial klt, $\Ff$ is terminal, and $C=C^\nu$ is non-singular.

Let $f:Y\to X$ be a resolution such that $\Mm_k$ descends to $Y$ for each $k$. By the negativity lemma, $h^*\Mm_{k,X}-\Mm_{k,Y}\geq 0$ for each $k$.

Since $\Mm_{k,X}$ is integral for each $k$, $I\mult_x \Mm_{k,X}|_C$ is an integer. In particular,
 $$v_k:=I\left(\mult_x \Mm_{k,X}|_C-\mult_x \Mm_{k,C}^C\right)=I\mult_x\left(\left(h|_{C_Y})_*(h^*\Mm_{k,X}-\Mm_{k,Y}\right)|_{C_Y}\right)\geq 0$$
is a non-negative integer for each $k$. Then the lemma follows from Lemma \ref{lem: surface pia terminal}. 

The moreover part is clear.
%By Lemma \ref{lem: surface pia terminal}, the vanishing order $${\rm ord}_x\left(\left(K_{\Ff}+\sum_{j=1}^mb_j'B_j+\sum_{k=1}^nr_k'\Mm_{k,X}\right)\Biggm|_{C}-\sum_{k=1}^nr_k'\Mm_{k,C}^C\right)=\frac{I-1+\sum_{j=1}^mw_jb_j'+\sum_{k=1}^nv_kr_k'}{I}.$$ If $(X,\Ff,\sum_{j=1}^mb_j'B_j,\sum_{k=1}^nr_k'\Mm_{k,X})$ is lc, then $q\geq 0$. By Theorem \ref{thm: not precise adjunction}, $q\leq 1$. (5) follows.
\end{proof}

\begin{lem}\label{lem: surface pia not terminal}
     Let $(X,\Ff,B,\Mm)$ be an lc gfq such that $\dim X=2$, $\rk\Ff=1$, and $B_j$ are the irreducible components of $B$. Let $C$ be an $\Ff$-invariant curve with normalization $\nu: C^\nu\rightarrow C$. Let $x\in C$ be a closed point, such that $\Ff$ is not terminal near $x$. Then:
     \begin{enumerate}
         \item $x\not\in\Supp B$ and $\Mm$ descends to $X$ over a neighborhood of $x$.
         \item For any closed point $y\in\nu^{-1}(x)$, the vanishing order of $K_{\Ff}|_{C^\nu}$ at $y$ is a non-negative integer.
     \end{enumerate}
\end{lem}
\begin{proof}
  (1) Since $(X,\Ff,B,\Mm)$ is lc, by Lemma \ref{lem: anum same as a}, $(X,\Ff,B,\Mm)$ is num-lc. By Lemma \ref{lem: add component worse sing}, $(X,\Ff,B)$ is num-lc near $x$. Since $\Ff$ is not terminal near $x$, by \cite[Theorem 3.19]{LMX23a}, $B=0$ near $x$. By Lemma \ref{lem: add component worse sing} again, $\Mm$ descends to $X$ over a neighborhood of $x$. 

  (2) By considering a local analytic neighborhood of $x$ and separate $C$ into different analytic irreducible components, we may assume that $y=\nu^{-1}(x)$. (2) follows from \cite[Theorem 3.19]{LMX23a}. More precisely, we let $h: Y\rightarrow X$ be the minimal resolution of $\Ff$ near $x$ and let $C_Y:=h^{-1}_*C$, then we only need to show that
  $$K_{\Ff}\cdot C-K_{C^\nu}=h^*K_{\Ff}\cdot C_Y-K_{C_Y}$$
  is a positive integer over a neighborhood of $x$. (2) follows by checking all cases of \cite[Theorem 3.19]{LMX23a} and applying \cite[Proposition 2.16(3)]{CS20} to $C_Y$ for each case.
\end{proof}

\subsection{Precise adjunction formula when the foliation is induced by a morphism}

\begin{thm}\label{thm: adjunction foliation nonnqc}
Under the assumption of Theorem \ref{thm: not precise adjunction} and further assuming that $\Ff$ is induced by a morphism $f:X\to Z$, $(S^\nu,\Ff_S,B_S,\Mm^S)$ is lc.
\end{thm}
\begin{proof}
By Theorem \ref{thm: not precise adjunction}, $(S^\nu,\Ff_S,B_S,\Mm^S)$ is sub-lc. It suffices to show that $B_S\ge0.$
%The rest part of Theorem \ref{thm: adjunction foliation nonnqc} is only about the coefficients of divisors on $S^\nu$, which is a codimension 2 property on $X$. 
%By Lemma \ref{lem: stein induce same foliation}, we may assume that $f$ is a contraction. 
By Propositions \ref{prop: general hyperplane invariant} and \ref{prop: general hyperplane non-invariant}, we may cut $X$ by general elements in base-point-free linear systems and assume that $\dim X=2$.

If $\rk\Ff=0$, then since $(X,\Ff,B,\Mm)$ is lc, $B=0$ and $\Mm$ descends to $X$, and the theorem is trivial. If $\rk\Ff=2$ then the theorem follows from Lemma \ref{lem: lc adjunction keep lc}. In the following we assume that $\rk\Ff=1$. 

If $S$ is not $\Ff$-invariant, i.e., $S$ is transverse to $\Ff$. According to \cite[Proposition 3.4]{Spi20} and Theorem \ref{thm: not precise adjunction}, we can see that%there exists an $\Rr$-divisor $\Delta_S\geq 0$ on $S^\nu$ such that
$$(K_{\Ff}+S)|_{S^\nu}=K_{\Ff_S}.$$ 
%By Theorem \ref{thm: not precise adjunction}, $(S^\nu,\Ff_S,\Delta_S)$ is lc, so $\Delta_S=0$. 
It follows that $B_S\ge0$. %The theorem follows in this case. %Thus we may assume that $S$ is $\Ff$-invariant.

Assume that $S$ is $\Ff$-invariant. We only need to check the coefficients of $B_S$ near any closed point $y$ on $S^\nu$. Let $x$ be the image of $y$ in $S$. If $\Ff$ is terminal at $x$, then the theorem follows from Lemma \ref{lem: surface pia terminal}. If $\Ff$ is not terminal at $x$, then the theorem follows from Lemma \ref{lem: surface pia not terminal}.
\end{proof}

\begin{thm}\label{thm: precise adjunction when induced}
   Theorem \ref{thm: precise adj gfq} holds when $\Ff$ is induced by a morphism $f:X\rightarrow Z$. 
\end{thm}
\begin{proof}
Suppose that we are under the conditions of Theorem \ref{thm: precise adj gfq}. According to Theorem \ref{thm: adjunction foliation nonnqc}, we only need to prove that $B_S$ has the form as in Theorem \ref{thm: precise adj gfq}.
%By Lemma \ref{lem: stein induce same foliation}, we may assume that $f$ is a contraction. According to Theorem \ref{thm: adjunction foliation nonnqc}, if $(X,\Ff,\epsilon_{\Ff}(S)S+\sum_{j=1}^mb_j'B_j,\sum_{k=1}^nr_k'\Mm_k)$ is lc, then $$\left(S^\nu,\Ff_S,\frac{w_i-1+\sum_{j=1}^mw_{i,j}b_j'+\sum_{k=1}^nv_{i,k}r_k'}{w_i}C_i,\sum_{k=1}^nr_k'\Mm_{k,S^\nu}^S\right)$$ is lc. 
%The rest part of Theorem \ref{thm: precise adj gfq} is only about the coefficients of divisors on $S^\nu$, which is a codimension 2 property on $X$. Since $\Ff$ is induced by a contraction $X\rightarrow Z$, 
To this end, by Propositions \ref{prop: general hyperplane invariant} and \ref{prop: general hyperplane non-invariant}, we may assume that $\dim X=2$.

If $\rk\Ff=0$, then since $(X,\Ff,B,\Mm)$ is lc, $B=0$ and $\Mm$ descends to $X$, and the theorem is trivial. If $\rk\Ff=2$, then the theorem follows from the usual precise adjunction formula for lc g-pairs \cite[Page 306, Line 30]{BZ16}. In what follows we may assume that $\rk\Ff=1$. 

If $S$ is not $\Ff$-invariant, then by \cite[Proposition 3.4]{Spi20} and Theorem \ref{thm: not precise adjunction}, 
$$(K_{\Ff}+S)|_{S^\nu}=K_{\Ff_S}.$$ 
The theorem follows in this case. Thus we may assume that $S$ is $\Ff$-invariant. Then the statement follows from Lemmas \ref{lem: surface pia terminal} and \ref{lem: surface pia not terminal}.
%We only need to check the coefficient near any closed point $y$ on $S^\nu$. Let $x$ be the image of $y$ in $S$. If $\Ff$ is terminal at $x$, then the theorem follows from Lemma \ref{lem: surface pia terminal}. If $\Ff$ is not terminal at $x$, then the theorem follows from Lemma \ref{lem: surface pia not terminal}.
\end{proof}

\begin{rem}
(1) Theorem \ref{thm: precise adjunction when induced}, even without the control on the coefficients and with $\Mm=\bm{0}$, is already stronger than \cite[Proposition 3.2]{ACSS21} as the latter requires that $X$ is $\Qq$-factorial. 

(2) The complete versions of Theorem \ref{thm: precise adj gfq} will be proven after we establish the existence of Property $(*)$ modifications.
\end{rem}

\section{Property \texorpdfstring{$(*)$}{} and ACSS generalized foliated quadruples}\label{sec: acss gfq}

In this section, we introduce the concepts of Property $(*)$ and ACSS generalized foliated quadruples and study their basic properties.

\subsection{Qdlt generalized pairs}
\begin{defn}[Qdlt]\label{defn: qdlt}
Let $(X,B,\Mm)/U$ be an lc g-pair. We say that $(X,B,\Mm)$ is \emph{qdlt} if there exists an open (possibly empty) subset $V\subset X$ satisfying the following.
\begin{enumerate}
    \item $(V,B|_V)$ is $\Qq$-factorial toroidal. In particular, $B|_V$ is a reduced divisor. 
    \item $V$ contains the generic point of any lc center of $(X,B,\Mm)$.
    \item The generic point of any lc center of $(X,B,\Mm)$ is the generic point of an lc center of $(V,B|_V)$.
\end{enumerate}
\end{defn}

\begin{lem}\label{lem: equi def qdlt}
    Let $(X,B,\Mm)/U$ be an lc g-pair. Then the following conditions are equivalent
    \begin{enumerate}
        \item $(X,B,\Mm)$ is qdlt.
        \item For any lc center of $(X,B,\Mm)$ with generic point $\eta$, near $\eta$, $(X,B)$ is $\Qq$-factorial toroidal and $\Mm$ descends to $X$.
    \end{enumerate}
\end{lem}
\begin{proof}
It is clear that (2) implies (1). Thus we only need to prove (1) implies (2).
    
Let $\eta$ be the generic point of an lc center of $(X,B,\Mm)$. Since $(X,B,\Mm)$ is qdlt, there exists an open subset $V\subset X$ satisfying Definition \ref{defn: qdlt}. In particular, $\eta$ is an lc center of $(V,B|_V)$ and $\Mm_X|_V$ is $\Rr$-Cartier. We let $\Mm^V:=\Mm|_V$ be the restricted $\bb$-divisor of $\Mm$ on $V$, then $\Mm^V$ is nef$/V$ and $\Mm^V_V=\Mm_X|_V$. Suppose that $h: V'\rightarrow V$ is a resolution of $V$ such that $\Mm^V$ descends to $V'$ and there exists a prime divisor $E$ on $V'$ such that $\Center_{V}E=\bar\eta$ and $E$ is an lc place of $(V,B|_V)$. By the negativity lemma,
$$\Mm^V_{V'}=h^*\Mm^V_V-F$$
for some $F\geq 0$. Moreover, we have either $F=0$ over $\eta$ or $\Supp F=\Supp h^{-1}(\bar\eta)$. Since $(X,B,\Mm)$ is lc, $(V,B|_V,\Mm^V)$ is lc. Thus $F=0$ over $\eta$. Possibly shrinking $V$, we may assume that $\Mm$ descends to $V$. The lemma follows.
\end{proof}

\begin{lem}\label{lem: qdlt equivalent definition}
Let $(X,B,\Mm)$ be an lc g-pair and $x$ a (not necessarily closed) point of $X$ such that $\bar x$ is an lc center of $(X,B,\Mm)$. Let $d:=\dim X-\dim \bar x$. Then the following conditions are equivalent
\begin{enumerate}
  \item $(X,B,\Mm)$ is qdlt near $x$.
  \item There exist components $D_1,\dots,D_{d}$ of $\lfloor B\rfloor$, such that $K_X$ and each $D_i$ are $\Qq$-Cartier near $x$, and $x\in\Supp D_i$ for each $i$.
\end{enumerate}
\end{lem}
\begin{proof}
It is obvious that (1) implies (2). We prove (2) implies (1). Possibly shrinking $X$ around $x$, we may assume that $(X,\sum_{i=1}^{d}D_i)$ is a pair. Since $B\geq\sum_{i=1}^{d}D_i$, $(X,D)$ is lc near $x$. By \cite[Proposition 34]{dFKX17}, $B=\sum_{i=1}^d D_i$ near $x$, $(X,B)$ is qdlt near $x$, and $\bar x$ is an lc center of $(X,B)$. Since $(X,B,\Mm)$ is lc, $x$ is an lc center of $(X,B,\Mm)$, and $(X,B,\Mm)$ is qdlt near $x$.
\end{proof}

\begin{lem}\label{lem: qdlt perturbation}
Let $(X,B,\Mm)$ be a qdlt g-pair and $D\geq 0$ an $\Rr$-Cartier $\Rr$-divisor on $X$ such that $D\subset\Supp\{B\}$. Then there exists a positive real number $\delta$ such that $(X,B+\delta D,\Mm)$ is qdlt.
\end{lem}
\begin{proof}
By definition, $\Supp\{B\}$ does not contain any lc center of $(X,B,\Mm)$. Thus $(X,B+\epsilon D,\Mm)$ is lc for some positive real number $\epsilon$. Let $\delta:=\frac{\epsilon}{2}$, then $(X,B+\delta D,\Mm)$ is lc, and any lc center of $(X,B+\delta D,\Mm)$ is an lc center of $(X,B,\Mm)$. By definition, $(X,B+\delta D,\Mm)$ is qdlt.
\end{proof}

\begin{lem}\label{lem: mmp preserves qdlt}
Let $(X,B,\Mm)/U$ be an lc g-pair and $\phi: (X,B,\Mm)\dashrightarrow (X',B',\Mm)$ a sequence of steps of a $(K_X+B+\Mm_X)$-MMP. Suppose that $(X,B,\Mm)$ is qdlt. Then $(X',B',\Mm)$ is qdlt. 
\end{lem}
We remark here that $\phi$ may not be an MMP$/U$, so $(X',B',\Mm)/U$ may not be a g-pair, but $(X',B',\Mm)/X'$ is a g-pair.
\begin{proof}
By Lemma \ref{lem: qdlt equivalent definition}, we only need to show that, for any lc center $S'$ of $(X',B',\Mm)$ with generic point $\eta_{S'}$, near $\eta_{S'}$, we have that $(X',B')$ is $\Qq$-factorial toroidal and $S'$ is an lc center of $(X',B')$. Indeed, let $E$ be an lc place of $(X',B',\Mm)$ such that $\Center_{X'}E=S'$. By our assumption, $E$ is also an lc place of $(X,B,\Mm)$ and $\phi^{-1}$ is an isomorphism near $\eta_{S'}$. Then the lemma follows by Lemma \ref{lem: qdlt equivalent definition} again.
\end{proof}

\subsection{Property \texorpdfstring{$(*)$}{} generalized foliated quadruples}

\begin{defn}
Let $f: X\rightarrow Z$ be a morphism between normal varieties and $G$ an $\Rr$-divisor on $X$. We say that $G$ is \emph{super$/Z$}, or \emph{$f$-super}, if there exist ample Cartier divisors $H_1,\dots,H_{2\dim X+1}$ on $Z$ such that $G\geq\sum_{i=1}^{2\dim X+1}f^*H_i$.
\end{defn}

\begin{defn}[Property $(*)$ gfq]\label{defn: foliation property *}
Let $(X,\Ff,B,\Mm)/U$ be a sub-gfq. Let $G\geq 0$ be a reduced divisor on $X$ and let $f: X\rightarrow Z$ be a contraction. We say that $(X,\Ff,B,\Mm;G)/Z$ \emph{satisfies Property $(*)$} if the following holds
\begin{enumerate}
  \item $\Ff$ is induced by $f$.
  \item $G$ is an $\Ff$-invariant divisor.
  \item $f: (X,B+G,\Mm)\rightarrow Z$ satisfies Property $(*)$ (See Definition \ref{defn: property *}). 
\end{enumerate}
We say that $(X,\Ff,B,\Mm)/Z$ satisfies Property $(*)$ if $(X,\Ff,B,\Mm;G)/Z$ satisfies Property $(*)$ for some $G\geq 0$. We say that $(X,\Ff,B,\Mm)$ satisfies Property $(*)$ if $(X,\Ff,B,\Mm;G)/Z$ satisfies Property $(*)$ for some $G$ and $X\rightarrow Z$. In this case, we say that $f: X\rightarrow Z$ is an \emph{associated contraction} of $(X,\Ff,B,\Mm)$ and $Z$ an \emph{associated base} of $(X,\Ff,B,\Mm)$. 
We say that $G$ is \emph{associated to} $(X,\Ff,B,\Mm)/Z$, and if $G$ is \emph{super}$/Z$, then we say that $G$ is \emph{superbly associated to} $(X,\Ff,B,\Mm)/Z$.

It is clear that property $(*)$ is independent of the choice of $U$. We remark that the choice of $f$ and $G$ may not be unique. We also remark that $f$ may not be a morphism$/U$.
\end{defn}

\begin{defn}[ACSS gfq, {cf. \cite[Definition 4.3]{DLM23}}]\label{defn: ACSS f-triple}
Let $(X,\Ff,B,\Mm)/U$ be a gfq, $G\geq 0$ a reduced divisor on $X$, and $f: X\rightarrow Z$ a projective morphism. We say that $(X,\Ff,B,\Mm;G)/Z$ is \emph{weak ACSS} if 
\begin{enumerate}
    \item $(X,\Ff,B,\Mm;G)/Z$ satisfies Property $(*)$ and $(X,\Ff,B,\Mm)$ is lc, and
    \item $f$ is equidimensional.
\end{enumerate}
We say that $(X,\Ff,B,\Mm;G)/Z$ is \emph{ACSS} if the following additional conditions are satisfied
\begin{enumerate}
  \item[(3)] There exists an \(\Rr\)-divisor \(D\geq 0\) on \(X\) and a nef\(/X\) \(\bb\)-divisor \(\Nn\) such that
    \begin{enumerate}
      \item $\Supp\{B\}\subset\Supp D$,
      \item $\Nn-\alpha \Mm$ is nef$/X$ for some $\alpha>1$,
      \item $D+\Nn_X-\Mm_X$ is $\Rr$-Cartier, and
      \item $(X,B+D+G+f^*(\Sigma-f(G)),\Nn)$ is qdlt, where $\Sigma\ge f(G)$ is a reduced divisor such that $(Z,\Sigma)$ is log smooth.
    \end{enumerate}
  \item[(4)] For any lc center of $(X,\Ff,B,\Mm)$ with generic point $\eta$, over a neighborhood of $\eta$
    \begin{enumerate}
      \item $\Mm$ descends to $X$,
      \item $\eta$ is the generic point of an lc center of $(X,\Ff,\lfloor B\rfloor)$, and
       \item $f: (X,B+G)\rightarrow (Z,f(G))$ is a toroidal morphism; in particular, $(X,B)$ is toroidal and $B=\lfloor B\rfloor$.
    \end{enumerate}
\end{enumerate}
  In this case, we say that the divisor $G$ and the variety $Z$ are \emph{properly associated to} $(X,\Ff,B,\Mm)$. If additionally $G$ is super$/Z$, then we say that $(X,\Ff,B,\Mm;G)/Z$ is \emph{super ACSS}.

We say that $(X,\Ff,B,\Mm)/Z$ is \emph{weak ACSS} (resp. \emph{ACSS}, \emph{super ACSS}) if $(X,\Ff,B,\Mm;G)/Z$ is \emph{weak ACSS} (resp. \emph{ACSS}, \emph{super ACSS}) for some $G$. We say that $(X,\Ff,B,\Mm)$ is \emph{weak ACSS} (resp. \emph{ACSS}, \emph{super ACSS}) if $(X,\Ff,B,\Mm;G)/Z$ is \emph{weak ACSS} (resp. \emph{ACSS}, \emph{super ACSS}) for some $G$ and $f: X\rightarrow Z$.
\end{defn}

\begin{lem}\label{lem: weak acss can be super}
    Let $(X,\Ff,B,\Mm)/U$ be a gfq and $f: X\rightarrow Z$ a contraction. Then there exists a super$/Z$ divisor $G$ on $X$ such that if $(X,\Ff,B,\Mm)/Z$ satisfies Property $(*)$ (resp. is weak ACSS), then $(X,\Ff,B,\Mm;G)/Z$ satisfies Property $(*)$ (resp. is weak ACSS).
\end{lem}
\begin{proof}
If $(X,\Ff,B,\Mm)/Z$ satisfies Property $(*)$ (resp. is weak ACSS), then there exists a divisor $G_0\geq 0$ on $X$ such that $(X,\Ff,B,\Mm;G_0)/Z$ satisfies Property $(*)$ (resp. is weak ACSS). We let $H_1,\dots,H_{2\dim X+1}$ be general elements of a very ample linear system on $Z$, and let $$G:=G_0+\sum_{i} f^*H_i.$$ 
Then by definition, $(X,\Ff,B,\Mm;G)/Z$ satisfies Property $(*)$ (resp. is weak ACSS).
\end{proof}

The following lemmas will be very useful when applying to the minimal model program for algebraically integrable foliations.

\begin{lem}\label{lem: acss smaller coefficient}
Assume that $(X,\Ff,B,\Mm)/U$ and $(X,\Ff,B',\Mm')/U$ are two gfqs such that $B\geq B'$ and $\Mm-\Mm'$ is nef$/U$. If $(X,\Ff,B,\Mm)/U$ satisfies Property $(*)$ (resp. is weak ACSS, ACSS, super ACSS) with associated contraction $f:X\rightarrow Z$ and divisor $G$, then so does (resp. is) $(X,\Ff,B',\Mm';G)/Z$. 
\end{lem}
\begin{proof}
The lemma follows by checking the definitions.
\end{proof}

\begin{lem}\label{lem: fls imply acss}
If $(X,\Ff,B,\Mm)$ is foliated log smooth, then $(X,\Ff,B^\Ff,\Mm)$ is $\Qq$-factorial super ACSS. 
\end{lem}
\begin{proof}
By definition, $X$ is $\Qq$-factorial and there exists an equidimensional toroidal morphism $f: (X,\Sigma_X,\Mm)\rightarrow (Z,\Sigma_Z)$ such that $\Supp B\subset\Sigma_X$. We let $G$ be the vertical$/Z$ part of $\Sigma_X$. By Lemma \ref{lem: acss smaller coefficient}, it suffices to show that $(X,\Ff,B, \Mm)$ is ACSS. To this end, we may assume that $B=B^\Ff=\Sigma_X-G$. By \cite[Proposition 3.2]{AK00}, possibly adding general hyperplane sections to $\Sigma_Z$ and adding their pullbacks to $\Sigma_X$, we may assume that $G$ is super$/Z$. We only need to show that $(X,\Ff,B,\Mm;G)/Z$ is ACSS.  By Proposition \ref{prop: weak ss satisfies *} and Lemma \ref{lem: foliated log smooth imply lc}, we can see that $(X,\Ff,B,\Mm;G)/Z$ satisfies Property $(*)$. Since $\{B\}=0$, $\Mm$ descends to $X$, and $f^*(\Sigma-f(G))=f^{-1}(\Sigma-f(G))$, by \cite[Proposition 3.2]{AK00}, 
$$f: (X,B+G+D+f^*(\Sigma-f(G)),\Mm)\rightarrow (Z,\Sigma)$$
is toroidal. In particular, $(X,B+G+D+f^*(\Sigma-f(G)),\Nn)$ is qdlt. Finally, Definition \ref{defn: ACSS f-triple}(4) immediately follows from the definition of foliated log smooth.
\end{proof}

\begin{lem}\label{lem: alg int foliation lct achieved}
Let $(X,\Ff,B,\Mm)/U$ be a sub-gfq, $D$ an $\Rr$-divisor on $X$, and $\Nn$ a \(\bb\)-divisor on \(X\) such that \(D+\Nn_X\) is \(\Rr\)-Cartier and \(\Nn\) descends to a birational model of \(X\). Suppose that \(\Ff\) is algebraically integrable. Let
$$t:=\sup\{s\mid s\geq 0, \Mm+s\Nn\text{ is nef}/U,\text{ and } (X,\Ff,B+sD,\Mm+s\Nn)/X\text{ is sub-lc}\}.$$
Then either $t=+\infty$, or 
$$t:=\max\{s\mid s\geq 0, \Mm+s\Nn\text{ is nef}/U,\text{ and } (X,\Ff,B+sD,\Mm+s\Nn)/X\text{ is sub-lc}\}.$$
In particular, $(X,\Ff,B+tD)$ is sub-lc and $\Mm+t\Nn$ is nef$/U$ if $t<+\infty$.
\end{lem}
\begin{proof}
Let $$l:=\sup\{s\mid s\geq 0, \Mm+s\Nn\text{ is nef}/U\}=\max\{s\mid s\geq 0, \Mm+s\Nn\text{ is nef}/U\}$$ since nef is a closed condition. Moreover, $t\leq l$. 

Let $f: X\dashrightarrow Z$ be a dominant map that induces $\Ff$. Let $g: \bar X\rightarrow Z$ be a birational morphism such that $f\circ g$ is a morphism and both $\Mm$ and $\Nn$ descend to $\bar X$. By Lemma \ref{lem: stein induce same foliation}, we may assume that $f\circ g$ is a contraction. Let $\bar B$ be the reduced divisor supported on $$\Supp(f^{-1}_*B)\cup\Supp(D)\cup\Supp\Exc(g).$$ 
By Definition-Theorem \ref{defthm: weak ss reduction}, there exists an equidimensional model$/U$ $f': (X', \Sigma_{X'})\rightarrow (Z', \Sigma_{Z'})$ of $f\circ g: (\bar X, \bar B)\rightarrow Z$ associated with $h: X'\rightarrow\bar X$ and $h_Z: Z'\rightarrow Z$. Since $\Mm$ and $\Nn$ descend to $X'$, let $\phi:=g\circ h$, $\Ff':=\phi^{-1}\Ff$, and 
$$K_{\Ff'}+B'+\Mm_{X'}:=\phi^*(K_{\Ff}+B+\Mm_X).$$
Then $\phi$ is a foliated log resolution of $(X,\Ff,\Supp B\cup\Supp D,\Mm)$, and $f'$ induces $\Ff'$. Since $\Mm$ and $\Nn$ descend to $X'$, by Lemma \ref{lem: foliated log smooth imply lc},
\begin{align*}
    t&=\sup\{s\mid 0\leq s\leq l, a(E,X,\Ff,B+sD,\Mm+s\Nn)\geq-\epsilon_{\Ff}(E)\text{ for any prime divisor } E\text{ on }X\}\\
    &=\sup\{s\mid 0\leq s\leq l, a(E,X,\Ff,B+sD,\Mm+s\Nn)\geq-\epsilon_{\Ff}(E)\text{ for any prime divisor } E\subset\Supp \bar B\}.
\end{align*}
Since there are only finitely many components of $\Supp \bar B$, the lemma follows.
\end{proof}

\begin{lem}\label{lem: acss f-triple perturb coefficient}
Let $(X,\Ff,B,\Mm)/U$ be a gfq, $G\geq 0$ a reduced divisor on $X$, and $f: X\rightarrow Z$ a contraction, such that $(X,\Ff,B,\Mm;G)/Z$ is ACSS. Let $D$ be an $\Rr$-divisor on $X$ and $\Nn$ a $\bb$-divisor on $X$ such that
\begin{itemize}
  \item both $D$ and $\Nn_X$ are $\Rr$-Cartier,
  \item $\Supp D\subset\Supp\{B\}$, and $\Nn$ descends to a birational model of $X$, and
  \item $\Mm+\Nn$ is nef$/U$, and $\Mm-\delta\Nn$ is nef$/U$ for some $\delta\in(0,1)$.
\end{itemize}
Then there is a positive real number $\gamma$ such that $(X,\Ff,B+\alpha D,\Mm+\beta\Nn;G)/Z$ is ACSS for any $\alpha,\beta\in [0,\gamma]$.
\end{lem}
\begin{proof}
By assumption, $\Supp D$ does not contain any lc center of $(X,\Ff,B,\Mm)$, and $\Mm$ descends to $X$ near the generic point of any lc center of $(X,\Ff,B,\Mm)$. Thus $\Nn$ descends to $X$ near the generic point of any lc center of $(X,\Ff,B,\Mm)$. Since $\Mm+\Nn$ is nef$/U$, by Lemma \ref{lem: alg int foliation lct achieved}, there exists a real number $\gamma_0\in (0,1)$ such that both $(X,\Ff,B+\gamma_0 D,\Mm)$ and $(X,\Ff,B,\Mm+\gamma_0\Nn)$ are lc. Possibly replacing \(\gamma_0\) with \(\frac{1}{2}\gamma_0\), we may assume that \((X,\Ff,B+\gamma_0 D,\Mm)\), \((X,\Ff,B,\Mm+\gamma_0\Nn)\), and \((X,\Ff,B,\Mm)\) have the same lc centers. 

By convexity of discrepancies, for any $\alpha,\beta\in [0,\frac{1}{2}\gamma_0]$, $(X,\Ff,B+\alpha D, \Mm+\beta\Nn)$ is lc with the same lc centers as $(X,\Ff,B,\Mm)$. In particular, for any lc center of $(X,\Ff,B+\alpha D, \Mm+\beta\Nn)$ with generic point $\eta$, near $\eta$, $\Mm+\beta\Nn$ descends to $X$, and $B=\lfloor B+\alpha D\rfloor$. Thus $f: (X,B+\alpha D+G)\rightarrow (Z,f(G))$ is a toroidal morphism near $\eta$, and $\eta$ is the generic point of an lc center of $(X,\Ff,\lfloor B+\alpha D\rfloor)$ over a neighborhood of $\eta$.

Since $(X,\Ff,B,\Mm;G)/Z$ is ACSS, there exists an $\Rr$-divisor $D'$ on $X$ such that $\Supp\{B\}\subset\Supp D'$, and a nef$/X$ $\bb$-divisor $\Nn'$, such that $\Nn'-\alpha'\Mm$ is nef$/X$ for some $\alpha'>1$, and for any $\Sigma$ on $X$ such that $\Sigma\geq f(G)$ and $(Z,\Sigma)$ is log smooth, $(X,B+G+D'+f^*(\Sigma-f(G)),\Nn')$ is qdlt. 

We will show that 
$$\gamma:=\min\left\{\frac{\gamma_0}{3},\left(\alpha'-1\right)\delta\right\}$$ 
has the required properties. Indeed, As $\Supp D\subset\Supp\{B\}\subset\Supp D'$, we see that for any $\alpha\in(0,\gamma)$, it holds that 
$$D'-\alpha D\ge0\text{ and }\Supp (D'-\alpha D)=\Supp D'\supset\Supp\{B+\alpha D\}.$$
By definition, it is enough to prove that $D'-\alpha D$ and $\Nn'$ satisfy Definition \ref{defn: ACSS f-triple}(3). In particular, we only need to show that there exists a positive real number $\alpha''>1$ such that $\Nn'-\alpha''(\Mm+\beta\Nn)$ is nef$/X$. Take 
$$\alpha'':=\frac{\alpha'}{1+\frac{\beta}{\delta}},$$
then it is clear that $\alpha''>1$ as $\beta\le\gamma<(\alpha'-1)\delta$.
Then
$$\Nn'-\alpha''(\Mm+\beta\Nn)=\frac{\alpha''}{\alpha'}(1+\frac{\beta}{\delta})(\Nn'-\alpha'\Mm)+\frac{\alpha''\beta}{\delta}(\Mm-\delta\Nn)$$
is nef$/X$. We may finish the proof.
\end{proof}

\begin{prop}[cf. {\cite[Proposition 3.6]{ACSS21}}]\label{prop: weak cbf gfq}
Suppose that $f: (X,B+G,\Mm)\rightarrow Z$ satisfies Property $(*)$ and is equidimensional. Assume that $B$ is horizontal$/Z$ and $G$ is vertical$/Z$. Let $\Ff$ be the foliation induced by $f$, and let $\Nn$ be the moduli part of $f: (X,B+G,\Mm)\rightarrow Z$. Then
\begin{enumerate}
  \item $K_{\Ff}+B+\Mm_X\sim \Nn_X$, and
  \item $K_{\Ff}+B+\Mm_X\sim_{Z}K_X+B+G+\Mm_X$.
\end{enumerate}
In particular, $K_{\Ff}+B+\Mm_X$ is $\Rr$-Cartier.
\end{prop}
\begin{proof}
By assumption, $Z$ is smooth and $B_Z=f(G)$ is the discriminant part of $f: (X,B+G,\Mm)\rightarrow Z$. Moreover, $B_Z$ is reduced by Lemma \ref{lem: basic property (*) gpair}. Then (2) follows from (1) immediately. Thus we only need to prove (1). 

Since $f$ is equidimensional, we have
$$K_{\Ff}=K_{X/Z}-R$$
where $R:=\sum(f^*D-f^{-1}(D))$ and $D$ runs over the prime divisors on $Z$. We claim that $f^*B_Z=R+G$. Indeed, let $D$ be a prime divisor on $X$ such that $D$ is vertical$/Z$. Since $f$ is equidimensional, $D_Z:=f(D)$ is a prime divisor. If $D_Z\subset\Supp B_Z$, then $D\subset\Supp G$ and $\mult_DG=1$. Therefore, 
\begin{align*}
    \mult_Df^*B_Z&=\mult_Df^*D_Z=\mult_Df^{-1}(D_Z)+\mult_D(f^*D_Z-f^{-1}(D_Z))=\mult_DG+\mult_DR.
\end{align*}
If $D_Z\not\subset \Supp B_Z$, then $\mult_Df^*B_Z=\mult_DG=0$. By assumption, $(X,B+G+f^*D_Z,\Mm)$ is sub-lc over the generic point of $D_Z$. In particular, $\mult_Df^*D_Z=1$ and hence $\mult_DR=0$. Since $f^*B_Z$ and $R+G$ are both vertical$/Z$, the claim follows.

Since $f^*B_Z=R+G$, one can see that
\begin{align*}
    \Nn_X&\sim K_X+B+G+\Mm_X-f^*(K_Z+B_Z)=K_{X/Z}+B+G+\Mm_X-f^*B_Z\\
    &=K_{\Ff}+R+B+G-f^*B_Z=K_{\Ff}+B+\Mm_X.
\end{align*}
This completes the proof.
\end{proof}

\subsection{\texorpdfstring{$(*)$-}{ }models and ACSS models}

\begin{defn}\label{defn: acss model}
Let $(X,\Ff,B,\Mm)/U$ be a gfq such that $\Ff$ is algebraically integrable. A \emph{$(*)$-modification} of $(X,\Ff,B, \Mm)$ is a projective birational morphism $h:Y\rightarrow X$ such that
\begin{enumerate}
   \item[(a)] $Y$ is klt,
   \item[(b)] $\left(Y,\Ff_Y:=h^{-1}\Ff,B_Y:=h^{-1}_*(B\wedge\Supp B)+(\Supp\Exc(h))^{\Ff_Y},\Mm\right)$ is weak ACSS, and
  \item[(c)] $a(E,\Ff,B,\Mm)\leq-\epsilon(E)$ for any $h$-exceptional prime divisor $E$.
  In particular, if $(X,\Ff,B,\Mm)$ is lc, then 
  $$K_{\Ff_Y}+B_Y+\Mm_Y=h^*(K_{\Ff}+B+\Mm_X).$$
\end{enumerate}
In this case, we say that $(Y,\Ff_Y,B_Y,\Mm)$ is a \emph{$(*)$-model} of $(X,\Ff,B,\Mm)$. If additionally, $(Y,\Ff_Y,B_Y,\Mm)$ is $\Qq$-factorial ACSS (resp. $\Qq$-factorial super ACSS), then we say that $h:Y\to X$ is an \emph{ACSS modification} (resp. \emph{super ACSS modification}) of $(X,\Ff,B,\Mm)$, and $(Y,\Ff_Y,B_Y,\Mm)$ is an \emph{ACSS model} (resp. \emph{super ACSS model}) of $(X,\Ff,B,\Mm)$. 

Let $f: Y\rightarrow X$ be a $\Qq$-factorial $(*)$-modification (resp. ACSS modification, super ACSS modification) of $(X,\Ff,B,\Mm)$, let $\Ff_Y:=f^{-1}\Ff$, and let $$B_Y:=f^{-1}_*(B\wedge\Supp B)+(\Supp\Exc(f))^{\Ff_Y}.$$ 
We say that $f$ is a \emph{proper $(*)$-modification} (resp. \emph{proper ACSS modification}, \emph{great ACSS modification}) of $(X,\Ff,B,\Mm)$ if there exists a contraction $\pi: Y\rightarrow Z$ and a reduced divisor $G$ on $Y$, such that $(Y,\Ff_Y,B_Y,\Mm;G)/Z$ is weak ACSS (resp. ACSS, super ACSS), and for any $f$-exceptional $\Ff_Y$-invariant divisor $D$, $D\subset\Supp G$. We call $(Y,\Ff_Y,B_Y,\Mm)$, $(Y,\Ff_Y,B_Y,\Mm)/Z$, and $(Y,\Ff_Y,B_Y,\Mm;G)/Z$ \emph{proper Property $(*)$ models} (resp. \emph{proper ACSS models}, \emph{great ACSS models}) of $(X,\Ff,B,\Mm)$.
\end{defn}

\begin{nota}
Let $(X_0,\Ff_0,B_0,\Mm)/U$ be a Property $(*)$ gfq, $Z$ a base associated to $(X_0,\Ff_0,B_0,\Mm)$, and $G_0$ a divisor associated to $(X_0,\Ff_0,B_0,\Mm)/Z$. When we say the following
\begin{center}$\xymatrix{
(X_0,\Ff_0,B_0,\Mm;G_0)\ar@{-->}[r]^{f_0} & (X_1,\Ff_1,B_1,\Mm;G_1)\ar@{-->}[r]^{\ \ \ \ \ \ \ \ \ \ f_1} & \dots\ar@{-->}[r] & (X_n,\Ff_n,B_n,\Mm;G_n)\ar@{-->}[r]^{\ \ \ \ \ \ \ \ \ \ f_n} & \dots 
}$
\end{center}
is a (possibly infinite) sequence of steps of a $(K_{\Ff_0}+B_0+\Mm_{X_0})$-MMP$/U$, we mean the following: for any $i$, $f_i: X_{i}\dashrightarrow X_{i+1}$ is a step of a $(K_{\Ff_i}+B_i+\Mm_{X_i})$-MMP$/U$ that is not a Mori fiber space, $\Ff_{i+1}:=(f_i)_*\Ff_i$, $B_{i+1}:=(f_i)_*B_i$, and $G_{i+1}:=(f_i)_*G_i$.
\end{nota}

\section{Cone theorem and ACSS modifications}\label{sec: cone}

In this section, we prove the cone theorem (Theorem \ref{thm: cone theorem gfq}) and the existence of ACSS modifications (Theorem \ref{thm: ACSS model}). As an immediate corollary, we will prove the precise adjunction formula (Theorem \ref{thm: precise adj gfq}) in full generality, without assuming that $\Ff$ is induced by a contraction.

\subsection{Bend and break} 
It is important to notice that we will work under the relative setting, so the following relative bend and break theorem is crucial for our proofs.

\begin{thm}[Relative bend and break]\label{thm: relative bb}
Let $d$ be a positive integer, $\pi: X\rightarrow U$ a contraction from a normal quasi-projective variety to a variety such that $\dim X-\dim U=d$, $M,D_1,\dots,D_d$ $\Rr$-divisors on $X$ that are nef along general fibers of $\pi$, $B\geq 0$ an $\Rr$-divisor on $X$, and $\Ff$ a foliation on $X$. Suppose that for any general fiber $F$ of $\pi$,
\begin{enumerate}
    \item $(D_1|_F)\cdot (D_2|_F)\cdot\dots\cdot (D_d|_F)=0$, and
    \item $-(K_{\Ff}+B)|_F\cdot (D_2|_F)\cdot\dots\cdot (D_d|_F)>0$.
\end{enumerate}
Then for any general closed point $x\in X$, there exists a rational curve $C_x$ satisfying the following.
\begin{enumerate}
    \item $x\in C_x$, 
    \item $\pi(C_x)$ is a point, and
    \item $D_1\cdot C_x=0$ and
    $$M\cdot C_x\leq 2d\frac{M|_F\cdot  (D_2|_F)\cdot\dots\cdot (D_d|_F)}{-K_{\Ff}|_F\cdot  (D_2|_F)\cdot\dots\cdot (D_d|_F)}.$$
\end{enumerate}
\end{thm}
\begin{proof}
Since statement (3) is a closed condition and $M$ is a limit of $\Qq$-divisors that are nef along general fibers of $\pi$, we may assume that $M$ is a $\Qq$-divisor. Possibly replacing $M$ with a multiple, we may assume that $M$ is a Weil divisor.

We let $X^c$ and $U^c$ be compactifications of $X$ and $U$, such that $X^c$ and $U^c$ are normal projective, $X$ is a dense open subset of $X^c$, $U$ is a dense open subset of $U^c$, and there exists a contraction $\pi^c: X^c\rightarrow U^c$ such that $\pi^C|_{X}=\pi$. Let $M^c,D^c_1,\dots,D^c_d,B^c$ be the closures of $M,D_1,\dots,D_d,B$ in $X^c$ respectively, and let $\Ff^c$ be the natural extension of $\Ff$ in $X^c$. Then the general fibers of $\pi^c$ are general fibers of $\pi$, and $M^c,D^c_1,\dots,D^c_d$ are $\Rr$-divisors that are nef along general fibers of $\pi$. Since we only care about properties of general fibers of $\pi$ and properties near a general closed point $x\in X$, we may replace $\pi: X\rightarrow U$ with $\pi^c: X^c\rightarrow U^c$, $M,D_1,\dots,D_d,B$ with $M^c,D^c_1,\dots,D^c_d,B^c$, and $\Ff$ with $\Ff^c$, and assume that $\pi$ is a projective morphism between normal projective varieties.

Let $x\in X$ be a general closed point. Then $x$ is contained in a general fiber $F$ of $\pi$. Let $q:=\dim U$. Then there exist general hyperplane sections $H_1,\dots,H_q$ with $A_i:=\pi^*H_i$, such that $F=\cap_{i=1}^q\pi^*A_i$. Let $V_k:=X\cap_{i=1}^kA_i$ and $W_k:=U\cap_{i=1}^kH_i$ for each $0\leq k\leq q$. Then we have
$$F=V_q\subset V_{q-1}\subset\dots\subset V_0=X$$
and
$$z:=W_q\subset W_{q-1}\subset\dots\subset W_0=U$$
where $z$ is a closed point. We may inductively define $\Ff_{k}$ to be the restricted foliation of $\Ff$ on $V_k$ for each $k$, and let $\Ff_F:=\Ff_{q}$. We let $M_k:=M|_{V_k}$, $B_k:=B|_{V_k}$, $M_F:=M|_F$, and $B_F:=B|_F$. Then it is clear that $M_k|_F=M|_F$, $B_k|_F=B_F$ for each $k$, and $B_{V_k}\geq 0$ for each $k$. Moreover, since $H_1,\dots,H_q$ are general hyperplane sections, $M_k$ is a Weil divisor for each $k$.

\begin{claim}\label{claim: induction bend and break}
There exists a rational curve $C_x$ such that $x\in C_x$, $\pi(C_x)$ is a closed point, $D_1\cdot C_x=0$, and
$$M|_F\cdot C_x\leq2d\frac{M|_F\cdot  (D_2|_F)\cdot\dots\cdot (D_d|_F)}{-K_{\Ff_k}|_F\cdot  (D_2|_F)\cdot\dots\cdot (D_d|_F)}$$
for each $k$.
\end{claim}
\begin{proof}
    We apply induction on $q-k$. When $q-k=0$, the existence of $C_x$ follows from \cite[Corollary 2.28]{Spi20}. We will show that this $C_x$ satisfies our requirement for all $q-k$ as well. In the following, we may assume that $q>k$.

    We let $\pi_k: V_k\rightarrow W_k$ be the restricted contraction of $\pi$ to $V_k$ for each $k$. We consider $W_{k+1}$ as a divisor on $W_k$ and $V_{k+1}$ as a divisor on $V_k$. There are two possibilities.

\medskip

\noindent\textbf{Case 1}. $V_{k+1}$ is $\Ff_k$-invariant. In this case, the general fibers of $\pi_k$ are tangent to $\Ff_k$, so 
$$K_F=K_{\Ff_F}=K_{\Ff_k}|_F.$$
Thus by the $q-k=0$ case,
$$M|_F\cdot C_x\leq2d\frac{M|_F\cdot  (D_2|_F)\cdot\dots\cdot (D_d|_F)}{-K_{\Ff_F}\cdot  (D_2|_F)\cdot\dots\cdot (D_d|_F)}=2d\frac{M|_F\cdot  (D_2|_F)\cdot\dots\cdot (D_d|_F)}{-K_{\Ff_k}|_F\cdot  (D_2|_F)\cdot\dots\cdot (D_d|_F)}.$$

\medskip

\noindent\textbf{Case 2}. $V_{k+1}$ is not $\Ff_k$-invariant. In this case, by \cite[3.6(1)]{Dru21}, we have
$$(K_{\Ff_k}+V_{k+1})|_{V_{k+1}}\sim K_{\Ff_{k+1}}+D_{k+1}$$
for some $\Qq$-divisor $D_{k+1}\geq 0$. 

Since $H_{k+1}$ is a general hyperplane section, there exists $H_{k+1}'\sim H_{k+1}$ such that $H_{k+1}'$ does not contain $z$. Thus
$$V_{k+1}|_F=(H_{k+1}|_{V_k})|_F=H_{k+1}|_F\sim H'_{k+1}|_F=0.$$
Since $H_{k+2},\dots,H_q$ are general hyperplane sections, $D_{k+1}|_F\geq 0$. Therefore,
\begin{align*}
    &-K_{\Ff_k}|_F\cdot  (D_2|_F)\cdot\dots\cdot (D_d|_F)\\
    =&-(K_{\Ff_k}+V_{k+1})|_F\cdot  (D_2|_F)\cdot\dots\cdot (D_d|_F)\\
    =&-((K_{\Ff_k}+V_{k+1})|_{V_{k+1}})|_F\cdot  (D_2|_F)\cdot\dots\cdot (D_d|_F)\\
    =&-(K_{\Ff_{k+1}}+D_{k+1})|_F\cdot  (D_2|_F)\cdot\dots\cdot (D_d|_F)\\
    \leq&-K_{\Ff_{k+1}}|_F\cdot  (D_2|_F)\cdot\dots\cdot (D_d|_F).
\end{align*}
By induction hypothesis,
$$M|_F\cdot C_x\leq2d\frac{M|_F\cdot  (D_2|_F)\cdot\dots\cdot (D_d|_F)}{-K_{\Ff_{k+1}}|_F\cdot  (D_2|_F)\cdot\dots\cdot (D_d|_F)}=2d\frac{M|_F\cdot  (D_2|_F)\cdot\dots\cdot (D_d|_F)}{-K_{\Ff_k}|_F\cdot  (D_2|_F)\cdot\dots\cdot (D_d|_F)}.$$
\end{proof}
\noindent\textit{Proof of Theorem \ref{thm: relative bb} continued}. It immediately follows from Claim \ref{claim: induction bend and break} by letting $k=0$.
\end{proof}

\subsection{Inductive approach to cone theorem}

Similar to \cite[Theorems 3.9, 3.19]{ACSS21}, the cone theorem for generalized foliated quadruples is closely related to the existence of Property $(*)$ models for generalized foliated quadruples, and their proofs are done inductively. Moreover, we shall directly establish the existence of proper $\Qq$-factorial ACSS models with controlled extraction of divisors. This kind of model is more technically constructed but also more useful in practice.

\begin{thm}[Cone theorem for induction, cf. {\cite[Theorem 3.9]{ACSS21}}]\label{thm: cone theorem induction}
Let $d$ be a positive integer. Let $(X,\Ff,B,\Mm)/U$ be a gfq of dimension $d$ such that $\Ff$ is algebraically integrable. Let $\{R_j\}_{j\in\Lambda}$ be the set of all $(K_{\Ff}+B+\Mm_X)$-negative extremal rays$/U$ that are not contained in the non-lc locus of $(X,\Ff,B,\Mm)$. Then
$$\overline{NE}(X/U)=\overline{NE}(X/U)_{K_{\Ff}+B+\Mm_X\geq 0}+\overline{NE}(X/U)_{\Nlc(X,\Ff,B,\Mm)}+\sum_{j\in\Lambda} R_j$$
and for any $j\in\Lambda$, $R_j$ is spanned by a rational curve $C_j$ such that $C_j$ is tangent to $\Ff$ and $0<-(K_{\Ff}+B+\Mm_X)\cdot C_j\leq 2d$.
\end{thm}

\begin{thm}[Existence of ACSS models, cf. {\cite[Theorem 3.10]{ACSS21}, \cite[Proposition 4.14]{DLM23}}]\label{thm: property * induction}
Let $d$ be a positive integer. Let $(X,\Ff,B,\Mm)/U$ be a gfq of dimension $d$ such that $\Ff$ is algebraically integrable, and $E_1,\dots,E_s$ lc places of $(X,\Ff,B,\Mm)$, such that $(X,\Ff,B,\Mm)$ is lc near the generic point of $\Center_X{E_i}$ for each $i$. Then $(X,\Ff,B,\Mm)$ has a great ACSS model $(Y,\Ff_Y,B_Y,\Mm)$ such that $E_1,\dots,E_s$ are on $Y$ if $(X,\Ff,B,\Mm)$ is lc.
\end{thm}

\begin{lem}\label{lem: induction cone 1}
Let $d$ be a positive integer. Assume that Theorem \ref{thm: cone theorem induction} holds in dimension $\leq d-1$.

Let $(X,\Ff,B,\Mm)/U$ be an lc gfq of dimension $d$ satisfying Property $(*)$ with an associated contraction $f: X\rightarrow Z$. Suppose that for any $(K_{\Ff}+B+\Mm_X)$-negative extremal ray$/U$ $R$, there exists a prime divisor $E$ on $X$ such that $R$ is contained in the image of $\overline{NE}(E/U)\rightarrow\overline{NE}(X/U)$ and $\mult_EB=\epsilon_{\Ff}(E)$. Let $\{R_j\}_{j\in\Lambda}$ be the set of $(K_{\Ff}+B+\Mm_X)$-negative extremal rays$/U$. Then
\begin{enumerate}
    \item $\overline{NE}(X/U)=\overline{NE}(X/U)_{K_{\Ff}+B+\Mm_X\geq 0}+\sum_{j\in\Lambda} R_j$.
  \item Each $R_j$ is spanned by a rational curve $C_j$ such that $C_j$ is tangent to $\Ff$ and $0\leq -(K_{\Ff}+B+\Mm_X)\cdot C_j\leq 2(d-1)$.
  \item For any curve $C_j'$ such that $[C_j']\in R_i$, $C_j'$ is contracted by $f$.
  \item Assume that $f$ is equi-dimensional and either $X$ is $\Qq$-factorial klt or $\Mm$ is NQC$/U$. Then
\begin{enumerate}
    \item $\Lambda$ is a countable set.
    \item For any ample$/U$ $\Rr$-divisor $A$ on $X$, there exists a finite set $\Lambda_A\subset\Lambda$ such that
    $$\overline{NE}(X/U)=\overline{NE}(X/U)_{K_{\Ff}+B+A+\Mm_X\geq 0}+\sum_{j\in\Lambda_A}R_j.$$
    \item For any $j\in\Lambda$, there exists a contraction $\phi_j: X\rightarrow X_j'$ of $R_j$ such that 
    \begin{enumerate}
    \item $\phi_j$ is a contraction$/U$ as well as a contraction$/Z$, and
    \item If $\phi_j$ is small, then there exists a small contraction $\phi_j^+: X_j^+\rightarrow X_j'$ such that the induced birational map $\psi_j: X\dashrightarrow X_j^+$ is both a $(K_{\Ff}+B+\Mm_X)$-flip$/U$ and a $(K_{\Ff}+B+\Mm_X)$-flip$/Z$.
    \end{enumerate}
    \item Let $G$ be any divisor associated to $(X,\Ff,B,\Mm)/U$. For any $j$,
    $$(K_{\Ff}+B+\Mm_X)\cdot R_j=(K_X+B+G+\Mm_X)\cdot R_j.$$
    In particular,
    \begin{enumerate}
        \item each $R_j$ is a $(K_X+B+G+\Mm_X)$-negative extremal ray, and
        \item $\phi_j$ is a $(K_X+B+G+\Mm_X)$-negative extremal contraction, and if $\phi_j$ is small, then $\psi_j$ is a $(K_X+B+G+\Mm_X)$-flip.
    \end{enumerate}
\end{enumerate}
\end{enumerate}
\end{lem}
\begin{proof}
(1) is obvious.

Pick a $(K_{\Ff}+B+\Mm_X)$-negative extremal ray $R$. By our assumption, there exists a prime divisor $E$ on $X$ such that $R$ is contained in the image of $\overline{NE}(E/U)\rightarrow\overline{NE}(X/U)$ and $\mult_EB=\epsilon_{\Ff}(E)$. We let $S$ be the normalization of $E$. Then there exists a natural surjection $\overline{NE}(S/U)\rightarrow\overline{NE}(E/U)$, and thus $R$ is contained in the image of $\iota: \overline{NE}(S/U)\rightarrow\overline{NE}(E/U)\rightarrow\overline{NE}(X/U)$. We claim that there exists an extremal ray $R_S$ in $\overline{NE}(S/U)$ such that $R=\iota(R_S)$. In fact, there exist extremal rays $R_i'$ in $\overline{NE}(S/U)$ such that $\iota(\sum a_iR_i')=R$ for some $a_i>0$. Since $R$ is extremal$/U$, either $\iota(R_i')=R$ or $\iota(R_i')=0$ for each $i$. Since $R\not=0$, there exists $j$ such that $\iota(R_j')\not=0$. We may take $R_S:=R_j'$ and the claim holds.

Let $\Ff_S$ be the restricted foliation of $\Ff$ on $S$ which is algebraically integrable by Proposition \ref{prop: a.i preserved adjunction}, $\Mm^S:=\Mm|_S$, and
$$K_{\Ff_S}+B_S+\Mm^S_S:=(K_{\Ff}+B+\Mm_X)|_S.$$
Then $R_S$ is a $(K_{\Ff}+B+\Mm_X)|_S$-negative extremal ray. By Theorem \ref{thm: adjunction foliation nonnqc}, $(S,\Ff_S,B_S,\Mm^S)/U$ is an lc gfq. Since we assume Theorem \ref{thm: cone theorem induction} in dimension $\leq d-1$, $R_S$ is spanned by a rational curve $C$ such that $C$ is tangent to $\Ff_S$ and
$$0<-\left(K_{\Ff_S}+B_S+\Mm^S_S\right)\cdot C\leq 2(d-1).$$
We identify $C$ with its image in $X$ under the natural inclusion $S\rightarrow E\rightarrow X$. Then $C$ spans $R$ and
$$0<-\left(K_{\Ff_S}+B_S+\Mm^S_S\right)\cdot C=-\left(K_{\Ff}+B+\Mm_X\right)\cdot C\leq 2(d-1).$$
Moreover, by \cite[Lemma 3.3]{ACSS21}, $C$ is tangent to $\Ff$ and is contracted by $f$. This implies (2).

For any curve $C_j'$ such that $[C_j']\in R$, we let $C''$ be any irreducible component of $C_j'$. Since $R$ is extremal, $[C'']\in R$ which implies that $C\equiv\lambda C''$ for some positive rational number $\lambda$. Thus $C''$ is contracted by $f$ and therefore $C_j'$ is contracted by $f$, and we get (3).

Now we may assume that $f$ is equi-dimensional, $A$ is an ample$/U$ $\Rr$-divisor on $X$, and $G$ is a divisor associated to $(X,\Ff,B,\Mm)/Z$. By Lemma \ref{lem: basic property (*) gpair} and Proposition \ref{prop: weak cbf gfq}, $(X,B+G,\Mm)$ is lc and
$$K_{\Ff}+B+\Mm_X\sim_{\mathbb R,Z}K_{X}+B+G+\Mm_X.$$
By (3), $R_j$ is a $(K_{X}+B+G+\Mm_X)$-negative extremal ray$/U$ for any $j\in\Lambda$. 

If $\Mm$ is NQC$/U$, then by \cite[Theorem 1.1(3)]{HL21a}, $\Lambda_A$ is a finite set, and thus $\Lambda$ is a countable set. This implies (4.a) and (4.b). (4.c.i) follows from \cite[Theorem 1.5]{Xie22} (see also \cite[Theorem 1.7]{CLX23}), and (4.c.ii) follows from \cite[Theorem 1.2]{LX23b}. 

Suppose that $X$ is $\Qq$-factorial klt. Possibly adding the pull-back of some very ample divisor to $A$, we may assume that $A$ is ample. By \cite[Lemma 3.4]{HL22}, there exists an $\Rr$-divisor $\Delta_A\ge0$ such that $\Delta_A\sim_{\mathbb R}B+G+A+\Mm_X$ and $(X,\Delta_A)$ is klt. It follows that
$$\Lambda_A=\{j\in\Lambda\mid(K_{\Ff}+\Delta)\cdot R_j<0\}$$
is a finite set by the classical cone theorem (cf. \cite[Theorem 4-2-1]{KMM87}, \cite[Theorem 4.5.2]{Fuj17}), and $\Lambda=\cup_{n=1}^{+\infty}\Lambda_{\frac{1}{n}A}$ is a countable set. This implies (4.a) and (4.b). (4.c.i) follows from the classical contraction theorem (cf. \cite[Theorem 3-2-1]{KMM87}, \cite[Theorem 4.5.2]{Fuj17}) and (4.c.ii) follows from the existence of flips \cite[Corollary 1.4.1]{BCHM10}.

(4.d) follows immediately from (4.c).
\end{proof}

\begin{prop}\label{prop: cone d-1 imply * dim d part 1}
Let $d$ be a positive integer. Assume that Theorem \ref{thm: cone theorem induction} holds in dimension $\leq d-1$. 

Let $(X,\Ff,B,\Mm)/U$ a gfq of dimension $d$ such that $\Ff$ is algebraically integrable. Let $E_1,\dots,E_s$ be lc places of $(X,\Ff,B,\Mm)$ and $T$ a reduced $\Ff$-invariant divisor. Then $(X,\Ff,B,\Mm)$ has a great ACSS model $(Y,\Ff_Y,B_Y,\Mm;G_Y)$ such that
\begin{enumerate}
\item $G_Y$ contains the strict transform of $T$ on $Y$, and
\item $E_1,\dots,E_s$ are on $Y$ if $(X,\Ff,B,\Mm)$ is lc.
\end{enumerate}
\end{prop}
\begin{proof}
Let $f: X\dashrightarrow \tilde Z$ be a dominant map which induces $\Ff$ and $g: \bar X\rightarrow X$ a birational morphism such that the induced map $f\circ g$ is a morphism. By Lemma \ref{lem: stein induce same foliation}, possibly replacing $\tilde Z$, we may assume that $f\circ g$ is a contraction. Let $\bar B:=\Supp(g^{-1}_*B)\cup\Supp\Exc(g)\cup\Supp\Exc(g^{-1}_*T)$.
    
By Definition-Theorem \ref{defthm: weak ss reduction} there exists an equi-dimensional model $f': (X',\Sigma_{X'},\Mm)\rightarrow (Z,\Sigma_{Z})$ of $f\circ g: (\bar X,\bar B,\Mm)\rightarrow \tilde Z$ associated with $h: X'\rightarrow \bar X$ and $h_Z: Z\rightarrow \tilde Z$, such that $E_1,\dots,E_s$ are on $X'$. Let $\phi:=g\circ h$, $\Ff':=\phi^{-1}\Ff$, and 
$$B':=\phi^{-1}_*(B\wedge\Supp B)+(\Supp\Exc(\phi))^{\Ff'}.$$
Then $\phi$ is a foliated log resolution of $(X,\Ff,B,\Mm)$, and any component of $B'$ is a component of the horizontal$/Z$ part of $\Sigma_{X'}$. We let $H$ be a very ample divisor on $Z$, $H_1,\dots,H_{2d+1}\in |H|$ general elements, $\Sigma':=\Sigma_{X'}+\sum_{i=1}^{2d+1}f'^*H_i$, $G'$ the vertical$/Z$ part of $\Sigma'$, and $B_Z:=\Sigma_Z+\sum_{i=1}^{2d+1}H_i$. By \cite[Proposition 3.2]{AK00}, $f': (X',\Sigma',\Mm)\rightarrow (Z,B_Z)$ is toroidal. By Lemma \ref{lem: fls imply acss}, $(X',\Ff',B',\Mm;G')/Z$ is $\Qq$-factorial ACSS.
        
Since $G'\geq\sum_{i=1}^{2d+1}f'^*H_i$, $G'$ is super$/Z$. Moreover, any $\Ff'$-invariant $\phi$-exceptional prime divisor is contained in $G'$, and the strict transform of $T$ on $Y$ is contained in $G'$.

\begin{claim}\label{claim: induction run mmp with scaling}
Let $A$ be an ample  $\Rr$-divisor on $X$. Then we may run a $(K_{\Ff'}+B'+\Mm_{X'})$-MMP$/X$ 
  \begin{center}
 $\xymatrix{(X_0,\Ff_0,B_0,\Mm;G_0)\ar@{-->}[r]^{\psi_0} & (X_1,\Ff_1,B_1,\Mm;G_1)\ar@{-->}[r]^{\ \ \ \ \ \ \ \ \ \ \psi_1} & \dots\ar@{-->}[r] & (X_n,\Ff_n,B_n,\Mm;G_n)\ar@{-->}[r]^{\ \ \ \ \ \ \ \ \ \ \psi_n} & \dots 
}$
  \end{center}
where $(X_0,\Ff_0,B_0,\Mm;G_0):=(X',\Ff',B',\Mm;G')$, so that the following conditions are satisfied for each $i$. Let $A_i$ be the strict transform of $A$ on $X_i$.
\begin{enumerate}
    \item There exists a contraction $f_i: X_i\rightarrow Z$ such that $f_{i+1}=f_i\circ\psi_i$.
    \item There exists a contraction $\phi_i: X_i\rightarrow X$ such that $\phi_{i+1}=\phi_i\circ\psi_i$.
    \item $(X_i,\Ff_i,B_i,\Mm;G_i)/Z$ is $\Qq$-factorial 
 super ACSS.
    \item For any $(K_{\Ff_i}+B_i+\Mm_{X_i})$-negative extremal ray$/X$ $R$, there exists a prime divisor $F$ on $X_i$ such that $R$ is contained in ${\rm Im}\left(\overline{NE}(F/U)\rightarrow\overline{NE}(X/U)\right)$ and $\mult_FB_i=\epsilon_{\Ff_i}(F)$.
    \item For any extremal ray$/X$ $R$ on $X_i$ such that $R$ is either a $(K_{\Ff_i}+B_i+\Mm_{X_i})$-negative extremal ray or a $(K_{X_i}+B_i+G_i+\Mm_{X_i})$-negative extremal ray, 
    \begin{enumerate}
      \item $R$ is an extremal ray$/Z$,
      \item $(K_{\Ff_i}+B_i+\Mm_{X_i})\cdot R=(K_{X_i}+B_i+G_i+\Mm_{X_i})\cdot R$, and
      \item $R$ is a $(K_{\Ff_i}+B_i+\Mm_{X_i})$-negative extremal ray if and only if $R$ is a $(K_{X_i}+B_i+G_i+\Mm_{X_i})$-negative extremal ray.
    \end{enumerate}
    \item $\psi_i$ is a step of a $(K_{X_i}+B_i+G_i+\Mm_{X_i})$-MMP$/X$ with scaling of $A_i$ as well as a $(K_{\Ff_i}+B_i+\Mm_{X_i})$-MMP$/X$ with scaling of $A_i$.
    \item $\psi_i$ is a step of a $(K_{\Ff_i}+B_i+\Mm_{X_i})$-MMP$/Z$ as well as a step of a $(K_{X_i}+B_i+G_i+\Mm_{X_i})$-MMP$/Z$.
\end{enumerate}
Moreover, there exists a positive integer $m$ satisfying the following.
\begin{enumerate}
    \item[(8)] The induced birational map $X_0\dashrightarrow X_m$ contracts any $\phi$-exceptional prime divisor $F$ such that $a(F,\Ff,B,\Mm)>-\epsilon_{\Ff}(F)$. 
    \item[(9)] If $(X,\Ff,B,\Mm)$ is lc, then any divisor $F$ contracted by $X_0\dashrightarrow X_m$ satisfies that $a(F,\Ff,B,\Mm)>-\epsilon_{\Ff}(F)$. 
\end{enumerate}
    \end{claim}
\begin{proof}
 
\noindent\textbf{Step 1}. In this step, we prove (1-4) for $i=0$. (1) We have $f_0:=f$. (2) We have $\phi_0:\phi$. (3) It follows from our construction. (4) The image of $R$ on $X$ is a closed point, so $R$ is contained in a $\phi$-exceptional divisor $F$. By our construction, $\mult_FB_0=\epsilon_{\Ff_0}(F)$. 

\medskip

\noindent\textbf{Step 2}. In this step, we prove that (1-4) for $i=n$ implies (5) for $i=n$.

First we prove (5.a). Assume that $R$ is a $(K_{\Ff_n}+B_n+\Mm_{X_n})$-negative extremal ray$/X$. Then by Lemma \ref{lem: induction cone 1}(2), $R$ is a $(K_{\Ff_n}+B_n+\Mm_{X_n})$-negative extremal ray$/Z$. Now assume that $R$ is a $(K_{X_n}+B_n+G_n+\Mm_{X_n})$-negative extremal ray$/X$. Since $G_0\geq  \sum_{j=1}^{2d+1} f_0^*H_j$, $L_n:=G_n-\sum_{j=1}^{2d+1} f_n^*H_j\geq 0$. 
By (3), $(X_n,B_n+G_n,\Mm)$ is $\Qq$-factorial lc and $X$ is klt, so $(X_n,B_n+L_n,\Mm)$ is $\Qq$-factorial lc. By the length of extremal rays for lc g-pairs over $\Qq$-factorial klt varieties (cf. \cite[Proposition 3.17]{HL22}), $R$ is spanned by a rational curve $C$ such that
$$0>(K_{X_n}+B_n+G_n+\Mm_{X_n})\cdot C=(K_{X_n}+B_n+L_n+\Mm_{X_n})\cdot C+\left(\sum_{j=1}^{2d+1}f_n^*H_j\right)\cdot C\geq -2d.$$
Therefore, $f_n^*H_j\cdot C=0$ for each $j$, so $R$ is an extremal ray$/Z$. This implies (5.a). 

(5.b) follows from (6.a) and Proposition \ref{prop: weak cbf gfq}, and (5.c) follows from (5.b). Thus (5) holds.

\medskip

\noindent\textbf{Step 4}. In this step, we prove that (1-5) for $i=n$ and (1-7) for $i\leq n-1$ imply (6) and (7) for $i=n$, and also imply (1)(2) for $i=n+1$.

By induction hypothesis, the induced birational map $X_0\dashrightarrow X_n$ is a sequence of steps of a $(K_{X_0}+B_0+G_0+\Mm_{X_0})$-MMP$/X$ with scaling of $A$. By Lemma \ref{lem: scaling number go to 0}, either this MMP already terminates at $X_n$ and we are done, or we may run the next step of this $(K_{X_0}+B_0+G_0+\Mm_{X_0})$-MMP$/X$ with scaling of $A$, which is a step of a $(K_{X_n}+B_n+G_n+\Mm_{X_n})$-MMP$/X$ with scaling of $A_n$. (6) and (7) for $i=n$ now follow from (5) for $i=n$. (1) for $i=n+1$ follows from (7) for $i=n$ and (2) for $i=n+1$ follows from (6) for $i=n$.

\medskip

\noindent\textbf{Step 5}. In this step, we prove that (1-7) for $i\leq n-1$ and (1)(2) for $i=n$ imply (3) for $i=n$. 

By (3)(7) for $i=n-1$, $X_n$ is $\Qq$-factorial. By (1) for $i=n$ and (3) for $i=n-1$, $G_n$ is super$/Z$. So we only need to show that $(X_n,\Ff_n,B_n,\Mm;G_n)/Z$ is ACSS. We check conditions (1-4) of Definition \ref{defn: ACSS f-triple} for $(X_n,\Ff_n,B_n,\Mm;G_n)/Z$. 

Definition \ref{defn: ACSS f-triple}(1) for $(X_n,\Ff_n,B_n,\Mm;G_n)/Z$: By (6) for $i=n-1$ and Proposition \ref{prop: MMP preserves *}, $(X_n,B_n+G_n,\Mm)/Z$ satisfies Property $(*)$. Since $\Ff_{n-1}$ is induced by $f_{n-1}$, $\Ff_n$ is induced by $f_n$. Since $G_{n-1}\geq 0$ is $\Ff_{n-1}$-invariant, $G_n\geq 0$ if $\Ff_n$-invariant. Thus $(X_n,\Ff_n,B_n,\Mm;G_n)/Z$ satisfies Property $(*)$. By (3)(7) for $i=n-1$, $(X_n,\Ff_n,B_n,\Mm)$ is lc, so Definition \ref{defn: ACSS f-triple}(1) holds for $(X_n,\Ff_n,B_n,\Mm;G_n)/Z.$

Definition \ref{defn: ACSS f-triple}(2) for $(X_n,\Ff_n,B_n,\Mm;G_n)/Z$: it immediately follows from (3)(6) for $i=n-1$ and Proposition \ref{prop: MMP preserves *}.

Definition \ref{defn: ACSS f-triple}(3) for $(X_n,\Ff_n,B_n,\Mm;G_n)/Z$: For any divisor $\Sigma$ on $Z$ such that $\Sigma\geq B_Z$ and $(Z,\Sigma)$ is log smooth, there exists $D$ and $\Nn$ such that
$(X_{n-1},B_{n-1}+G_{n-1}+D+f_{n-1}^*(\Sigma-B_Z),\Nn)$ is qdlt, $\Supp\{B_{n-1}\}\subset\Supp D$, and $\Nn-\alpha\Mm$ is nef$/X$ for some $\alpha>1.$
Let $\Pp:=\Nn-\Mm$. By (7) for $i=n-1$, $\psi_{n-1}$ is also a step of a $(K_{\Ff_{n-1}}+B_{n-1}+f_{n-1}^*(\Sigma-B_Z)+\Mm_{X_{n-1}})\text{-MMP}/Z,$ hence a step of a $(K_{\Ff_{n-1}}+B_{n-1}+\delta D+f_{n-1}^*(\Sigma-B_Z)+\Mm_{X_{n-1}}+\delta\Pp_{X_{n-1}})\text{-MMP}/Z$ for any $0<\delta\le 1$. 
Note that $(X_{n-1},B_{n-1}+G_{n-1}+\delta D+f_{n-1}^*(\Sigma-B_Z),\Mm+\delta\Pp)$ is qdlt. By Lemma \ref{lem: mmp preserves qdlt}, $(X_{n},B_{n}+G_{n}+\delta(\psi_{n-1})_*D+f_{n}^*(\Sigma-B_Z),\Mm+\delta\Pp)$ is also qdlt.

Definition \ref{defn: ACSS f-triple}(4) for $(X_n,\Ff_n,B_n,\Mm;G_n)/Z$: For any lc place $S$ of $(X_n,\Ff_n,B_n,\Mm)$ with generic point $\eta_W$, we have
      $$-\epsilon_{\Ff}(S)=a(S,\Ff_n,B_n,\Mm)\geq a(S,\Ff_{n-1},B_{n-1},\Mm)\geq -\epsilon_{\Ff}(S).$$
      Therefore, $S$ is an lc place of $a(S,\Ff_{n-1},B_{n-1},\Mm)$. Thus $\psi_{n-1}$ is an isomorphism near the generic point of $\Center_{X_{n-1}}S$. Since Definition \ref{defn: ACSS f-triple}(4) is a property near the generic point of lc places,       
Definition \ref{defn: ACSS f-triple}(4) holds for $(X_n,\Ff_n,B_n,\Mm;G_n)/Z$.

Therefore, $(X_n,\Ff_n,B_n,\Mm;G_n)/Z$ is $\Qq$-factorial ACSS.

\medskip

\noindent\textbf{Step 6}. In this step, we prove (4) for $i=n$ assuming that (1-7) hold for $i=n-1$, hence conclude the proof of (1-7). For any $(K_{\Ff_n}+B_n+\Mm_{X_n})$-negative extremal ray$/X$ $R$, $R$ is contained in a prime $\phi_n$-exceptional divisor $F$. Let $F'$ be the strict transform of $F$ on $X'$. Then $F'$ is a prime $\phi$-exceptional divisor, so 
$$\mult_FB_n=\mult_{F'}B_0=\epsilon_{\Ff'}(E)=\epsilon_{\Ff_n}B_n.$$
This implies (4). By Lemma \ref{lem: induction cone 1}(4.c), we may construct $\psi_n$ and get (5)(6).

By induction, (1-7) hold.

\medskip

\noindent\textbf{Step 7}. In this step, we prove (8) and (9) and conclude the proof of the claim.

If this MMP terminates, then we let $m$ be the index such that $(X_m,\Ff_m,B_m,\Mm;G_m)$ is the last output of this MMP. In particular, $K_{\Ff_m}+B_m+\Mm_{X_m}$ is nef$/X$. In particular, it is movable$/X$. If this MMP does not terminate, then we let $m$ be the index such that $\psi_i$ is a flip for any $i\geq m$. By (5),
$$\lambda_i:=\inf\{t\geq 0\mid K_{\Ff_i}+B_i+\Mm_{X_i}+tA_i\text{ is nef}/X\}=\inf\{t\geq 0\mid K_{X_i}+B_i+G_i+\Mm_{X_i}+tA_i\text{ is nef}/X\}.$$
Then by Lemma \ref{lem: scaling number go to 0}, $\lim_{i\rightarrow+\infty}\lambda_i=0$. Therefore,
$$K_{\Ff_m}+B_m+\Mm_{X_m}=\lim_{i\rightarrow+\infty}(\psi_i)^{-1}_*(K_{\Ff_i}+B_i+\Mm_{X_i}+\lambda_iA_i)$$
is a movable$/X$, where $\psi_i: X_m\dashrightarrow X_i$ is the induced birational map.
Let $F_1,\dots,F_l$ be the $\phi_m$-exceptional prime divisors and let $a_k:=a(F_k,\Ff,B,\Mm)+\epsilon_{\Ff_m}(F_k)$ for each $k$. Then
\begin{align*}
   &K_{\Ff_m}+B_m+\Mm_{X_m}\\
   &=\phi_m^*(K_{\Ff}+B+\Mm_X)+\sum_{k}a_kF_k-\sum_{D|\mult_DB>1}(\mult_DB-1)(\phi_m^{-1})_*D\\&
    \sim_{\mathbb R,X}\sum_{k}a_kF_k-\sum_{D|\mult_DB>1}(\mult_DB-1)(\phi_m^{-1})_*D.
\end{align*}
By \cite[Lemma 3.3]{Bir12}, $a_k\leq 0$ for any $k$. This implies (8). 

If $(X,\Ff,B,\Mm)$ is lc, then
$$K_{\Ff'}+B'+\Mm_{X'}\sim_{\mathbb R,X}\sum_{F\mid F\subset\Exc(\phi)}(\epsilon_{\Ff}(F)+a(F,\Ff,X,\Mm))F\geq 0$$
and (9) follows. This completes the proof of the claim.
\end{proof}
\noindent\textit{Proof of Proposition \ref{prop: cone d-1 imply * dim d part 1} continued}. 
Now $(X_m,\Ff_m,B_m,\Mm;G_m)/Z$ is a super ACSS model of $(X,\Ff,B,\Mm)$ such that $E_1,\dots,E_s$ are on $X_m$ if $(X,\Ff,B,\Mm)$ is lc. Since any $\Ff'$-invariant exceptional$/X$ prime divisor is contained in $G'$, any $\Ff_m$-invariant exceptional$/X$ prime divisor is contained in $G_m$. Therefore, $(X_m,\Ff_m,B_m,\Mm;G_m)/Z$ is a great ACSS model of $(X,\Ff,B,\Mm)$. Finally, since the strict transform of $T$ on $X'$ is contained in $G'$, the strict transform of $T$ on $X_m$ is contained in $G_m$. The proposition follows by taking $(Y,\Ff_Y,B_Y,\Mm;G_Y):=(X_m,\Ff_m,B_m,\Mm;G_m).$
\end{proof}

\begin{prop}\label{prop: cone d-1 imply * dim d part 2}
Let $d$ be a positive integer. Assume that Theorem \ref{thm: cone theorem induction} holds in dimension $\leq d-1$. Then Theorem \ref{thm: property * induction} holds in dimension $d$.
\end{prop}
\begin{proof}
Notation and assumptions as in Theorem \ref{thm: property * induction}. By Proposition \ref{prop: cone d-1 imply * dim d part 1}, $(X,\Ff,B,\Mm)$ has a great ACSS model $(Y',\Ff_{Y'},B_{Y'},\Mm;G_{Y'})/Z$. 
Moreover, $E_1,\dots,E_s$ are also lc places of  $(Y',\Ff_{Y'},B_{Y'},\Mm)$. By Proposition \ref{prop: cone d-1 imply * dim d part 1} again, $(Y',\Ff_{Y'},B_{Y'},\Mm)$ has a great ACSS model $(Y,\Ff_{Y},B_{Y},\Mm;G_Y)$ such that $E_1,\dots,E_s$ are on $Y$ and $G_Y$ contains the strict transform of $G_{Y'}$ on $Y$. Therefore, $G_Y$ contains all $\Ff_Y$-exceptional prime divisors. Since 
 $$g^*(K_{\Ff}+B+\Mm_X)\geq  K_{\Ff_{Y'}}+B_{Y'}+\Mm_{Y'}$$
 the induced birational morphism $Y\rightarrow X$ is a great ACSS modification $(X,\Ff,B,\Mm)$.
\end{proof}

The following lemma is well-known to experts. For the reader's convenience, we conclude a proof here.

\begin{lem}\label{lem: supporting function are +A}
    Let $X\rightarrow U$ be a projective morphism from a normal quasi-projective variety to a variety. Let $D$ be an $\Rr$-Cartier $\Rr$-divisor on $X$ and $R$ a $D$-negative extremal ray in $\overline{NE}(X/U)$. Then there exists an ample$/U$ $\Rr$-divisor $A$ on $X$ such that $H:=D+A$ is the supporting function of $R$.
\end{lem}
\begin{proof}
   Let $H_R$ be a supporting function of $R$. Then $H_R\cdot R=0$ and $H_R\cdot R'>0$ for any $R'\not=R$ in $\overline{NE}(X/U)$. Let 
   $$C:=\left\{D\in N^1(X/U)\mid D\cdot z\geq 0 \text{ for any }z\in\overline{NE}(X/U)_{D\geq0}\right\}.$$
   Then $C$ is the dual cone of $\overline{NE}(X/U)_{D\geq 0}$ and is generated by nef$/U$ divisors and $D$. Since $H_R$ is positive on $\overline{NE}(X/U)_{D\geq 0}\backslash\{0\}$, $H_R$ is contained in the interior of $C$. Therefore, there exists an ample$/U$ $\Rr$-divisor $\tilde A$ such that $H_R-\tilde A=L+pD$ in $N^1(X/U)$, where $L$ is a nef$/U$ $\Rr$-divisor and $p$ is a non-negative real number. Let $A':=\tilde A+L$. Then $A'$ is ample$/U$. We may let $H:=\frac{1}{p}H_R=\frac{1}{p}\tilde A'+D$ and $A:=\frac{1}{p} A'$.
\end{proof}

\begin{prop}\label{prop: * to cone}
    Let $d$ be a positive integer. Assume that Theorem \ref{thm: cone theorem induction} holds in dimension $\leq d-1$ and Theorem \ref{thm: property * induction} holds in dimension $d$. Then Theorem \ref{thm: cone theorem induction} holds in dimension $d$.
\end{prop}
\begin{proof}
\medskip

\noindent\textbf{Step 1}. In this step, we construct a supporting function$/U$ of $R$.

Notation and assumptions as in Theorem \ref{thm: cone theorem induction}. It is clear that
    $$\overline{NE}(X/U)=\overline{NE}(X/U)_{K_{\Ff}+B+\Mm_X\geq 0}+\overline{NE}(X/U)_{\Nlc(X,\Ff,B,\Mm)}+\sum_{j\in\Lambda} R_j$$
Let $R$ be a $(K_{\Ff}+B+\Mm_X)$-negative extremal ray$/U$ such that $R\not\subset\overline{NE}(X/U)_{\Nlc(X,\Ff,B,\Mm)}$. By Lemma \ref{lem: supporting function are +A}, there exists an ample$/U$ $\Rr$-divisor $A$ on $X$ such that
$$H_R:=K_{\Ff}+B+A+\Mm_X$$
is a supporting function$/U$ of $R$. In particular, $H_R\not\equiv_U 0$ is nef, $H_R\cdot R=0$, and $H_R\cdot R'>0$ for any $R'\in\overline{NE}(X)\backslash R$. 

\medskip

\noindent\textbf{Step 2}. In this step, we deal with the case when $H_R$ is not big$/U$. 

Let $F$ be the Stein factorization of a general fiber of the morphism $\pi: X\rightarrow U$. Then $H_F:=H_R|_F$ is nef, not big, and is not numerically trivial. Let $q:=\dim F$ and $A_F:=A|_F$. Then there exists an integer $1\leq k\leq q-1$ such that
$$H_F^k\cdot A_F^{q-k}>H_F^{k+1}\cdot A_F^{q-k-1}=0.$$
Let $D_i:=H_R$ for any $1\leq i\leq k+1$, and let $D_i:=A$ for any $k+2\leq i\leq q$. Then
$$(D_1|_F)\cdot (D_2|_F)\cdots\dots\cdot (D_q|_F)=H_F^{k+1}\cdot A_F^{q-k-1}=0$$
and
$$-(K_{\Ff}+B)|_F\cdot (D_2|_F)\cdots\dots\cdot (D_q|_F)=(A_F-H_F)\cdot H_F^{k}\cdot A_F^{q-k-1}=H_F^{k}\cdot A_F^{q-k}>0.$$
Let $M:=H_R+A=K_{\Ff}+B+2A+\Mm_X$, which is ample$/U$. By Theorem \ref{thm: relative bb}, for any general closed point $x\in X$, there exists a rational curve $C_x$ such that $x\in C_x$, $\pi(C_x)$ is a closed point, $0=D_1\cdot C_x=H_R\cdot C_x,$ and
\begin{align*}
    0<&-(K_{\Ff}+B+\Mm_X)\cdot C_x=M\cdot C_x\\
    \leq& 2d\cdot\frac{M|_F\cdot  (D_2|_F)\cdot\dots\cdot (D_q|_F)}{-K_{\Ff}|_F\cdot  (D_2|_F)\cdot\dots\cdot (D_q|_F)}=2d\cdot\frac{-(K_{\Ff}+B+\Mm_X)|_F\cdot H_F^k\cdot A_F^{q-k-1}}{-K_{\Ff}|_F\cdot H_F^k\cdot A_F^{q-k-1}}.
\end{align*}
Let $\Mm^F:=\Mm|_F$ and $B_F:=B|_F$. Since $F$ is a general fiber of $\pi$, $B_F\geq 0$ and $\Mm^F$ is nef. Thus $\Mm^F_F$ is pseudo-effective and $(B+\Mm_X)|_F\cdot H_F^k\cdot A_F^{q-k-1}\geq 0$. Therefore 
$$0<-(K_{\Ff}+B+\Mm_X)\cdot C_x\leq 2d.$$

\medskip

\noindent\textbf{Step 3}. From now on, we may assume that $H_R$ is big$/U$. In this step, we construct a set $\Ii$ of tuples $(W,\lambda)$ and show that it contains a minimal element. 

Since $H_R$ is big$/U$, $H_R=A'+P$ for some ample$/U$ $\Rr$-divisor $A'$ and some $\Rr$-divisor $P\geq 0$. In particular, $P$ is $\Rr$-Cartier and $P\cdot R<0$. Let $S$ be the normalization of $\Supp P$. Then $R$ is contained in the image of $\overline{NE}(S/U)\rightarrow\overline{NE}(X/U)$ induced by the natural inclusion $S\rightarrow \Supp P\rightarrow X.$ Consider the set $\Ii$ of all $(W,\lambda)$ such that
\begin{enumerate}
    \item $\lambda$ is a non-negative real number,
    \item $W$ is an lc center of $(X,\Ff,B+\lambda P,\Mm)$ with normalization $W^\nu$, and
    \item $R$ is contained in ${\rm Im}\left(\overline{NE}(W^\nu/U)\rightarrow\overline{NE}(X/U)\right)$  induced by the natural inclusion $W^\nu\rightarrow W\rightarrow X.$
\end{enumerate}
By construction, there exists a component $L$ of $S$ such that $(L,1)\in\Ii$. Thus $\Ii\not=\emptyset$. 

In the rest of this step, we show that there exists $(W_0,\lambda_0)\in\Ii$ that is minimal in the following way: for any $(W,\lambda)\in\Ii$, one of the following cases holds.
\begin{itemize}
\item $\lambda_0<\lambda$.
\item $\lambda_0=\lambda$ and $W_0\subsetneq W$.
\item  $(W,\lambda)=(W_0,\lambda_0)$.
\end{itemize}
Let $h: X'\rightarrow X$ be a foliated log resolution of $(X,\Ff,\Supp B\cup\Supp P,\Mm)$ with associated morphism $f: X'\rightarrow Z$. Then there exists a toroidal morphism $f: (X',\Sigma)\rightarrow (Z,\Sigma_Z)$ such that $h^{-1}(\Supp B\cup\Supp P)\cup\Supp\Exc(h)$ is contained in $\Sigma$. By Lemma \ref{lem: foliated log smooth imply lc}, for any $(W,\lambda)\in\Ii$, either $W$ is the image of a stratum of $(X',\Sigma)$ on $X$ or $\lambda=0$. Therefore, the set
\begin{align*}
    \Ii':=\{\lambda\mid& \text{ there exists an lc center of }(X,\Ff,B+\lambda P,\Mm)\\
    & \text{ that is not an lc center of }(X,\Ff,B+(\lambda-\delta)P,\Mm) \text{ for any }0<\delta\ll 1\}
\end{align*}
is finite, so we may let
$$\lambda_0:=\min\{\lambda\mid \text{ there exists }(W,\lambda)\in\Ii\}.$$
Now by the Noetherian property, there exists $(W_0,\lambda_0)\in\Ii$ such that $W_0\subset W$ for any $(W,\lambda_0)\in\Ii$. 

\medskip

\noindent\textbf{Step 4}. In this step, we construct an $\Rr$-divisor $\tilde B$ on $X$ and a $\Qq$-factorial ACSS model $(Y,\Ff_Y,\tilde B_Y,\Mm)$ of $(X,\Ff,B,\Mm)$, so that $R$ is the image of a $(K_{\Ff_Y}+\tilde B_Y+\Mm_Y)$-negative extremal ray$/U$ in $X$.

Let $\tilde B:=B+\lambda_0P$ and let $E$ be an lc place of $(X,\Ff,\tilde B,\Mm)$ such that $\Center_EX=W_0$. By Theorem \ref{thm: property * induction} in dimension $d$, there exists a great ACSS model $(Y,\Ff_Y,\tilde B_Y,\Mm;G)$ of $(X,\Ff,\tilde B,\Mm)$ such that $E$ is on $Y$ with induced birational morphism $g: Y\rightarrow X$. We have
$$K_{\Ff_Y}+\tilde B_Y+\Mm_Y+F=g^*\left(K_{\Ff}+\tilde B+\Mm_X\right)$$
for some $F\geq 0$. Let $V:=g(\Supp F)$. Then $V\subset\Nlc(X,\Ff,\tilde B,\Mm)$ is a reduced subscheme of $X$. 

    Let $C_i$ be a sequence of curves such that $[C_i]\in\overline{NE}(X/U)$ and $\lim [C_i]=R$. Then there exist curves $C_{Y,i}$ on $Y$ such that $g(C_{Y,i})=C_i$ for each $i$. Let $R':=\lim [C_{Y,i}]$. Then $g(R')=R$. Let $R_1,\dots,R_l$ be extremal rays in $\overline{NE}(Y/U)$ such that $R'=\sum a_iR_i$ for some $a_i>0$. Since $R$ is extremal, there exists $i$ such that $g(R_i)=R$, and let $R_Y:=R_i$. By the projection formula,
\begin{align*}
\left(K_{\Ff_Y}+\tilde B_Y+F+\Mm_Y\right)\cdot R_Y&=\left(K_{\Ff_Y}+\tilde B_Y+F+\Mm_Y\right)\cdot R'\\
        &=\lim \left(K_{\Ff_Y}+\tilde B_Y+F+\Mm_Y\right)\cdot C_i\\
        &=\lim \left(K_{\Ff}+\tilde B+\Mm_X\right)\cdot C_i=\left(K_{\Ff}+\tilde B+\Mm_X\right)\cdot R<0.
\end{align*}
    Thus $R_Y$ is a $(K_{\Ff_Y}+\tilde B_Y+F+\Mm_Y)$-negative extremal ray. 

    If $F\cdot R_Y<0$, then $R_Y$ is contained in the image of $\overline{NE}(\Supp F/U)\rightarrow\overline{NE}(Y/U)$. Then $R=g(R_Y)$ is contained in the image of $\overline{NE}(V/U)\rightarrow\overline{NE}(X/U)$. Thus there exists an irreducible component $V_0$ of $V$ such that $R$ is contained in the image of $\overline{NE}(V_0/U)\rightarrow\overline{NE}(X/U)$. Since $R$ is not contained in $\overline{NE}(X/U)_{\Nlc(X,\Ff,B,\Mm)}$, $V_0$ is not an lc center of $(X,\Ff,B,\Mm)$. Since $V_0\subset V=g(F)\subset\Nlc(X,\Ff,\tilde B,\Mm)$ and $\tilde B=B+\lambda_0P$, there exists a real number $0<\lambda_1<\lambda_0$ such that $V_0$ is an lc center of $\Nlc(X,\Ff,B+\lambda_1P,\Mm)$. This contradicts the minimality of $(W_0,\lambda_0)$, as $(V_0,\lambda_1)\in\Ii$ and $\lambda_1<\lambda_0$. Hence $F\cdot R_Y\geq 0$ and $R_Y$ is a  $(K_{\Ff_Y}+\tilde B_Y+\Mm_Y)$-negative extremal ray.

\medskip

\noindent\textbf{Step 5}. In this step, we prove the proposition assuming $X$ is $\Qq$-factorial. 

Assume that $X$ is $\Qq$-factorial. By \cite[Lemma 3.6.2]{BCHM10}, $\Exc(g)$ is a divisor, so $g^{-1}(W_0)$ is a divisor. Since $R\subset{\rm Im}(\overline{NE}(W/U)\rightarrow\overline{NE}(X/U))$ and $g(R_Y)=R$, there exists a divisor $E_0$ on $Y$ such that $R_Y\subset{\rm Im}(\overline{NE}(E_0/U)\rightarrow\overline{NE}(Y/U))$. Since $g$ is an ACSS modification of $(X,\Ff,\tilde B,\Mm)$, $E_0$ is an lc place of $(X,\Ff,\tilde B,\Mm)$ and an lc place of $(Y,\Ff_Y,\tilde B_Y,\Mm)$. 

Let $T$ be the normalization of $E_0$, let $\Ff_{T}:=\Ff_Y|_{T}$ be the restricted foliation, let $\Mm^{T}:=\Mm|_{T}$, and
$$K_{\Ff_T}+\tilde B_T+\Mm^T_T:=(K_{\Ff_Y}+\tilde B_Y+\Mm_Y)|_{T}.$$
Since $R_Y\subset{\rm Im}(\overline{NE}(E_0/U)\rightarrow\overline{NE}(Y/U))$, $R_Y$ is contained in the image of 
$$\iota: \overline{NE}(T/U)\rightarrow \overline{NE}(E_0/U)\rightarrow\overline{NE}(Y/U).$$ Therefore, there exists $\tilde R\in\overline{NE}(T/U)$ such that $\iota(\tilde R)=R_Y$. We write $\tilde R=\sum\tilde a_i\tilde R_i$ where each $\tilde R_i$ is an extremal ray in $\overline{NE}(T/U)$. Since $R_Y$ is extremal$/U$, for each $i$, either $\iota(\tilde R_i)=0$ or $\iota(\tilde R_i)=R_Y$. Since $R_Y\not=0$, there exists $i_0$ such that $\iota(\tilde R_{i_0})=R_Y$. We let $R_T:=\tilde R_{i_0}$. Thus $R_T$ is a $(K_{\Ff_T}+B_T+\Mm^T_T)$-negative extremal ray$/U$. By Theorem \ref{thm: adjunction foliation nonnqc} and Theorem \ref{thm: cone theorem induction} in dimension $\leq d-1$, $R_T$ is spanned by a rational curve $C_T$ such that $C_T$ is tangent to $\Ff_T$ and
$$0<-(K_{\Ff_T}+\tilde B_T+\Mm^T_T)\cdot C_T\leq 2(d-1).$$
Let $C_Y$ be the image of $C_T$ in $Y$. Then $C_Y$ spans $R_Y$,
$$0<-(K_{\Ff_T}+\tilde B_T+\Mm^T_T)\cdot C_T=-(K_{\Ff_Y}+\tilde B_Y+\Mm_Y)\cdot C_Y\leq 2(d-1),$$
and by \cite[Lemma 3.3(4)]{ACSS21}, $C_Y$ is tangent to $\Ff_Y$. Let $C:=g(C_Y)$. Then $C$ is tangent to $\Ff$. By \textbf{Step 4}, $F\cdot C_Y\geq 0$, so
$$2d\geq -(K_{\Ff_Y}+\tilde B_Y+\Mm_Y)\cdot C_Y\geq -(K_{\Ff_Y}+\tilde B_Y+F+\Mm_Y)\cdot C_Y=-(K_{\Ff}+B+\Mm_X)\cdot C>0.$$
We are done for the case when $X$ is $\Qq$-factorial.

\medskip

\noindent\textbf{Step 6}. In this step, we conclude the proof of the theorem. 

Since $Y$ is $\Qq$-factorial and $R_Y$ is a $(K_{\Ff_Y}+\tilde B_Y+\Mm_Y)$-negative extremal ray, by \textbf{Step 5}, $R_Y$ is spanned by a rational curve $C_Y$ that is tangent to $\Ff_Y$ and 
$$0<-(K_{\Ff_Y}+\tilde B_Y+\Mm_Y)\cdot C_Y\leq 2d.$$
 Let $C:=g(C_Y)$. Then $C$ is tangent to $\Ff$. Since $F\cdot C_Y\geq 0$,
$$2d\geq -(K_{\Ff_Y}+\tilde B_Y+\Mm_Y)\cdot C_Y\geq -(K_{\Ff_Y}+\tilde B_Y+F+\Mm_Y)\cdot C_Y=-(K_{\Ff}+B+\Mm_X)\cdot C>0.$$
This concludes the proof of the theorem.
\end{proof}

\subsection{Proofs of  Theorems \ref{thm:  ACSS model}, \ref{thm: precise adj gfq}, \ref{thm: dcc adjunction is dcc}, and \ref{thm: lc adjunction foliation nonnqc}}\label{subsec: proof of adj}
\begin{proof}[Proofs of Theorems \ref{thm: cone theorem induction} and \ref{thm: property * induction}]
    It is obvious that Theorems \ref{thm: cone theorem induction} and \ref{thm: property * induction} hold trivially when $d=1$, so by Propositions \ref{prop: cone d-1 imply * dim d part 2} and \ref{prop: * to cone}, Theorems \ref{thm: cone theorem induction} and \ref{thm: property * induction} hold.
\end{proof}

\begin{proof}[Proof of Theorem \ref{thm:  ACSS model}]
    It is a special case of Theorem \ref{thm: property * induction}.
\end{proof}

\begin{proof}[Proof of Theorems \ref{thm: lc adjunction foliation nonnqc} and \ref{thm: precise adj gfq}]
By Theorem \ref{thm:  ACSS model}, there exists a $\Qq$-factorial ACSS model $(X',\Ff',B',\Mm)$ of $(X,\Ff,B,\Mm)$ with induced birational morphism $h: X'\rightarrow X$ and properly associated $f: X'\rightarrow Z$. Since  $(X,\Ff,B,\Mm)$ is lc,
$$K_{\Ff'}+B'+\Mm_{X'}=h^*(K_{\Ff}+B+\Mm_X).$$
Let $S'$ be the normalization of $h^{-1}_*S$, let $\Ff_{S'}$ be the restricted foliation of $\Ff$ on $S'$, and let $\Mm^S:=\Mm|_{S^\nu}$. Then there exists an induced birational morphism $h_S: S'\rightarrow S^\nu$. Let
$$K_{\Ff_{S'}}+B_{S'}+\Mm^S_{S'}:=(K_{\Ff'}+B'+\Mm_{X'})|_{S'},$$
then $(S',\Ff_{S'},B_{S'},\Mm^S)$ is lc and
$$K_{\Ff_{S'}}+B_{S'}+\Mm^S_{S'}=h_S^*(K_{\Ff_S}+B_S+\Mm^S_{S^\nu}).$$
Then Theorems \ref{thm: lc adjunction foliation nonnqc} and \ref{thm: precise adj gfq} hold by Theorem \ref{thm: adjunction foliation nonnqc} and Theorem \ref{thm: precise adjunction when induced}.
\end{proof}

\begin{proof}[Proof of Theorem \ref{thm: dcc adjunction is dcc}]
    It is an immediate corollary of Theorem \ref{thm: precise adj gfq}.
\end{proof}

\subsection{Proof of Theorem \ref{thm: cone theorem gfq}}\label{subsec: proof of cone}

\begin{lem}\label{lem: cone finiteness rays}
    Let $(X,\Ff,B,\Mm)/U$ be a gfq and $A$ an ample$/U$ $\Rr$-divisor on $X$. Then there are finitely many $(K_{\Ff}+B+A+\Mm_X)$-negative extremal rays$/U$ that are not contained in $\overline{NE}(X/U)_{\Nlc(X,\Ff,B,\Mm)}$.
\end{lem}
\begin{proof}
Let $d:=\dim X$, let $\omega:=K_{\Ff}+B+\Mm_X$, let $\rho:=\rho(X/U)$, and let $A_1,\dots,A_{\rho-1}$ be ample$/U$ Cartier divisors on $X$ such that $\omega,A_1,\dots,A_{\rho-1}$ form a basis of $N^1_{\mathbb R}(X/U)$. Let $0<\epsilon\ll 1$ be a rational number such that $A-\epsilon\sum_{i=1}^{\rho-1}A_i$ is ample$/U$. Then we only need to show that there are finitely many $(K_{\Ff}+B+\epsilon\sum_{i=1}^{\rho-1}A_i+\Mm_X)$-negative extremal rays$/U$ that are not contained in $\overline{NE}(X/U)_{\Nlc(X,\Ff,B,\Mm)}$. Possibly replacing $A$, we may assume that $A=\epsilon\sum_{i=1}^{\rho-1}A_i$.

Suppose that the lemma does not hold. Then there exist an infinite set $\Lambda$ and an infinite set $\{R_j\}_{j\in\Lambda}$ of $(K_{\Ff}+B+A+\Mm_X)$-negative extremal rays$/U$ that are not contained in $\overline{NE}(X/U)_{\Nlc(X,\Ff,B,\Mm)}$. By Theorem \ref{thm: cone theorem induction}, for any $j\in\Lambda$, there exists a rational curve $C_j$ on $X$ that is tangent to $\Ff$ and $R_j=[C_j]$, such that $-2d\leq \omega\cdot C_j<0$. For each $j\in\Lambda$, by Lemma \ref{lem: supporting function are +A}, there exists an ample$/U$ $\Rr$-divisor $L_j$ and a nef$/U$ $\Rr$-divisor $H_j$ such that
$$H_j=L_j+(K_{\Ff}+B+A+\Mm_X)=L_j+\epsilon\sum_{i=1}^{\rho-1}A_i+\omega$$
and $H_j$ is the supporting function of $R_j$. We have
$$0=H_j\cdot C_j=L_j\cdot C_j+\epsilon\sum_{i=1}^{\rho-1}A_i\cdot C_j+\omega\cdot C_j\geq-2d+\epsilon\sum_{i=1}^{\rho-1}A_i\cdot C_j.$$
Therefore, $A_i\cdot C_j\leq\frac{2d}{\epsilon}$ for any $i,j$. Since $A_i\cdot C_j\in\mathbb N^+$, there are finitely many possibilities of $A_i\cdot C_j$. Possibly replacing $\Lambda$ with an infinite subset, we may assume that $A_i\cdot C_j=A_i\cdot C_{j'}$ for any $i$ and any $j,j'\in\Lambda$.

We may write $\omega=\sum_{i=1}^c r_iD_i$ such that $r_1,\dots,r_c$ are linearly independent over $\Qq$ and $D_i$ are Weil divisors. By \cite[Lemma 5.3]{HLS19}, each $D_i$ is a $\Qq$-Cartier divisor. Thus there exist real numbers $a_{i,k}$ and $b_i$ such that
$$D_i\equiv_U\sum_{k=1}^{\rho-1}a_{i,k}A_k+b_i\omega$$
for each $i$. 
    
We let $\delta_1,\dots,\delta_c$ be real numbers such that $\sum_{i=1}^cb_i\delta_i>-1$ and $r_i':=\delta_i+r_i\in\mathbb Q.$ Let $\omega':=\sum_{i=1}^cr_i'D_i$. Then 
$$\omega'=\omega+\sum_{i=1}^c\delta_iD_i=\left(\sum_{k=1}^{\rho-1}\left(\sum_{i=1}^c\delta_ia_{i,k}\right)A_k\right)+\left(1+\sum_{i=1}^c\delta_ib_i\right)\omega.$$
Since $\sum_{i=1}^cb_i\delta_i>-1$, $\omega'$ and $A_1,\dots,A_{\rho-1}$ form a basis of $N^1_\Rr(X/U)$. Moreover,
$$\omega'\cdot C_j=\left(\sum_{i=1}^c\sum_{k=1}^{\rho-1}\delta_ia_{i,k}\cdot(A_k\cdot C_j)\right)+\left(1+\sum_{i=1}^c\delta_ib_i\right)(\omega\cdot C_j).$$
By our assumptions,
$$\alpha:=\sum_{i=1}^c\sum_{k=1}^{\rho-1}\delta_ia_{i,k}\cdot(A_k\cdot C_j)\text{ and }\beta:=1+\sum_{i=1}^c\delta_ib_i>0$$
are constants which do not depend on $j$, and $\omega\cdot C_j\in [-2d,0)$. Therefore,
$$\omega'\cdot C_j\in [-2d\beta+\alpha,\alpha)$$
for any $j$.
Note that $\omega'$ is a $\Qq$-Cartier $\Qq$-divisor. Let $I$ be the Cartier index of $\omega'$, then 
$$\omega'\cdot C_j\in [-2d\beta+\alpha,\alpha)\cap \frac{1}{I}\mathbb Z$$
for any $j$. Therefore, there are only finitely many possibilities of $\omega'\cdot C_j$. Possibly replacing \(\Lambda\) with an infinite subset, we may assume that \(\omega'\cdot C_j=\omega'\cdot C_{j'}\) for any \(j,j'\in\Lambda\). Since \(\omega',A_1,\dots,A_{\rho-1}\) form a basis of \(N^1_\Rr(X/U)\), \(C_j\equiv_U C_{j'}\), which is absurd as \(R_j\) and \(R_j'\) are different rays in \(\overline{NE}(X/U)\).
\end{proof}

\begin{lem}\label{lem: gfq extremal ray rational}
Let $(X,\Ff,B,\Mm)/U$ be a gfq such that $\Ff$ is algebraically integrable. Assume that $R$ is a $(K_{\Ff}+B+\Mm_X)$-negative extremal ray in $\overline{NE}(X/U)$ that is not contained in $\overline{NE}(X/U)_{\Nlc(X,\Ff,B,\Mm)}$. Then $R$ is a rational extremal ray in $\overline{NE}(X/U)$.
\end{lem}
\begin{proof}
By Lemma \ref{lem: supporting function are +A}, there exists an ample$/U$ $\Rr$-divisor $A$ such that $H_R:=K_{\Ff}+B+A+\Mm_X$ is a supporting function of $R$. We let $\delta\in (0,1)$ be a rational number such that $R$ is a $(K_{\Ff}+B+\delta A+\Mm_X)$-negative extremal ray$/U$ that is not contained in $\overline{NE}(X/U)_{\Nlc(X,\Ff,B,\Mm)}$. Let $\Lambda$ be the set of all $(K_{\Ff}+B+\delta A+\Mm_X)$-negative extremal rays$/U$. By Lemma \ref{lem: cone finiteness rays}, we may write $\Lambda=\{R,R_1,\dots,R_l\}$. Then 
$$V:=\overline{NE}(X/U)_{K_{\Ff}+B+\delta A+\Mm_X\geq 0}+\overline{NE}(X/U)_{\Nlc(X,\Ff,B,\Mm)}+\sum_{i=1}^lR_i$$
is a closed sub-cone of $\overline{NE}(X/U)$ and $R\not\in V$. Let $C$ be the dual cone of $V$ in $N^1(X/U)$. Then since $H_R\cdot R'>0$ for any $R'\in V$, $H_R$ is contained in the interior of $C$. Therefore, there exists a positive real number $\epsilon\in (0,1)$ such that $H_R-\epsilon A$ is contained in the interior of $C$. In particular, $(H_R-\epsilon A)\cdot R'>0$ for any $R'\in V$.

We write $H_R=\sum_{i=1}^cr_iD_i$, where $r_1,\dots,r_c$ are real numbers that are linearly independent over $\Qq$ and $D_i$ are $\Qq$-Cartier Weil divisors. By Theorem \ref{thm: cone theorem induction}, $R$ is spanned by a rational curve $C$. Since $H_R\cdot C=0$, $D_i\cdot C=0$ for each $i$.

Pick rational numbers $r_1',\dots,r_c'$ such that $\sum_{i=1}^c(r_i'-r_i)D_i+\epsilon A$ is ample$/U$. Let $H_R':=\sum_{i=1}^cr_i'D_i$. Then $H_R'\cdot R=0$. For any extremal ray $R'\not=R\in\overline{NE}(X/U)$, then $R'\in V$ and
$$H_R'\cdot R'=H_R\cdot R'+\sum_{i=1}^c(r_i'-r_i)D_i\cdot R'=(H_R-\epsilon A)\cdot R'+\left(\sum_{i=1}^c(r_i'-r_i)D_i+\epsilon A\right)\cdot R'>0.$$
Thus $H_R'$ is a supporting function of $R$ and the lemma holds. 
\end{proof}

\begin{proof}[Proof of Theorem \ref{thm: cone theorem gfq}]
    Theorem \ref{thm: cone theorem gfq}(1) follows from Lemma \ref{lem: gfq extremal ray rational} and Theorem \ref{thm: cone theorem induction}. Theorem \ref{thm: cone theorem gfq}(2) follows from Theorem \ref{thm: cone theorem induction}. Theorem \ref{thm: cone theorem gfq}(3) follows from Lemma \ref{lem: cone finiteness rays} and that $\Lambda=\cup_{n=1}^{+\infty}\Lambda_{\frac{1}{n}A}$ for any ample$/U$ $\Rr$-divisor $A$. We have yet to prove (4).

     For any $(K_{\Ff}+B+\Mm_X)$-negative extremal face $F$ in $\overline{NE}(X/U)$ that is relatively ample at infinity with respect to $(X,\Ff,B,\Mm)$, $F$ is also a $(K_{\Ff}+B+\Mm_X+A)$-negative extremal face for some ample$/U$ $\Rr$-divisor $A$ on $X$. Let $V:=F^\bot\subset N^1(X/U)$. Then since $F$ is spanned by a subset of $\{R_j\}_{j\in\Lambda_A}$ and each $R_j$ is rational, $V$ is defined over $\Qq$. We let
$$W_F:=\overline{NE}(X/U)_{K_X+B+\Mm_X+A\geq 0}+\overline{NE}(X/U)_{\Nlc(X,\Ff,B,\Mm)}+\sum_{j\mid j\in\Lambda_A,R_j\not\subset F}R_j.$$
Then $W_F$ is a closed cone, $\overline{NE}(X/U)=W_F+F$, and $W_F\cap F=\{0\}$. The supporting functions of $F$ are the elements in $V$ that are positive on $W_F\backslash\{0\}$, which is a non-empty open subset of $V$, and hence contains a rational element $H$. In particular, $F=H^\bot\cap \overline{NE}(X/U)$, hence $F$ is rational, and we get (4). This concludes the proof of Theorem \ref{thm: cone theorem gfq}.
\end{proof}

\section{Minimal model program for ACSS generalized foliated quadruples}\label{sec: mmp gfq}

With the establishment of the cone theorem, we are ready to study the minimal model program for algebraically integrable generalized foliated quadruples. Unfortunately for us, we cannot prove the contraction theorem and the cone theorem for the time being due to technical reasons. However, we are able to run some special types of the minimal model program for foliations.

First, we need to introduce the concept of different models of generalized foliated quadruples, similar to those of usual pairs, generalized pairs, and foliated triples.

\subsection{Models}

\begin{defn}[Models, II]\label{defn: models ii}
Let $(X,\Ff,B,\Mm)/U$ be an lc gfq and $(X',\Ff',B',\Mm)/U$ a log birational model of $(X,\Ff,B,\Mm)/U$. We say that $(X',\Ff',B',\Mm)/U$ is a \emph{log minimal model} of $(X,\Ff,B,\Mm)/U$ if 
\begin{enumerate}
    \item $(X',\Ff',B',\Mm)/U$ is a weak lc model of $(X,\Ff,B)/U$,
    \item $(X',\Ff',B',\Mm)$ is $\Qq$-factorial ACSS, and
    \item for any prime divisor $D$ on $X$ which is exceptional over $X'$, 
    $$a(D,\Ff,B,\Mm)<a(D,\Ff',B',\Mm).$$
\end{enumerate}
We say that $(X',\Ff',B',\Mm)/U$ is a \emph{good minimal model} of $(X,\Ff,B,\Mm)/U$ if $(X',\Ff',B',\Mm)/U$ is a log minimal model of $(X,\Ff,B,\Mm)/U$ and a semi-good minimal model of $(X,\Ff,B,\Mm)/U$.
\end{defn}

The following lemma is straightforward but also convenient for us to apply in some scenarios.
\begin{lem}\label{lem: acss model is gmm}
Let $(X,\Ff,B,\Mm)/U$ be an lc gfq and $(X',\Ff',B',\Mm)$ a $\Qq$-factorial ACSS model of $(X,\Ff,B,\Mm)$. Then $(X',\Ff',B',\Mm)/X$ is a good minimal model of $(X,\Ff,B,\Mm)/X$.
\end{lem}
\begin{proof}
    It immediately follows from the definitions.
\end{proof}

\begin{lem}\label{lem: ACSS mmp can run}
Let $(X,\Ff,B,\Mm)/U$ be a gfq, $\Delta\geq 0$ an $\Rr$-divisor on $X$, and $\Nn$ a nef$/U$ $\bb$-divisor on $X$. Assume that $\Ff$ is induced by a contraction $f: X\rightarrow Z$ and 
$$K_{\Ff}+B+\Mm_X\sim_{\mathbb R,Z}K_X+\Delta+\Nn_X.$$
 Then the following hold.
\begin{enumerate}
    \item Any $(K_{\Ff}+B+\Mm_X)$-negative extremal ray$/U$ $R$ is a $(K_X+\Delta+\Nn_X)$-negative extremal ray$/Z$, and $(K_{\Ff}+B+\Mm_X)\cdot R=(K_X+\Delta+\Nn_X)\cdot R.$
    \item Any step of a $(K_{\Ff}+B+\Mm_X)$-MMP$/U$ is a step of a $(K_X+\Delta+\Nn_X)$-MMP$/Z$. Moreover, assume that $(X,\Delta,\Nn)$ is lc and either $X$ is $\Qq$-factorial klt or $\Nn$ is NQC$/U$. Then we may run a step of a $(K_{\Ff}+B+\Mm_X)$-MMP$/U$.
    \item Assume that $(X,\Ff,B,\Mm;G)/Z$ is weak ACSS (resp. ACSS, super ACSS) with an associated divisor $G$, $\Delta=B+G$, and $\Nn=\Mm$. Assume that
    $$\phi: (X,\Ff,B,\Mm;G)\dashrightarrow (X',\Ff',B',\Mm;G')$$ 
    is a sequence of steps of a $(K_{\Ff}+B+\Mm_X)$-MMP$/U$ whose last step is not a Mori fiber space$/U$. Then $(X',\Ff',B',\Mm;G')/Z$ is weak ACSS (resp. ACSS, super ACSS). Moreover, if $X$ is $\Qq$-factorial and $\Mm$ is NQC$/U$ (resp. $X$ is klt, $(X,\Ff,B,\Mm;G)/Z$ is (super) ACSS), then $X'$ is $\Qq$-factorial (resp. and klt, and $(X,\Ff,B,\Mm;G)/Z$ is (super) ACSS).
    \item Any sequence of steps of a $(K_{\Ff}+B+\Mm_X)$-MMP$/U$ is a sequence of steps of a $(K_X+\Delta+\Nn_X)$-MMP$/Z$.
\end{enumerate}
\end{lem}
\begin{proof}
(1) By Theorem \ref{thm: cone theorem gfq}, any $(K_{\Ff}+B+\Mm_X)$-negative extremal ray$/U$ is tangent to $\Ff$, hence is an extremal ray$/Z$. We get (1).

(2) By (1), any $(K_{\Ff}+B+\Mm_X)$-negative extremal ray$/U$ $R$ is a $(K_X+\Delta+\Nn_X)$-negative extremal ray$/U$. If $X$ is $\Qq$-factorial klt, then by \cite[Lemma 3.4]{HL22} and the cone/contraction theorem/existence of flips for usual klt pairs, we get a step of a $(K_X+\Delta+\Nn_X)$-MMP$/U$ associated to $R$, which is also a step of a  $(K_{\Ff}+B+\Mm_X)$-MMP$/U$ associated to $R$. If $\Nn$ is NQC$/U$, then by the cone theorem (\cite[Theorem 1.1]{HL21a}, Theorem \ref{thm: cone theorem gfq}), the contraction theorem (\cite[Theorem 1.5]{Xie22}, \cite[Theorem 1.7]{CLX23}), and the existence of flips (\cite[Theorem 1.2]{LX23b}), we get a step of a $(K_{X}+\Delta+\Nn_X)$-MMP$/U$ associated to $R$, which is also a step of a  $(K_{\Ff}+B+\Mm_X)$-MMP$/U$ associated to $R$. Moreover, by (1), $R$ is a negative extremal ray$/Z$, so this step of the MMP is also a step of an MMP$/Z$.

(3) Without loss of generality, we may assume that $\phi$ is a single step of a $(K_{\Ff}+B+\Mm_X)$-MMP$/U$. It is clear that $(X',\Ff',B',\Mm;G')/Z$ is weak ACSS by Proposition \ref{prop: MMP preserves *}. 

If $(X,\Ff,B,\Mm;G)/Z$ is ACSS, then there exist $D$, $\Mm'$, such that $\Supp\{B\}\subset\Supp D$, $\Mm'-\alpha \Mm$ is nef$/X$ for some $\alpha>1$, and for any divisor $\Sigma$ on $Z$ such that $\Sigma\geq f(G)$ and $(Z,\Sigma)$ is log smooth, $(X,B+G+D+f^*(\Sigma-f(G)),\Mm')$ is qdlt. Let $\Pp:=\Mm'-\Mm$, then $(X,B+G+\delta D+f^*(\Sigma-\pi(G)),\Mm+\delta\Pp)$ is qdlt for any $0\leq \delta\leq 1$, and 
$$\Mm+\delta\Pp-(1+\delta(\alpha-1))\Mm=\delta(\Mm'-\alpha\Mm)$$ 
is nef$/X$. By (2), $\phi$ is a step of a $(K_X+B+G+\Mm_X)$-MMP$/Z$, hence a step of a  $(K_X+B+G+\delta D+f^*(\Sigma-f(G))+\Mm_X+\delta\Pp_X)$-MMP$/Z$. It follows that $(K_{X'}+B'+G'+\delta\phi_*D+f'^*(\Sigma-\pi(G)),\Mm_{X'}+\delta\Pp)$ is qdlt, where $f': X'\rightarrow Z$ is the induced contraction. Moreover, for any lc place $E$ of $(X',\Ff',B',\Mm)$, since $$-\epsilon_{\Ff}(E)\leq a(E,\Ff,B,\Mm)\leq a(E,\Ff',B',\Mm)\leq -\epsilon_{\Ff'}(E)=-\epsilon_{\Ff}(E),$$
$E$ is also an lc place of $(X,\Ff,B,\Mm)$ and $\phi$ is an isomorphism near the generic point of $\Center_XE$. Therefore $(X',\Ff',B',\Mm;G')/Z$ is ACSS. Additionally if $G$ is super, then it is clear that $G'$ is super and hence $(X',\Ff',B',\Mm;G')/Z$ is super ACSS.

Now we prove the Moreover part. By (2), $\phi$ is a step of a $(K_X+B+G+\Mm_X)$-MMP$/U$, hence a step of a $(K_X+B+G+\Mm_X+A)$-MMP$/U$ for some ample$/U$ $\Rr$-divisor $A$. The statement follows from \cite[Lemma 3.4]{HL22} and \cite[Corollary 5.20, Theorem 6.3]{HL21a}.

(4) follows from (2) and (3).
\end{proof}

\subsection{MMP with super divisors}

\begin{lem}\label{lem: super mmp with scaling}
Let $(X,\Ff,B,\Mm)/U$ be an lc gfq, $(X,\Delta,\Nn)/U$ an lc g-pair, and $f: X\rightarrow Z$ a contraction, such that $\Ff$ is induced by $f$, $\Delta$ is super$/Z$, and
$$K_{\Ff}+B+\Mm_X\sim_{\mathbb R,Z}K_X+\Delta+\Nn_X.$$
Then the following hold.
\begin{enumerate}
  \item Any $(K_X+\Delta+\Nn_X)$-negative extremal ray$/U$ $R$ is a $(K_{\Ff}+B+\Mm_X)$-negative extremal ray$/Z$ and $(K_{\Ff}+B+\Mm_X)\cdot R=(K_X+\Delta+\Nn_X)\cdot R.$
  \item A step of a $(K_X+\Delta+\Nn_X)$-MMP$/U$ is a step of a $(K_{\Ff}+B+\Mm_X)$-MMP$/Z$.
  \item Any sequence of steps of a $(K_X+\Delta+\Nn_X)$-MMP$/U$ is a sequence of steps of a $(K_{\Ff}+B+\Mm_X)$-MMP$/Z$.
  \item Let $D\geq 0$ be an $\Rr$-divisor on $X$ and $\Nn'$ a nef$/U$ $\bb$-divisor on $X$ such that $D+\Nn'_X$ is $\Rr$-Cartier. Then any sequence of steps of a $(K_X+\Delta+\Nn_X)$-MMP$/U$ with scaling of $(D,\Nn')$ is a sequence of steps of a $(K_{\Ff}+B+\Mm_X)$-MMP$/U$ with scaling of $(D,\Nn')$, and any sequence of steps of a $(K_{\Ff}+B+\Mm_X)$-MMP$/U$ with scaling of $(D,\Nn')$ is a sequence of steps of a $(K_X+\Delta+\Nn_X)$-MMP$/U$ with scaling of $(D,\Nn')$.
\end{enumerate}
\end{lem}
\begin{proof}
(1) Let $d:=\dim X$. Since $\Delta$ is super, $\Delta\geq\sum_{i=1}^{2d+1}f^*H_i$ for some ample Cartier divisors $H_i$ on $Z$. Let $L:=\Delta-\sum_{i=1}^{2d+1}f^*H_i$, then $(X,L,\Mm)$ is lc and $R$ is a $(K_X+L+\Nn_X)$-negative extremal ray. By Theorem \ref{thm: cone theorem nonnqc gpair}, there exists a rational curve $C$ on $X$ such that $C$ spans $R$ and $-2d\leq (K_X+L+\Nn_X)\cdot C<0.$ Therefore,
$$0>(K_X+\Delta+\Nn_X)\cdot C=(K_X+L+\Nn_X)\cdot C+\left(\sum_{i=1}^{2d+1}f^*H_i\cdot C\right)\geq -2d+\left(\sum_{i=1}^{2d+1}f^*H_i\cdot C\right).$$
It implies that $f(C)$ is a point and $R$ is an extremal ray$/Z$. By Proposition \ref{prop: weak cbf gfq}, we get (1).

(2) It immediately follows from (1).

Finally, (3-4) follow from (1-2) and Lemma \ref{lem: ACSS mmp can run}.
\end{proof}

\begin{lem}\label{lem: equivalence over bases}
Let $(X,\Ff,B,\Mm)/U$ be a gfq and $(X,\Delta,\Nn)/U$ an lc g-pair such that $\Ff$ is induced by a contraction $X\to Z$ and
$$K_\Ff+B+\Mm_X\sim_{\Rr,Z}K_X+\Delta+\Nn_X.$$
\begin{enumerate}
    \item If $K_X+\Delta+\Nn_X$ is either nef$/Z$ or nef$/U$, then $K_{\Ff}+B+\Mm_X$ is nef$/U$.
    \item If $\Delta$ is super$/Z$ and either $K_X+\Delta+\Nn_X$ is nef$/Z$ or $K_{\Ff}+B+\Mm_X$ is nef$/U$, then $K_X+\Delta+\Nn_X$ is nef$/U$.
\end{enumerate}
\end{lem}
\begin{proof}
The lemma follows from Lemma \ref{lem: ACSS mmp can run}(1) and Lemma \ref{lem: super mmp with scaling}(1).
\end{proof}

\subsection{MMP with scaling and existence of Mori fiber spaces}\label{subsec: eomfs}

\begin{nota}
    In the following, we need to identify ``one special MMP satisfying certain properties" and ``all MMP satisfying certain properties". For example, there are some arguments which hold for ``all MMP with scaling of an ample divisor" while some other arguments hold for ``one MMP with scaling of an ample divisor". Due to this subtlety, we need to consider ``MMP" as objects, and usually denote them by $\mathcal{P}$ or similar notations. 
\end{nota}

\begin{prop}\label{prop: run mmp with scaling gfq}
Let $(X,\Ff,B,\Mm)/U$ be an lc gfq and $f: X\rightarrow Z$ a contraction, such that 
$$K_{\Ff}+B+\Mm_X\sim_{\mathbb R,Z}K_X+\Delta+\Nn_X$$
for some lc g-pair $(X,\Delta,\Nn)/U$. 
Assume that either $X$ is $\Qq$-factorial klt or $\Nn$ is NQC$/U$. Then for any ample$/U$ $\Rr$-divisor $A$, we can run a $(K_{\Ff}+B+\Mm_X)$-MMP$/U$ with scaling of $A$. 

Moreover, there exists a $(K_{\Ff}+B+\Mm_X)$-MMP$/U$ with scaling of $A$, say $\mathcal{P}_0$, satisfying the following. Let $\mathcal{P}=\mathcal{P}_0$ if $X$ is not $\Qq$-factorial, otherwise let $\mathcal{P}$ be any $(K_{\Ff}+B+\Mm_X)$-MMP$/U$ with scaling of $A$ with the scaling numbers $\lambda_i$. 
Then the following hold.
\begin{enumerate}
  \item Suppose that there exists an ample$/U$ $\Rr$-divisor $H\geq 0$, such that either $\Delta\geq H$ or $\Nn-\overline{H}$ is nef$/U$.  Then $\mathcal{P}$ terminates at a model $(X',\Ff',B',\Mm)/U$ such that either 
  \begin{enumerate}
    \item there exists a $(K_{\Ff'}+B'+\Mm_{X'})$-Mori fiber space$/U$ which is also a $(K_{\Ff'}+B'+\Mm_{X'})$-Mori fiber space$/Z$, or 
    \item $K_{\Ff'}+B'+\Mm_{X'}\sim_{\mathbb R,Z}D$ for some semi-ample$/U$ $\Rr$-divisor $D$.
  \end{enumerate}
  \item Either $\mathcal{P}$ terminates, or $\lim_{i\rightarrow+\infty}\lambda_i=0$.  
\end{enumerate}
\end{prop}
\begin{proof}
We first construct $\mathcal{P}_0$. Possibly replacing $\Delta$, we may assume that $\Delta$ is super$/Z$. By Lemma \ref{lem: scaling number go to 0} and \cite[Lemma 2.17]{TX23}, we may run a $(K_X+\Delta+\Nn_X)$-MMP$/U$ with scaling of $A$. By Lemmas \ref{lem: super mmp with scaling} and \ref{lem: equivalence over bases}, this MMP is also a $(K_{\Ff}+B+\Mm_{X})$-MMP$/U$ with scaling of $A$, and
$$\lambda_i=\inf\{t\geq 0\mid K_{X_i}+\Delta_i+tA_i+\Nn_{X_i}\text{ is nef}/U\}$$
for each $i$, where $\Delta_i$ is the strict transform of $\Delta$ on $X_i$. This shows the existence of $\mathcal{P}_0$. 

Suppose that $X$ is not $\Qq$-factorial. Then $\Nn$ is NQC$/U$. According to \cite[Theorem A, Theorem F, Lemma 4.3]{TX23}, there is a choice of $\mathcal{P}_0$ satisfying the properties in the Moreover part.
In the following we may assume that $X$ is $\Qq$-factorial. 

\smallskip

If $\Delta\geq H$, then we let $\Delta':=\Delta$ and $\Nn':=\Nn$, and if $\Nn-\bar H$ is nef$/U$, then we let $\Nn':=\Nn-\bar H$ and $\Delta':=\Delta+H'$, where $H'$ is a general element of $|H/U|_{\mathbb R}$. Possibly replacing $\Delta,\Nn,H$ with $\Delta',H',\Nn'$ respectively, we may assume that $\Delta\geq H$. Since $\Delta$ is super$/Z$, by Lemma \ref{lem: super mmp with scaling}(4), any $(K_{\Ff}+B+\Mm_X)$-MMP$/U$ with scaling of $A$ is a $(K_X+\Delta+\Nn_X)$-MMP$/U$ with scaling of $A$. By Lemma \ref{lem: gklt+ample terminate} and Proposition \ref{prop: qfact nqc any scaling terminate}, $\mathcal{P}$ terminates with a model $(X',\Ff',B',\Mm)/U$. Let $\Delta'$ be the strict transform of $\Delta$ on $X'$. Then $\Delta'$ is super$/Z$. If we have a $(K_{X'}+\Delta'+\Nn_{X'})$-Mori fiber space $X'\rightarrow T$ over $U$, then by Lemma \ref{lem: super mmp with scaling}(1), $X'\rightarrow T$ is a $(K_{\Ff'}+B'+\Mm_{X'})$-Mori fiber space$/Z$ and we get (1.a). Thus we may assume that $(K_{X'}+\Delta'+\Mm_{X'})$ is semi-ample$/U$. (1.b) immediately follows.

\smallskip

Suppose that $\mathcal{P}$ does not terminate and $\lambda:=\lim_{i\rightarrow+\infty}\lambda_i>0$. Then $\mathcal{P}$ is an infinite sequence of steps of a $(K_{\Ff}+B+\frac{\lambda}{2}A+\Mm_X)$-MMP$/U$. Then we get a contradiction by (1). Therefore (2) holds.
\end{proof}

\begin{prop}\label{prop: run mmp get mfs}
Let $(X,\Ff,B,\Mm)/U$ be a weak ACSS gfq. 
Assume that
\begin{itemize}
    \item either $X$ is $\Qq$-factorial klt or $\Mm$ is NQC$/U$, and
    \item $K_{\Ff}+B+\Mm_X$ is not pseudo-effective$/U$.
\end{itemize}
Then there exists $\mathcal{P}_0$, a $(K_{\Ff}+B+\Mm_X)$-MMP$/U$ with scaling of an ample$/U$ $\Rr$-divisor $A$, satisfying the following. 

Let $\mathcal{P}:=\mathcal{P}_0$ if $X$ is not $\Qq$-factorial, and let $\mathcal{P}$ be any $(K_{\Ff}+B+\Mm_X)$-MMP$/U$ with scaling of $A$ if  $X$ is $\Qq$-factorial. Then $\mathcal{P}$ terminates with a Mori fiber space$/U$.
\end{prop}
\begin{proof}
According to Proposition \ref{prop: run mmp with scaling gfq}, we have that either $\mathcal{P}$ terminates, or $\lim_{i\rightarrow+\infty}\lambda_i=0$, where $\lambda_i$ are the scaling numbers. We first show that $\mathcal{P}$ terminates. Otherwise $\lim_{i\rightarrow+\infty}\lambda_i=0$. We may pick a positive real number $\epsilon$ such that $K_{\Ff}+B+\epsilon A+\Mm_X$ is not pseudo-effective$/U$. Since $\lim_{i\rightarrow+\infty}\lambda_i=0$, there exists an integer $m$ such that $\lambda_m<\epsilon$ which implies that $K_{\Ff_m}+B_m+\epsilon A_m+\Mm_{X_m}$ is pseudo-effective$/U$, which is not possible. Thus $\mathcal{P}$ terminates.

Suppose that $\mathcal{P}$ terminates at $(X_m,\Ff_m,B_m,\Mm)$ for some $m\geq 0$. Since $K_{\Ff}+B+\Mm_X$ is not pseudo-effective$/U$, $K_{\Ff_m}+B_m+\Mm_{X_m}$ is not pseudo-effective$/U$. Thus $K_{\Ff_m}+B_m+\Mm_{X_m}$ is not nef$/U$, so there exists a $(K_{\Ff_m}+B_m+\Mm_{X_m})$-Mori fiber space$/U$. The proposition follows.
\end{proof}

\begin{thm}\label{thm: existence mfs}
Let $(X,\Ff,B,\Mm)/U$ be an lc gfq. Assume that $\Ff$ is algebraically integrable and $K_{\Ff}+B+\Mm_X$ is not pseudo-effective$/U$. Then:
\begin{enumerate}
  \item $(X,\Ff,B,\Mm)/U$ has a Mori fiber space.
  \item Suppose that $(X,\Ff,B,\Mm)$ is weak ACSS, and either $X$ is $\Qq$-factorial klt or $\Mm$ is NQC$/U$. Then:
  \begin{enumerate}
    \item We may run a $(K_{\Ff}+B+\Mm_X)$-MMP$/U$ with scaling of an ample$/U$ $\Rr$-divisor, which terminates with a Mori fiber space$/U$.
    \item If $X$ is $\Qq$-factorial, then any $(K_{\Ff}+B+\Mm_X)$-MMP$/U$ with scaling of an ample$/U$ $\Rr$-divisor terminates with a Mori fiber space$/U$.
   \end{enumerate}
   \end{enumerate}
\end{thm}
\begin{proof}
(2) follows from Proposition \ref{prop: run mmp get mfs} so we only need to show (1). 

According to Theorem \ref{thm:  ACSS model}, $(X,\Ff,B,\Mm)$ has a $\Qq$-factorial ACSS model $(Y,\Ff_Y,B_Y,\Mm)$. Let $g: Y\rightarrow X$ be the induced birational morphism, then $g$ only extracts divisors $E$ such that $-\epsilon_{\Ff}(E)=a(E,\Ff,B,\Mm)$, and
$$K_{\Ff_Y}+B_Y+\Mm_Y=g^*(K_{\Ff}+B+\Mm_X)$$
is not pseudo-effective$/U$. By Proposition \ref{prop: run mmp get mfs}, we may run a $(K_{\Ff_Y}+B_Y+\Mm_Y)$-MMP$/U$ which terminates with a Mori fiber space $(Y',\Ff_{Y'},B_{Y'},\Mm)/U$ associated with $Y'\rightarrow T$. Then $(Y',\Ff_{Y'},B_{Y'},\Mm)\rightarrow T$ is a Mori fiber space$/U$ of $(X,\Ff,B,\Mm)/U$.
\end{proof}

\subsection{MMP for very exceptional divisors}\label{subsec: very exceptional}
\begin{thm}\label{thm: mmp very exceptional alg int fol}
Let $(X,\Ff,B,\Mm)/U$ be a weak ACSS gfq. Let $E_1,E_2\geq 0$ be two $\Rr$-divisors on $X$ such that $E_1\wedge E_2=0$, $E_1$ is very exceptional$/U$, and
$$K_{\Ff}+B+\Mm_X\sim_{\mathbb R,U}\text{(resp. }\equiv_U,\sim_{\mathbb Q,U}\text{) }E_1-E_2.$$ 
Assume that either $X$ is $\Qq$-factorial klt or $\Mm$ is NQC$/U$. Let $A$ be an ample$/U$ $\Rr$-divisor.
\begin{enumerate}
    \item We may run a $(K_{\Ff}+B+\Mm_X)$-MMP$/U$ with scaling  of $A$.
    \item Let $\mathcal{P}$ be the $(K_{\Ff}+B+\Mm_X)$-MMP$/U$ constructed in (1) if $X$ is not $\Qq$-factorial, and let $\mathcal{P}$ be any $(K_{\Ff}+B+\Mm_X)$-MMP$/U$ with scaling of $A$ if $X$ is $\Qq$-factorial. Then:
    \begin{enumerate}
        \item  Either $\mathcal{P}$ terminates with a Mori fiber space or contracts $E_1$ after finitely many steps.
        \item Suppose that $E_2=0$. Then:
        \begin{enumerate}
        \item $\mathcal{P}$ terminates with a weak lc model $(X',\Ff',B',\Mm)/U$ such that
        $$K_{\Ff'}+B'+\Mm_{X'}\sim_{\mathbb R,U}\text{(resp. }\equiv_U,\sim_{\mathbb Q,U}\text{) }0.$$
         \item The divisors contracted by $X\dashrightarrow X'$ are exactly $\Supp E_1$.
        \item If $(X,\Ff,B,\Mm)$ is $\Qq$-factorial ACSS, then $(X',\Ff',B',\Mm)/U$ is a good minimal model of $(X,\Ff,B,\Mm)/U$. 
        \end{enumerate}
    \end{enumerate}
\end{enumerate}
\end{thm}
\begin{proof}
(1) is a direct corollary of Proposition \ref{prop: run mmp with scaling gfq}. 

(2.a) Suppose that $\mathcal{P}$ does not terminate with a Mori fiber space. Then by Proposition \ref{prop: run mmp with scaling gfq}, there is a model $(X',\Ff',B',\Mm)$ in $\mathcal{P}$ such that $K_{\Ff'}+B'+\Mm_{X'}$ is movable$/U$. Let $E_{1}'$ and $E_{2}'$ be the strict transforms of $E_1$ and $E_2$ on $X'$ respectively. Then $E_{1}'$ is very exceptional$/U$ and 
$$K_{\Ff'}+B'+\Mm_{X'}\sim_{\mathbb R,U}\text{(resp. }\equiv_U,\sim_{\mathbb Q,U}\text{) }E_{1,m}-E_{2,m}$$
is movable over $U$. By \cite[Lemma 3.3]{Bir12}, $E_{1}'=0$. This implies (2.a).

(2.b) Now we assume that $E_2=0$. By (2.a), $\mathcal{P}$ contracts $E_1$ after finitely many steps, i.e., there is a model $(X',\Ff',B',\Mm)$ in $\mathcal{P}$ such that $E_1$ is contracted via $X\dashrightarrow X'$. In particular, 
$$K_{\Ff'}+B'+\Mm_{X'}\sim_{\mathbb R,U}\text{(resp. }\equiv_U,\sim_{\mathbb Q,U}\text{) }0.$$
Since the induced birational map $X\dashrightarrow X'$ does not extract any divisor,  $(X',\Ff',B',\Mm)/U$ is a weak lc model of $(X,\Ff,B,\Mm)/U$, which implies  (2.b.i). Since $\mathcal{P}$ is also an $E_1$-MMP$/U$, we get (2.b.ii). (2.b.iii) follows from Lemma \ref{lem: ACSS mmp can run}.
\end{proof}

\section{ACC for lc thresholds and the global ACC}\label{sec: acc gfq}

\subsection{The global ACC}\label{subsec: global acc}

\begin{lem}\label{lem: trivial trace nef imply trivial}
    Let $X$ be a normal projective variety and $\Mm$ a nef $\bb$-divisor on $X$. If $\Mm_X\equiv 0$, then $\Mm\equiv\bm{0}$.
\end{lem}
\begin{proof}
Let $f: Y\rightarrow X$ be a birational morphism such that $\Mm$ descends to $Y$. By the negativity lemma, 
$\Mm_Y=f^*\Mm_X-E\equiv -E$ for some $E\geq 0$. Since $\Mm_Y$ is nef, $\Mm_Y$ is pseudo-effective, so $-E$ is pseudo-effective. Thus $E=0$ and $\Mm_Y\equiv 0$, so $\Mm\equiv\bm{0}$.
\end{proof}

\begin{proof}[Proof of Theorem \ref{thm: global acc alg int gfq}]
We may assume that $1\in\Gamma$. According to Theorem \ref{thm:  ACSS model}, $(X,\Ff,B,\Mm)$ has a $\Qq$-factorial ACSS model. Possibly replacing $(X,\Ff,B,\Mm)$ with an ACSS model, we may assume that there exists a contraction $f: X\rightarrow Z$ such that $(X,\Ff,B,\Mm)/Z$ is $\Qq$-factorial ACSS. Let $F$ be a general fiber of $f$, $B_F:=B|_F$, $\Mm^F:=\Mm|_F$, and $\Mm^F_j:=\Mm_j|_F$ for each $j$. Since $K_F=K_X|_F=K_{\Ff}|_F$,
$$(F,B_F,\Mm^F=\sum\gamma_j\Mm_j^F)$$
is an lc g-pair of dimension $r$ such that $K_F+B_F+\Mm^F_F\equiv 0$. Moreover, $B_F\in\Ii$. By \cite[Theorem 1.6]{BZ16}, there exists a finite set $\Ii_1\subset\Ii$ depending only on $r$ and $\Ii$ such that $B_F\in\Ii_1$. Since $(X,\Ff,B,\Mm)$ is lc, $B$ is horizontal$/Z$. Thus $B\in\Ii_1$.

Possibly rewriting $\Mm$, we may assume that $\Mm_j\not\equiv\bm{0}$ and $\gamma_j>0$ for any $j$. By Lemma \ref{lem: trivial trace nef imply trivial}, $\Mm_{j,X}\not\equiv 0$ for each $j$. For any $j$, we let $\delta_j\in (0,\gamma_j)$ be a real number, then
$$K_{\Ff}+B+\Mm_X-\delta_j\Mm_{j,X}\equiv-\delta_j\Mm_{j,X}$$
is not pseudo-effective$/Z$. We may run a $(K_{\Ff}+B+\Mm_X-\delta_j\Mm_{j,X})\text{-MMP}/Z$ and it terminates with a Mori fiber space$/Z$ $\pi_j: (X_j,\Ff_j,B_j,\Mm-\delta_j\Mm_j)\rightarrow T_j$ by Theorem \ref{thm: existence mfs}.
    
Note that $(X_j,\Ff_j,B_j,\Mm)$ is lc and $K_{\Ff_j}+B_j+\Mm_{X_j}\equiv 0$, as $K_{\Ff}+B+\Mm_X\equiv 0$. Since $K_{\Ff_j}+B_j+\Mm_{X_j}-\delta_j\Mm_{j,X_j}$ is anti-ample$/T_j$, $\Mm_{j,X_j}$ is ample$/T_j$. Let $F_j$ be a general fiber of $\pi_j$, $r_j:=\dim F_j$, $B_{F_j}:=B_j|_{F_j}$, $\Mm^{(j)}:=\Mm|_{F_j}$, and $\Mm^{(j)}_i:=\Mm_i|_{F_j}$. Then $r_j\leq r$. Since $K_{F_j}=K_{X_j}|_{F_j}=K_{\Ff_j}|_{F_j}$,
$(F_j,B_{F_j},\Mm^{(j)}=\sum_i\gamma_i\Mm^{(j)}_i)$
is an lc g-pair of dimension $r_j$, and $B_{F_j}\in\Ii$. Moreover, since $\Mm_{j,X_j}$ is ample$/T_j$, $\Mm^{(j)}_{j,X_j}$ is ample. Thus $\Mm^{(j)}_j\not\equiv\bm{0}$ by Lemma \ref{lem: trivial trace nef imply trivial}. By \cite[Theorem 1.6]{BZ16}, there exists a finite set $\Ii_2\subset\Ii$ depending only on $r$ and $\Ii$ such that $\gamma_j\in\Ii_2$. Since $j$ can be any index, we may take $\Ii_0:=\Ii_1\cup\Ii_2$. We finish the proof.
\end{proof}

\subsection{ACC for lc thresholds}\label{subsec: acc}

\begin{lem}\label{lem: find nontrivial divisor on ACSS model}
Let $(X,\Ff,B,\Mm)/X$ be an lc gfq, $D$ an $\Rr$-divisor on $X$, and $\Nn$ a $\bb$-divisor on $X$ satisfying the following.
\begin{itemize}
    \item[(i)] $\Ff$ is algebraically integrable.
    \item[(ii)] $\Supp B=\Supp(B+D)$.
    \item[(iii)] Both $\Mm+\Nn$ and $\Mm-\delta\Nn$ are nef$/X$ for some $\delta\in(0,1)$.
    \item[(iv)] $(X,\Ff,B+D,\Mm+\Nn)/X$ is an lc gfq. In particular, $D+\Nn_X$ is $\Rr$-Cartier.
    \item[(v)] $(X,\Ff,B+(1+\epsilon)D,\Mm+(1+\epsilon)\Nn)$ is not lc for any positive real number $\epsilon$.
    \item[(vi)] For any prime divisor $P$ on $X$ with $a(P,\Ff,B+D,\Mm+\Nn)=-\epsilon_{\Ff}(D)$, $\mult_PD=0$.
\end{itemize}
Then there are two projective birational morphisms $h: X'\rightarrow X$ and $g: Y'\rightarrow X'$ and a real number $t\in (0,1)$ satisfying the following. 
\begin{enumerate}
    \item $h$ is an ACSS modification of $(X,\Ff,B+tD,\Mm+t\Nn)$.
    \item For any prime $h$-exceptional divisor $P$, $a(P,\Ff,B,\Mm)=-\epsilon_{\Ff}(P)$. In particular, $$a(D,\Ff,B+sD,\Mm+s\Nn)=-\epsilon_{\Ff}(D)$$ 
    for any real number $s$.
    \item $g$ extracts a unique prime divisor $E$. In particular, $-E$ is ample over $X'$.
    \item $a(E,\Ff,B+D,\Mm+\Nn)=-\epsilon_{\Ff}(E)$ and $a(E,\Ff,B,\Mm)>-\epsilon_{\Ff}(E)$. In particular,  $$a(E,\Ff,B+sD,\Mm+s\Nn)>-\epsilon_{\Ff}(E)$$ for any real number $s<1$.
    \item Let $B_{Y'},D_{Y'}$ be the strict transforms of $B,D$ on $Y'$ respectively, $\Ff_{Y'}:=(h\circ g)^{-1}\Ff$, and $F_{Y'}:=(\Supp\Exc(h\circ g))^{\Ff_{Y'}}$. Then $(Y',\Ff_{Y'},B_{Y'}+tD_{Y'}+F_{Y'};\Mm+t\Nn)$ is $\Qq$-factorial ACSS.
\end{enumerate}
\end{lem}
\begin{proof}
By assumption, there exists a prime divisor $P$ which is exceptional$/X$ such that 
$$a(P,\Ff,B+D,\Mm+\Nn)=-\epsilon_{\Ff}(P)\text{ and }a(P,\Ff,B+tD,\Mm+\alpha\Nn)<-\epsilon_{\Ff}(P)$$
 for any $\alpha>1$. According to Theorem \ref{thm: property * induction}, there exists a proper ACSS model $(Y,\Ff_Y,\tilde B_Y,\Mm+\Nn;G_Y)/Z$ of $(X,\Ff,B+D,\Mm+\Nn)$ such that $P$ is a prime divisor on $Y$. Let $B_Y,D_Y$ be the strict transforms of $B,D$ on $Y$ respectively, and $F_{Y}:=(\Supp\Exc(f))^{\Ff_{Y}},$ where $f: Y\rightarrow X$ is the induced birational morphism. Then $\tilde B_Y=B_Y+D_Y+F_Y$. By conditions (ii) and (iii) and Lemma \ref{lem: acss f-triple perturb coefficient}, there exists a real number $t\in (0,1)$ such that 
$$(Y,\Ff_Y,B_Y+tD_Y+F_Y,\Mm+t\Nn;G_Y)/Z$$
is ACSS. Let $E_1,\dots,E_n$ be the prime $f$-exceptional divisors, then
$$K_{\Ff_Y}+B_Y+tD_Y+\Mm_Y+t\Nn_Y+F_Y\sim_{\mathbb R,X}\sum_{i}\left(\epsilon_{\Ff}(E_i)+a(E_i,\Ff,B+tD,\Mm+t\Nn)\right)E_i\geq 0.$$
By Theorem \ref{thm: mmp very exceptional alg int fol}, we may run an MMP$/X$ on $K_{\Ff_Y}+B_Y+tD_Y+\Mm_Y+t\Nn_Y+F_Y$ which terminates with a good minimal model $(X',\Ff',B'+tD'+F',\Mm+t\Nn)/X$ such that $$K_{\Ff'}+B'+tD'+F'+\Mm_{X'}+t\Nn_{X'}\sim_{\mathbb R,X}0,$$ where $B',D',F'$ are the strict transforms of $B_Y,D_Y,F_Y$ on $X'$ respectively. Moreover, $(X',\Ff',B'+tD'+F',\Mm+t\Nn;G')/Z$ is $\Qq$-factorial ACSS by Lemma \ref{lem: ACSS mmp can run}, where $G'$ is the strict transform of $G_Y$ on $X'$. In particular, the induced morphism $h: X'\rightarrow X$ is an ACSS modification of $(X,\Ff,B+tD,\Mm+t\Nn)$.

By construction, the divisors contracted by the induced birational map $Y\dashrightarrow X'$ are the divisors $E_i$ satisfying 
$$\epsilon_{\Ff}(E_i)+a(E_i,\Ff,B_Y+tD_Y,\Mm+t\Nn)>0.$$ 
It is clear that $Y\dashrightarrow X'$ contracts $P$, hence $Y\dashrightarrow X'$ contains a divisorial contraction. We let $g: Y'\dashrightarrow X'$ be the last step of this MMP. Since $X'$ is $\Qq$-factorial and $K_{\Ff'}+B'+tD'+F'+\Mm_{X'}+t\Nn_{X'}\sim_{\mathbb R,X}0$, $g$ is a divisorial contraction of a prime divisor $E$. 

We show that $h,$ $g$, and $t$ satisfy our requirements. (1) and (5) immediately follow from our construction.
For any prime divisor $Q$ on $X'$ that is exceptional over $X$, we have that
$$a(Q,\Ff,B+D,\Mm+\Nn)=-\epsilon_{\Ff}(Q)=a(Q,\Ff,B+tD,\Mm+t\Nn).$$
Hence $a(Q,\Ff,B+sD,\Mm+s\Nn)=-\epsilon_{\Ff}(Q)$ for any real number $s$. This implies (2). Since $g$ is a divisorial contraction of the prime divisor $E$, $-E$ is ample$/X'$ and
\begin{align*}
a(E,\Ff,B+tD,\Mm+t\Nn)&=a(E,\Ff',B'+tD'+F',\Mm+t\Nn)\\
&>a(E,\Ff_Y,B_Y+tD_Y+F_Y,\Mm+t\Nn)=-\epsilon_{\Ff}(E),
\end{align*}
so $a(E,\Ff,B,\Mm)>-\epsilon_{\Ff}(E)$ and (3) and (4) hold. This completes the proof. 
\end{proof}

\begin{proof}[Proof of Theorem \ref{thm: acc lct alg int gfq}]
Suppose that the theorem does not hold. Then there exists a sequence of NQC lc gfqs $(X_i,\Ff_i,B_i,\Mm_i)$, $\Rr$-Cartier $\Rr$-divisors $D_i$ on $X_i$, and $\bb$-divisors $\Nn_i$ on $X_i$, such that $\rk\Ff_i=r$, $B_i,D_i\in\Ii$, $\Mm_i,\Nn_i$ are $\Ii$-linear combinations of $\bb$-nef$/X$ $\bb$-divisors, and
$$t_i:=\lct(X_i,\Ff_i,B_i,\Mm_i;D_i,\Nn_i)$$ is strictly increasing. By Lemma \ref{lem: find nontrivial divisor on ACSS model}, possibly replacing $\Ii$ with $\Ii\cup\{1\}$, we may assume that
\begin{itemize}
    \item[(i)] $(X_i,\Ff_i,B_i+t_i'D_i,\Mm_i+t_i'\Nn_i)$ is $\Qq$-factorial ACSS for some $0<t_i'<t_i$,
    \item[(ii)] There exists a divisorial contraction $f_i: Y_i\rightarrow X_i$ of a prime divisor $E_i$, such that $-E_i$ is ample$/X_i$, $$a(E_i,\Ff_i,B_i+t_iD_i,\Mm_i+t_i\Nn_i)=-\epsilon_{\Ff_i}(E_i)$$ 
    and
    $$a(E_i,\Ff_i,B_i+sD_i,\Mm_i+s\Nn_i)\not=-\epsilon_{\Ff_i}(E_i)$$ 
    for any $s\not=t_i$, and
    \item[(iii)] Let $B_{Y_i},D_{Y_i}$ be the strict transforms of $B_i,D_i$ on $Y_i$ respectively, $\Ff_{Y_i}:=f_i^{-1}\Ff_i$, and $F_i:=(\Supp\Exc(f_i))^{\Ff_i}$. Then
$$(Y_i,\Ff_{Y_i},B_{Y_i}+t_i'D_{Y_i}+F_i,\Mm_i+t_i'\Nn_i)$$
    is $\Qq$-factorial ACSS.
\end{itemize}
We let $E_i^\nu$ be the normalization of $E_i$, $\Ff_{E_i}$ the restricted foliation of $\Ff_{Y_i}$ on $E_i$, $\Mm^{E}_i:=\Mm_i|_{E_i}$, and $\Nn^E_{i}:=\Nn_i|_{E_i}$.
For any real number $t$, we let $\Mm^E(t)_i:=\Mm^E_i+t\Nn^E_i$, and
$$K_{\Ff_{E_i}}+B_{E_i}(t)+\Mm^E(t)_{i,E_i^\nu}:=(K_{\Ff_{Y_i}}+B_{Y_i}+tD_{Y_i}+F_i+\Mm_{i,Y_i}+t\Nn_{i,Y_i})|_{E_i^\nu}.$$
Let $V_i$ be the center of $E_i$ on $X_i$. Then there exists an induced birational morphism $\phi_i: E_i^\nu\rightarrow V_i$ such that
$$K_{\Ff_{E_i}}+B_{E_i}(t_i)+\Mm^E(t_i)_{i,E_i^\nu}\sim_{\mathbb R,V_i}0.$$
Since $-E_i$ is ample$/X_i$,
$$K_{\Ff_{E_i}}+B_{E_i}(t_i')+\Mm^E(t_i')_{i,E_i^\nu}$$
is anti-ample$/V_i$.

By Proposition \ref{prop: a.i preserved adjunction}, $\Ff_{E_i}$ is algebraically integrable. By Theorem \ref{thm: precise adj gfq}, 
$$(E_i^\nu,\Ff_{E_i},B_{E_i}(t_i),\Mm^E(t_i)_i)/V_i$$
is lc, and
$$(E_i^\nu,\Ff_{E_i},B_{E_i}(t),\Mm^E(t)_i)/V_i$$
is lc for any $0\leq t\leq t_i$. By Theorem \ref{thm:  ACSS model}, we may let 
$(W_i,\Ff_{W_i},B_{W_i}(t_i),\Mm^E(t_i)_i;G_i)/Z_i$
be an ACSS model of 
$$(E_i^\nu,\Ff_{E_i},B_{E_i}(t_i),\Mm^E(t_i)_i)$$ with induced birational morphism $g_i: W_i\rightarrow E_i^\nu$, and let 
$$B_{W_i}(t):=(g_i^{-1})_*B_{E_i}(t)+(\Supp\Exc(g_i))^{\Ff_{W_i}}$$
for each $i$. Since $$K_{\Ff_{E_i}}+B_{E_i}(t_i')+\Mm^E(t_i')_{i,E_i^\nu}$$
is anti-ample$/V_i$,
$K_{\Ff_{W_i}}+B_{W_i}(t_i')+\Mm^E(t_i')_{i,W_i}$ is not pseudo-effective$/V_i$. Thus we may run a 
$$(K_{\Ff_{W_i}}+B_{W_i}(t_i')+\Mm^E(t_i')_{i,W_i})\text{-MMP}/V_i$$
with scaling of an ample$/V_i$ divisor. By Theorem \ref{thm: existence mfs}, this MMP terminates with a Mori fiber space$/V_i$ $\phi_i: (\bar W_i,\Ff_{\bar W_i},B_{\bar W_i}(t_i'),\Mm^E(t_i')_{i})\rightarrow T_i$ of $(W_i,\Ff_{W_i},B_{W_i}(t_i'),\Mm^E(t_i)_i)/V_i$. By Lemma \ref{lem: ACSS mmp can run}, $\phi_i$ is also a Mori fiber space$/Z_i$. Since 
$$K_{\Ff_{W_i}}+B_{W_i}(t_i)+\Mm^E(t_i)_{i,W_i}\sim_{\mathbb R,V_i}0,$$
$(\bar W_i,\Ff_{\bar W_i},B_{\bar W_i}(t_i),\Mm^E(t_i)_{i})$ and $(W_i,\Ff_{W_i},B_{W_i}(t_i),\Mm^E(t_i)_i)$ are crepant, where $B_{\bar W_i}(t)$ is the image of $B_{W_i}(t)$ on $\bar W_i$ for any $t$. Then
$$K_{\Ff_{\bar W_i}}+B_{\bar W_i}(t_i)+\Mm^E(t_i)_{i,\bar W_i}\sim_{\mathbb R,V_i}0,$$
so
$$K_{\Ff_{\bar W_i}}+B_{\bar W_i}(t_i)+\Mm^E(t_i)_{i,\bar W_i}\sim_{\mathbb R,T_i}0.$$
Let $L_i$ be a general fiber of $\phi_i$, $B_{L_i}(t):=B_{\bar W_i}(t)|_{L_i}$ for any $t$, and $\Mm^L(t)_i:=\Mm^E(t)_i|_{L_i}$ for any $t$. Then $K_{\Ff_{\bar W_i}}|_{L_i}=K_{L_i}$, $(L_i,B_{L_i}(t_i),\Mm^L(t_i)_i)$ is lc,
$$K_L+B_{L_i}(t_i)+\Mm^L(t_i)_{i,L_i}\equiv 0,$$
and
$$K_L+B_{L_i}(t_i')+\Mm^L(t_i')_{i,L_i}$$
is anti-ample. Moreover, since $\phi_i$ is a Mori fiber space$/Z_i$, by Proposition \ref{prop: a.i preserved adjunction}, 
$$\dim L_i\leq\rk\Ff_{\bar W_i}=\rk\Ff_{E_i}\leq\rk\Ff_i=r.$$
We get a contradiction to \cite[Theorem 1.6]{BZ16} by considering the coefficients of $B_{L_i}(t_i)$ and $\Mm^L(t_i)_{i,L_i}$, which can be precisely computed by Theorem \ref{thm: precise adj gfq}. Theorem \ref{thm: acc lct alg int gfq} is proved.
\end{proof}

\subsection{Uniform rational polytopes}\label{subsec: urp}

\begin{defn}\label{defn: r affine functional divisor}
Let $X$ be a normal variety, $D_i$ $\Rr$-divisors on $X$, $\Mm_i$ $\bb$-divisors on $X$, and $d_i(t):\mathbb R\rightarrow\mathbb R$ $\Rr$-affine functions. Then we call the formal finite sum $\sum d_i(t)D_i$ an \emph{$\Rr$-affine functional divisor}, and call the formal finite set $\sum d_i(t)\Mm_i$ an \emph{$\Rr$-affine $\bb$-divisor}.
\end{defn}

\begin{defn}
Let $c$ be a non-negative real number, and $\Ii\subset[0,+\infty)$ a set of real numbers. Let $X$ be a normal variety.

For any \(\Rr\)-affine functional divisor \(\Delta(t)\) on \(X\), we write \(\Delta(t)\in\mathcal{D}_c(\Ii)\) if we may write \(\Delta(t)=\sum_id_i(t)D_i\), where \(D_i\) are distinct prime divisors, and the following condition is satisfied. For any \(i\), either \(d_i(t)=1\), or 
$$d_i(t)=\frac{m-1+\gamma+kt}{m},$$
where $m\in\mathbb{N}^{+}$, $\gamma\in\Ii_+$, $k\in\mathbb{Z}$, and $f+kt=\sum_{j}(f_j+k_jt)$, where $f_j\in\Ii\cup\{0\}$, $k_j\in\mathbb{Z}$, and $f_j+k_jc\ge0$ for any $j$.

For any \(\Rr\)-affine functional \(\bb\)-divisor \(\Mm(t)\) on \(X\) and any projective morphism \(X\rightarrow Z\), we write \(\Mm(t)\in\mathcal{D}_c(\Ii/Z)\) if we can write \(\Mm(t)=\sum_i\mu_i(t)\Mm_i\), where \(\Mm_i\) are nef\(/Z\) \(\bb\)-Cartier \(\bb\)-divisors, and the following condition is satisfied. For any \(i\), either \(\mu_i(t)=1\), or
$$\mu_i(t)=v+nt=\sum_j(v_j+n_jt),$$
where $v_j\in\Ii$, $n_j\in\mathbb Z$, and $v_j+n_jc\geq 0$ for any $j$. Moreover, if $Z=\{pt\}$, then we may omit $Z$ and write $\Mm(t)\in\mathcal{D}_c(\Ii)$.
\end{defn}

\begin{defn}
Let $d$ be a positive integer and $\Ii\subset[0,+\infty)$ a set of real numbers. We define $\mathcal{B}_{d}(\Ii),\mathcal{B}'_{d}(\Ii)\subset [0,+\infty)$ as follows, $c\in\mathcal{B}_{d}(\Ii)$ (resp. $\mathcal{B}'_{d}(\Ii)$) if and only if there exist a normal projective variety $X$ (resp. a $\Qq$-factorial normal projective variety $X$), an $\Rr$-affine functional divisor $\Delta(t)$ on $X$, and an $\Rr$-affine functional $\bb$-divisor $\Mm(t)$ satisfying the following.
\begin{enumerate}
	\item $\dim X\le d$, 
	\item $\Delta(t)\in\mathcal{D}_{c}(\Ii)$, $\Mm(t)\in\mathcal{D}_c(\Ii)$,
	\item $(X,\Delta(c),\Mm(c))$ is lc,
	\item $K_X+\Delta(c)+\Mm(c)_X\equiv0$, and
	\item $K_X+\Delta(c')+\Mm(c')_X\not\equiv 0$ for any $c'\neq c$.
\end{enumerate}
\end{defn}

\begin{defn}
Let $r$ be a positive integer and $\Ii\subset[0,+\infty)$ a set of real numbers. We define $\mathcal{C}_{r}(\Ii),\mathcal{C}'_{r}(\Ii)\subset [0,+\infty)$ as follows, $c\in\mathcal{C}_{r}(\Ii)$ (resp. $c\in\mathcal{C}_{r}'(\Ii)$) if and only if there exist a normal projective variety $X$ (resp. a $\Qq$-factorial normal projective variety $X$), an algebraically integrable foliation $\Ff$ on $X$, an $\Rr$-affine functional divisor $\Delta(t)$ on $X$, and an $\Rr$-affine functional $\bb$-divisor $\Mm(t)$ on $X$ satisfying the following.
\begin{enumerate}
	\item $\rk\Ff\le r$, 
	\item $\Delta(t)\in\mathcal{D}_{c}(\Ii)$, $\Mm(t)\in\mathcal{D}_c(\Ii)$,
	\item $(X,\Ff,\Delta(c),\Mm(c))$ is lc,
	\item $K_\Ff+\Delta(c)+\Mm(c)_X\equiv0$, and
	\item $K_\Ff+\Delta(c')+\Mm(c')_X\not\equiv 0$ for any $c'\not=c$.
\end{enumerate}
\end{defn}

\begin{prop}\label{prop: nak special set equal for foliation}
Let $r$ be a positive integer and $\Ii\subset[0,+\infty)$ a set of real numbers. Then $\mathcal{B}_r(\Ii)=\mathcal{C}_r(\Ii)=\mathcal{B}'_r(\Ii)=\mathcal{C}'_r(\Ii)$.
\end{prop}

\begin{proof}
By considering the foliation $\Ff=T_X$, we have $\mathcal{B}_r(\Ii)\subset\mathcal{C}_r(\Ii)$. By the existence of dlt modifications, $\mathcal{B}_r(\Ii)=\mathcal{B}_r'(\Ii)$. By Theorem \ref{thm:  ACSS model}, $\mathcal{C}_r(\Ii)=\mathcal{C}_r'(\Ii)$. We only need to show that $\mathcal{C}_r(\Ii)\subset\mathcal{B}_r(\Ii)$.

Pick $c\in\mathcal{C}_r(\Ii)$. Then there exists a $\Qq$-factorial normal projective variety $X$, an algebraically integrable foliation $\Ff$ on $X$, an \(\Rr\)-affine functional divisor \(\Delta(t)\) on \(X\), and an \(\Rr\)-affine functional \(\bb\)-divisor \(\Mm(t)\) on \(X\), such that
\begin{enumerate}
	\item $\rk\Ff\le r$, 
	\item $\Delta(t)\in\mathcal{D}_{c}(\Ii), \Mm(t)\in\mathcal{D}_c(\Ii)$,
	\item $(X,\Ff,\Delta(c),\Mm(c))$ is lc,
	\item $K_{\Ff}+\Delta(c)+\Mm(c)_X\equiv0$, and
	\item $K_{\Ff}+\Delta(t)+\Mm(t)_X\not\equiv 0$ for any $t\not=c$.
\end{enumerate}
By Theorem \ref{thm:  ACSS model}, we may let $f: X'\rightarrow X$ be an ACSS modification of $(X,\Ff,\Delta(c),\Mm(c))$, $\Ff':=f^{-1}\Ff$, $E:=(\Supp\Exc(f))^{\Ff'}$, and $\Delta'(t):=f^{-1}_*\Delta(t)+E$ for any real number $t$. Then $\rk\Ff'\leq r$, $\Delta'(t)\in\mathcal{D}_c(\Ii)$, $(X',\Ff',\Delta'(c),\Mm(c))$ is lc, and $K_{\Ff'}+\Delta'(c)+\Mm(c)_{X'}\equiv 0$. Moreover, for any $t\not=c$, since $$0\not\equiv K_{\Ff}+\Delta(t)+\Mm(t)_X=f_*(K_{\Ff'}+\Delta'(t)+\Mm(t)_{X'}),$$ $K_{\Ff'}+\Delta'(t)+\Mm(t)_{X'}\not\equiv 0$. Therefore, we may replace $(X,\Ff,\Delta(t),\Mm(t))$ with $(X',\Ff',\Delta'(t),\Mm(t))$, and assume that $(X,\Ff,\Delta(c),\Mm(c))$ is $\Qq$-factorial ACSS. Thus there exists a contraction $\pi: X\rightarrow Z$ and a reduced divisor $G$ such that $(X,\Ff,\Delta(c),\Mm(c);G)/Z$ is ACSS.

Suppose that for any $0<\delta\ll 1$, $(X,\Ff,\Delta(c+\delta),\Mm(c+\delta);G)/Z$ or $(X,\Ff,\Delta(c-\delta),\Mm(c+\delta);G)/Z$ is not ACSS. By Lemmas \ref{lem: acss smaller coefficient} and \ref{lem: acss f-triple perturb coefficient},
\begin{itemize}
    \item either there exists a component $D$ of $\Delta(c)$, such that $\mult_D\Delta(c)=1$ and $\mult_D\Delta(t)\not=1$ for any $t\not=c$, or
    \item $\Mm(t)=\sum \mu_i(t)\Mm_i$, where each $\Mm_i$ is $\bb$-nef, and $\mu_i(t)=v_{i}+n_{i}t=\sum_i(v_{i,j}+n_{i,j}t)$ for any $v_{i,j}\in\Ii$, $n_{i,j}\in\mathbb Z$, $v_{i}+n_{i}c\geq 0$ for any $i$, and $v_i+n_ic=0$ for some $i$.
\end{itemize}
By \cite[Lemma 3.7]{Nak16}, $c\in\mathcal{B}_1(\Ii)\subset\mathcal{B}_r(\Ii)$. Therefore, we may assume that $(X,\Ff,\Delta(c+\delta),\Mm(c+\delta);G)/Z$ and $(X,\Ff,\Delta(c-\delta),\Mm(c-\delta);G)/Z$ are ACSS for any $0<\delta\ll 1$. 

Fix $0<\delta\ll 1$. Since $K_{\Ff}+\Delta(t)+\Mm(t)_X\not\equiv 0$ for any $t\not=c$ and $K_{\Ff}+\Delta(c)+\Mm(c)_X\equiv 0$, either $K_{\Ff}+\Delta(c+\delta)+\Mm(c+\delta)_X$ or $K_{\Ff}+\Delta(c-\delta)+\Mm(c-\delta)_X$ is not pseudo-effective. By Theorem \ref{thm: existence mfs}, we may run a $(K_{\Ff}+\Delta(c+\delta)+\Mm(c+\delta)_X)$-MMP (resp. $(K_{\Ff}+\Delta(c-\delta)+\Mm(c-\delta)_X)$-MMP) with scaling of an ample divisor if $(K_{\Ff}+\Delta(c+\delta)+\Mm(c+\delta)_X)$ (resp. $(K_{\Ff}+\Delta(c-\delta)+\Mm(c-\delta)_X)$) is not pseudo-effective, which terminates with a Mori fiber space $\phi: (X'',\Ff'',\Delta''(c+\delta),\Mm(c+\delta))\rightarrow T$ (resp. $\phi: (X'',\Ff'',\Delta''(c-\delta),\Mm(c-\delta))\rightarrow T$) of $(X,\Ff,\Delta(c+\delta),\Mm(c+\delta))$ (resp. $(X,\Ff,\Delta(c-\delta),\Mm(c-\delta))$), where $\Delta''(t)$ is the image of $\Delta(t)$ on $X''$ for any $t$. By Lemma \ref{lem: ACSS mmp can run}(4), this MMP is also an MMP$/Z$ and $\phi$ is a contraction$/Z$.

Since $K_{\Ff}+\Delta(c)+\Mm(c)_X\equiv 0$, $(X'',\Ff'',\Delta''(c),\Mm(c))$ and $(X,\Ff,\Delta(c),\Mm(c))$ are crepant, so $K_{\Ff''}+\Delta''(c)+\Mm(c)_{X''}\equiv 0$ and $(X'',\Ff'',\Delta''(c),\Mm(c))$ is lc.

Let $F$ be a general fiber of $\phi$. By Theorem \ref{thm: cone theorem gfq}, $F$ is tangent to $\Ff''$, so $K_{\Ff''}|_F=K_{X''}|_F=K_F$. Let $\Delta_{F}(t):=\Delta''(t)|_F$ and $\Mm^F(t):=\Mm(t)|_F$. Then
\begin{itemize}
   \item $\dim F\leq \dim X-\dim Z=\rk\Ff\leq r$,
    \item $\Delta_F(t)\in\mathcal{D}_c(\Ii)$ and $\Mm^F(t)\in\mathcal{D}_c(\Ii)$.
    \item $(F,\Delta_F(c),\Mm^F(c))$ is lc,
    \item $K_F+\Delta_F(c)+\Mm^F(c)_F\equiv 0$, and
    \item $K_F+\Delta_F(c+\delta)+\Mm^F(c+\delta)_F$ or $K_F+\Delta_F(c-\delta)+\Mm^F(c-\delta)_F$ is anti-ample.
\end{itemize}
Thus $c\in\mathcal{B}_r(\Ii)$.
\end{proof}

\begin{thm}\label{thm: uniform rational polytope foliation one variable}
Let $d,c,m,n$ be positive integers, $r_1,\dots,r_c$ real numbers such that $1,r_1,\dots,r_c$ are linearly independent over $\mathbb Q$, $\bm{r}:=(r_1,\dots,r_c)$, and $s_1,\dots,s_m, \mu_1,\dots,\mu_n: \mathbb R^{c+1}\rightarrow\mathbb R$ $\mathbb Q$-linear functions. Then there exists a positive real number $\delta$ depending only on $d, \bm{r}$ and $s_1,\dots,s_m, \mu_1,\dots,\mu_n$ satisfying the following. Assume that
\begin{enumerate}
    \item 
    $$\left(X,\Ff,B=\sum_{i=1}^ms_i(1,r_1,\dots,r_{c-1},t)B_i, \Mm=\sum_{i=1}^n\mu_i(1,r_1,\dots,r_{c-1},t)\Mm_i\right)/X$$ is an lc gfq such that $\Ff$ is algebraically integrable and $\rk\Ff\leq d$,
    \item $B_i\geq 0$ are distinct Weil divisors (possibly $0$) and $s_i(1,\bm{r})\geq 0$ for each $i$,
    \item $\Mm_i$ are nef$/X$ $\bb$-Cartier $\bb$-divisors and $\mu_i(1,\bm{r})\geq 0$ for each $i$, and
    \item $B(t):=\sum_{i=1}^ms_i(1,r_1,\dots,r_{c-1},t)B_i$ and $\Mm(t):=\sum_{i=1}^n\mu_i(1,r_1,\dots,r_{c-1},t)\Mm_i$ for any $t\in\mathbb R$.
\end{enumerate}
Then $(X,\Ff,B(t),\Mm(t))$ is lc for any $t\in (r_c-\delta,r_c+\delta)$.
\end{thm}
\begin{proof}
We let $s_i(t):=s_i(1,r_1,\dots,r_{c-1},t)$ and $\mu_i(t):=\mu_i(1,r_1,\dots,r_{c-1},t)$ for any $t\in\mathbb R$. If $s_i(r_c)=0$, then $s_i(t)=0$ for any $i$, so we may assume that $s_i(r_c)\not=0$ for any $i$. Let $(X',\Ff',B'(r_c),\Mm(r_c))$ be an ACSS model of $(X,\Ff,B(r_c),\Mm(r_c))$, $f: X'\rightarrow X$ the induced birational morphism, $E:=(\Supp(\Exc(f)))^{\Ff'}$, and $B'(t):=f^{-1}_*B(t)+E$ for any $t$. Then
$$K_{\Ff'}+B'(r_c)+\Mm(r_c)_{X'}=f^*(K_{\Ff}+B(r_c)+\Mm(r_c)_{X}).$$
Since $1,r_1,\dots,r_c$ are linearly independent over $\mathbb Q$, for $B'(t):=f^{-1}_*B(t)+E$, we have
$$K_{\Ff'}+B'(t)+\Mm(t)_{X'}=f^*(K_{\Ff}+B(t)+\Mm(t)_{X})$$
for any $t\in\mathbb R$. Thus possibly replacing $(X,\Ff,B(t),\Mm(t))$ with $(X',\Ff',B'(t),\Mm(t))$, we may assume that $(X,\Ff,B(r_c),\Mm(r_c))$ is $\Qq$-factorial ACSS.

Let $$t_1:=\inf\{t\geq r_c\mid (X,\Ff,B(r_c),\Mm(r_c)\text{ is lc}\}$$
and 
$$t_2:=\sup\{t\leq r_c\mid (X,\Ff,B(r_c),\Mm(r_c)\text{ is lc}\}.$$
If $t_1\leq t_2$, then we let $t_0:=t_1$. Otherwise, we let $t_0:=t_2$. We only need to show that there exists a positive real number $\epsilon$ depending only on $d,\bm{r}$, and $s_1,\dots,s_m,\mu_1,\dots,\mu_n$, such that $|t_0-r_c|\geq\epsilon$. 

Since $1,r_1,\dots,r_c$ are linearly independent over $\mathbb Q$, there exists a positive real number $\delta_1$ depending only on $\bm{r}$ and $s_1,\dots,s_m,\mu_1,\dots,\mu_n$, such that $s_i(t)>0$ and $\mu_i(t)>0$ for any $t\in (r_c-\delta_1,r_c+\delta_1)$. We may assume that $|t_0-r_c|<\delta_1$. In particular, for any $0<\delta\ll 1$, $B(t_0+\delta(t_0-r_c))\geq 0$, and $\mu_i(t_0+\delta(t_0-r_c))>0$. Thus $(X,\Ff,B(t_0))$ has an lc center $V_0$ such that $\dim V_0\leq \dim X-2$, and $V_0$ is not an lc center of $(X,\Ff,B(r_c),\Mm(r_c))$. 

By Lemma \ref{lem: find nontrivial divisor on ACSS model}, possibly replacing $(X,\Ff,B(t),\Mm(t))$, we may assume that there exists a divisorial contraction $g: Y\rightarrow X$ of a prime divisor \(\tilde E\) and a real number \(s\) satisfying the following, let \(B_Y(t)\) be the strict transform of \(B(t)\) on \(Y\) for any \(t\) and \(\Ff_Y:=g^{-1}\Ff\).
\begin{itemize}
\item[(i)] $s\in (r_c,t_0)$ if $r_c>t_0$, and $s\in (t_0,r_c)$ if $t_0<r_c$.
\item[(ii)] $(X,\Ff,B(s),\Mm(s))$ is $\Qq$-factorial ACSS, $(X,\Ff,B(r_c),\Mm(r_c))$ is lc, and $(X,\Ff,B(t_0),\Mm(t_0))$ is lc.
\item[(iii)] $-\tilde E$ is ample over $X$.
\item[(iv)] $(Y,\Ff_Y,B_Y(s)+\epsilon_{\Ff}(\tilde E),\Mm(s))$ is $\Qq$-factorial ACSS.
\item[(v)] $a(\tilde E,\Ff,B(t_0),\Mm(t_0))=-\epsilon_{\Ff}(\tilde E)$ and $a(\tilde E,\Ff,B(r_c),\Mm(r_c))>-\epsilon_{\Ff}(\tilde E)$. In particular, $(Y,\Ff_Y,B_Y(t_0)+\epsilon_{\Ff}(\tilde E),\Mm(t_0))$ is lc and $a(\tilde E,\Ff,B(s),\Mm(s))>-\epsilon_{\Ff}(\tilde E)$.
\end{itemize}
We let $E$ be the normalization of $\tilde E$, $\Ff_E$ the restricted foliation of $\Ff_Y$ on $E$, $V:=\Center_XE$,
$\Mm^E(t):=\Mm(t)|_E$, and
$$K_{\Ff_E}+B_E(t)+\Mm^E(t)_E:=(K_{\Ff_Y}+B_Y(t)+\epsilon_{\Ff}(\tilde E)+\Mm(t)_Y)|_E$$
for any real number $t$. By Theorem \ref{thm: precise adj gfq}, $B_E(t)$ is an \(\Rr\)-affine functional divisor, \(\Mm^E(t)\) is an \(\Rr\)-affine functional \(\bb\)-divisor, and
$$(E,\Ff_E,B_E(t_0),\Mm^E(t_0)),(E,\Ff_E,B_E(s),\Mm^E(s))$$
are lc gfqs. By Proposition \ref{prop: a.i preserved adjunction}, $\Ff_E$ is algebraically integrable and $\rk\Ff\leq d$. 

Let $E\rightarrow V$ be the induced projective surjective morphism. Since $-\tilde E$ is ample$/X$,
$$K_{\Ff_E}+B_E(t_0)+\Mm^E(t_0)_E\sim_{\mathbb R,V}0$$
and
$$K_{\Ff_E}+B_E(s)+\Mm^E(s)_E$$
is anti-ample$/V$. Thus $$K_{\Ff_E}+B_E(t)+\Mm^E(t)_E$$ is anti-ample$/V$ for any $t\in (t_0,s)$ if $t_0<r_c$, and for any $t\in (s,t_0)$ if $t_0>r_c$.

By Theorem \ref{thm:  ACSS model}, we may let $$(W,\Ff_W,B_W(t_0),\Mm^E(t_0);G)/Z$$ be an ACSS model of $(E,\Ff_E,B_E(t_0),\Mm^E(t_0))$ with induced birational morphism $g: W\rightarrow E$. Let $F_W:=(\Supp\Exc(g))^{\Ff_W}$ and let $B_W(t):=g^{-1}B_E(t_0)+F_W$ for any $t\in\mathbb R$. Since $s\in (r_c-\delta_1,r_c+\delta_1)$, $s_i(s)>0$ and $\mu_i(s)>0$. By Theorem \ref{thm: precise adj gfq} and Lemma \ref{lem: acss f-triple perturb coefficient}, there exists a real number $u$, such that $u\in (t_0,s)$ if $t_0<r_c$, $u\in (s,t_0)$ if $t_0>r_c$, and 
$(W,\Ff_W,B_W(u),\Mm^E(u);G)/Z$
is ACSS. Since
$$K_{\Ff_E}+B_E(u)+\Mm^E(u)_E$$ is anti-ample$/V$,
$K_{\Ff_W}+B_W(u)+\Mm^E(u)_W$ is not pseudo-effective$/V$. Thus we may run a 
$$(K_{\Ff_W}+B_W(u)+\Mm^E(u)_W)\text{-MMP}/V$$
with scaling of an ample$/V$ divisor. By Theorem \ref{thm: existence mfs}, this MMP terminates with a Mori fiber space$/V$ $\phi: (\bar W,\Ff_{\bar W},B_{\bar W}(u),\Mm^E(u))\rightarrow T$ of $(W,\Ff_{W},B_{W}(u),\Mm^E(u))$. By Lemma \ref{lem: ACSS mmp can run}, $\phi$ is also a Mori fiber space$/Z$. 

Let $B_{\bar W}(t)$ be the image of $B_W(t)$ on $\bar W$ for any $t$. Since 
$$K_{\Ff_E}+B_E(t_0)+\Mm^E(t_0)_E\sim_{\mathbb R,V}0,$$
we have
$$K_{\Ff_W}+B_W(t_0)+\Mm^E(t_0)_W\sim_{\mathbb R,V}0,$$
so $(\bar W,\Ff_{\bar W},B_{\bar W}(u),\Mm^E(u))$ and $(W,\Ff_{W},B_{W}(u),\Mm^E(u))$ are crepant, and
$$K_{\Ff_{\bar W}}+B_{\bar W}(t_0)+\Mm^E(t_0)_{\bar W}\sim_{\mathbb R,V}0.$$
Let $L$ be a general fiber of $\phi$, $B_L(t):=B_{\bar W}(t)|_L$ for any $t$, and $\Mm^L(t):=\Mm^E(t)|_L$ for any $t$. Since $\phi$ is a Mori fiber space$/Z$, the general fibers of $\phi$ are tangent to $\Ff_{\bar W}$. Thus $K_{\Ff_{\bar W}}|_L=K_L$, $(L,B_L(t_0),\Mm^L(t_0))$ is lc,
$$K_L+B_L(t_0)+\Mm^L(t_0)_L\equiv 0,$$
and
$$K_L+B_L(u)+\Mm^L(u)_L$$
is anti-ample. Moreover, since $\phi$ is a Mori fiber space$/Z$, 
$$\dim L\leq\rk\Ff_{\bar W}=\rk\Ff_{E}\leq d.$$
By \cite[Theorem 3.6]{Che20} and considering the coefficients of $B_L(t_0)$ and $\Mm^L(u)$, which can be precisely computed by Theorem \ref{thm: precise adj gfq}, there exists a positive real number $\epsilon$ depending only on $d, \bm{r}, s_1,\dots,s_m, \mu_1,\dots,\mu_n$, such that $|t_0-r_c|\geq\epsilon$. This concludes the proof of the theorem.
\end{proof}

\begin{thm}\label{thm: uniform rational polytope}
Let $d,c,m,n$ be positive integers, $r_1,\dots,r_c$ real numbers such that $1,r_1,\dots,r_c$ are linearly independent over $\mathbb Q$, $\bm{r}:=(r_1,\dots,r_c)$, and $s_1,\dots,s_m,\mu_1,\dots,\mu_n: \mathbb R^{c+1}\rightarrow\mathbb R$ $\mathbb Q$-linear functions. Then there exists an open subset $U\ni\bm{r}$ depending only on $d,\bm{r}$ and $s_1,\dots,s_m,\mu_1,\dots\mu_n$ satisfying the following. Assume that
\begin{enumerate}
    \item $$\left(X,\Ff,B(\bm{r}):=\sum_{i=1}^ms_i(1,\bm{r})B_i, \Mm(\bm{r}):=\sum_{i=1}^n\mu_i(1,\bm{r})\Mm_i\right)/X$$ is an lc gfq such that $\Ff$ is algebraically integrable and $\rk\Ff\leq d$,
    \item $B_i\geq 0$ are distinct Weil divisors (possibly $0$) and $s_i(1,\bm{r})\geq 0$,
    \item $\Mm_i$ are nef$/X$ $\bb$-Cartier $\bb$-divisors and $\mu_i(1,\bm{r})\geq 0$, and
    \item $B(\bm{v}):=\sum_{i=1}^ms_i(1, \bm{v})B_i$ and \(\Mm(\bm{v}):=\sum_{i=1}^n\mu_i(1,\bm{v})\Mm_i\) for any \(t\in\mathbb R\).
\end{enumerate}
Then $(X,\Ff,B(\bm{v}),\Mm(\bm{v}))$ is lc for any $\bm{v}\in U$.
\end{thm}
\begin{proof}
We apply induction on $c$. When $c=1$, Theorem \ref{thm: uniform rational polytope} directly follows from Theorem \ref{thm: uniform rational polytope foliation one variable}. When $c\geq 2$, by Theorem \ref{thm: uniform rational polytope foliation one variable}, there exists a positive integer $\delta$ depending only on $r_1,\dots,r_c,s_1,\dots,s_m$, such that for any $t\in (r_c-\delta,r_c+\delta)$, $$\left(X,\Ff,\sum_{i=1}^ms_i(1,r_1,\dots,r_{c-1},t)B_i,\sum_{i=1}^n\mu_i(1,r_1,\dots,r_{c-1},t)\Mm_i\right)$$ is lc. We pick rational numbers $r_{c,1}\in (r_c-\delta,r_c)$ and $r_{c,2}\in (r_c,r_c+\delta)$ depending only on $r_1,\dots,r_c,s_1,\dots,s_m$. By induction on $c$, there exists an open subset $U_0$ of $\mathbb R^{c-1}$ containing $(r_1,\dots,r_{c-1})$, such that for any $\bm{v}\in U_0$, $$\left(X,\Ff,\sum_{i=1}^ms_i(1,\bm{v},r_{c,1})B_i,\sum_{i=1}^n\mu_i(1,\bm{v},r_{c,1})\Mm_i\right)$$ and $$\left(X,\Ff,\sum_{i=1}^ms_i(1,\bm{v},r_{c,2})B_i,\sum_{i=1}^n\mu_i(1,\bm{v},r_{c,2})\Mm_i\right)$$ are lc. We may pick $U:=U_0\times (r_{c,1},r_{c,2})$.
\end{proof}

\part{Canonical bundle formula and MMP for generalized pairs}\label{part:cbf}

\section{Canonical bundle formula for lc-trivial fibrations}\label{sec: cbf}

\subsection{Stability of generalized foliated quadruples}\label{subsec: stability gfq}

\begin{prop}[{cf. \cite[Proposition 3.7]{ACSS21}}]\label{prop: bp semistable foliation lc}
    Let $(X,\Ff,B,\Mm)/U$ be a sub-gfq satisfying Property $(*)$ with an associated contraction $f: X\rightarrow Z$. Let $G$ be a divisor associated to $(X,\Ff,B,\Mm)/Z$. Assume that $f$ is equidimensional and $B$ is horizontal$/Z$. 

    Then $(X,B+G,\Mm)$ is BP semistable$/Z$ if and only if $(X,\Ff,B,\Mm)$ is sub-lc.
\end{prop}
\begin{proof}
Since $(X,\Ff,B,\Mm)/U$ satisfies Property $(*)$, $f: (X,B+G,\Mm)\rightarrow Z$ satisfies Property $(*)$. By Lemma \ref{lem: basic property (*) gpair}(1), $(X,B+G,\Mm)$ is sub-lc. Let $f': (X',\Sigma_{X'},\Mm)\rightarrow (Z',\Sigma_{Z'})$ be any equidimensional model of $(X,B+G,\Mm)$ associated with $h: X'\rightarrow X$ and $h_Z: Z'\rightarrow Z$. By Proposition \ref{prop: weak ss imply *}, there exists an $\Rr$-divisor $\bar B$ on $X'$ satisfying the following.
\begin{itemize}
    \item $\Supp \bar B\subset\Supp\Sigma_{X'}$.
    \item $K_{X'}+\bar B+\Mm_{X'}=h^*(K_X+B+G+\Mm_X)+F$ for some $F\geq 0$ that is vertical$/Z'$.
\item $(X',\bar B,\Mm)$ and $(X,B+G,\Mm)$ are crepant over the generic point of $Z$. In particular, $(X',\bar B,\Mm)$ and $(X,B,\Mm)$ are crepant over the generic point of $Z$.
\item $f': (X',\bar B,\Mm)\rightarrow Z'$ satisfies Property $(*)$. By Lemma \ref{lem: basic property (*) gpair}(1), $(X',\bar B,\Mm)$ is sub-lc.
\end{itemize}
Let $\Ff' := h^{-1}\Ff$, $G'$ the vertical$/Z'$ part of $\bar B$, and $B'$ the horizontal$/Z'$ part of $\bar B$. Then $(X',\Ff',B',\Mm;G')/Z'$ satisfies Property $(*)$. Since $f'$ is vertical$/Z'$, and $\Ff'$ is induced by $f'$, $F$ is $\Ff'$-invariant.

We let $B_Z$ and $\Nn$ be the discriminant part and moduli part of $f: (X,B+G,\Mm)\rightarrow Z$ respectively, and let $\bar B_{Z'}$ and $\Nn'$ be the discriminant part and moduli part of $f': (X',\bar B,\Mm)\rightarrow Z'$ respectively. Since $(X,\Ff,B,\Mm)/Z$ satisfies Property $(*)$, $Z$ is smooth, so $K_Z+B_Z$ is $\Rr$-Cartier, and we may define
$$K_{Z'}+B_{Z'} := h_Z^*(K_Z+B_Z).$$
Since $f': (X',\bar B,\Mm)\rightarrow Z'$ satisfies Property $(*)$, $\bar B_{Z'}=f'(G')$. Since $f: (X,B+G,\Mm)\rightarrow Z$ satisfies Property $(*)$, $B_Z=f(G)$.

By Proposition \ref{prop: weak cbf gfq}, $K_{\Ff}+B+\Mm_X\sim\Nn_X$ and $K_{\Ff'}+B'+\Mm_{X'}\sim\Nn'_{X'}$. In particular, $\Nn_X$ is $\Rr$-Cartier. Let $A:=\Nn'_{X'}-h^*\Nn_{X'}$. Then 
\begin{align*}
    A&=K_{X'}+\bar B+\Mm_{X'}-f'^*(K_{Z'}+\bar B_{Z'})-h^*(K_X+B+G+\Mm_X-f^*(K_Z+B_Z))\\
    &=F-f'^*(K_{Z'}+\bar B_{Z'})+f'^*(K_{Z'}+B_{Z'})=F-f'^*(\bar B_{Z'}-B_{Z'}).
\end{align*}
In particular, $A$ is vertical$/Z'$.
\begin{claim}\label{claim: a' for semistable}
    $(X,\Ff,B,\Mm)$ is sub-lc if and only if $A\geq 0$.
\end{claim}
\begin{proof}
    Let 
    $$A':=K_{\Ff'}+h^{-1}_*B+\Mm_{X'}+(\Supp\Exc(h))^{\Ff'}-h^*(K_{\Ff}+B+\Mm_X).$$
    Since $h$ is a foliated log resolution of $(X,\Ff,B,\Mm)$, $(X,\Ff,B,\Mm)$ is sub-lc if and only if $A'\geq 0$. Since
    $$A\sim K_{\Ff'}+B'+\Mm_{X'}-h^*(K_{\Ff}+B+\Mm_X),$$
    for suitable choices of $K_{\Ff}$ and $K_{\Ff'}$, we have
    $$A'-A=h^{-1}_*B-B'+(\Supp\Exc(h))^{\Ff'}.$$
    For any horizontal$/Z$ prime divisor $P$ on $X'$, if $P$ is not exceptional$/X$, then
    $$\mult_PA'=\mult_P(h^{-1}_*B-B')=0$$
    as $G$ and $F$ are vertical$/Z$. If $P$ is exceptional$/X$, then
    $$\mult_PA'=1+\mult_P(h^{-1}_*B-B')\geq 1-\mult_P B'.$$
    Since $f': (X',\bar B,\Mm)\rightarrow Z'$ satisfies Property $(*)$, $\mult_P \bar B\leq 1$, so $\mult_PA'\geq 0$.

    For any vertical$/Z'$ prime divisor $P$ on $X'$, since $B$ is horizontal$/Z$ and $B'$ is horizontal$/Z'$, $\mult_PA'=\mult_PA$. The claim follows.
\end{proof}
 Let $B'_{Z'}$ be the discriminant part of $f': (X',\bar B-F,\Mm)\rightarrow Z'$. Then for any prime divisor $D$ on $Z'$,
     $$\mult_D\bar B_{Z'}=1-\sup\{t\mid (X',\bar B+tf'^*D,\Mm)\text{ is sub-lc over the generic point of }D\}$$
     and
    $$\mult_D B'_{Z'}=1-\sup\{t\mid (X',\bar B-F+tf'^*D,\Mm)\text{ is sub-lc over the generic point of }D\}.$$
    Since $\Supp\bar B\subset\Supp\Sigma_{X'}$, by the definition of an equidimensional model, $\Supp(\bar B-F)\subset\Supp\Sigma_{X'}$. Therefore, if $D\subset\Supp\Sigma_{Z'}$, then
    $$\mult_D B'_{Z'}=\mult_D\bar B_{Z'}=\mult_{f'^*D}F=0.$$
    Otherwise,
    $$\mult_D(\bar B_{Z'}-B'_{Z'})=\sup\{t\geq 0\mid F-tf'^*D\geq 0\}.$$
    Therefore, 
    \begin{itemize}
        \item  $F-f'^*(\bar B_{Z'}-B'_{Z'})\geq 0$, and
        \item  $F-f'^*(\bar B_{Z'}-B'_{Z'})-\delta f'^*D\not\geq 0$ for any prime divisor $D$ on $Z'$ and any $\delta>0$.
    \end{itemize}
Since
$$A=f'^*(B_{Z'}-B'_{Z'})+(F-f'^*(\bar B_{Z'}-B'_{Z'})),$$
we have that $A\geq 0$ if and only if $B_{Z'}-B'_{Z'}\geq 0$. The proposition follows from Claim \ref{claim: a' for semistable}.
\end{proof}

\begin{prop}\label{prop: bp stable nef}
    Let $(X,B,\Mm)/U$ be an lc g-pair and $f: X\rightarrow Z$ a contraction such that
    \begin{itemize}
        \item $f: (X,B,\Mm)\rightarrow Z$ satisfies Property $(*)$,
        \item $(X,B,\Mm)$ is BP semistable$/Z$,
        \item $f$ is equidimensional, and
        \item $K_X+B+\Mm_X$ is nef$/Z$.
    \end{itemize}
    Let $\Nn$ be the moduli part of $f: X\rightarrow Z$. Then $\Nn_X$ is nef$/U$. Moreover, if $(X,B,\Mm)$ is BP stable$/Z$, then $\Nn$ descends to $X$, and in particular, $\Nn$ is nef$/U$.
\end{prop}
\begin{proof}
Let $\Ff$ be the foliation induced by $f$ and $B^h$ the horizontal$/Z$ part of $B$. By Proposition \ref{prop: bp semistable foliation lc}, $(X,B^h,\Mm)$ is lc and thus $(X,\Ff,B^h,\Mm;G)/Z$ is weak ACSS by definition. Since $K_X+B+\Mm_X$ is nef$/Z$, by Lemma \ref{lem: equivalence over bases}, $K_{\Ff}+B^h+\Mm_X$ is nef$/U$. By Proposition \ref{prop: weak cbf gfq}, $\Nn_X\sim K_{\Ff}+B^h+\Mm_X$ is nef$/U$. Moreover, if $(X,B,\Mm)$ is BP stable$/Z$, then $\Nn$ descends to $X$ by Lemma \ref{lem: bp stable implies descend}.
\end{proof}

\begin{lem}\label{lem: special acss model}
Let $(X,\Ff,B,\Mm)/U$ be an lc gfq such that $\Ff$ is induced by a contraction $f: X\rightarrow Z$. Let $D_Z$ be a divisor over $Z$. Then there exists an ACSS model $(X',\Ff',B',\Mm)/Z'$ of $(X,\Ff,B,\Mm)$ with induced morphisms $f': X'\rightarrow Z'$ and $g: X'\rightarrow X$, and a birational morphism $h_Z: Z'\rightarrow Z$ such that $h_Z\circ f'=f\circ g$ and $D_Z$ is on $Z'$.
\end{lem}
\begin{proof}
By Definition-Theorem \ref{defthm: weak ss reduction}, there exists an equidimensional model $(Y,\Sigma_Y,\Mm)\rightarrow Z'$ of $f: (X,B,\Mm)\rightarrow Z$ associated with $h: Y\rightarrow X$ and $h_Z: Z'\rightarrow Z$, such that $D_Z$ is on $Z'$. Set $\Ff_Y:=h^{-1}\Ff$ and $B_Y:=h^{-1}_*B+(\Supp\Exc(h))^{\Ff_Y}$, then $(Y,F_Y,B_Y,\Mm)$ is foliated log smooth, and
$$K_{\Ff_Y}+B_Y+\Mm_Y\sim_{\mathbb R,X}\sum_{E\subset\Exc(h)}(\epsilon_{\Ff}(E)-a(E,\Ff,B,\Mm))E\geq 0.$$
By Theorem \ref{thm: mmp very exceptional alg int fol}, we may run a $(K_{\Ff_Y}+B_Y+\Mm_Y)$-MMP$/X$ which terminates with a good minimal model$/X$ $(X',\Ff',B',\Mm)$. By Lemma \ref{lem: ACSS mmp can run}, this MMP is also a $(K_{\Ff_Y}+B_Y+\Mm_Y)$-MMP$/Z'$ and $(X',\Ff',B',\Mm)/Z'$ satisfies our requirements.
\end{proof}

\begin{thm}\label{thm: lc+weak acc=bpstable}
    Let $(X,\Ff,B,\Mm)/U$ be an lc gfq, $f: X\rightarrow Z$ a contraction, and $G$ a reduced divisor on $X$ such that $(X,\Ff,B,\Mm;G)/Z$ is weak ACSS. Then $(X,B+G,\Mm)$ is BP stable$/Z$.
\end{thm}
\begin{proof}
For any prime divisor $D_Z$ over $Z$, by Lemma \ref{lem: special acss model}, there exist two birational morphisms $h_Z: Z'\rightarrow Z$ and $h: X'\rightarrow X$, and an ACSS model $(X',\Ff',B',\Mm)/Z'$ of $(X,\Ff,B,\Mm)$ with induced morphism $f': X'\rightarrow Z'$, such that $f\circ h=h_Z\circ f'$ and $D_Z$ is on $Z'$. 

We let $G'$ be a divisor on $X'$ such that $(X',\Ff',B',\Mm;G')/Z'$ is ACSS. Let $B_Z$ and $\Nn$ be the discriminant and moduli parts of $f: (X,B+G,\Mm)\rightarrow Z$ respectively, $K_{Z'}+B_{Z'}:=h_Z^*(K_Z+B_Z)$, and let $B'_{Z'}$ and $\Nn'$ be the discriminant part of $f': (X',B'+G',\Mm)\rightarrow Z'$ respectively. 

By Proposition \ref{prop: weak cbf gfq}, we have
$$\Nn'_{X'}\sim K_{\Ff'}+B'+\Mm_{X'}=h^*(K_{\Ff}+B+\Mm_X)\sim h^*\Nn_X.$$
Thus for suitable choices of $\Nn'$ and $\Nn$, we may assume that $\Nn'_{X'}=h^*\Nn_X$. Let
$$K_{X'}+\tilde B'+\Mm_{X'}:=h^*(K_X+B+G+\Mm_X).$$
Let $\tilde B_{Z'}$ and $\tilde\Nn$ be the discriminant and moduli parts of $f': (X',\tilde B',\Mm)\rightarrow Z'$ respectively. Since
\begin{align*}
    &\Nn'_{X'}-h^*\Nn_X\\
    =&\left(\left(K_{X'}+B'+G'+\Mm_{X'}\right)-f'^*\left(K_{Z'}+B'_{Z'}\right)\right)-h^*\left(K_X+B+G+\Mm_X-f^*(K_Z+B_Z)\right)\\
    =&B'+G'-\tilde B'-f'^*(B'_{Z'}-B_{Z'}),
\end{align*}
it follows that $B'+G'-\tilde B'-f'^*(B'_{Z'}-B_{Z'})=0$. Moreover, by Proposition \ref{prop: bp semistable foliation lc}, $(X,B+G,\Mm)$ is BP semistable$/Z$. Thus $\mult_{D_Z}\tilde B_{Z'}\leq\mult_{D_Z}B_{Z'}.$ For any prime divisor $D$ on $X'$ with $f'(D)=D_Z$, we have $\mult_DB'=0$ and
$$\mult_DG'=\mult_D\left(\tilde B'+f'^*\left(B'_{Z'}-B_{Z'}\right)\right).$$
There are two cases.

\medskip

\noindent\textbf{Case 1}. $D_Z$ is not a component of $B'_{Z'}$. 

In this case, $\mult_DG'=0$, $D_Z\not\subset f'(G')$, and
$$1=\sup\left\{t\mid \left(X',B'+G'+f'^*D_Z,\Mm\right)\text{ is lc over the generic point of }D_Z\right\}.$$
For any prime divisor $D$ on $X'$ that is vertical$/Z'$, since $f'$ is equidimensional, $f'(D)$ is a prime divisor $D_Z$ on $Z'$. Moreover, $\mult_DB'=0$ since $B'$ is horizontal$/Z'$.

If $D_Z$ is not a component of $B'_{Z'}$, then $\mult_DG'=0$, and
$$\mult_D\tilde B'=\mult_Df'^*B_{Z'}=\mult_{D_Z}B_{Z'}\cdot\mult_Df'^*D_Z.$$
Thus $\mult_D\left(\tilde B'-\mult_{D_Z}B_{Z'}f'^*D_Z\right)=0,$ and
$$\mult_D\left(\tilde B'+(1-\mult_{D_Z}B_{Z'})f'^*D_Z\right)=\mult_Df'^*D_Z\geq 1.$$
It follows that
\begin{align*}
1-\mult_{D_Z}\tilde B_{Z'}&=\sup\left\{t\mid \left(X',\tilde B'+tf'^*D_Z,\Mm\right)\text{ is lc over the generic point of }D_Z\right\}\\
&\leq 1-\mult_{D_Z}B_{Z'}
\end{align*}
and hence
$$\mult_{D_Z}B_{Z'}=\mult_{D_Z}\tilde B_{Z'}.$$

\medskip

\noindent\textbf{Case 2}. $D_Z$ is a component of $B'_{Z'}$. 

In this case, $\mult_DG'=1$ and $\mult_DB'_{Z'}=1$. Therefore,
$$1=\mult_D\left(\tilde B'+\left(\mult_DB'_{Z'}-\mult_DB_{Z'}\right)f'^*D_Z\right)=\mult_D\left(\tilde B'+\left(1-\mult_DB_{Z'}\right)f'^*D_Z\right).$$
Thus one can see that
\begin{align*}
1-\mult_{D_Z}\tilde B_{Z'}&=\sup\{t\mid (X',\tilde B'+tf'^*D_Z,\Mm)\text{ is lc over the generic point of }D_Z\}\\
&\leq 1-\mult_{D_Z}B_{Z'}
\end{align*}
which implies that
$$\mult_{D_Z}B_{Z'}=\mult_{D_Z}\tilde B_{Z'}.$$

\medskip

In either case, we have $\mult_{D_Z}B_{Z'}=\mult_{D_Z}\tilde B_{Z'}$. Since $D_Z$ can be any prime divisor over $Z$, $(X,B+G,\Mm)$ is BP stable$/Z$. We finish the proof.
\end{proof}

\subsection{Numerical dimension zero generalized foliated quadruples}\label{sub: num 0 mmp}

\begin{prop}\label{prop: weak ss num 0 mmp}
    Let $(X,\Ff,B,\Mm)/U$ be a $\Qq$-factorial lc gfq such that
    \begin{itemize}
        \item $(X,\Ff,B,\Mm)$ is weak ACSS,
        \item $\kappa_{\sigma}(X/U,K_{\Ff}+B+\Mm_X)=0$, and
        \item either $X$ is klt or $\Mm$ is NQC$/U$.
    \end{itemize}
    Let $\mathcal{P}$ be any $(K_{\Ff}+B+\Mm_X)$-MMP$/U$ with scaling of an ample$/U$ $\Rr$-divisor $A$, whose existence is guaranteed by Proposition \ref{prop: run mmp with scaling gfq}. Then after a sequence of steps in $\mathcal{P}$, we get a log birational model $(X',\Ff',B',\Mm)/U$ satisfying the following.
    \begin{enumerate}
        \item For any very general fiber $F'$ of $X'\rightarrow U$, $(K_{\Ff'}+B'+\Mm_{X'})|_{F'}\equiv 0$, and if $\kappa_{\iota}(X/U,K_{\Ff}+B+\Mm_X)=0$, then $(K_{\Ff'}+B'+\Mm_{X'})|_{F'}\sim_{\mathbb R}0.$
        \item Suppose that $\pi:X\to U$ is an equidimensional contraction and $U$ is $\Qq$-factorial. 
        \begin{enumerate}
            \item If $K_{\Ff}+B+\Mm_X\equiv_U\text{(resp. }\sim_{\mathbb R,U}\text{) }E_1+E_2$ for some $\Rr$-divisors $E_1$ and $E_2$ such that $E_1\geq 0$ and $E_2$ is vertical$/U$. Then
            $$K_{\Ff'}+B'+\Mm_{X'}\equiv_U\text{(resp. }\sim_{\mathbb R,U}\text{) }0.$$
            In particular, $(X',\Ff',B',\Mm)/U$ is a weak lc model of $(X,\Ff,B,\Mm)/U$.
            \item If $\kappa_{\iota}(X/U,K_{\Ff}+B+\Mm_X)=0$, then
            \begin{enumerate}
                \item $K_{\Ff'}+B'+\Mm_{X'}\sim_{\mathbb R,U}0.$
                \item $(X',\Ff',B',\Mm)/U$ is a weak lc model of $(X,\Ff,B,\Mm)/U$.
                \item If $(X,\Ff,B,\Mm)/U$ is ACSS, then $(X',\Ff',B',\Mm)/U$ is a good minimal model of $(X,\Ff,B,\Mm)/U$.
            \end{enumerate}
        \end{enumerate}
    \end{enumerate}
\end{prop}
\begin{proof}
Let $(X_0,\Ff_0,B_0,\Mm):=(X,\Ff,B,\Mm)$.
We denote $\mathcal{P}$ by
\begin{center}$\xymatrix{
(X_0,\Ff_0,B_0,\Mm)\ar@{-->}[r]^{f_0} & (X_1,\Ff_1,B_1,\Mm)\ar@{-->}[r]^{\ \ \ \ \ \ \ \ f_1} & \dots\ar@{-->}[r] & (X_n,\Ff_n,B_n,\Mm)\ar@{-->}[r]^{\ \ \ \ \ \ \ \ \ f_n} & \dots 
}.$
\end{center}
Let $A_i$ be the strict transform of $A$ on $X_i$, $\pi_i: X_i\rightarrow U$ the induced contraction for each $i$, and
$$\lambda_i:=\inf\{t\geq 0\mid K_{\Ff_i}+B_i+\Mm_{X_i}+tA_i\text{ is nef}/U\}$$
the scaling numbers. According to Proposition \ref{prop: run mmp with scaling gfq}, either $\mathcal{P}$ terminates, or $\lim_{i\rightarrow+\infty}\lambda_i=0$. 

If $\mathcal{P}$ terminates, then we let $(X_m,\Ff_m,B_m,\Mm)/U$ be the output of $\mathcal{P}$. If $\mathcal{P}$ does not terminate, then we let $m$ be a positive integer such that $f_i$ is a flip for any $i\geq m$. We let $\phi_i: X_m\dashrightarrow X_i$ be the induced birational map for any $i\geq m$. Since $\mathcal{P}$ contains countably many steps, there are countably many closed points $z\in Z$ such that for some $i\geq m$, $\pi_i^{-1}(z)$ is contained in either the flipping locus of $f_i$ or the flipped locus of $f_{i-1}$. Therefore, for a very general point $z\in Z$ and any $i\geq m$, $\pi_i^{-1}(z)$ is neither contained in the flipping locus of $f_i$ nor the flipped locus of $f_{i-1}$. We let $F_m$ be a very general fiber of $\pi_m$, $z_0:=\pi_m(F_m)$, and let $F_{i}$ be the fiber of $\pi_i$ over $z_0$ for each $i$. Then the induced birational map 
$$\phi_{F,i}:=\phi_i|_{F_m}: F_m\dashrightarrow F_i$$
is small for any $i\geq m$. Let $\Mm^F:=\Mm|_{F_m}$, $B_{F_i}:=B_i|_{F_i}$, and $A_{F_i}:=A_i|_{F_i}$ for each $i\geq m$. Note that for each $i\geq m$, we have $K_{F_i}=K_{\Ff_i}|_{F_i}$. 

We will show that $(X',\Ff',B',\Mm)/U=(X_m,\Ff_m,B_m,\Mm)/U$ satisfies our requirements.

\begin{claim}\label{claim: movable along very general fiber}
$K_{\Ff_m}+B_m+\Mm_{X_m}$ is movable$/U$ and $(K_{\Ff_m}+B_m+\Mm_{X_m})|_{F_m}$ is movable.
\end{claim}
\begin{proof}
  If $K_{\Ff_m}+B_m+\Mm_{X_m}$ is nef$/U$ then the claim is obvious, so we may assume that $K_{\Ff_m}+B_m+\Mm_{X_m}$ is not nef$/U$. In particular, $\mathcal{P}$ does not terminate. Since $K_{\Ff_i}+B_i+\Mm_{X_i}+tA_i$ is nef$/U$ for any $i\geq m$,
  $$K_{X_m}+B_{m}+\Mm_{X_m}=\lim_{i\rightarrow+\infty}(\phi_{i}^{-1})_*(K_{X_i}+B_{i}+\Mm_{X_i}+tA_{i})$$
  is movable$/U$, and
  $$K_{F_i}+B_{F_i}+\Mm^F_{F_i}+tA_{F_i}=(K_{\Ff_i}+B_i+\Mm_{X_i}+tA_i)|_{F_i}$$
is nef for each $i\geq m$. Thus
$$K_{F_m}+B_{F_m}+\Mm^F_{F_m}=\lim_{i\rightarrow+\infty}(\phi_{F,i}^{-1})_*(K_{F_i}+B_{F_i}+\Mm^F_{F_i}+tA_{F_i})$$
is movable, and the claim follows.
\end{proof}
\noindent\textit{Proof of Proposition \ref{prop: weak ss num 0 mmp} continued}. Since $\kappa_{\sigma}(X/U,K_{\Ff}+B+\Mm_X)=0$, $\kappa_{\sigma}(X_m/U,K_{\Ff_m}+B_m+\Mm_{X_m})=0$ and $\kappa_{\sigma}(K_{F_m}+B_{F_m}+\Mm^F_{F_m})=0$. By Claim \ref{claim: movable along very general fiber} and Lemma \ref{lem: movable num 0 is 0}, $K_{F_m}+B_{F_m}+\Mm^F_{F_m}\equiv 0$. Suppose that $\kappa_{\iota}(X/U,K_{\Ff}+B+\Mm_X)=0$, then $\kappa_{\iota}(X_m/U,K_{\Ff_m}+B_m+\Mm_{X_m})=0$ and hence $\kappa_{\iota}(K_{F_m}+B_{F_m}+\Mm^F_{F_m})=0$. It implies that $K_{F_m}+B_{F_m}+\Mm^F_{F_m}\sim_{\mathbb R}0$. This implies (1).

From now on, we may assume that $\pi$ is equidimensional and $U$ is $\Qq$-factorial. Suppose that $K_{\Ff}+B+\Mm_X\equiv_U\text{(resp. }\sim_{\mathbb R,U}\text{) }E_1+E_2$ where $E_1\geq 0$ and $E_2$ is vertical$/U$. Since $U$ is $\Qq$-factorial, for any prime divisor $D$ on $U$, $\pi^*D$ is $\Rr$-Cartier, and we may define
$$t_D:=\sup\{s\mid E_2-s\pi^*D\geq 0\text{ over the generic point of }D\}.$$
Let
$$\tilde E_2:=E_2-\sum t_D\pi^*D,$$
where $D$ runs over the prime divisors on $U$. Since $\pi$ is equidimensional, $\tilde E_2\geq 0$ and $\tilde E_2$ is very exceptional$/U$. Possibly replacing $E_2$ with $\tilde E_2$, we may assume that $E_2 \ge 0$ is very exceptional$/U$. Let $E_{1,m}$ and $E_{2,m}$ be the strict transforms of $E_1$ and $E_2$ on $X'=X_m$ respectively. Then $E_{1,m}|_{F_m}\equiv 0$, and thus $E_{1,m}|_{F_m}=0$ and
$$E_{2,m}\equiv_U\text{(resp. }\sim_{\mathbb R,U}\text{) }K_{\Ff_m}+B_m+\Mm_{X_m}$$
is movable$/U$. For any prime divisor $S$ on $X_m$ and very general curves $C$ on $S$ over $U$, $E_{2,m}\cdot C\geq 0$. By \cite[Lemma 3.3]{Bir12}, we see that $E_{2,m}=0$. This implies (2.a).

If $\kappa_{\iota}(X/U,K_{\Ff}+B+\Mm_X)=0$, then $K_{\Ff}+B+\Mm_X\sim_{\mathbb R,U}E\geq 0$ for some $\Rr$-divisor $E$ on $X$. We let $E^h$ be the horizontal$/U$ part of $E$ and let $E^v$ be the vertical$/U$ part of $E$. Then (2.b.i) and (2.b.ii) follow from (2.a). (2.b.iii) follows from (2.b.i), (2.b.ii), and Lemma \ref{lem: ACSS mmp can run}.
\end{proof}

The following proposition is similar to Proposition \ref{prop: weak ss num 0 mmp}, but does not immediately follow from Proposition \ref{prop: weak ss num 0 mmp}.

\begin{prop}\label{prop: projective num 0 mmp}
Let $(X,\Ff,B,\Mm)$ be a projective lc gfq such that
\begin{itemize}
   \item $(X,\Ff,B,\Mm)$ is weak ACSS,
   \item $\kappa_{\sigma}(K_{\Ff}+B+\Mm_X)=0$, and
   \item either $X$ is $\Qq$-factorial klt or $\Mm$ is NQC.
\end{itemize} Then for any ample $\Rr$-divisor $A$, there exists a $(K_{\Ff}+B+\Mm_X)$-MMP with scaling of $A$, say $\mathcal{P}_0$, satisfying the following. Let $\mathcal{P}=\mathcal{P}_0$ if $X$ is not $\Qq$-factorial, and let $\mathcal{P}$ be any $(K_{\Ff}+B+\Mm_X)$-MMP with scaling of an ample $\Rr$-divisor if $X$ is $\Qq$-factorial. Then
\begin{enumerate}
   \item $\mathcal{P}$ terminates with a weak lc model $(X',\Ff',B',\Mm)$ of $(X,\Ff,B,\Mm)$ such that
   $$K_{\Ff'}+B'+\Mm_{X'}\equiv 0.$$
   \item Suppose that $\kappa_{\iota}(K_{\Ff}+B+\Mm_X)=0$. Then
   \begin{enumerate}
       \item $K_{\Ff'}+B'+\Mm_{X'}\sim_{\mathbb R}0.$
       \item If $(X,\Ff,B,\Mm)$ is $\Qq$-factorial ACSS, then $(X',\Ff',B',\Mm)$ is a good minimal model of $(X,\Ff,B,\Mm)$.
   \end{enumerate}
\end{enumerate}
\end{prop}
\begin{proof}
Let $(X_0,\Ff_0,B_0,\Mm):=(X,\Ff,B,\Mm)$. According to Proposition \ref{prop: run mmp with scaling gfq}, we may suppose that $\mathcal{P}$ is an MMP with scaling of $A$
     \begin{center}$\xymatrix{
(X_0,\Ff_0,B_0,\Mm)\ar@{-->}[r]^{f_0} & (X_1,\Ff_1,B_1,\Mm)\ar@{-->}[r]^{\ \ \ \ \ \ \ \ f_1} & \dots\ar@{-->}[r] & (X_n,\Ff_n,B_n,\Mm)\ar@{-->}[r]^{\ \ \ \ \ \ \ \ \ f_n} & \dots 
}$
\end{center}
such that either this MMP terminates, or $\lim_{i\rightarrow+\infty}\lambda_i=0$, where $\lambda_i$ is the scaling number and $A_i$ is the strict transforms of $A$ on $X_i$ for each $i$. We first show that $\mathcal{P}$ terminates. Otherwise, there is an integer $m$ such that $f_i$ is a flip for any $i\geq m$. Since $\lim_{i\rightarrow+\infty}\lambda_i=0$, we can see that 
$$K_{X_m}+B_m+\Mm_{X_m}=\lim_{i\rightarrow+\infty}(\phi_i^{-1})_*(K_{X_i}+B_i+\lambda_iA_i+\Mm_{X_i})$$
is movable, where $\phi_i: X_m\dashrightarrow X_i$ is the induced birational map for any $i\geq m$. By Lemma \ref{lem: movable num 0 is 0}, $K_{X_m}+B_m+\Mm_{X_m}\equiv 0$, which is absurd. Thus $\mathcal{P}$ terminates with a weak lc model $(X',\Ff',B',\Mm)$ such that $K_{\Ff'}+B'+\Mm_{X'}\equiv 0$ as $\kappa_{\sigma}(K_{\Ff'}+B'+\Mm_{X'})=0$. This implies (1).

If $\kappa_{\iota}(K_{\Ff}+B+\Mm_X)=0$, then $\kappa_{\iota}(K_{\Ff'}+B'+\Mm_{X'})=0$, so (2.a) follows from (1). (2.b) follows from (1), (2.a), and Lemma \ref{lem: ACSS mmp can run}.
\end{proof}

\subsection{Refined definition of lc-trivial fibrations}

\begin{defn}\label{defn: lc trivial fibration gfq}
Let $(X,\Ff,B,\Mm)/U$ be a sub-gfq and $f: X\rightarrow Z$ a contraction$/U$, such that the general fibers of $f$ are tangent to $\Ff$. We say that $f: (X,\Ff,B,\Mm)\rightarrow Z$ is an \emph{lc-trivial fibration} if
\begin{enumerate}
\item $(X,\Ff,B,\Mm)$ is sub-lc over the generic point of $Z$,  
\item $K_{\Ff}+B+\Mm_X\sim_{\mathbb R,Z}0$, and
\item there exists a birational morphism $h: Y\rightarrow X$ with $\Ff_Y := h^{-1}\Ff$ and $K_{\Ff_Y}+B_Y+\Mm_Y=h^*(K_{\Ff}+B+\Mm_X)$, such that $-B_Y^{\leq 0}$ is $\Rr$-Cartier and 
$$\kappa_{\sigma}(Y/Z,-B_Y^{\leq 0})=0.$$
\end{enumerate}
It is clear that an lc-trivial fibration does not depend on the choice of $U$.
\end{defn}

\begin{rem}\label{rem: lc trivial fibration definition}
    It is very important to notice that our definition of an lc-trivial fibration is different from the one defined in the classical way, even when $\Mm=\bm{0}$ and $\Ff=T_X$. This is with good reason. For simplicity, in the following, we shall assume that $\Ff=T_X$. 
    
 In the classical definition, condition (3) is replaced with the condition
    \begin{enumerate}
        \item[(3')] $\rk f_*\mathcal{O}_X(\lceil\Aa^*(X,B,\Mm)\rceil)=1.$
    \end{enumerate} 
    The condition (3') was used in the earliest version of the canonical bundle formula \cite[Condition (3) of Theorem 2]{Kaw98}. It has also been used in later versions of the canonical bundle formula, e.g., \cite[Theorem 0.2]{Amb05} for sub-klt sub-pairs, \cite[Theorem 8.3.7]{Kol07} (see also \cite[Theorem 3.6]{FG14}) for lc-trivial fibrations for sub-lc sub-pairs. 
    
    However, for generalized sub-pairs, when condition (3') is being applied, we cannot get a complete version of the canonical bundle formula. More precisely, for lc-trivial fibrations for NQC generalized pairs defined under condition (3') instead of (3), we have to add one of the following two conditions to get the canonical bundle formula.
    \begin{itemize}
        \item[(4.1)] $B\geq 0$ over the generic point of $Z$ (rational coefficient case \cite[Theorem 2.20]{FS23}; real coefficient case \cite[Theorem 2.23]{JLX22}).
        \item[(4.2)] $\Mm$ is $\bb$-semi-ample$/Z$ (rational coefficient case \cite[Chapter 6, Theorem 7]{Fil19}; real coefficient case \cite[Theorem 2.23]{JLX22}).
    \end{itemize}
 The canonical bundle formula for NQC generalized pairs under these additional conditions (4.1) or (4.2) is usually enough for us to apply, as it guarantees that the structure of NQC generalized pairs is preserved under the canonical bundle formula. However, this causes big trouble when we want to consider the canonical bundle formula for generalized foliated quadruples. This is because, in the construction of the canonical bundle formula for generalized foliated quadruples, \cite[Definition-Theorem 6.12]{LLM23}, we need to pass through an equidimensional model. But then we will potentially get a sub-lc g-sub-pair with some negative coefficients, which in general do not satisfy (4.1) or (4.2). In this case, we cannot define the canonical bundle formula for generalized foliated quadruples in general. Condition (3), on the other hand, is introduced to resolve this issue.

 In fact, the most important cases of the canonical bundle formula are the cases when $B\geq 0$ over the generic point of $Z$. However, since we need to check the coefficients of the discriminant part along any high model of the base, we need to do base change in many scenarios. Therefore, we have to consider the crepant pullbacks. Moreover, since we may want to run minimal model programs over the base to get new structures, we need to consider crepant transformations over the generic point of $Z$. This will inevitably introduce sub-pairs or g-sub-pairs and force us to consider a larger category of structures so that the canonical bundle formula can be applied. More precisely, we want to find a category $\mathcal{D}$ of structures
 $$f: (X,B,\Mm)\rightarrow Z,$$
such that
\begin{enumerate}
\item[(i)] For any g-sub-pair $(X,B,\Mm)/U$ and contraction$/U$ $f: X\rightarrow Z$ such that $K_X+B+\Mm_X\sim_{\mathbb R,Z}0$ and $(X,B,\Mm)$ is lc over the generic point of $Z$, $f: (X,B,\Mm)\rightarrow Z$ belongs to $\mathcal{D}$.
\item[(ii)] For any g-sub-pairs $(X,B,\Mm)/U$ and $(X',B',\Mm')/U$ and birationally equivalent contractions$/U$ $f: X\rightarrow Z$, $f': X'\rightarrow Z'$ such that $K_X+B+\Mm_X\sim_{\mathbb R,Z}0$, $K_{X'}+B'+\Mm_{X'}\sim_{\mathbb R,Z'}0$, and $(X,B,\Mm)$ and $(X',B',\Mm')$ are crepant over the generic point of $Z$, $f: (X,B,\Mm)\rightarrow Z$ belongs to $\mathcal{D}$ if and only if $f': (X',B',\Mm')\rightarrow Z'$ belongs to $\mathcal{D}$.
\end{enumerate}
 Condition (3') is actually one natural condition to add in order to form the category $\mathcal{D}$. For generalized pairs, however, the category $\mathcal{D}$ constructed by adding condition (3') is a category that is too large to prove the canonical bundle formula in general. In fact, comparing our condition (3) with condition (3'), we can easily see that (3') can be roughly interpreted as
 $$\kappa(Y/Z,-B_Y^{\leq 0})=0$$
 (cf. \cite[Definitions 8.4.1, 8.4.2]{Kol07}). This actually means that some kind of existence of good minimal models should hold for generalized pairs with Kodaira dimension $0$. But this is absurd due to numerous counterexamples (cf. \cite[1.1 Example]{Sho00}). In fact, even for usual pairs, since the existence of good minimal models is unknown for pairs with Kodaira dimension $0$, the use of mixed Hodge structures was essentially used in all literature on the canonical bundle formula except \cite{ACSS21}, while \cite{ACSS21} does not deal with lc-trivial fibrations in general.

Therefore, instead of considering the classical category defined using (3'), we turn to consider a new category $\mathcal{D}$ of g-sub-pairs which satisfies (i) and (ii) but does not rely on condition (3'). It turns out that (3) is also a natural condition (see Lemmas \ref{lem: lc trivial preserved crepant} and \ref{lem: lc trivial holds for lc gpair} below) for us to add, and it turns out that we can completely bypass the abundance conjecture or the mixed Hodge structure by using (3) instead of (3'). This will eventually lead to the canonical bundle formula for non-NQC g-pairs and gfqs in full generality.
\end{rem}

The following lemmas are analogues of Lemmas \ref{lem: lc trivial (3) invariant under pullback}, \ref{lem: lc trivial preserved crepant}, and \ref{lem: lc trivial holds for lc gpair} for foliations, and their proofs are similar.

\begin{lem}\label{lem: lc trivial (3) invariant under pullback}
    Let $(X,\Ff,B,\Mm)/U$ be a sub-gfq. Assume that $-B^{\leq 0}$ is $\Rr$-Cartier and $\kappa_{\sigma}(X/U,-B^{\leq 0})=0$. Then for any birational morphism $g: W\rightarrow X$, such that
    \begin{enumerate}
        \item $K_{g^{-1}\Ff}+B_W+\Mm_W:=g^*(K_{\Ff}+B+\Mm_X)$ satisfies that $-B_W^{\leq 0}$ is $\Rr$-Cartier, and
        \item there exists an $\Rr$-Cartier $\Rr$-divisor $0\leq F\subset\Supp\Exc(g)$,
    \end{enumerate}
we have that
$$\kappa_{\sigma}(W/U,-B_W^{\leq 0})=0.$$
\end{lem}
\begin{proof}
Let $D:=-B^{\leq 0}$, $D_W:=-B_W^{\leq 0}$, and $m\gg 0$ an integer. Then we have
$$D_W=g^{-1}_*D+E$$
for some $E\geq 0$ that is exceptional$/X$. Thus
$$0=\kappa_{\sigma}(X/U,D)=\kappa_{\sigma}(W/U,g^*D+mF)\geq\kappa_{\sigma}(W/U,g^{-1}_*D+E)=\kappa_{\sigma}(W/U,D_W)\geq 0.$$
So $\kappa_{\sigma}(W/U,D_W)=0$.
\end{proof}

\begin{lem}\label{lem: lc trivial preserved crepant}
    Let $(X,\Ff,B,\Mm)/U$ and $(X',\Ff',B',\Mm')/U$ be two sub-gfqs. Let $f: X\rightarrow Z$ and $f': X'\rightarrow Z'$ be two birationally equivalent contractions$/U$ such that $(X,\Ff,B,\Mm)$ and $(X',\Ff',B',\Mm')$ are crepant over the generic point of $Z$. Assume that $K_{\Ff}+B+\Mm_X\sim_{\mathbb R,Z}0$, and $K_{\Ff'}+B'+\Mm_{X'}\sim_{\mathbb R,Z'}0$.

    Then $f: (X,\Ff,B,\Mm)\rightarrow Z$ is an lc-trivial fibration if and only if $f': (X',\Ff',B',\Mm')\rightarrow Z'$ is an lc-trivial fibration.
\end{lem}
\begin{proof}
By symmetry, we only need to prove the only if part, and we may assume that $f: (X,\Ff,B,\Mm)\rightarrow Z$ is an lc-trivial fibration.

Let $p: W\rightarrow X$ and $q: W\rightarrow X'$ be a resolution of indeterminacy of the induced birational map $\phi: X\dashrightarrow X'$ such that $\Mm$ descends to $W$, $\Ff_W:=p^{-1}\Ff=q^{-1}\Ff'$,
$K_{\Ff_W}+B_W+\Mm_{W}:=p^*(K_{\Ff}+B+\Mm_X)$, and $K_{\Ff_W}+B'_W+\Mm'_{W}:=q^*(K_{\Ff'}+B'+\Mm'_{X'})$. Moreover, by Lemma \ref{lem: lc trivial (3) invariant under pullback}, possibly replacing $W$ with a higher resolution, we may assume that $W$ is smooth and $\kappa_{\sigma}(W/Z,-B_W^{\leq 0})=0$. 

Since $(X,\Ff,B,\Mm)$ and $(X',\Ff',B',\Mm')$ are crepant over the generic point of $Z$, over the generic point of $Z$, we have $B_W=B'_W$, $\Ff_W$ is the common transform of $\Ff, \Ff'$, and $\Mm_W=\Mm'_W$. Thus $\kappa_{\sigma}(W/Z,-B_W'^{\leq 0})=0$. Moreover, since $(X,\Ff,B,\Mm)$ is sub-lc over the generic point of $Z$, $(W,\Ff_W,B_W,\Mm)$ is sub-lc over the generic point of $Z$, so $(W,\Ff_W,B_W',\Mm')$ is sub-lc over the generic point of $Z$, so $(W,\Ff_W,B_W',\Mm')$ is sub-lc over the generic point of $Z'$, so $(X',\Ff',B',\Mm')$ is sub-lc over the generic point of $Z'$. The lemma follows.
\end{proof}

\begin{lem}\label{lem: lc trivial holds for lc gpair}
    Let $(X,\Ff,B,\Mm)/U$ be a sub-gfq and $f: X\rightarrow Z$ a contraction$/U$ such that $(X,\Ff,B,\Mm)$ is lc over the generic point of $Z$ and $K_{\Ff}+B+\Mm_X\sim_{\mathbb R,Z}0$. Then $f: (X,\Ff,B,\Mm)\rightarrow Z$ is an lc-trivial fibration.
\end{lem}
\begin{proof}
    Over the generic point of $Z$, $B^{\leq 0}=0$, so $\kappa_{\sigma}(X/Z,B^{\leq 0})=0$. The lemma follows from the definition.
\end{proof}

\subsection{Canonical bundle formula for generalized pairs}\label{subsec: cbf gpair}

\begin{defn}\label{defn: cbf gpair}
Let $(X,B,\Mm)/U$ be a g-sub-pair and $f: X\rightarrow Z$ a contraction$/U$ such that $f: (X,B,\Mm)\rightarrow Z$ is an lc-trivial fibration$/U$. Then there exists an $\Rr$-divisor $L$ on $Z$ such that $K_X+B+\Mm_X\sim_{\mathbb R}f^*L$. There exists a unique (up to $\Rr$-linear equivalence) $\bb$-divisor $\Mm^Z$ on $Z$ satisfying the following.
    
Let $f': X'\rightarrow Z'$ be any contraction that is birationally equivalent to $f$ such that the induced birational maps $h: X'\dashrightarrow X$ and $h_Z: Z'\dashrightarrow Z$ are morphisms. We let $$K_{X'}+B'+\Mm_{X'}:=h^*(K_X+B+\Mm_X)$$
and let $B_{Z'}$ be the discriminant part of $f': (X',B',\Mm)\rightarrow Z'$. Then
$$\Mm^Z_{Z'}=h_Z^*L-K_{Z'}-B_{Z'}.$$
We call $\Mm^Z$ the \emph{base moduli part} of $f: (X,B,\Mm)\rightarrow Z$. If there is no confusion, we may also call $\Mm^Z$ the \emph{moduli part} of $f: (X,B,\Mm)\rightarrow Z$. It is clear that $\Mm^Z$ only depends on the choice of $L$, which is unique up to $\Rr$-linear equivalence.
\end{defn}

\begin{lem}\label{lem: order along generic fiber cbf}
Let $(X,B,\Mm)/U$ be a g-sub-pair and $f: X\rightarrow Z$ a contraction such that $f: (X,B,\Mm)\rightarrow Z$ is an lc-trivial fibration$/U$. Suppose that $n(K_X+B+\Mm_X)\sim 0$ over the generic point of $Z$ for some positive integer $n$. Then there exists a choice $\Mm^Z$ of the base moduli part of $f: (X,B,\Mm)\rightarrow Z$ such that 
$$n(K_X+B+\Mm_X)\sim nf^*\left(K_Z+B_Z+\Mm^Z_Z\right),$$
where $B_Z$ is the discriminant part of $f: (X,B,\Mm)\rightarrow Z$.
\end{lem}
\begin{proof}
By assumption, there exists a rational function $\psi\in K(X)$ such that $n(K_X+B+\Mm_X)+(\psi)$ is vertical$/Z$. Then by \cite[Lemma 2.5]{CHL23}, there exists an $\Rr$-Cartier $\Rr$-divisor $L$ on $Z$ such that
$$n(K_X+B+\Mm_X)+(\psi)=nf^*L.$$
The lemma follows from our construction of $\Mm^Z$ as in Definition \ref{defn: cbf gpair}.
\end{proof}

\begin{lem}\label{lem: m preserved under crepant}
    Let $(X,B,\Mm)/U$ and $(X',B',\Mm')/U$ be two g-sub-pairs. Let $f: (X,B,\Mm)\rightarrow Z$ and $f': (X',B',\Mm')\rightarrow Z'$ be two lc-trivial fibrations$/U$ such that $f$ and $f'$ are birationally equivalent, and $(X,B,\Mm)$ and $(X',B',\Mm')$ are crepant over the generic point of $Z$. Let $\Mm^Z$ be the base moduli part of $f: (X,B,\Mm)\rightarrow Z$ and let $\Mm^{Z'}$ be the base moduli part of $Z'$. Then $\Mm^Z\sim_{\mathbb R}\Mm^{Z'}$.
\end{lem}
\begin{proof}
Possibly passing to a common base and resolving the indeterminacy of the induced birational map $X\dashrightarrow X'$, we may assume that $f=f'$, $X=X'$, and $Z=Z'$. Now $K_X+B+\Mm_X=K_{X}+B'+\Mm_{X}$ over the generic point of $Z$, so $B-B'$ is vertical$/Z$. Since $K_X+B+\Mm_X\sim_{\mathbb R,Z}0$ and $K_{X}+B'+\Mm_{X}\sim_{\mathbb R,Z}0$, $B-B'\sim_{\mathbb R,Z}0$, so $B-B'=f^*P$ for some $\Rr$-divisor $P$ on $Z$ by \cite[Lemma 2.5]{CHL23}. 

Let $B_Z$ and $B_Z'$ be the discriminant parts of $f: (X,B,\Mm)\rightarrow Z$ and $f: (X,B',\Mm)\rightarrow Z$ respectively. By the definition of the discriminant part, $B_Z=B_Z'+P$. Since
$$K_Z+B_Z'+P+\Mm^{Z'}_{Z}\sim_{\mathbb R}K_Z+B_Z+\Mm^Z_Z,$$
$\Mm^{Z'}_{Z}\sim_{\mathbb R}\Mm^Z_Z$. Since we may pass to an arbitrarily high base change, we have $\Mm^Z\sim_{\mathbb R}\Mm^{Z'}$.
\end{proof}

\begin{thm}\label{thm: cbf gpair nonnqc}
Let $(X,B,\Mm)/U$ be a g-sub-pair and $f: X\rightarrow Z$ a contraction$/U$ such that $f: (X,B,\Mm)\rightarrow Z$ is an lc-trivial fibration. Let $B_Z$ and $\Mm^Z$ be the discriminant part and a base moduli part of $f: (X,B,\Mm)\rightarrow Z$ respectively. Then $\Mm^Z$ is nef$/U$. Moreover
\begin{enumerate}
    \item $(Z,B_Z,\Mm^Z)/U$ is a g-sub-pair.
    \item If the vertical$/Z$ part of $B$ is $\geq 0$, then $(Z,B_Z,\Mm^Z)/U$ is a g-pair.
    \item If $(X,B,\Mm)$ is sub-lc (resp. lc, sub-klt, klt), then $(Z,B_Z,\Mm^Z)$ is sub-lc (resp. lc, sub-klt, klt).
    \item Any lc center of $(Z,B_Z,\Mm^Z)$ is the image of an lc center of $(X,B,\Mm)$.
    \item The image of any lc center of $(X,B,\Mm)$ on $Z$ is an lc center of $(Z,B_Z,\Mm^Z)$.
    \item If $\Mm$ is NQC$/U$, then $\Mm^Z$ is NQC$/U$.
\end{enumerate}
\end{thm}
\begin{proof}
According to Definition-Theorem \ref{defthm: weak ss reduction}, $f: (X,B,\Mm)\rightarrow Z$ has an equidimensional model $f': (X',\Sigma_{X'},\Mm)\rightarrow Z'$ with associated $h:X'\to X$. 
Let
$$K_{X'}+\tilde B'+\Mm_{X'}:=h^*(K_X+B+\Mm_X),$$
$\tilde B'^h$ the horizontal$/Z'$ part of $\tilde B'$, and $B':=(\tilde B'^h)^{\geq 0}$. Let $G'$ be the vertical$/Z'$ part of $\Sigma_{X'}$, $\tilde B'^v$ the vertical$/Z'$ part of $\tilde B'$, $E^h:=-(\tilde B'^h)^{\leq 0}$, and $E^v:=G'-\tilde B'^v$. Then $E^h\geq 0$ and $E^v$ is vertical$/Z'$. By Lemma \ref{lem: lc trivial (3) invariant under pullback}, $\kappa_{\sigma}(X'/Z,E^h)=0$. We have
\begin{align*}
    K_{X'}+B'+G'+\Mm_{X'}&=h^*(K_X+B+\Mm_X)+B'+G'-\tilde B'\\
    &\sim_{\mathbb R,Z'}\left(B'-\tilde B'^h\right)+G'-\tilde B'^v=-\left(\tilde B'^h\right)^{\leq 0}+\left(G'-\tilde B'^v\right)=E^h+E^v.
\end{align*}
Since $(X,B,\Mm)$ is sub-lc over the generic point of $Z$, $\Sigma_{X'}\geq B'\geq 0$. Let $\Ff'$ be the foliation induced by $f': X'\rightarrow Z'$, then $(X',\Ff',B',\Mm;G')/Z'$ is ACSS. By Proposition \ref{prop: weak cbf gfq}, 
$$K_{\Ff'}+B'+\Mm_{X'}\sim_{\mathbb R,Z'}K_{X'}+B'+G'+\Mm_{X'}\sim_{\mathbb R,Z'}E^h+E^v.$$
Thus
$$\kappa_{\sigma}(X'/Z',K_{\Ff'}+B'+\Mm_{X'})=\kappa_{\sigma}(X'/Z',E^h)=0.$$
By Proposition \ref{prop: weak ss num 0 mmp}, we may run a $(K_{\Ff'}+B'+\Mm_{X'})$-MMP$/Z'$ which terminates with a good minimal model $(X'',\Ff'',B'',\Mm)/Z'$. By Lemma \ref{lem: ACSS mmp can run}, $(X'',\Ff'',B'',\Mm;G'')$ is ACSS, where $G''$ is the strict transform of $G'$ on $X''$.

Since $X'\rightarrow U$ factors through $Z'$, $X'\dashrightarrow X''$ is a sequences of steps of a $(K_{\Ff'}+B'+\Mm_{X'})$-MMP$/U$. By Lemma \ref{lem: equivalence over bases}, $K_{\Ff''}+B''+\Mm_{X''}$ is nef$/U$. By Theorem \ref{thm: lc+weak acc=bpstable}, $(X'',B''+G'',\Mm)$ is BP stable$/Z'$. Let $f'': X''\rightarrow Z'$ be the induced contraction and let $\Nn$ be the moduli part of $f'': (X'',B''+G'',\Mm)\rightarrow Z'$. By Proposition \ref{prop: bp stable nef}, $\Nn$ is nef$/U$ and $\Nn$ descends to $X$. By Proposition \ref{prop: weak cbf gfq}, $$K_{X''}+B''+G''+\Mm_{X''}\sim_{\mathbb R,Z'}0.$$
Let $\Mm'$ be the base moduli part of $f'': (X'',B''+G'',\Mm)\rightarrow Z'$, then by the definition of base moduli part, $\Mm'$ descends to $Z'$ and $f''^*\Mm'_{X'}=\Nn_{X'}$ is nef, so $\Mm'_{X'}$ is nef, hence $\Mm'$ is nef.

Let $\tilde B''^h$ be the image of $\tilde B'^h$ on $X''$. Since $K_{X'}+\tilde B'+\Mm_{X'}\sim_{\mathbb R,Z'}0$, $K_{X'}+\tilde B'^h+\Mm_{X'}\sim_{\mathbb R}0$ over the generic point of $Z'$. Thus $K_{X''}+\tilde B''^h+\Mm_{X''}\sim_{\mathbb R}0$ over the generic point of $Z'$. Since $B''\geq\tilde B''^h$ and $K_{X''}+B''+\Mm_{X''}\sim_{\mathbb R,Z'}0$, $B''=\tilde B''^h$ over the generic point of $Z'$. Since $(X',\tilde B'^h,\Mm)$ and $(X'',\tilde B''^h,\Mm)$ are crepant over the generic point of $Z'$, $(X',\tilde B',\Mm)$ and $(X'',B''+G'',\Mm)$ are crepant over the generic point of $Z'$. Thus $(X,B,M)$ and $(X'',B''+G'',\Mm)$ are crepant over the generic point of $Z$. By Lemma \ref{lem: m preserved under crepant}, $\Mm^Z=\Mm'$. The main part of the theorem follows. (1) immediately follows.

(2-4) follow immediately from the definition of the discriminant part. (5) follows from the definition of the discriminant part and Lemma \ref{lem: alg int foliation lct achieved}. By \cite[Theorem 2.23]{JLX22}, if $\Mm$ is NQC$/U$, then $\Mm'$ is NQC$/U$, hence $\Mm^Z$ is NQC$/U$. (6) follows.
\end{proof}

\subsection{Canonical bundle formula for generalized foliated quadruples}\label{subsec: cbf gfq}

\begin{deflem}\label{deflem: cbf gfq}
    Let $(X,\Ff,B,\Mm)/U$ be a sub-gfq and $f: X\rightarrow Z$ a contraction$/U$ such that the general fibers of $f$ are tangent to $\Ff$ and $f: (X,\Ff,B,\Mm)\rightarrow Z$ is an lc-trivial fibration$/U$. We define two $\bb$-divisors $\Bb$ and $\Mm^Z$ on $Z$ in the following way.

 By Lemma \ref{lem: gen fiber tangent mean induce}, there exists a foliation $\Ff_Z$ on $Z$ such that $\Ff=f^{-1}\Ff_Z$. Let $f': (X',\Sigma_{X'},\Mm)\rightarrow (Z',\Sigma_{Z'})$ be any equidimensional model of $f: (X,B,\Mm)\rightarrow Z$ with associated morphisms $h: X'\rightarrow X$ and $h_Z: Z'\rightarrow Z$. Let $\Ff_{Z'}:=h_Z^{-1}\Ff_Z$ and $\Ff':=h^*\Ff$, then $\Ff'=f'^{-1}\Ff_{Z'}$. We define
 $$R':=\sum\left(f'^*D-f'^{-1}\left(D\right)\right),$$
where $D$ runs over all $\Ff_{Z'}$-invariant prime divisors on $Z'$. By \cite[2.9]{Dru17}, we have 
 $$K_{\Ff'/\Ff_{Z'}}=K_{X'/Z'}-R'.$$
Let $K_{\Ff'}+B'+\Mm_{X'}:=h^*(K_{\Ff}+B+\Mm_X)$. Then $K_{\Ff'}+B'+\Mm_{X'}\sim_{\mathbb R,Z'}0$, so 
$$K_{X'}+B'-R'+\Mm_{X'}\sim_{\mathbb R,Z'}0.$$
Since $R'=0$ and $K_{X'}=K_{\Ff'}$ over the generic point of $Z'$, $f': (X',B'-R',\Mm)\rightarrow Z'$ is an lc-trivial fibration. By Theorem \ref{thm: cbf gpair nonnqc}, there exist two $\bb$-divisors $\Bb$ and $\Mm^Z$ on $Z$, such that $\Bb$ is uniquely determined and $\Mm^Z$ is uniquely determined up to $\Rr$-linear equivalence, and the following conditions are satisfied
    \begin{enumerate}
        \item[(i)] $K_{X'}+B'-R'+\Mm_{X'}\sim_{\mathbb R}f'^*(K_{Z'}+\Bb_{Z'}+\Mm^Z_{Z'})$.
        \item[(ii)] $\Mm^Z$ is nef$/U$.
        \item[(iii)] For any birational morphism $g_Z: Z''\rightarrow Z'$ and $g: X''\rightarrow X'$ such that the induced map $f'': X''\dashrightarrow Z''$ is a morphism, we let 
        $$K_{X''}+\tilde B''+\Mm_{X''}:=g^*(K_{X'}+B'-R'+\Mm_{X'}),$$ 
        then $\Bb_{Z''}$ is the discriminant part of $f'': (X'',\tilde B'',\Mm)\rightarrow Z''$.
    \end{enumerate}
 We call $\Bb$ the \emph{discriminant $\bb$-divisor} of $f: (X,\Ff,B,\Mm)\rightarrow Z$ and call $\Mm^Z$ the \emph{base moduli part} of $f: (X,\Ff,B,\Mm)\rightarrow Z$. We also call $\Bb_Z$ the \emph{discriminant part} of $f: (X,\Ff,B,\Mm)\rightarrow Z$. Then
    \begin{enumerate}
      \item $\Bb$ and $\Mm^Z$ are well-defined, i.e., $\Bb$ and the $\Rr$-linear equivalence class of $\Mm^Z$ are independent of the choices of the equidimensional model$/U$ of $f: (X,B,\Mm)\rightarrow Z$. 
 \item $(Z,\Ff_Z,B_Z:=\Bb_Z,\Mm^Z)/U$ is a sub-gfq. 
 \item If $\Mm$ is NQC$/U$, then $\Mm^Z$ is NQC$/U$.
    \end{enumerate}
 We say that $(Z,\Ff_Z,B_Z,\Mm^Z)/U$ is a sub-gfq \emph{induced by a canonical bundle formula$/U$} of $f: (X,\Ff,B,\Mm)\rightarrow Z$.
\end{deflem}

\begin{proof}
By \cite[Definition-Lemma 6.11]{LLM23}, $\Bb$ is independent of the choices of the equidimensional model$/U$ of $f: (X,B,\Mm)\rightarrow Z$. 

Since $K_{\Ff}+B+\Mm_X\sim_{\mathbb R,Z}0$, there exists an $\Rr$-divisor $L$ on $Z$ which is uniquely determined up to $\Rr$-linear equivalence, such that
$$K_{\Ff}+B+\Mm_X\sim_{\mathbb R}f^*L.$$
By condition (i), we have
$$K_{\Ff'}+B'+\Mm_{X'}\sim_{\mathbb R}f'^*\left(K_{\Ff_{Z'}}+\Bb_{Z'}+\Mm^Z_{Z'}\right).$$
Therefore, for any birational morphism $g_Z: Z''\rightarrow Z'$ with $\Ff_{Z''}:=g_Z^{-1}\Ff_{Z'}$, we have
$$\Mm^Z_{Z''}\sim_{\mathbb R}(h_Z\circ g_Z)^*L-K_{\Ff_{Z''}}-\Bb_{Z''}.$$
Thus $\Mm^Z_{Z''}$ is uniquely determined up to the choices of $L$ in its $\Rr$-linear equivalence class. Thus $\Mm^Z$ is uniquely determined up to $\Rr$-linear equivalence. This implies (1).

Moreover, we have
$$L=(h_Z)_*h_Z^*L\sim_{\mathbb R}(h_Z)_*\left(K_{\Ff_{Z'}}+\Bb_{Z'}+\Mm^Z_{Z'}\right)=K_{\Ff_Z}+B_Z+\Mm^Z_Z,$$
so $K_{\Ff_Z}+B_Z+\Mm^Z_Z$ is $\Rr$-Cartier. By our condition (ii), $(Z,\Ff_Z,B_Z:=\Bb_Z,\Mm^Z)/U$ is a sub-gfq. This implies (2).

(3) follows from Theorem \ref{thm: cbf gpair nonnqc}(6).
\end{proof}

\begin{lem}\label{lem: order along generic fiber cbf gfq}
Let $(X,\Ff,B,\Mm)/U$ be a sub-gfq and $f: X\rightarrow Z$ a contraction$/U$ such that $f: (X,\Ff,B,\Mm)\rightarrow Z$ is an lc-trivial fibration. Suppose that $n(K_{\Ff}+B+\Mm_X)\sim 0$ over the generic point of $Z$ for some positive integer $n$ and there is a foliation $\Ff_Z$ on $Z$ such that $\Ff=f^{-1}\Ff_Z$. Then there is a choice $\Mm^Z$ of the base moduli part of $f: (X,\Ff,B,\Mm)\rightarrow Z$, such that 
$$n(K_{\Ff}+B+\Mm_X)\sim nf^*\left(K_{\Ff_Z}+B_Z+\Mm^Z_Z\right),$$
where $B_Z$ is the discriminant part of $f: (X,\Ff,B,\Mm)\rightarrow Z$.
\end{lem}
\begin{proof}
Let $f': (X',\Sigma_{X'},\Mm)\rightarrow (Z',\Sigma_{Z'})$ be a sufficiently high equidimensional model of $f: (X,B,\Mm)\rightarrow Z$ with associated morphisms $h: X'\rightarrow X$ and $h_Z: Z'\rightarrow Z$. Let $\Ff_{Z'}:=h_Z^{-1}\Ff_Z$ and let 
$$R':=\sum_{D\mid D\text{ is an }\Ff_{Z'}\text{-invariant prime divisor}}(f'^*D-f'^{-1}(D)).$$
Then $f': (X',B'-R',\Mm)\rightarrow Z'$ is an lc-trivial fibration.
Since $R'$ is vertical$/Z'$, $n(K_{X'}+B'-R'+\Mm_{X'})\sim 0$ over the generic point of $Z$. The lemma follows from Lemma \ref{lem: order along generic fiber cbf}.
\end{proof}

\begin{lem}\label{lem: m preserved under crepant gfq}
Let $(X,\Ff,B,\Mm)/U$ and $(X',\Ff',B',\Mm')/U$ be two sub-gfqs. Let $f: (X,\Ff,B,\Mm)\rightarrow Z$ and $f': (X',\Ff',B',\Mm')\rightarrow Z'$ be two lc-trivial fibrations$/U$ such that $f$ and $f'$ are birationally equivalent, and $(X,\Ff,B,\Mm)$ and $(X',\Ff',B',\Mm')$ are crepant over the generic point of $Z$. Let $\Mm^Z$ be the base moduli part of $f: (X,\Ff,B,\Mm)\rightarrow Z$ and let $\Mm^{Z'}$ be the base moduli part of $f': (X',\Ff',B',\Mm')\rightarrow Z'$. Then $\Mm^Z\sim_{\mathbb R}\Mm^{Z'}$.
\end{lem}
\begin{proof}
Possibly passing to a common base and resolving the indeterminacy of the induced birational map $X\dashrightarrow X'$, we may assume that $f=f'$, $X=X'$, $Z=Z'$, and $\Ff=\Ff'$ over the generic point of $Z$, $\Mm$ and $\Mm'$ descend to $X$, $f: (X,\Sigma)\rightarrow (Z,\Sigma_Z)$ is equidimensional toroidal for some $\Sigma\supset\Supp B\cup\Supp B'$, and $(Z,\Sigma_Z)$ is log smooth. Let $\Ff_Z$ and $\Ff_Z'$ be two foliations on $Z$ such that $\Ff=f^{-1}\Ff_Z$ and $\Ff'=f'^{-1}\Ff_Z'$,
 $$R:=\sum_{D\mid D\text{ is an }\Ff_{Z}\text{-invariant prime divisor}}(f^*D-f^{-1}(D)),$$
 and
$$R':=\sum_{D\mid D\text{ is an }\Ff'_{Z}\text{-invariant prime divisor}}(f^*D-f^{-1}(D)).$$
Then $\Mm^Z$ and $\Mm^{Z'}$ are the moduli parts of $f: (X,B-R,\Mm)\rightarrow Z$ and $f': (X,B'-R',\Mm)\rightarrow Z$ respectively. Since $(X,B-R,\Mm)$ and $(X',B'-R',\Mm)$ are crepant over the generic point of $Z$, by Lemma \ref{lem: m preserved under crepant}, $\Mm^Z\sim_{\mathbb R}\Mm^{Z'}$.
\end{proof}

\begin{lem}\label{lem: td=bd}
 Let $(X,\Ff,B,\Mm)/U$ be a sub-gfq and $f: X\rightarrow Z$ a contraction$/U$ such that $f: (X,\Ff,B,\Mm)\rightarrow Z$ is an lc-trivial fibration$/U$ with the discriminant $\bb$-divisor $\Bb$. Let $\Ff_Z$ be a foliation on $Z$ such that $\Ff=f^{-1}\Ff_Z$. Then for any prime divisor $D$ on $Z$, 
$$\mult_D\Bb_Z=\epsilon_{\Ff_Z}(D)-\sup\{t\geq 0\mid (X,\Ff,B+tf^*D,\Mm)\text{ is sub-lc over the generic point of } D\}.$$
Moreover, there exists an lc center of $(X,\Ff,B+(\epsilon_{\Ff_Z}(D)-\mult_D\Bb_Z)f^*D,\Mm)$ over the generic point of $D$.
\end{lem}
\begin{proof}
Set $B_Z:=\Bb_Z$. By Definition-Lemma \ref{deflem: cbf gfq}, possibly replacing $f: X\rightarrow Z$ with an equidimensional model of $f: (X,B,\Mm)\rightarrow Z$, we may assume that $X$ is $\Qq$-factorial klt with at most toric quotient singularities, $f$ is equidimensional, $\Mm$ descends to $X$, and there exists a toroidal morphism $f: (X,\Sigma_X,\Mm)\rightarrow (Z,\Sigma_Z)$ such that $\Supp B\subset\Sigma_X$. We may define $R:=\sum(f^*D-f^{-1}(D))$, where $D$ runs over $\Ff_Z$-invariant prime divisors on $Z$. For any prime divisor $D$ on $Z$, we define 
$$t_D:=1-\sup\{t\geq 0\mid (X,B-R+tf^*D,\Mm)\text{ is lc over the generic point of } D\}$$
and
$$b_D:=\epsilon_{\Ff_Z}(D)-\sup\{t\geq 0\mid (X,\Ff,B+tf^*D,\Mm)\text{ is lc over the generic point of } D\}.$$
For any prime divisor $D$ on $Z$, $\mult_DB_Z=t_D$ by definition. Denote by $\eta_D$ the generic point of $D$. There are three cases.

\medskip

\noindent\textbf{Case 1}. $D$ is not $\Ff_Z$-invariant. 

In this case, $R=0$ and $K_{\Ff}=K_X$ over $\eta_D$. It immediately implies that $b_D=t_D$. Moreover, any lc center of $(X,B-R+(1-b_D)f^*D,\Mm)$ over $\eta_D$ is an lc center of $(X,\Ff,B+(1-b_D)f^*D,\Mm)$ over $\eta_D$. Then we may conclude the moreover part.

\medskip

\noindent\textbf{Case 2}. $D$ is $\Ff_Z$-invariant and $D\not\subset\Sigma_Z$. 

Let $B^h$ be the horizontal$/Z$ part of $B$, then $B=B^h$ over $\eta_D$. Since $(X,\Ff,B,\Mm)$ is sub-lc over the generic point of $Z$, $\Sigma_X\geq B^h$. By \cite[Lemma 6.6]{LLM23}, $(X,B^h+f^{-1}(D),\Mm)$ is sub-lc over $\eta_D$. Since
$$(X,B-R+f^*D,\Mm)=(X,B^h+f^{-1}(D),\Mm)$$
over $\eta_D$, $t_D=0$. Thus $\mult_DB_Z=0$. Since $D$ is $\Ff_Z$-invariant, any component of $f^{-1}(D)$ is $\Ff$-invariant and hence is an lc center of $(X,\Ff,B,\Mm)$. In particular, $b_D=0=t_D.$

\medskip

\noindent\textbf{Case 3}. $D$ is $\Ff_Z$-invariant and $D\subset\Sigma_Z$. 

In this case, 
$$-t_D=\sup\{t\mid (X,B+f^{-1}(D)+tf^*D,\Mm)\text{ is sub-lc over } \eta_D\}.$$
Since $f: (X,\Sigma_X,\Mm)\rightarrow (Z,\Sigma_Z)$ is toroidal, there exists a component $S$ of $f^*D$ such that 
$$\mult_S\left(B+f^{-1}(D)-t_Df^*D\right)=1.$$ 
Moreover, as $\lfloor B+f^{-1}(D)+tf^*D\rfloor\leq 0$ over $\eta_D$ for any $t<-t_D$, we can see that 
$$\mult_S(B-t_Df^*D)=0\text{ and }B-t_Df^*D\leq 0\text{ over }\eta_D.$$ 
Note that any component of $f^{-1}(D)$ is $\Ff$-invariant, we have
$$-t_D\geq \sup\{t\geq 0\mid (X,\Ff,B+tf^*D,\Mm)\text{ is sub-lc over } \eta_D\}= -b_D.$$
Suppose that $-t_D>-b_D$. Let $s\in (-b_D,-t_D)$ be a real number, then over $\eta_D$, $(X,B+f^{-1}(D)+sf^*D,\Mm)$ is sub-lc, and $(X,\Ff,B+sf^*D,\Mm)$ is not sub-lc. In particular, there exists a prime divisor $D_X$ over $X$ such that the image of $D_X$ on $Z$ is $D$, and 
$$a(D_X,\Ff,B+sf^*D,\Mm)<-\epsilon_{\Ff}(D_X).$$ 
By Definition-Theorem \ref{defthm: weak ss reduction}, there exists an equidimensional model $f': (X',\Sigma_{X'},\Mm)\rightarrow (Z',\Sigma_{Z'})$ of $f: (X,\Supp B+\Supp f^*D,\Mm)\rightarrow Z$ associated with $h: X'\rightarrow X$ and $h_Z: Z'\rightarrow Z$, such that $D_X$ is on $X'$. Let $\Ff':=h^{-1}\Ff$, $\Ff_{Z'}:=h_Z^{-1}\Ff_Z$, $K_{\Ff'}+B'+\Mm_{X'}:=h^*(K_{\Ff}+B+\Mm_X)$, $D':=(h_Z^{-1})_*D$, and 
$$R':=\sum_{L\mid L\text{ is an }\Ff_{Z'}\text{-invariant prime divisor}}\left(f'^*L-f'^{-1}(L)\right).$$
Note that $D_X$ is a component of $f'^{-1}(D')$ and is $\Ff'$-invariant. As $a(D_X,\Ff,B+sf^*D,\Mm)<-\epsilon_{\Ff}(D_X)=0$, we see that $\mult_{D_X}(B'+sf'^*D')>0$. By Definition-Lemma \ref{deflem: cbf gfq}(1),
\begin{align*}
   -t_D&=\sup\left\{t\geq 0\mid (X',B'-R'+tf'^*D',\Mm)\text{ is lc over the generic point $\eta_{D'}$ of } D'\right\}-1\\
   &=\sup\left\{t\geq 0\mid (X',B'+f'^{-1}(D')+tf'^*D',\Mm)\text{ is lc over } \eta_{D'}\right\}<s<-t_D,
\end{align*}
a contradiction. Therefore $b_D=t_D$. Since $\mult_S(B-t_Df^*D)=0$, $S$ is an lc center of $(X,B-b_Df^*D,\Mm)$ over $\eta_D$. The lemma follows in this case.
\end{proof}

\begin{prop}\label{prop: gfq cbf preserve sing}
Let $(X,\Ff,B,\Mm)/U$ be a sub-gfq and $f: X\rightarrow Z$ a contraction$/U$ such that $f: (X,\Ff,B,\Mm)\rightarrow Z$ is an lc-trivial fibration. Let $\Bb$ be the discriminant $\bb$-divisor of $f: (X,\Ff,B,\Mm)\rightarrow Z$, $B_Z:=\Bb_Z$, and $\Mm^Z$ the base moduli part of $f: (X,\Ff,B,\Mm)\rightarrow Z$. Let $\Ff_Z$ be a foliation on $Z$ such that $\Ff=f^{-1}\Ff_Z$. Then
\begin{enumerate}
    \item If the vertical$/Z$ part of $B$ is $\geq 0$, then $B_Z\geq 0$.
    \item If $(X,\Ff,B,\Mm)$ is sub-lc (resp. lc), then $(Z,\Ff_Z,B_Z,\Mm^Z)$ is sub-lc (resp. lc).
    \item Any lc center of $(Z,\Ff_Z,B_Z,\Mm^Z)$ is the image of an lc center of $(X,\Ff,B,\Mm)$.
    \item The image of any lc center of $(X,\Ff,B,\Mm)$ on $Z$ is an lc center of $(Z,\Ff_Z,B_Z,\Mm^Z)$.
\end{enumerate}
\end{prop}
\begin{proof}
The proposition immediately follows from Lemma \ref{lem: td=bd}.
\end{proof}

Finally, we state the following proposition that can be useful for inductive purposes.

\begin{prop}\label{prop: composition lc trivial fibration}
    Let $(X,\Ff,B,\Mm)$ be a sub-gfq and $X\xrightarrow{f}Y\xrightarrow{g}Z$ two contractions$/U$. Let $h:=g\circ f$. Suppose that $h: (X,\Ff,B,\Mm)\rightarrow Z$ is an lc-trivial fibration. Let $(Z,\Ff_Z,B_Z,\Mm^Z)$ be the sub-gfq induced by $h: (X,\Ff,B,\Mm)\rightarrow Z$. Then
    \begin{enumerate}
        \item $f: (X,\Ff,B,\Mm)\rightarrow Y$ is an lc-trivial fibration.
        \item Let $(Y,\Ff_Y,B_Y,\Mm^Y)$ be a sub-gfq induced by $f: (X,\Ff,B,\Mm)\rightarrow Y$. Then
        \begin{enumerate}
            \item $g: (Y,\Ff_Y,B_Y,\Mm^Y)\rightarrow Z$ is an lc-trivial fibration.
            \item The discriminant part of $g: (Y,\Ff_Y,B_Y,\Mm^Y)\rightarrow Z$ is $B_Z$.
            \item $(Z,\Ff_Z,B_Z,\Mm^Z)$ is a sub-gfq induced by $g: (Y,\Ff_Y,B_Y,\Mm^Y)\rightarrow Z$.
        \end{enumerate}
    \end{enumerate}
\end{prop}
\begin{proof}
Possibly replacing $X$ and $Y$ with higher models, we may assume that $X$ and $Y$ are smooth, and
$\kappa_{\sigma}(X/Z,-B^{\leq 0})=0.$

(1) Since $(X,\Ff,B,\Mm)$ is sub-lc over the generic point of $Z$, $(X,\Ff,B,\Mm)$ is sub-lc over the generic point of $Y$. Since $K_{\Ff}+B+\Mm_X\sim_{\mathbb R,Z}0,$ $K_{\Ff}+B+\Mm_X\sim_{\mathbb R,Y}0$. Since $\kappa_{\sigma}(X/Z,-B^{\leq 0})=0$, $\kappa_{\sigma}(X/Y,-B^{\leq 0})=0$. This implies (1).

(2.a) Since $(X,\Ff,B,\Mm)$ is sub-lc over the generic point of $Z$, by Theorem \ref{thm: cbf gpair nonnqc}, $(Y,\Ff_Y,B_Y,\Mm^Y)$ is sub-lc over the generic point of $Z$. Since 
$$f^*(K_{\Ff_Y}+B_Y+\Mm_Y)\sim_{\mathbb R}K_{\Ff}+B+\Mm_X\sim_{\mathbb R,Z}0,$$
    $K_{\Ff_Y}+B_Y+\Mm_Y\sim_{\mathbb R,Z}0$. By Lemma \ref{lem: td=bd}, for any component $D$ of $B_Y^{\leq 0}$ and any irreducible component $D_X$ of $f^{-1}(D)$ over the generic point of $D$, $D_X$ is a component of $B^{\leq 0}$. Therefore, over the generic point of $Z$, there exists a positive real number $\epsilon$ such that
    $$-B^{\leq 0}\geq \epsilon f^*(-B_Y^{\leq 0}).$$
    Thus
    $$0\leq \kappa_{\sigma}(X/Z,f^*(-B_Y^{\leq 0}))=\kappa_{\sigma}(X/Z,\epsilon f^*(-B_Y^{\leq 0}))\leq \kappa_{\sigma}(X/Z,-B^{\leq 0})=0,$$
    so
    $$\kappa_{\sigma}(Y/Z,-B_Y^{\leq 0})=\kappa_{\sigma}(X/Z,f^*(-B_Y^{\leq 0}))=0.$$
    Therefore, $g: (Y,\Ff_Y,B_Y,\Mm^Y)\rightarrow Z$ is an lc-trivial fibration.

(2.b) Let $B_{Z}'$ be the discriminant part of $g: (Y,\Ff_Y,B_Y,\Mm^Y)\rightarrow Z$. For any prime divisor $D$ over $Z$, let $s_D:=\epsilon_{\Ff_Z}(D)-\mult_DB_Z$ and $s'_D:=\epsilon_{\Ff_Z}(D)-\mult_DB_Z'$. 

By Lemma \ref{lem: td=bd}, for any positive real number $t$ and any prime divisor $D$ on $Z$, $(Y,\Ff_Y,B_Y+tg^*D,\Mm)$ is the sub-gfq induced by $f: (X,\Ff,B+th^*D,\Mm)\rightarrow Y$ over the generic point of $D$. By Proposition \ref{prop: gfq cbf preserve sing}(3)(4),
\begin{align*}
s'_D&=\sup\{t\geq 0\mid (Y,\Ff_Y,B_Y+tg^*D,\Mm^Y)\text{ is sub-lc over the generic point of }D\}\\
&=\sup\{t\geq 0\mid (X,\Ff,B+th^*D,\Mm)\text{ is sub-lc over the generic point of }D\}=s_D.
\end{align*}
Thus $B_Z=B_Z'$. 

(2.c) By applying (2.b) to all high models of $Z$, we get (2.c).
\end{proof}

\section{Canonical bundle formula for lc-trivial morphisms and subadjunction formula}\label{sec: subadj}

\subsection{Canonical bundle formula for lc-trivial morphisms}

\begin{deflem}[{\cite[Proposition 3.4]{Dru21}; cf. \cite[Proposition 3.7]{Spi20}}]\label{deflem: hurwitz foliation}
    Let $f: X'\rightarrow X$ be a surjective finite morphism between normal varieties and $\Ff$ a foliation on $X$. Assume that $K_{\Ff}$ is $\Qq$-Cartier and $\Ff':=f^{-1}\Ff$. For any prime divisor on $X$, we let $r_D$ be the ramification index of $f$ along $D$. We call
    $$R:=\sum_{D\mid D\text{ is a non-}\Ff\text{-invariant prime divisor}}(r_D-1)D$$
    the \emph{ramification divisor} of $f$ with respect to $\Ff$. Then we have
    $$K_{\Ff'}=f^*K_{\Ff}+R.$$
\end{deflem}

\begin{deflem}\label{deflem: cbf finite}
Let $(X,\Ff,B,\Mm)/U$ be a sub-gfq and $f: X\rightarrow Z$ a finite morphism $/U$. Suppose that there exists a foliation $\Ff_Z$ on $Z$ such that $\Ff=f^{-1}\Ff_Z$, and suppose that $K_{\Ff}+B+\Mm_X\sim_{\mathbb R,Z}0$.

We define two $\bb$-divisors, $\Bb$ and $\Mm^Z$ on $Z$, in the following way. Let $h_Z: Z'\rightarrow Z$ be any birational morphism, $X'$ the main component of $Z'\times_{Z}X$, $f': X'\rightarrow Z'$ and $h: X'\rightarrow X$ the induced morphisms, $\Ff':=h^{-1}\Ff$, and $\Ff_{Z'}:=h_Z^{-1}\Ff_Z$. We let $$K_{\Ff'}+B'+\Mm_{X'}:=h^*(K_{\Ff}+B+\Mm_X).$$
Let $Z'^0$ be the largest open subset of $Z'$ which does not contain $\Sing(\Ff_{Z'})\cup\Sing(Z')$ and let $X'^0:=f'^{-1}(Z'^0)$. By Definition-Lemma \ref{deflem: hurwitz foliation}, 
$$K_{\Ff'|_{X'^0}}=(f'|_{X'^0})^*K_{\Ff_{Z'}|_{Z'^0}}+R'^0$$
where $R'^0$ is the ramification divisor of $f'|_{X'^0}$ with respect to $\Ff_{Z'}|_{Z'^0}$. We let $R'$ be the closure of $R'^0$ in $X'$. We let $\Bb$ and $\Mm^Z$ be the $\bb$-divisors such that $\Bb_{Z'}=\frac{1}{\deg f}f'_*(R'+B')$ and $\Mm^Z_{Z'}=\frac{1}{\deg f}f'_*\Mm_{X'}$ for any choice of $Z'$. Then
\begin{enumerate}
   \item $\Bb$ and $\Mm^Z$ are well-defined and uniquely determined.
    \item For any choice of $Z'$,
    $$K_{\Ff'}+B'+\Mm_{X'}\sim_{\mathbb R}f'^*(K_{\Ff_{Z'}}+\Bb_{Z'}+\Mm^Z_{Z'}).$$
    \item $\Mm^Z$ is nef $/U$.
    \item If $B\geq 0$, then $\Bb_Z\geq 0$.
    \item If $(X,\Ff,B,\Mm)$ is (sub-)lc, then $(Z,\Ff_Z,\Bb_Z,\Mm^Z)$ is (sub-)lc, and for any lc center $T$ of $(Z,\Ff_Z,B_Z,\Mm^Z)$, any component of $f^{-1}(T)$ is an lc center of $(X,\Ff,B,\Mm)$.
    \item If $\Mm$ is NQC $/U$, then $\Mm^Z$ is NQC $/U$.
\end{enumerate}
We call $\Bb$ the \emph{discriminant $\bb$-divisor} of $f: (X,B,\Mm)\rightarrow Z$, and call $B_Z:=\Bb_Z$ the \emph{discriminant part} of $f: (X,B,\Mm)\rightarrow Z$. We call $\Mm^Z$ the \emph{base moduli part} of $f: (X,B,\Mm)\rightarrow Z$. We say that $(Z,\Ff_Z,B_Z,\Mm^Z)/U$ is the sub-gfq \emph{induced} by $f: (X,\Ff,B,\Mm)\rightarrow Z$.
\end{deflem}
\begin{proof}
(1) We only need to show that for any birational morphism $g_Z: Z''\rightarrow Z'$, $(g_Z)_*\Bb_{Z''}=\Bb_{Z'}$ and $(g_Z)_*\Mm^Z_{Z''}=\Mm^Z_{Z'}$. We let $X''$ be the main component of $X'\times_{Z'}Z''$ and $g: X''\rightarrow X'$, $f'': X''\rightarrow Z''$ the induced morphisms. Let $\Ff'':=g^{-1}\Ff', \Ff_{Z''}:=g^{-1}_Z\Ff_{Z'}$, $Z''^0$ be the largest open subset of $Z''$ which does not contain $\Sing(\Ff_{Z''})\cup\Sing(Z'')$, $X''^0:=f''^{-1}(Z''^0)$, $R''^0$ the ramification divisor of $f''|_{X''^0}$ with respect to $\Ff_{Z''}|_{Z''^0}$, and $R''$ the closure of $R''^0$ in $X''$. Then
$$\Bb_{Z'}=\frac{1}{\deg f}f'_*(B'+R')=\frac{1}{\deg f}f'_*g_*(B''+R'')=\frac{1}{\deg f}(g_Z)_*f''_*(B''+R'')=(g_Z)_*\Bb_{Z''}$$
and
$$\Mm^Z_{Z'}=\frac{1}{\deg f}f'_*\Mm_{X'}=\frac{1}{\deg f}f'_*g_*\Mm_{X''}=\frac{1}{\deg f}(g_Z)_*f''_*\Mm_{X''}=(g_Z)_*\Mm^Z_{Z''}.$$

(2) By (1), we only need to prove (2) for any sufficiently high model $Z'$ of $Z$. In particular, we may assume that $Z'$ is $\Qq$-factorial. Then $f'^*(\frac{1}{\deg f}f'_*R')=R'$, $f'^*(\frac{1}{\deg f}f'_*B')=B'$, and $f'^*(\frac{1}{\deg f}f'_*\Mm_{X'})=\Mm_{X'}$, so (2) immediately follows.

(3)(6) By \cite[Lemma 4.2]{HL21b}, there exists a birational morphism $h_Z: Z''\rightarrow Z$ satisfying the following. Let $X''$ be the main component of $Z''\times_Z X$, then $\Mm$ descends to $X''$. By definition, $\Mm^Z$ descends to $Z''$. Since $\Mm_{X''}$ is nef, $\Mm^Z_{Z''}$ is nef. Thus $\Mm^Z$ is nef. This implies (3). Moreover, if $\Mm$ is NQC $/U$, then $\Mm_{X''}$ is NQC $/U$, so $\Mm^Z_{Z''}$ is NQC $/U$, hence $\Mm^Z$ is NQC $/U$. This implies (6).

(4) It is obvious from the definition.

(5) By (4) we only need to prove the sub-lc case. Suppose that $(X,\Ff,B,\Mm)$ is sub-lc, then $(X',\Ff',B',\Mm)$ is sub-lc. Let $D$ be a prime divisor on $Z'$. Let $E_1,\dots,E_m$ be all components of $f'^{-1}(D)$ and let $r_i$ be the ramification index along $E_i$. 

If $D$ is $\Ff_{Z'}$-invariant, then each $E_i$ is $\Ff'$-invariant, and $E_i\not\subset\Supp R'$. Since $(X',\Ff',B',\Mm)$ is sub-lc, $\mult_{E_i}B'\leq 0$ for any $i$. Thus
$$\mult_D\Bb_{Z'}=\mult_D\frac{1}{\deg f}f'_*(B'+R')=\sum_{i=1}^m\frac{1}{\deg f}(\mult_{E_i}B')\leq 0=\epsilon_{\Ff_{Z'}}(D).$$
Moreover, if $D$ is an lc place of $(Z,\Ff_Z,B_Z,\Mm^Z)$, then $\mult_D\Bb_{Z'}=0$, so $\mult_{E_i}B'=0$ for each $i$. Therefore, each $E_i$ is an lc place of $(X',\Ff',B',\Mm)$, hence an lc place of $(X,\Ff,B,\Mm)$.

If $D$ is not $\Ff_Z$-invariant, then each $E_i$ is not $\Ff'$-invariant, and $\sum_{i=1}^m r_i \leq \deg f$. Since $(X',\Ff',B',\Mm)$ is sub-lc, $\mult_{E_i}B'\leq 1$ for any $i$. Thus
$$\mult_D\Bb_{Z'}=\mult_D\frac{1}{\deg f}f'_*(B'+R')=\sum_{i=1}^m\frac{1}{\deg f}(r_i-1+\mult_{E_i}B')\leq\frac{\sum_{i=1}^m r_i}{\deg f}\leq 1=\epsilon_{\Ff_{Z'}}(D).$$
Moreover, if $D$ is an lc place of $(Z,\Ff_Z,B_Z,\Mm^Z)$, then $\mult_D\Bb_{Z'}=1$, so $\mult_{E_i}B'=1$ for each $i$. Therefore, each $E_i$ is an lc place of $(X',\Ff',B',\Mm)$, hence an lc place of $(X,\Ff,B,\Mm)$.

Since $h_Z: Z'\rightarrow Z$ can be any birational morphism, we get (5).
\end{proof}

\begin{defn}[lc-trivial morphism]\label{defn: lc trivial morphism}
Let $(X,\Ff,B,\Mm)/U$ be a sub-gfq and $f: X\rightarrow Z$ a projective surjective morphism over $U$. Let $X\xrightarrow{\tau}\tilde Z\xrightarrow{\gamma}Z$ be the Stein factorization of $f$. We say that $f: (X,\Ff,B,\Mm)\rightarrow Z$ is an \emph{lc-trivial morphism} if
\begin{enumerate}
\item $K_{\Ff}+B+\Mm_X\sim_{\Rr,Z}0$,
\item $\tau: (X,\Ff,B,\Mm)\rightarrow\tilde Z$ is an lc-trivial fibration, and
\item there exists a foliation $\Ff_{Z}$ on $Z$ such that $\Ff=f^{-1}\Ff_Z$. 
\end{enumerate}
\end{defn}

\begin{defthm}[Canonical bundle formula for lc-trivial morphisms]\label{defthm: cbf lctrivial morphism}
    Let $$(X,\Ff,B,\Mm)/U$$ be a sub-gfq and $f: X\rightarrow Z$ an lc-trivial morphism $/U$, and $\Ff_Z$ a foliation on $Z$ such that $\Ff=f^{-1}\Ff_Z$. Then there is a sub-gfq $(Z,\Ff_Z,B_Z,\Mm^Z)/U$, such that $B_Z$ is uniquely determined and $\Mm^Z$ is determined up to $\Rr$-linear equivalence, defined in the following way.

    Let $X\xrightarrow{\tau}\tilde Z\xrightarrow{\gamma}Z$ be the Stein factorization of $f$. By Definition-Lemma \ref{deflem: cbf gfq}, there exists a sub-gfq $$(\tilde Z,\Ff_{\tilde Z},B_{\tilde Z},\tilde\Mm^Z)/U$$
    induced by $\tau: (X,\Ff,B,\Mm)\rightarrow\tilde Z$, such that $B_{\tilde Z}$ is uniquely determined, and $\tilde\Mm^Z$ is uniquely determined up to $\Rr$-linear equivalence. Moreover, we have $\Ff_{\tilde Z}=\gamma^{-1}\Ff_Z$ and
    $$K_{\Ff_{\tilde Z}}+B_{\tilde Z}+\tilde\Mm^Z_{\tilde Z}\sim_{\mathbb R,Z}0.$$
    By Definition-Lemma \ref{deflem: cbf finite}, there exists a sub-gfq 
  $$(Z,\Ff_Z,B_Z,\Mm^Z)/U$$
  induced by $\gamma: (\tilde Z,\Ff_{\tilde Z},B_{\tilde Z},\tilde\Mm^Z)\rightarrow Z$, such that $B_Z$ is uniquely determined, and $\Mm^Z$ is uniquely determined up to $\Rr$-linear equivalence. We say that $B_Z$ is the \emph{discriminant part} of $f: (X,\Ff,B,\Mm)\rightarrow Z$, $\Mm^Z$ is the \emph{base moduli part} of $f: (X,\Ff,B,\Mm)\rightarrow Z$, and say that $(Z,\Ff_Z,B_Z,\Mm^Z)$ is a sub-gfq induced by $f: (X,\Ff,B,\Mm)\rightarrow Z$.

 Moreover, we have the following
\begin{enumerate}
\item If the vertical $/Z$ part of $B$ is $\geq 0$, then $B_Z\geq 0$.
\item If $(X,\Ff,B,\Mm)$ is (sub-)lc, then $(Z,\Ff_Z,B_Z,\Mm^Z)$ is (sub-)lc.
\item $B_Z$ is uniquely determined, and $\Mm^Z$ is uniquely determined up to $\Rr$-linear equivalence.
\item Suppose that $(X,\Ff,B,\Mm)$ is sub-lc. Then for any lc center $T$ of $(Z,\Ff_Z,B_Z,\Mm^Z)$, $T$ is the image of an lc center of $(X,\Ff,B,\Mm)$ on $Z$.
\end{enumerate}
\end{defthm}
\begin{proof}
(1) It follows from Definition-Lemma \ref{deflem: cbf finite}(4) and Proposition \ref{prop: gfq cbf preserve sing}(1).

(2) It follows from Definition-Lemma \ref{deflem: cbf finite}(5) and Proposition \ref{prop: gfq cbf preserve sing}(2).

(3) It follows from Definition-Lemma \ref{deflem: cbf gfq}(1) and Definition-Lemma \ref{deflem: cbf finite}(1).

(4) It follows from Proposition \ref{prop: gfq cbf preserve sing}(3) and Definition-Lemma \ref{deflem: cbf finite}(5).
\end{proof}

\subsection{Subadjunction formula for g-pairs}

In this section, we shall introduce and discuss the subadjunction formula for lc g-pairs. Since the canonical bundle formula for lc-trivial fibrations for gfqs requires that the general fibers are tangent to the foliation, the subadjunction formula for foliations is more subtle and we will leave it to a later part of the paper.

\begin{defthm}[Subadjunction formula via log resolutions]\label{defthm: subadjun}
    Let $(X,B,\Mm)/U$ be a g-sub-pair and $V$ an lc center of $(X,B,\Mm)$ with normalization $\nu: W\rightarrow V$, such that $B\geq 0$ near the generic point of $V$. Then there exists a naturally defined g-sub-pair $(W,B_W,\Mm^W)/U$ defined in the following way. 
    
    Let $S$ be an lc place of $(X,B,\Mm)$ so that $\Center_X S=V$. Let $h: Y\rightarrow X$ be a log resolution of $(X,\Supp B)$ such that $\Mm$ descends to $Y$ and $S$ is on $Y$. We let 
    $$K_Y+B_Y+\Mm_Y:=h^*(K_X+B+M_X)$$
    and let $(S,B_S,\Mm^S)/U$ be the g-sub-pair induced by the adjunction
    $$K_S+B_S+\Mm^S_S:=(K_Y+B_Y+\Mm_Y)|_S.$$
    Then there exists an induced projective surjective morphism $h_S: S\rightarrow W$ such that $\nu\circ f_S=h|_S$. By construction, we have
    $$K_S+B_S+\Mm^S_S\sim_{\mathbb R,W}0.$$
    Since $B\geq 0$ near the generic point of $V$, $B_W\geq 0$ near the generic point of $S$. Therefore, $h_S: (S,B_S,\Mm^S)\rightarrow W$ is an lc-trivial morphism. By Definition-Theorem \ref{defthm: cbf lctrivial morphism}, there exists a g-sub-pair $(W,B_W,\Mm^W)/U$ induced by $h_S: (S,B_S,\Mm^S)\rightarrow W$. Moreover, we have the following
    \begin{enumerate}
        \item For a fixed choice of $S$, $B_W$ is uniquely determined, and $\Mm^W$ is uniquely determined up to $\Rr$-linear equivalence. In particular, $B_W$ and the $\Rr$-linear equivalence class of $\Mm^W$ are independent of the choice of $h$.
        \item $K_W+B_W+\Mm_W\sim_{\mathbb R}(K_X+B+\Mm_X)|_{W}$.
        \item If $(X,B,M)$ is sub-lc near $V$, then $(W,B_W,\Mm^W)$ is sub-lc.
         \item Suppose that $(X,B,M)$ is sub-lc near $V$. Then for any lc center $T$ of $(W,B_W,\Mm^W)$, $\nu(T)$ is an lc center of $(X,B,\Mm)$.
    \end{enumerate}
    We say that $(W,B_W,\Mm^W)/U$ is a g-sub-pair induced by subadjunction
    $$K_W+B_W+\Mm^W_W:=(K_X+B+\Mm_X)|_W$$
    and say that $(W,B_W,\Mm^W)$ is associated with $S$.
\end{defthm}
\begin{proof}
    The construction is clear so we only need to prove (1-5). 

    (1) We let $h': Y'\rightarrow X$ be a log resolution of $(X,\Supp B)$ such that $\Mm$ descends to $Y'$ and $S$ is on $Y'$, so that the induced birational map $g: Y'\rightarrow Y$ is a morphism. Let $S':=g^{-1}_* S$, 
    $$K_{Y'}+B_{Y'}+\Mm_{Y'}:=h'^*(K_X+B+\Mm_X)$$
    and let $(S',B_{S'},\Mm^S)/U$ be the g-sub-pair induced by the adjunction
    $$K_{S'}+B_{S'}+\Mm^S_{S'}:=(K_{Y'}+B_{Y'}+\Mm_{Y'})|_{S'}.$$
    Then $g|_{S'}: S'\rightarrow S$ is a morphism, and we have
    \begin{align*}
       K_{S'}+B_{S'}+\Mm^S_{S'}&=(K_{Y'}+B_{Y'}+\Mm_{Y'})|_{S'}=g^*(K_Y+B_Y+\Mm_Y)|_{S'}\\
       &=g|_{S'}^*((K_Y+B_Y+\Mm_Y)|_S)=g|_{S'}^*(K_S+B_S+\Mm^S_S).
    \end{align*}
    By our construction, the g-sub-pair induced by $h_S\circ g|_{S'}: (S',B_{S'},\Mm^S)\rightarrow W$ is equal to the g-sub-pair induced by $h_S: (S,B_{S},\Mm^S)\rightarrow W$ modulo $\Rr$-linear equivalence of the base moduli part. Since $h'$ can be any high log resolution of $(X,\Supp B)$, (1) follows.

    (2) It immediately follows from the definition. 

    (3) Since $(X,B,\Mm)$ is sub-lc near $V$, $(S,B_S,\Mm^S)$ is sub-lc. By Definition-Theorem \ref{defthm: cbf lctrivial morphism}(2), we get (3).

    (4) By Definition-Theorem \ref{defthm: cbf lctrivial morphism}, $T$ is the image of an lc center $T_S$ of $(S,B_S,\Mm^S)$ on $W$. Since $(S,B_S,\Mm^S)$ is induced by adjunction from a log smooth g-sub-pair, $T_S$ is also an lc center of $(Y,B_Y,\Mm)$. Thus $h(T_S)$ is an lc center of $(X,B,\Mm)$. By construction, $\nu(T)=h(T_S)$.
\end{proof}

\begin{prop}[Subadjunction formula via dlt models]\label{prop: lc subadj is lc}
       Let $(X,B,\Mm)/U$ be a g-sub-pair and $V$ an lc center of $(X,B,\Mm)$ with normalization $\nu: W\rightarrow V$, such that $(X,B,\Mm)$ is lc near $W$. Let $S$ be an lc place of $(X,B,\Mm)$ such that $\Center_X S=V$. Let $(W,B_W,\Mm^W)/U$ be a g-sub-pair induced by subadjunction
    $$K_W+B_W+\Mm^W_W:=(K_X+B+\Mm_X)|_W$$
    and associated with $S$.

    Suppose that $f: Y\rightarrow X$ is a dlt modification of $(X,B,\Mm)$ near $W$ such that $S$ is on $Y$. Let 
    $$K_Y+B_Y+\Mm_Y:=f^*(K_X+B+\Mm_X),$$
     $(S,B_{S},\Mm^S)/U$ the g-sub-pair induced by the adjunction
    $$K_{S}+B_{S}+\Mm^S_{S}:=(K_Y+B_Y+\Mm_Y)|_{S},$$
    and $f_{S}: S\rightarrow W$ the induced projective surjective morphism such that $\nu\circ f_S=f|_{S}$. Then
    \begin{enumerate}
        \item $(W,B_W,\Mm^W)$ is the g-pair induced by $f_S: (S,B_S,\Mm^S)\rightarrow W$.
        \item $(W,B_W,\Mm^W)$ is lc.
    \end{enumerate}
\end{prop}
\begin{proof}
Let $g: Y'\rightarrow Y$ be a log resolution of $(Y,\Supp B_Y)$ such that $\Mm$ descends to $Y'$, 
    $$K_{Y'}+B_{Y'}+\Mm_{Y'}:=g^*(K_Y+B_Y+\Mm_Y),$$
    $S':=g^{-1}_* S$, and let $(S',B_{S'},\Mm^S)/U$ be the g-sub-pair induced by the adjunction
    $$K_{S'}+B_{S'}+\Mm^S_{S'}:=(K_{Y'}+B_{Y'}+\Mm_{Y'})|_{S'}.$$
    Then $g|_{S'}: S'\rightarrow S$ is a morphism, and we have
    \begin{align*}
       K_{S'}+B_{S'}+\Mm^S_{S'}&=(K_{Y'}+B_{Y'}+\Mm_{Y'})|_{S'}=g^*(K_Y+B_Y+\Mm_Y)|_{S'}\\
       &=g|_{S'}^*((K_Y+B_Y+\Mm_Y)|_{S})=g|_{S'}^*(K_{S}+B_{S}+\Mm^S_{S}).
    \end{align*}
    By our construction, $(W,B_W,\Mm^W)/U$ is the g-sub-pair induced by $f_S\circ g|_{S'}: (S',B_{S'},\Mm^S)\rightarrow W$, which is equal to the g-sub-pair induced by $f_S: (S,B_{S},\Mm^S)\rightarrow W$ modulo $\Rr$-linear equivalence of the base moduli part. By Definition-Theorem \ref{defthm: cbf lctrivial morphism}(2), $(W,B_W,\Mm^W)$ is lc.
\end{proof}

\begin{defthm}\label{thm: spring and source for glc crepant log structure}
Let $(X,B,\Mm)/U$ be a dlt g-pair and $f: (X,B,\Mm)\rightarrow Y$ a dlt crepant log structure $/U$. Let $Z\subset Y$ be an lc center of $f: (X,B,\Mm)\rightarrow Y$ with normalization $\nu: Z^n\rightarrow Z$. Let $\mathcal{S}$ be the set of all lc centers of $(X,B,\Mm)$ which dominate $Z$ and let $S\in\mathcal{S}$ be an element that is minimal under inclusion. Let $(S,B_S,\Mm^S)$ be the g-pair induced by adjunction
$$K_S+B_S+\Mm^S_S:=(K_X+B+\Mm_X)|_S,$$
$f_S: S\rightarrow Z^n$ the induced morphism such that $\nu\circ f_S=f|_S$, and let $f^n_S: S\xrightarrow{\tau} V\xrightarrow{\gamma} Z^n$ be the Stein factorization of $f|_S: S\rightarrow Z^n$.
Then
\begin{enumerate}
\item[(1)] (Crepant log structure) $(S,B_S,\Mm^S)$ is dlt, $K_S+B_S+\Mm^S_S\sim_{\Rr,Z^n}0$, and $(S,B_S,\Mm^S)$ is klt over the generic point of $Z^n$. In particular, $f|_S: (S,B_S,\Mm^S)\rightarrow Z^n$ is a dlt crepant log structure and an lc-trivial morphism.
\end{enumerate}
We let
$$(V,B_{V},\Mm^V)/U$$
be the g-pair induced by the lc-trivial fibration $\tau: (S,B_S,\Mm^S)\rightarrow V$. Then 
\begin{enumerate}
    \item[(2)] (Uniqueness of sources) The crepant birational equivalence class of $(S,B_S,\Mm^S)$ does not depend on the choice of $S$. We call the crepant birational equivalence class of $(S,B_S,\Mm^S)$ the \emph{source} of $Z$ with respect to $f: (X,B,\Mm)\rightarrow Y$, and it is denoted by $\Src(Z,X,B,\Mm)$.
    \item[(3)] (Uniqueness of springs) $(V,B_V,\Mm^V)$ modulo the $\Rr$-linear equivalence class of $\Mm^V$ is unique up to isomorphism. We call $(V,B_V,\Mm^V)$ the \emph{spring} of $Z$ with respect to $f: (X,B,\Mm)\rightarrow Y$, and it is denoted by $\Spr(Z,X,B,\Mm)$.
    \item[(4)] (Adjunction) Let $W\subset X$ be an lc center such that $Z\subset Y_W:=f(W)$, and let $(W,B_W,\Mm^W)/U$ be the lc g-pair induced by repeatedly applying adjunction
    $$K_W+B_W+\Mm^W_W:=(K_X+B+\Mm_X)|_W.$$  
    Let $\nu_Y: Y_W^n\rightarrow Y_W$ be the normalization of $Y_W$, $f_W: W\rightarrow Y_W^n$ the induced morphism such that $\nu_Y\circ f_W=f|_W$, and let
    $$W\xrightarrow{\tau_W} V_W\xrightarrow{\gamma_W}Y_W^n$$
    be the Stein factorization of $f_W$. Let $Z_W\subset V_W$ be an irreducible subvariety such that $(\nu_Y\circ\gamma_W)(Z_W)=Z^n$, and $(V_W,B_{V_W},\Mm^{V_W})/U$ a g-pair induced by the lc-trivial fibration $\tau_W: (W,B_W,\Mm^W)\rightarrow V_W$. Then
    \begin{enumerate}
    \item $Z_W$ is an lc center of $(V_W,B_{V_W},\Mm^{V_W})$.
    \item $\Src(Z,X,B,\Mm)=\Src(Z_W,W,B_W,\Mm^W)$.
    \item $\Spr(Z,X,B,\Mm)=\Spr(Z_W,W,B_W,\Mm^W)$.
    \end{enumerate}
\end{enumerate}
\end{defthm}
\begin{proof}
(1) By \cite[Lemma 2.9]{HL22}, $(S,B_S,\Mm^S)$ is dlt. Since $K_X+B+\Mm_X\sim_{\mathbb R,Z}0$, $K_S+B_S+\Mm^S_S\sim_{\mathbb R,Z}0$. By Lemma \ref{lem: inversion of adjunction gdlt} and since $S$ is minimal in $\mathcal{S}$, $(S,B_S,\Mm^S)$ is klt over the generic point of $Z^n$. (1) follows.

(2) By Theorem \ref{thm: P1 link for gdlt crepant log structure}, different choices of $S$ are $\mathbb P^1$-linked to each other, hence they are crepantly equivalent to each other by Definition \ref{defn: p1 link}(3). 

(3) It follows from (2) and Definition \ref{defn: cbf gpair}.

(4) By Lemma \ref{lem: gdlt crepant log structure is compatible under subadjunction}(3) and Theorem \ref{thm: cbf gpair nonnqc}, $Z_W$ is an lc center of $(V_W,B_{V_W},\Mm^{V_W})$ and an lc center of $\tau_W: (W,B_W,\Mm^W)\rightarrow V_W$. This implies (4.a).

Let $S'$ be a minimal lc center of $(W,B_W,\Mm^W)$ which dominates $Z_W$, then $S'$ is also an lc center of $(X,B,\Mm)$ which dominates $Z_W$. In particular, $S'$ dominates $Z$. If $S'$ is not minimal in $\mathcal{S}$, then there exists $S''\subsetneq S'$ such that $S''$ dominates $Z$, so $\tau_W(S'')\subset Z_W$ and $\tau_W(S'')$ dominates $Z$. This is not possible as $Z_W$ is irreducible and $\gamma_W$ is finite. Therefore, $S'$ is minimal in $\mathcal{S}$. This implies (4.b). (4.c) follows from (4.b) and (3).
\end{proof}

\begin{deflem}[Subadjunction formula via minimal lc centers]\label{deflem: subadj minimal lc center}
    Let $(X,B,\Mm)/U$ be a g-sub-pair and $V$ an lc center of $(X,B,\Mm)$ with normalization $\nu: W\rightarrow V$, such that $(X,B,\Mm)$ is lc near $W$. 

    Suppose that $f: Y\rightarrow X$ is a dlt modification of $(X,B,\Mm)$ near $W$ and let
    $$K_Y+B_Y+\Mm_Y:=f^*(K_X+B+\Mm_X).$$
    Let $\mathcal{S}$ be the set of all lc centers of $(Y,B_Y,\Mm)$ whose image on $X$ is $V$, and let $S$ be a minimal element of $\mathcal{S}$ up to inclusion. Let $(S,B_S,\Mm^S)/U$ be the g-pair induced by repeatedly applying adjunction
    $$K_S+B_S+\Mm^S_S:=(K_Y+B_Y+\Mm_Y)|_S,$$
    and let $f_S: S\rightarrow W$ be the induced projective surjective morphism such that $\nu\circ f_S:=f|_S$. 

    We let $(W,B_W,\Mm^W)/U$ be a g-pair induced by a canonical bundle formula of $f_S: (S,B_S,\Mm^S)\rightarrow W$. Then
    \begin{enumerate}
        \item There exists an lc place $S'$ of $(X,B,\Mm)$ such that $\Center_X S'=V$ and $(W,B_W,\Mm^W)$ is a g-pair induced by subadjunction
    $$K_W+B_W+\Mm^W_W:=(K_X+B+\Mm_X)|_W$$
    and $(W,B_W,\Mm^W)$ is associated with $S'$. 
        \item $K_W+B_W+\Mm_W\sim_{\mathbb R}(K_X+B+\Mm_X)|_W$.
            \item $(W,B_W,\Mm^W)$ is lc.
            \item For any lc center $T$ of $(W,B_W,\Mm^W)$, $\nu(T)$ is an lc center of $(X,B,\Mm)$.
            \item $(W,B_W,\Mm^W)$ does not depend on the choice of $S$ (but may depend on the choice of $f$).
    \end{enumerate}
    We say that $(W,B_W,\Mm^W)/U$ is associated to $f$.
\end{deflem}
\begin{proof}
(1) We let $g: Y'\rightarrow Y$ be the blow-up of the generic point of $S$ and let $S'$ be the reduced exceptional divisor. Let 
$$K_{Y'}+B_{Y'}+\Mm_{Y'}=g^*(K_Y+B_Y+\Mm_Y).$$
Then $(Y',B_{Y'},\Mm)$ is dlt over a neighborhood of $W$. Let $(S',B_{S'},\Mm^{S'})/U$ be the g-pair induced by adjunction
$$K_{S'}+B_{S'}+\Mm^{S'}_{S'}:=(K_{Y'}+B_{Y'}+\Mm_{Y'})|_{S'}.$$
Since $(Y,B_Y)$ is log smooth near the generic point of $S$ and $\Mm$ descends to $Y$ near the generic point of $S$, $g|_{S'}: S'\rightarrow S$ is a contraction, and $(S,B_S,\Mm^S)$ is induced by $g|_{S'}: (S',B_{S'},\Mm^{S'})\rightarrow S$.

Thus the Stein factorization of the induced morphism $S'\rightarrow W$ factors through $S$. By Proposition \ref{prop: composition lc trivial fibration}, we get (1).

(2) It follows from (1) and Definition-Theorem \ref{defthm: subadjun}(2).

(3) It follows from (1) and Proposition \ref{prop: lc subadj is lc}(2).

(4) It follows from (1) and Definition-Theorem \ref{defthm: subadjun}(4).

(5) It follows from Definition-Theorem \ref{thm: spring and source for glc crepant log structure}.
\end{proof}

\section{Stratification of generalized pairs and Du Bois property}\label{sec: du bois}

The goal of this section is to study the stratification properties of lc generalized pairs and prove Theorem \ref{thm: glc sings are Du Bois}.

\subsection{Stratification}

In this subsection we recall some basic definitions of stratifications.

\begin{defn}[{\cite[Definition 9.15]{Kol13}}] 
Let $X$ be a scheme. A {\it stratification} of $X$ is a decomposition of $X$ into a finite disjoint union of reduced locally closed subschemes. We will consider stratifications where the strata are of pure dimension
and are indexed by their dimensions. We write $X=\cup_{i}S_iX$ where $S_iX\subset X$ is the $i$-dimensional stratum. Such a stratified scheme is denoted by $(X,S_*)$. We also
assume that $\cup_{i\le j}S_iX$ is closed for every $j$. The {\it boundary} of $(X,S_*)$ is the closed subscheme
$$
B(X,S_*):=\cup_{i<\dim X}S_iX=X\backslash S_{\dim X}X,
$$
and is denoted by $B(X)$ if the stratification $S_*$ is clear. 

Let $(X, S_*)$ and $(Y, S_*)$ be stratified schemes. We say that $f: X\to Y$ is a {\it stratified morphism} if $f(S_iX)\subset S_iY$ for every $i$. Since the strata $S_iX$ are disjoint from each other, $f: X\to Y$ is a stratified morphism if and only if $S_iX=f^{-1}(S_iY)$.

Let $(Y, S_*)$ be a stratified scheme and $f:X\to Y$ a quasi-finite morphism such that $f^{-1} (S_iY)$ has pure dimension $i$ for every $i$. Then $S_iX:=f^{-1}(S_iY)$ defines a stratification of $X$. We denote it by $(X,f^{-1}S_*)$, and we say that $f:X\to(Y,S_*)$ is \emph{stratifiable}.
\end{defn}

\begin{defn}[{\cite[Definition 9.16]{Kol13}}]
Let $(X, S_*)$ be a stratified variety. A relation $(\sigma_1,\sigma_2): R\rightrightarrows (X,S_*)$ is {\it stratified} if each $\sigma_i$ is stratifiable and $\sigma_1^{-1}S_*=\sigma_2^{-1}S_*$. Equivalently,
there exists a stratification $(R,\sigma^{-1}S_i)$, such that $r\in\sigma^{-1}S_iR$ if and only if $\sigma_1(r)\in S_iX$ and if and only if $\sigma_2(r)\in S_iX$.
\end{defn}

\begin{defn}[{\cite[Definition 9.18]{Kol13}}]
Let $(X,S_*)$ be a stratified scheme such that $X$ is an excellent scheme. The normality conditions (N), (SN), (HN), and (HSN) are defined in the following ways.
\begin{enumerate}
    \item[(N)] We say that $(X,S_*)$ has {\it normal strata}, or that it satisfies (N), if each $S_iX$ is normal.
    \item[(SN)] We say that $(X,S_*)$ has {\it semi-normal boundary}, or that it satisfies (SN), if $X$ and $B(X,S_*)$ are both semi-normal.
    \item[(HN)] We say that $(X,S_*)$ has {\it hereditarily normal strata}, or that it satisfies (HN), if \begin{enumerate}
            \item the normalization $\pi: (X^n,\pi^{-1}S_*)\to (X,S_*)$ is stratifiable,
            \item $(X^n,S_*^n)$ satisfies (N), and
            \item $B(X^n,\pi^{-1}S_*)$ satisfies (HN).
    \end{enumerate}       
    \item[(HSN)] We say that $(X,S_*)$ has {\it hereditarily semi-normal boundary}, or that it
    satisfies (HSN), if \begin{enumerate}
            \item the normalization $\pi: (X^n,\pi^{-1}S_*)\to (X,S_*)$ is stratifiable,
            \item $(X,S_*)$ satisfies (SN), and
            \item $B(X^n,\pi^{-1}S_*)$ satisfies (HSN).
    \end{enumerate}
\end{enumerate}
\end{defn}

Next we give a special stratification that is induced by the lc crepant log structure. 

\begin{defn}[Lc stratification for generalized pairs]
Let $f:(X,\Delta,\Mm)\to Z$ be an lc crepant log structure. Let $S^*_i(Z,X,\Delta,\Mm)\subset Z$ be the union of all $\le i$-dimensional lc centers of $f:(X,\Delta,\Mm)\to Z$, and
$$
S_i(Z,X,\Delta,\Mm):=S^*_i(Z,X,\Delta,\Mm)~\backslash ~S^*_{i-1}(Z,X,\Delta,\Mm).
$$
If the lc crepant log structure $f:(X,\Delta,\Mm)\to Z$ is clear from the context, we will use $S_i(Z)$ as an abbreviation. It is clear that each $S_i(Z)$ is a locally closed subspace of $Z$ of pure dimension $i$, and $Z$ is the disjoint union of all $S_i(Z)$. 

The stratification of $Z$ induced by $S_i(Z)$ is called the \emph{lc stratification} of $Z$ induced by $f:(X,\Delta,\Mm)\to Z$. Since this is the only stratification we are going to use in the rest of this paper, we usually will not emphasize the lc crepant structure $f:(X,\Delta,\Mm)\to Z$, and we will denote the corresponding stratified scheme by $(Z,S_*)$. The \emph{boundary} of $(Z,S_*)$ is the closed subspace
$$B(Z,S_*):=Z\backslash S_{\dim Z}(Z)=\cup_{i<\dim Z}S_i(Z).$$
\end{defn}

\begin{defn}\label{defn: of glc origin}
We say that a semi-normal stratified space $(Y,S_*)$ is \textit{of lc origin} if $S_i(Y)$ is unibranch (see \cite[Definition 1.44]{Kol13}) for any $i$, and there are lc crepant log structures $f_j:(X_j,\Delta_j,\Mm^j)\to Z_j$ with lc stratifications $(Z_j,S_{*}^j)$ and a finite surjective stratified morphism $\pi: \amalg_j(Z_j,S_{*}^j)\to (Y,S_*)$.
\end{defn}

\subsection{Semi-normality of lc centers and lc origin}
In this subsection we show that lc centers of lc generalized pairs are semi-normal.

\begin{thm}\label{thm: glc locus is semi-normal}
Let $f:(X,\Delta,\Mm)\to Z$ be an lc crepant log structure. Let $W\subset Z$ be the union of all lc centers of $f:(X,\Delta,\Mm)\to Z$ except $Z$, and $B(W)\subset W$ the union of all non-maximal (with respect to inclusion) lc centers that are contained in $W$. Then
\begin{enumerate}
    \item $W$ is semi-normal, and
    \item $W\backslash B(W)$ is normal.
\end{enumerate}
\end{thm}
\begin{proof}
Let $(Z,\Delta_Z,\NN)/U$ be an lc g-pair induced by the canonical bundle formula $/U$ of $f: (X,\Delta,\Mm)\rightarrow Z$. By Theorem \ref{thm: cbf gpair nonnqc}, the lc centers of $(Z,\Delta_Z,\NN)$ are exactly the lc centers of $f: (X,\Delta,\Mm)\rightarrow Z$. Possibly replacing $(X,\Delta,\Mm)$ with a dlt model of $(Z,\Delta_Z,\NN)$, we may assume that $f$ is birational and $(X,\Delta,\Mm)$ is $\Qq$-factorial dlt. We have $W=f(\lf\Delta\rf)$. Let $\Delta':=\{\Delta\}$. We consider the exact sequence 
$$
0\to\Oo_X(-\lf\Delta\rf)\to\Oo_X\to\Oo_{\lf\Delta\rf}
$$
and its push-forward
$$
\Oo_Z=f_*\Oo_X\to f_*\Oo_{\lf\Delta\rf}\stackrel{\delta}{\longrightarrow}R^1f_*\Oo_X(-\lf\Delta\rf).
$$
By \cite[Lemma 3.4]{HL22}, we can find an $\Rr$-divisor $\Delta''\ge 0$ such that $$-\lf\Delta\rf\sim_{\Rr,Z}K_X+\Delta'+\Mm_X\sim_{\Rr,Z}K_X+\Delta''$$ and $(X,\Delta'')$ is klt. Since $-\lfloor\Delta\rfloor$ is a Weil divisor, by \cite[Lemma 5.3, Theorem 5.6]{HLS19}, possibly perturbing $\Delta''$, we may assume that $\Delta''$ is a $\Qq$-divisor and 
$$-\lf\Delta\rf\sim_{\mathbb Q,Z}K_X+\Delta''.$$
By \cite[Corollary 10.40]{Kol13}, $R^if_*\Oo_{X}(-\lf\Delta\rf)$ is torsion-free for every $i$. On the other hand, $f_*\Oo_{\lf\Delta\rf}$ is supported on $W$, hence it is a torsion sheaf. Thus the connecting map $\delta$ is zero, hence $\Oo_Z\twoheadrightarrow f_*\Oo_{\lf\Delta\rf}$ is surjective. Since this map factors through $\Oo_W$, and we conclude that $\Oo_W\twoheadrightarrow f_*\Oo_{\lf\Delta\rf}$ is also surjective, hence an isomorphism.

Note that $\lf\Delta\rf$ has only nodes at codimension 1 points and it is $S_2$ by \cite[Corollary 2.88]{Kol13}. By \cite[Lemma 10.14]{Kol13}, $\lf\Delta\rf$ is semi-normal. By \cite[Lemma 10.15]{Kol13}, $W$ is semi-normal. This is (1).

To prove (2), let $V\subset\lf\Delta\rf$ be an irreducible component of its non-normal locus. Then $V$ is an lc center of $(X,\Delta)$, hence an lc center of $(X,\Delta,\Mm)$. Thus $f(V)$ is an lc center of $f: (X,\Delta,\Mm)\rightarrow Z$. Hence either $f(V)$ is an irreducible component of $W$, or $f(V)\subset B(W)$. Thus \cite[Complement 10.15.1]{Kol13} implies that $W \backslash B(W)$ is normal.
\end{proof}

\begin{cor}
Let $(X,\Delta,\Mm)$ be an lc g-pair. Then $\Nklt(X,\Delta,\Mm)$ is semi-normal.
\end{cor}
\begin{proof}
It follows from Theorem \ref{thm: glc locus is semi-normal} when $f$ is the identity morphism.
\end{proof}

\begin{lem}\label{lem: (Z,S) is U and SN} (cf. \cite[Lemma 5.26]{Kol13})
Let $f:(X,\Delta,\Mm)\to Z$ be an lc crepant log structure and $(Z,S_*)$ the induced lc stratification. Then
\begin{enumerate}
    \item  $S_i(Z)$ is unibranch for every $i$, and 
    \item  $B(Z,S_*)$ is semi-normal.
\end{enumerate}
\end{lem}

\begin{proof}
(1) follows from Lemma \ref{lem: intersection of lc center gpair}(2) and (2) follows from Theorem \ref{thm: glc locus is semi-normal}.
\end{proof}

\begin{lem}\label{lem: stratification is compatible under adjunction} (cf. \cite[Proposition 4.42]{Kol13})
Let $f: (X,\Delta,\Mm)\to Z$ be a dlt crepant log structure, $(Z,S_*)$ its induced lc stratification, and $Y\subset X$ an lc center of $(X,\Delta,\Mm)$. Let $(Y,\Delta,\Mm^Y)/Z$ be the dlt g-pair induced by adjunction to the higher-codimensional lc center $Y$, i.e. 
$$K_Y+\Delta_Y+\Mm^Y_Y:=(K_X+\Delta+\Mm_X)|_Y.$$
We consider the Stein factorization of $f|_Y$
$$(Y,\Delta_Y,\Mm^Y)\stackrel{f_Y}{\longrightarrow}W\stackrel{\pi}{\longrightarrow}Z.$$
Then:
\begin{enumerate}
    \item $f_Y:(Y,\Delta_Y,\Mm^Y)\to W$ is a dlt crepant log structure which induces an lc stratification $(W,S_*)$.
    \item $S_i(W)=\pi^{-1}(S_i(Z))$ for every $i$.
\end{enumerate}
\end{lem}
\begin{proof}
It follows from Lemma \ref{lem: intersection of lc center gpair}.
\end{proof}

\begin{thm}\label{thm: (Z,S) is HN and HSN}
Let $f:(X,\Delta,\Mm)\to Z$ be an lc crepant log structure and $(Z,S_*)$ the induced lc stratification. Then $(Z,S_*)$ satisfies (HN) and (HSN).
\end{thm}
\begin{proof}
By Lemma \ref{lem: (Z,S) is U and SN} and \cite[Definitions 9.18,~9.19]{Kol13}, $(Z,S_*)$ satisfies (HU) and (HSN). By \cite[Theorem 9.21]{Kol13}, $(Z,S_*)$ satisfies (HN).
\end{proof}

\begin{lem}(cf. \cite[5.29]{Kol13})\label{lem: glc stratification is of glc origin}
Every lc stratification is of lc origin. More precisely, let $f:(X,\Delta,\Mm)\to W $ be an lc crepant log structure and $Y\subset W$ any union of lc centers. Then $(Y, S_*)$ is of lc origin, where $S_i(Y)=Y\cap S_i(W)$ for each $i$.
\end{lem}

\begin{proof}
By Theorem \ref{thm: (Z,S) is HN and HSN} and \cite[Theorem 9.26]{Kol13} we know that $Y$ is semi-normal and $S_i(Y)$ is unibranch for each $i$. Then we can apply Lemma \ref{lem: stratification is compatible under adjunction} to each lc center of $f: (X,\Delta,\Mm)$ contained in $Y$ to conclude that $(Y,S_*)$ is of lc origin.
\end{proof}

\subsection{Du Bois property}\label{subsec: du bois}

In this subsection, we show that lc generalized pairs have Du Bois singularities. This subsection is parallel to \cite[Section 6]{LX23b}. 

We recall the following definition in \cite{Kov11} (cf. \cite[Definition 6.10]{Kol13}). 

\begin{defn}
A \emph{DB pair} $(X,\Sigma)$ consists of a reduced scheme $X$ of finite type and a closed reduced subscheme $\Sigma$ in $X$ such that the natural morphism
$$
\mathcal{I}_{\Sigma\subset X}\to \underline{\Omega}_{X,\Sigma}^0
$$
is a quasi-isomorphism. We will also say $(X,\Sigma)$ is \emph{DB} in this case.
\end{defn}

The definition of DB pairs is subtle but what really matters here is the following lemma:

\begin{lem}[{\cite[Proposition 6.15]{Kol13}}]\label{lem: property of DB pairs}
Let $(X,\Sigma)$ be a DB pair. Then $X$ has Du Bois singularities if and only if $\Sigma$ has Du Bois singularities.
\end{lem}

The following theorems are analogues of \cite[Theorems 6.31, 6.33]{Kol13} for g-pairs and the proofs are similar. For the reader's convenience, we provide full proofs here.

\begin{thm}\label{thm: (Z,W) is DB for glc crepant log structure}
Let $(X,B,\Mm)/U$ be an lc g-pair and $f: (X,B,\Mm)\rightarrow Z$ an lc-trivial fibration $/U$.
Let $W\subset Z$ be the union of lc centers of $f: (X,B,\Mm)\rightarrow Z$ except $Z$. Then $(Z,W)$ is a DB pair.
\end{thm}

\begin{proof}
Let $(Z,B_Z,\Mm^Z)/U$ be a g-pair induced by $f: (X,B,\Mm)\rightarrow Z$. By Theorem \ref{thm: cbf gpair nonnqc}, the lc centers of $(Z,B_Z,\Mm^Z)$ are exactly the lc centers of $f: (X,B,\Mm)\rightarrow Z$. Thus we can assume that $f$ is the identity, $(X,B,\Mm)=(Z,B_Z,\Mm^Z)$, and $W=\Nklt(X,B,\Mm)$.

Let $g: Y\to X$ be a log resolution of $(X,\Supp B)$ such that $\Mm$ descends to $Y$ and $F:=g^{-1}(W)$ is an snc divisor. Let
$$K_Y+B_Y+\Mm_Y:=g^*(K_X+B+\Mm_X)$$
and $D:=B_Y^{=1}$. Since $\Mm_Y$ is nef $/X$ and big $/X$, there exists $0\le B'_Y\sim_{\Rr,X}\Mm_Y$ such that $(Y,B_Y-D+B'_Y)$ is sub-klt. Possibly replacing $Y$ with a higher resolution, we may assume that $(Y,\Supp B_Y\cup\Supp D\cup\Supp B_Y')$ is log smooth. Let 
$$\bar B_Y:=(B_Y-D+B'_Y)^{\ge0}+\{(B_Y-D+B'_Y)^{\le 0}\}$$ and 
$$E:=\lfloor (B_Y-D+B'_Y)^{\le 0}\rfloor$$, then $\lf\bar B_Y\rf=0$ and $E$ is a g-exceptional Weil divisor. In particular, $(Y,\bar B_Y)$ is klt.

Since $E-D\ge-F$, we have natural maps
$$
g_*\Oo_Y(-F)\to Rg_*\Oo_Y(-F)\to Rg_*\Oo_Y(E-D).
$$
Since $E-D\sim_{\Rr,X}K_Y+\bar B_Y$ and $E-D$ is a Weil divisor, by \cite[Lemma 5.3, Theorem 5.6]{HLS19}, $E-D\sim_{\Qq,X}K_Y+\bar B_Y'$ for some klt $\Qq$-pair $(Y,\bar B_Y')$. By \cite[Theorem 10.41]{Kol13},
$$
Rg_*\Oo_Y(E-D)\simeq_{qis}\sum_{i}R^ig_*\Oo_Y(E-D)[i].
$$
Thus we get a morphism 
$$
g_*\Oo_Y(-F)\to Rg_*\Oo_Y(-F)\to Rg_*\Oo_Y(E-D)\to g_*\Oo_Y(E-D).
$$
Note that 
$$
g_*\Oo_Y(E-D)=g_*\Oo_Y(E-D)\cap g_*\Oo_Y(E)=g_*\Oo_Y(E-D)\cap g_*\Oo_Y=g_*\Oo_Y(-D).
$$
Since $D$ is reduced and $g(D)=W$, we have $g_*\Oo_Y(-D)=\mathcal{I}_W$, the ideal sheaf of $W$ in $Z=X$. Moreover, $g_*\Oo_Y(-F)=\mathcal{I}_W$ since $F$ is also reduced. Therefore, we get an isomorphism $\mathcal{I}_W=g_*\Oo_Y(-F)\to g_*\Oo_Y(E-D)$, which implies that 
$$
\rho: \mathcal{I}_W\simeq g_*\mathcal{I}_F\to Rg_*\mathcal{I}_F
$$
has a left inverse. Since $Y$ is smooth and $F$ is an snc divisor, we see that $(Y,F)$ is a DB pair, thus by \cite[Theorem 3.3]{Kov12} (cf. \cite[Theorem 6.27]{Kol13}), $(Z,W)$ is also a DB pair. \end{proof}

\begin{thm}\label{thm: of glc origin implies DB}
Let $(X,S_*)$ be a stratified scheme of lc origin (Definition \ref{defn: of glc origin}). Then $X$ is Du Bois.
\end{thm}

\begin{proof}
We use induction on the dimension.

Let $\pi: (X^n,S^n_*)\to (X,S_*)$ denote the normalization. Let $B(X)\subset X$ and
$B(X^n)\subset X^n$ denote the corresponding boundaries. By \cite[9.15.1]{Kol13}, we have a universal push-out diagram
\begin{center}
$\xymatrix{
B(X^n)\ar@{^(->}[r]\ar@{->}[d] & X^n\ar@{->}[d]^{\pi}\\
 B(X)\ar@{^(->}[r]& X\\
}$
\end{center}
Notice that $B(X)$ and $B(X^n)$ are of lc origin by Lemma \ref{lem: glc stratification is of glc origin}, hence Du Bois by induction.

Since $\pi$ is finite, it follows that $R\pi_*\mathcal{I}_{B(X^n)\subset X^n}=\pi_*\mathcal{I}_{B(X^n)\subseteq X^n}$. Furthermore, $\pi_*\mathcal{I}_{B(X^n)\subseteq X^n}=\mathcal{I}_{B(X)\subseteq X}$ by \cite[Theorem 9.30]{Kol13}. By \cite[Theorem 3.3]{Kov12} and Lemma \ref{lem: property of DB pairs}, we only need to show that $X^n$ is Du Bois. By assumption, for each irreducible component $X_i^n\subset X^n$, there exists an lc crepant log structure $f_i:(Y_i,\Delta_i,\Mm)\to Z_i$ and a finite surjection $Z_i\to X_i^n$. By \cite[Corollary 2.5]{Kov99}, we only need to show that $Z_i$ is Du Bois for each $i$. Let $B(Z_i)\subset Z_i$ be the boundary of the lc stratification of $Z_i$. Then $B(Z_i)$ is of lc origin by Lemma \ref{lem: glc stratification is of glc origin}, hence Du Bois by induction. By Theorem \ref{thm: (Z,W) is DB for glc crepant log structure}, $(Z_i,B(Z_i))$ is a DB pair, hence $Z_i$ is Du Bois and we are done.
\end{proof}

\begin{proof}[Proof of Theorem \ref{thm: glc sings are Du Bois}]
   Let $W$ be any union of the glc centers, then by Lemma \ref{lem: glc stratification is of glc origin} the induced stratified space $(W,S_*)$ is of lc origin. Theorem \ref{thm: glc sings are Du Bois} follows from Theorem \ref{thm: of glc origin implies DB}.
\end{proof}

\section{Vanishing and contraction theorems for lc generalized pairs}\label{sec: vanishing gpair}

The goal of this section is to prove the vanishing theorems and contraction theorems for lc generalized pairs. This section is parallel to \cite{CLX23}, except that the canonical bundle formula and the subadjunction formulas are replaced by the ones established in Sections \ref{sec: cbf} and \ref{sec: subadj}.

\subsection{Vanishing theorems}

\begin{thm}\label{thm: kod vanishing with lc strata}
Let $(X,B,\Mm)/U$ be an lc g-pair associated with projective morphism $\pi: X\rightarrow U$, $D$ a Cartier divisor on $X$ such that $D-(K_X+B+\Mm_X)$ is nef$/U$ and log big$/U$ with respect to $(X,B,\Mm)$, and $Y$ a union of lc centers of $(X,B,\Mm)$ such that $Y\not=X$. Then:
\begin{enumerate}
	\item $R^i\pi_*\mathcal{O}_Y(D)=0$ for any positive integer $i$.
	\item $R^i\pi_*\mathcal{O}_X(D)=0$ for any positive integer $i$.
    \item The map $\pi_*\mathcal{O}_X(D)\rightarrow \pi_*\mathcal{O}_Y(D)$ is surjective.
    \item $R^i\pi_*(\mathcal{I}_Y\otimes\mathcal{O}_X(D))=0$ for any positive integer $i$, where $\mathcal{I}_Y$ is the defining ideal sheaf of $Y$ on $X$.
\end{enumerate}
\end{thm}
\begin{proof}
This follows from exactly the same lines of the proof of \cite[Theorem 1.1]{CLX23}, using the log canonical stratifications and induction on dimensions with the help of the push-out diagram (giving us some short exact sequences), except that we replace the canonical bundle formula and the subadjunction formula therein with Definition-Theorems \ref{thm: cbf gpair nonnqc} and \ref{prop: lc subadj is lc} respectively.
\end{proof}

\subsection{Base-point-freeness theorem and contraction theorem}\label{subsec: bpf nonnqc}

\begin{lem}\label{lem: non-vanishing of lc gpair}
 Let $(X,B,\Mm)/U$ be an lc g-pair and $D$ a nef$/U$ Cartier divisor on $X$ such that $aD-(K_X+B+\Mm_X)$ is ample$/U$ for some positive real number $a$. Let $Y$ be a minimal lc center of $(X,B,\Mm)$ if $(X,B,\Mm)$ is not klt, and let $Y:=X$ if $(X,B,\Mm)$ is klt. Let $D_Y:=D|_Y$. Then for any integer $m\gg 0$
\begin{enumerate}
    \item $\mathcal{O}_{Y}(mD_Y)$ is globally generated over $U$,
    \item $|mD/U|\not=\emptyset$, and
    \item $Y$ is not contained in $\Bs|mD/U|$.
\end{enumerate}
\end{lem}
\begin{proof}
When $(X,B,\Mm)$ is klt, by \cite[Lemma 3.4]{HL22}, there exists a klt pair $(X,\Delta)$ such that $D-(K_X+\Delta)$ is ample$/U$. By the usual base-point-freeness theorem (cf. \cite[Theorem 3-1-1]{KMM87}), the lemma follows.

When $(X,B,\Mm)$ is not klt, by Theorem \ref{thm: (Z,S) is HN and HSN}, $Y$ is normal. By Theorem \ref{thm: kod vanishing with lc strata}(3), the map $f_*\mathcal{O}_X(mD)\rightarrow f_*\mathcal{O}_Y(mD_Y)$ is surjective for any positive integer $m\geq a$. Thus (2)(3) follow from (1) and we only need to prove (1). If $\dim Y=0$ then there is nothing left to prove. If $\dim Y>0$, then by Definition-Lemma \ref{deflem: subadj minimal lc center}, there exists a klt g-pair $(Y,B_Y,\Mm^{Y})/U$ such that $K_{Y}+B_{Y}+\Mm^{Y}_{Y}\sim_{\Rr,U}(K_X+B+\Mm_X)|_{Y}$. Thus $D_Y-(K_{Y}+B_{Y}+\Mm^{Y}_{Y})$ is nef$/U$ and log big$/U$ with respect to $(Y,B_Y,\Mm^Y)$. By \cite[Lemma 3.4]{HL22}, there exists a klt pair $(Y,\Delta_Y)$ such that $D_Y-(K_Y+\Delta_Y)$ is ample$/U$. By the usual base-point-freeness theorem (cf. \cite[Theorem 3-1-1]{KMM87}), the lemma follows.
\end{proof}

\begin{proof}[Proof of Theorem \ref{thm:base-point-freeness intro}]
By Lemma \ref{lem: non-vanishing of lc gpair}, we may let $m_0$ be the minimal positive integer such that $|mD|\not=\emptyset$ for any integer $m\geq m_0$.

\begin{claim}\label{claim: induction bs}
Let $\{p_i\}_{i=1}^{+\infty}$ be a strictly increasing sequence of positive integers. There exist a non-negative integer $M$ and integers $i_1<i_2<\dots<i_{M+1}$ satisfying the following. Let $s_k:=\prod_{l=1}^kp_{i_l}$ for any $1\leq k\leq M+1$, then
\begin{enumerate}
    \item $|s_1D/U|\not=\emptyset$,
    \item $\Bs|s_kD/U|\supsetneq\Bs|s_{k+1}D/U|$ for any $1\leq k\leq M$, and
    \item $\Bs|s_{M+1}D/U|=\emptyset$.
\end{enumerate}
\end{claim}
\begin{proof}
We may take $i_1$ to be any integer such that $p_{i_1}\geq m_0$, then (1) holds. 

Suppose that we have already found $i_1,\dots,i_k$ for some positive integer $k$. Let $d:=\dim X$, let $H_{1},\dots,H_{d+1}$ be $d+1$ general elements in $|s_kD/U|$, and let $H:=H_{1}+\dots+H_{d+1}$. Then $(X,B+H,\Mm)$ is lc outside $\Bs|s_kD/U|$. If $\Bs|s_kD/U|=\emptyset$, then we may let $M:=k-1$ and we are done. Thus we may assume that $\Bs|s_kD/U|\not=\emptyset$.

Since every $H_{j}$ contains $\Bs|s_kD/U|$, by \cite[Theorem 18.22]{Kol+92}, $(X,B+H,\Mm)$ is not lc near $\Bs|s_kD/U|$. Let 
$$c:=\sup\{t\mid t\geq 0, (X,B+tH,\Mm)\text{ is lc}\},$$
then $c\in [0,1)$, and there exists at least one lc center of $(X,B+cH,\Mm)$ which is contained in $\Bs|s_kD/U|$. Let
$\mathcal{S}$ be the set of all lc centers of $(X,B+cH,\Mm)$ that are contained in $\Bs|s_kD/U|$, and let $Y$ be a minimal lc center in $\mathcal{S}$. Since
$$(a+s_k(d+1))D-(K_X+B+cH+\Mm_X)\sim_{\mathbb R}s_k(d+1)(1-c)D+(aD-(K_X+B+\Mm_X))$$
is ample$/U$, by Lemma \ref{lem: non-vanishing of lc gpair}, there exists a positive integer $n$ such that for any integer $m\geq n$, $|ms_kD/U|\not=\emptyset$ and $\Bs|ms_kD/U|$ does not contain $Y$. In particular, $\Bs|ms_kD/U|\subsetneq\Bs|s_kD/U|$. We may let $i_{k+1}$ be any integer such that $i_{k+1}>i_{k}$ and $p_{i_{k+1}}\geq n$. This construction implies (2). (3) follows from (2) and the Noetherian property.
\end{proof}

\noindent\textit{Proof of Theorem \ref{thm:base-point-freeness intro} continued}. We let $p$ and $q$ be two different prime numbers. By Claim \ref{claim: induction bs}, there exist two non-negative integers $M,N$ such that $\mathcal{O}_X(p^MD)$ and $\mathcal{O}_X(q^ND)$ are globally generated$/U$. Since $p^M$ and $q^N$ are coprime, for any integer $m\gg 0$, we may write $m=bp^M+cq^N$ for some non-negative integers $b,c$, hence
$$\Bs|mD/U|\subset\Bs|p^MD/U|\cup\Bs|q^ND/U|=\emptyset.$$
Therefore, $\mathcal{O}_X(mD)$ is globally generated over $U$ for any integer $m\gg 0$.
\end{proof}

\begin{thm}[Contraction theorem for lc generalized pairs, cf. {\cite[Theorem 1.5]{Xie22}}]\label{thm: cont thm gpair}
Let $(X,B,\Mm)/U$ be an lc generalized pair and $F$ a $(K_X+B+\Mm_X)$-negative extremal face$/U$. Then there exists a contraction$/U$ $\cont_F: X\rightarrow Z$ of $F$ satisfying the following.
\begin{enumerate}
    \item For any integral curve $C$ on $X$ such that the image of $C$ in $U$ is a closed point, $\cont_F(C)$ is a point if and only if $[C]\in F$.
    \item $\mathcal{O}_Z=(\cont_F)_*\mathcal{O}_X$. In other words, $\cont_F$ is a contraction.
    \item For any Cartier divisor $D$ on $X$ such that $D\cdot C=0$ for any curve $C$ contracted by $\cont_F$, there exists a Cartier divisor $D_Z$ on $Z$ such that $D=\cont_F^*D_Z$.
\end{enumerate}
\end{thm}
\begin{proof}%[Proof of Theorem \ref{thm: cont thm gpair}]
(1)(2) By Theorem \ref{thm: cone theorem gfq}, $F$ is a finite-dimensional rational $(K_X+B+\Mm_X)$-negative extremal face$/U$. Thus there exists a nef Cartier divisor $L$ on $X$ that is the supporting function of $F$. Then $L-(K_X+B+\Mm_X)$ is ample. By Theorem \ref{thm:base-point-freeness intro}, $mL$ is base-point-free$/U$, hence defines a contraction$/U$. Denote this contraction by $\cont_F$. Then $\cont_F$ satisfies (1) and (2).

(3) Since $D-(K_X+B+\Mm_X)$ is ample$/Z$, by Theorem \ref{thm:base-point-freeness intro}, $\mathcal{O}_X(mD)$ is globally generated over $Z$ for any integer $m\gg 0$. Since $D\cdot C=0$ for any curve $C$ contracted by $\cont_F$, $\cont_F$ is defined by $|mD|$ for any integer $m\gg 0$. Thus $mD=\cont_F^*D_{Z,m}$ and $(m+1)D=\cont_F^*D_{Z,m+1}$ for any integer $m\gg 0$. We may let $D_Z:=D_{Z,m+1}-D_{Z,m}$.
\end{proof}

\part{Good minimal model and the proofs of the main theorems}\label{part:gmm}

\section{Existence of good minimal models and \texorpdfstring{$\bb$}{}-semi-ampleness}\label{sec: gmm fdlt}

\subsection{Good minimal models for polarized foliations}

We remark that in the following lemma, in general, we do not require $\Ff$ to be algebraically integrable, so it may be applicable to other scenarios.
\begin{lem}\label{lem: +a keep under mmp}
Let $(X,\Ff,B,\Mm)/U$ be an lc gfq and $A$ an ample$/U$ $\Rr$-divisor on $X$. Let 
$$\phi: (X,\Ff,B+A,\Mm)\dashrightarrow (X',\Ff',B'+A',\Mm)$$
be a sequence of finite steps of an $(K_{\Ff}+B+A+\Mm_X)$-MMP$/U$. Then there exist a nef$/U$ $\bb$-divisor $\Nn$ and an ample$/U$ $\Rr$-divisor $\tilde A'$ on $X'$, such that 
\begin{enumerate}
    \item $(X',\Ff',B',\Nn)/U$ is lc,
    \item $\Nn_{X'}+\tilde A'\sim_{\mathbb R,U}\Mm_{X'}+A'$, and
    \item $\Nn-\Mm$ descends to an ample$/U$ $\Rr$-divisor on $X$.
\end{enumerate}
Moreover, if $(X,\Ff,B,\Mm;G)/Z$ satisfies Property $(*)$ (resp. is (weak) ACSS) with associated $X\to Z$ and $G$, then $(X',\Ff',B',\Nn;G':=\phi_*G)/Z$ satisfies Property $(*)$ (resp. is (weak) ACSS).
\end{lem}
\begin{proof}
    We may assume that $\phi$ is a single step of an MMP$/U$, and we have the following diagram$/U$
    \begin{center}$\xymatrix{
   X\ar@{->}[rd]^{g}\ar@{-->}[rr]^{\phi} & & X'\ar@{->}[dl]^{h}\\
    & T &
}$
\end{center}
such that either $\phi=g$ is a divisorial contraction, or $\phi$ is a flip, $g$ is the flipping contraction, and $h$ is the flipped contraction. Then there exists an ample$/T$ divisor $H$ such that $K_{\Ff}+B+A+H+\Mm_X\sim_{\mathbb R,T}0$. Let $H':=\phi_*H$, then $-H'$ is ample$/T$. We may choose an ample$/U$ divisor $C$ on $T$ and a positive real number $\epsilon$ such that both $\tilde A':=h^*C-\epsilon H'$ and $L:=A-g^*C+\epsilon H$ are ample$/U$, and 
$\phi$ is also a step of a $(K_\Ff+B+L+\Mm_X)$-MMP$/U$. 
Let $\Nn:=\Mm+\overline{L}$. By our construction, $\Nn$ and $\tilde A'$ satisfy (1-3). 

Moreover, if $(X,\Ff,B,\Mm;G)/Z$ satisfies Property $(*)$ (resp. is (weak) ACSS), then $(X,\Ff,B,\Nn;G)/Z$ satisfies Property $(*)$ (resp. is (weak) ACSS) as $\Nn-\Mm$ descends to $X$. Furthermore, since $\phi$ is a step of a $(K_{\Ff}+B+\Nn_X)$-MMP$/U$, $(X',\Ff',B',\Nn;G':=\phi_*G)/Z$ satisfies Property $(*)$ (resp. is (weak) ACSS) by Lemma \ref{lem: ACSS mmp can run}.
\end{proof}

\begin{lem}\label{lem:superundercbf}
Let $(X,B,\Mm)/U$ be an lc g-pair and $f:X\to Z$ a contraction such that $B$ is super over $Z$. Assume that $\phi:X\to T$ is a contraction$/U$ such that $K_X+B+\Mm_X\sim_{\Rr,T}0$ and $\phi$ is also a contraction$/Z$. Then $B_T$, the discriminant part of $(X,B,\Mm)$ over $T$, is super$/Z$.
\end{lem}
\begin{proof}
By assumption, there exist ample Cartier divisors $H_1,\dots,H_{2\dim X+1}$ on $Z$ such that $B\geq\sum f^*H_i$. In particular, $B-\sum f^*H_i\ge0$, and 
$$K_X+B-\sum f^*H_i+\Mm_X\sim_{\Rr,T}0.$$
Let $B_{T}'$ be the discriminant part of $(X,B-\sum f^*H_i,\Mm)$ over $T$ and $\psi:T\to Z$. Then we have $B_T=B_T'+\psi^*H_i$. The lemma follows.
\end{proof}

\begin{thm}\label{thm: bpf induced gfq}
Let $(X,\Ff,B,\Mm)/U$ be an lc gfq of dimension $d$ and $A\geq 0$ an ample$/U$ $\Rr$-divisor on $X$ such that
\begin{itemize}
  \item $\Ff$ is induced by a contraction $f: X\rightarrow Z$, 
  \item $K_{\Ff}+B+A+\Mm_X$ is nef$/U$, and
  \item $K_{\Ff}+B+\Mm_X\sim_{\mathbb R,Z}K_X+\Delta+\Nn_X$ for some lc g-pair $(X,\Delta,\Nn)/U$.
\end{itemize}
Then the following hold.
\begin{enumerate}
  \item $K_{\Ff}+B+A+\Mm_X$ is semi-ample$/U$.
  \item The contraction$/U$ defined by $K_{\Ff}+B+A+\Mm_X$ is a contraction$/Z$.
  \item Suppose that 
  $$m(K_{\Ff}+B+\Mm_X)\sim_{Z}m(K_X+\Delta+\Nn_X)$$ 
  and $m(K_{\Ff}+B+A+\Mm_X)$ is Cartier for some positive integer $m$. Then $$\mathcal{O}_X(nm(K_{\Ff}+B+A+\Mm_X))$$ is globally generated$/U$ for any integer $n\gg 0$.
\end{enumerate}
\end{thm}
\begin{proof}
Let $\pi: X\rightarrow U$ be the induced morphism and let $H'$ be a sufficiently ample Cartier divisor on $U$. Possibly replacing $A$ with $A+\pi^*H'$, we may assume that $A$ is ample. Let $H_Z'$ be a sufficiently ample Cartier divisor on $Z$. Possibly replacing $\Delta$, we may assume that $\Delta$ is super$/Z$. Then by Lemma \ref{lem: equivalence over bases}, $K_{X}+\Delta+A+\Nn_X$ is nef$/U$.
By Theorem \ref{thm: semi-ampleness intro}, $K_{X}+\Delta+A+\Nn_X$ is semi-ample$/U$, so $K_{X}+\Delta+A+\Nn_X$ defines a contraction$/U$ $\phi: X\rightarrow T$. 
Since $\phi$ only contracts $(K_{X}+\Delta+\Nn_X)$-negative extremal rays$/U$, $\phi$ is a contraction$/Z$ and thus $K_{\Ff}+B+A+\Mm_X\sim_{\mathbb R,T}0$. Let $\Ff_T$ be the foliation induced by the contraction $\psi: T\rightarrow Z$, then $\Ff=\phi^{-1}\Ff_T$. 

We let $H_T$ be a general ample $\Rr$-divisor on $T$ such that $H:=A-\phi^*H_T$ is ample$/U$. By Theorem \ref{thm: cbf gfq}, there exist an lc gfq $(T,\Ff_T,B_T,\Mm^T)/U$ induced by $\phi: (X,\Ff,B,\overline{H}+\Mm)\rightarrow T$, and an lc g-pair $(T,\Delta_T,\Nn^T)/U$ induced by $\phi: (X,\Delta,\overline{H}+\Nn)\rightarrow T$. Note that $\Delta_T$ is super$/Z$ by Lemma \ref{lem:superundercbf}, and
$$K_{\Ff_T}+B_T+\Mm^T_T\sim_{\mathbb R,Z}K_T+\Delta_T+\Nn^T_T.$$ 
Since $\phi$ is the morphism$/U$ defined by $K_{X}+\Delta+A+\Nn_X$, $K_T+\Delta_T+H_T+\Nn^T_T$ is ample$/U$. Thus $K_T+\Delta_T+(1-\delta)H_T+\Nn^T_T$ is ample$/U$ for any $0<\delta\ll 1$. 
Then by Theorem \ref{thm: cone theorem gfq},
$(K_{\Ff_T}+B_T+(1-\delta)H_T+\Mm^T_T)$ is nef$/U$. It immediately implies that $K_{\Ff_T}+B_T+H_T+\Mm^T_T$ is ample$/U$, and thus $K_{\Ff}+B+A+\Mm_X$ is semi-ample$/U$, and $\phi$ is the contraction$/U$ defined by $K_{\Ff}+B+A+\Mm_X$. This implies (1)(2). 

We prove (3). Since $\phi$ is a contraction defined by $K_X+\Delta+A+\Nn_X$, $m(K_{\Ff}+B+A+\Mm_X)$ is Cartier, and $K_{\Ff}+B+A+\Mm_X\sim_{\mathbb Q,T}0$, by Theorem \ref{thm: cont thm gpair}(3), there exists a Cartier divisor $L$ on $T$ such that
$$m(K_{\Ff}+B+A+\Mm_X)=\phi^*L.$$
Since $L\sim_{\mathbb R}K_{\Ff_T}+B_T+H_T+\Mm^T_T$, $L$ is ample$/U$. Thus $nL$ is very ample$/U$ for any integer $n\gg 0$, so 
$$\mathcal{O}_X(nm(K_{\Ff}+B+A+\Mm_X))=\mathcal{O}_X(\phi^*(nL))$$
is globally generated$/U$ for any integer $n\gg 0$.
\end{proof}

\begin{thm}\label{thm: gmm polarized gfq}
    Let $(X,\Ff,B,\Mm)/U$ be an lc gfq, and $A,H$ two ample$/U$ $\Rr$-divisors on $X$. Assume that
    \begin{itemize}
        \item $\Ff$ is induced by a contraction $\pi: X\rightarrow Z$,
        \item $K_{\Ff}+B+A+\Mm_X$ is pseudo-effective$/U$, 
        \item $K_{\Ff}+B+\Mm_X\sim_{\mathbb R,Z}K_X+\Delta+\Nn_X$ for some lc g-pair $(X,\Delta,\Nn)/U$, and
        \item either $X$ is $\Qq$-factorial klt or $\Mm$ is NQC$/U$.
    \end{itemize}
   Then there exists a $(K_{\Ff}+B+A+\Mm_X)$-MMP$/U$ with scaling of $H$, say $\mathcal{P}_0$, satisfying the following. Let $\mathcal{P}=\mathcal{P}_0$ if $X$ is not $\Qq$-factorial, and let $\mathcal{P}$ be any $(K_{\Ff}+B+A+\Mm_X)$-MMP$/U$ with scaling of an ample$/U$ $\Rr$-divisor if $X$ is $\Qq$-factorial. Then 
   \begin{enumerate}
       \item $\mathcal{P}$ terminates at a model $(X',\Ff', B'+A', \Mm)$ such that $K_{\Ff'}+B'+A'+\Mm_X$ is semi-ample$/U$, and
       \item the contraction$/U$ defined by $K_{\Ff'}+B'+A'+\Mm_{X'}$ is a contraction$/Z$.
   \end{enumerate}
\end{thm}
\begin{proof}
Let $\Mm':=\Mm+\bar A$. Then $\mathcal{P}$ is a $(K_{\Ff}+B+\Mm'_X)$-MMP$/U$ and $(X,\Ff,B,\Mm')$ is lc. By Proposition \ref{prop: run mmp with scaling gfq}, $\mathcal{P}$ terminates with a log birational model $(X',\Ff',B',\Mm')/U$ of $(X,\Ff,B,\Mm')/U$ such that $K_{\Ff'}+B'+\Mm'_{X'}$ is nef$/U$. By Lemma \ref{lem: +a keep under mmp}, there exists a nef$/U$ $\bb$-divisor $\Mm''$ and an ample$/U$ $\Rr$-divisor $A''$ such that $\Mm''_{X'}+A''\sim_{\mathbb R,U}\Mm_{X'}+A'$ and $(X',\Ff',B',\Mm'')/Z$ is lc, where $A'$ is the image of $A$ on $X'$. Moreover, by Lemma \ref{lem: ACSS mmp can run}, $\mathcal{P}$ is also a $(K_X+\Delta+A+\Nn_X)$-MMP$/U$. Since $(X,\Delta,\Nn+\bar A)$ is lc, $(X',\Delta',\Nn+\bar A)$ is lc, where $\Delta'$ is the strict transform of $\Delta$ on $X'$. Moreover,
$$K_{\Ff'}+B'+A''+\Mm''_{X'}\sim_{\mathbb R,Z}K_{X'}+\Delta'+\Nn_{X'}+\bar A_{X'}.$$
The theorem follows from Theorem \ref{thm: bpf induced gfq}.
\end{proof}

\subsection{A special case of Prokhorov-Shokurov's effective \texorpdfstring{$\bb$}{}-semi-ampleness conjecture}

\begin{thm}\label{thm: a special b-semiampleness}
Let $(X,B,\Mm)/U$ be an lc g-pair and $f:(X,B,\Mm)\rightarrow Z$ a contraction satisfying Property $(*)$. Let $\Nn$ be the moduli part of $f: (X,B,\Mm)\rightarrow Z$. Assume that
\begin{enumerate}
  \item $f$ is equi-dimensional,
  \item $K_X+B+\Mm_X$ is nef$/Z$, 
  \item $(X,B,\Mm)$ is BP semi-stable$/Z$, and
  \item there exists an ample$/U$ $\Rr$-divisor $H$ such that either $B^h\geq H$ or $\Mm-\bar H$ is nef$/U$, where $B^h$ is the horizontal$/Z$ part of $B$.
\end{enumerate}
Then $\Nn$ descends to $X$ and $\Nn_X$ is semi-ample$/U$.
\end{thm}
\begin{proof}
    Let $\Ff$ be the foliation induced by $f$. By Proposition \ref{prop: bp semistable foliation lc}, $(X,\Ff,B^h,\Mm)$ is lc and thus $(X,\Ff,B^h,\Mm;B-B^h)/Z$ is weak ACSS. Then it follows from Theorem \ref{thm: lc+weak acc=bpstable} that $(X,B,\Mm)$ is BP stable$/Z$. According to Proposition \ref{prop: bp stable nef}, $\Nn$ descends to $X$ and is nef$/U$. By Proposition \ref{prop: weak cbf gfq}, $K_{\Ff}+B^h+\Mm_X\sim\Nn_X$ is nef$/U$. We may conclude that $\Nn_X=K_{\Ff}+B^h+\Mm_X$ is semi-ample$/U$ by Theorem \ref{thm: bpf induced gfq}.
\end{proof}

When we have an lc-trivial fibration, we can prove stronger $\bb$-semi-ampleness.

\begin{thm}\label{thm: a special effective b-semiampleness}
Let $d$ and $m$ be two positive integers. Then there exists a positive integer $I$ depending only on $d$ and $m$ satisfying the following.

Assume that $(X,B,\Mm)/U$ is an lc g-pair and $f: (X,B,\Mm)\rightarrow Z$ is a contraction$/U$ satisfying Property $(*)$. Let $\Nn$ be the moduli part of $f: (X,B,\Mm)\rightarrow Z$. Assume that
 \begin{enumerate}
        \item $f$ is equi-dimensional,
        \item $X$ is of Fano type over $Z$,
        \item $K_X+B+\Mm_X\sim_{\mathbb Q,Z}0$,
        \item $(X,B,\Mm)$ is BP semi-stable$/Z$,
        \item $mB$ is a Weil divisor and $m\Mm$ is $\bb$-base-point-free$/U$, and
        \item there exists an ample$/U$ $\Rr$-divisor $H$ such that either $B^h\geq H$ or $\Mm-\bar H$ is nef$/U$, where $B^h$ is the horizontal$/Z$ part of $B$.
    \end{enumerate}
Then $\Nn$ descends to $X$, $I\Nn_X$ is Cartier, and $\mathcal{O}_X(nI\Nn_X)$ is globally generated$/U$ for any integer $n\gg 0$.
\end{thm}
\begin{proof}
According to \cite[Lemma 4.2]{Has22a} (cf. \cite[Proposition 6.3]{Bir19}), there exist a positive integer $q$ depending only on $d$ and $m$ and a choice $\Mm^Z$ of a moduli part of a canonical bundle formula for $(X,B,\Mm)\to Z$ such that 
$$q(K_X+B+\Mm_X)\sim qf^*(K_Z+B_Z+\Mm^Z_Z)$$ 
and $q\Mm^Z$ is nef$/U$, where $B_Z$ is the discriminant part. By the assumption that $f: (X,B,\Mm)\rightarrow Z$ satisfies Property $(*)$, $Z$ is smooth and $B_Z$ is reduced. Hence $q(K_X+B+\Mm_X)$ is Cartier.

By Theorem \ref{thm: a special b-semiampleness}, $\Nn$ descends to $X$ and $\Nn_X$ is semi-ample$/U$. By Proposition \ref{prop: weak cbf gfq}, 
$$\Nn_X\sim K_{\Ff}+B^h+\Mm_X\sim_Z K_X+B+\Mm_X$$ 
is semi-ample$/U$, where $\Ff$ is the foliation induced by $f$. Since $Z$ is smooth and $q(K_X+B+\Mm_X)$ is Cartier, $q\Nn_X$ and $q(K_{\Ff}+B^h+\Mm_X)$ are Cartier, and we may let $I:=q$. By Theorem \ref{thm: bpf induced gfq}(3), $\mathcal{O}_X(nI(K_{\Ff}+B^h+\Mm_X))=\mathcal{O}_X(nI\Nn_X)$ is globally generated$/U$ for any integer $n\gg 0$.
\end{proof}

\begin{thm}\label{thm: special b-semiampleness lc trivial}
Let $(X,B,\Mm)/U$ be an lc g-pair, $G\geq 0$ an $\Rr$-Cartier $\Rr$-divisor on $X$, and $\pi: X\rightarrow Z$ an equi-dimensional contraction$/U$. Assume that
\begin{itemize}
    \item $G$ is vertical$/Z$,
    \item $\pi: (X,B+G,\Mm)\rightarrow Z$ satisfies Property $(*)$,
    \item $(X,B+G,\Mm)$ is BP semi-stable$/Z$,
    \item $K_X+B+\Mm_X\sim_{\mathbb R,Z}0$,
    \item there exists an ample$/U$ $\Rr$-divisor $H$ such that either $B^h\geq H$ or $\Mm-\bar H$ is nef$/U$, where $B^h$ is the horizontal$/Z$ part of $B$, and
    \item either $X$ is $\Qq$-factorial klt or $\Mm$ is NQC$/U$.
\end{itemize}
Let $\Nn$ be the moduli part of $\pi: (X,B,\Mm)\rightarrow Z$. Then:
\begin{enumerate}
    \item $\Nn$ descends to $X$ and $\Nn_X$ is semi-ample$/U$.
    \item Suppose that there exists a positive integer $m$ such that $mB^h$ is a Weil divisor and $m\Mm$ is $\bb$-base-point-free$/U$, and $X$ is of Fano type over $Z$. Then there exists a positive integer $I$ depending only on $\dim X$ and $m$, such that $I\Nn_X$ is Cartier and $\mathcal{O}_X(nI\Nn_X)$ is globally generated$/U$ for any integer $n\gg 0$.
\end{enumerate}
\end{thm}
\begin{proof}
For any prime divisor $D$ on $Z$, we let 
$$t_D:=\sup\{t\geq 0\mid G-t\pi^*D\geq 0\}$$
and let 
$$G_0:=G-\sum_{D\mid D\text{ is a prime divisor on }Z}t_D\pi^*D.$$ 
Then $G_0\geq 0$ and $G_0$ is very exceptional$/Z$.

Let $\Ff$ be the foliation induced by $\pi$. By Proposition \ref{prop: bp semistable foliation lc}, $(X,\Ff,B^h,\Mm)$ is lc, so $(X,\Ff,B^h,\Mm;G+B-B^h)/Z$ is weak ACSS. By Proposition \ref{prop: weak cbf gfq}, $$K_{\Ff}+B^h+\Mm_X\sim_{\mathbb R,Z}K_X+B+G+\Mm_X\sim_{\mathbb R,Z}G\sim_{\mathbb R,Z}G_0.$$
By Theorem \ref{thm: mmp very exceptional alg int fol}, we may run a $(K_{\Ff}+B^h+\Mm_X)$-MMP$/Z$ which terminates with a weak lc model $(X',\Ff',(B^h)',\Mm)/Z$ of $(X,\Ff,B^h,\Mm)/Z$, such that $K_{\Ff'}+(B^h)'+\Mm_{X'}\sim_{\mathbb R,Z}0$. 

Let $B'$ and $G'$ be the images of $B$ and $G$ on $X'$ respectively, $\pi': X'\rightarrow Z$ the induced contraction, and let $\Nn'$ be the moduli part of $\pi': (X',B'+G',\Mm)\rightarrow Z$. By Lemma \ref{lem: ACSS mmp can run}, $(X',B'+G',\Mm)/Z$ satisfies Property $(*)$. By Proposition \ref{prop: bp semistable foliation lc}, $(X',B'+G',\Mm)/Z$ is BP semi-stable. By Proposition \ref{prop: weak cbf gfq}, $$K_{X'}+B'+G'+\Mm_{X'}\sim_{\mathbb R,Z}K_{\Ff'}+B'+\Mm_{X'}\sim_{\mathbb R,Z}0.$$
By Lemma \ref{lem: +a keep under mmp}, there exists an ample$/U$ $\Rr$-divisor $H'\geq 0$ on $X'$ such that either $(B^h)'\geq H'\geq 0$ or $\Mm-\bar H'$ is nef$/U$. By Theorem \ref{thm: a special b-semiampleness}, $\Nn'$ descends to $X'$ and $\Nn'_{X'}$ is semi-ample$/U$. Moreover, under the condition of (2), there exists a positive integer $I$ depending only on $d$ and $m$ such that $I\Nn'_{X'}$ is Cartier and $\mathcal{O}_X(nI\Nn'_{X'})$ is globally generated$/U$ for any integer $n\gg 0$.

Let $\Mm^Z$ and $\Mm'^Z$ be the base moduli parts of $\pi: (X,B,\Mm)\rightarrow Z$ and $\pi': (X',B'+G',\Mm)\rightarrow Z$ respectively. Since $\Nn'$ descends to $X'$ and $\Nn'_{X'}$ is semi-ample$/U$, $\Mm'^Z$ descends to $Z$ and $\Mm'^Z_Z$ is semi-ample$/U$. Since the induced birational map $\phi$ is a $G'$-MMP and $G'$ is vertical$/Z$, $\phi$ is an isomorphism over the generic point of $Z$. Thus $\pi: (X,B,\Mm)\rightarrow Z$ and $\pi': (X',B'+G',\Mm)\rightarrow Z$ are crepant over the generic point of $Z$. By Lemma \ref{lem: m preserved under crepant}, $\Mm^Z=\Mm'^Z$. Thus $\Nn=\Nn'$, and the theorem follows.
\end{proof}

\section{Proofs of the main theorems}\label{sec: proof of the main theorems}

In this section, we prove all the theorems that are listed in Sections \ref{sec:Introduction} and \ref{sec: statement of main results}. We remark that
\begin{enumerate}
    \item Theorem \ref{thm:base-point-freeness intro} was proven in Subsection \ref{subsec: bpf nonnqc}.
    \item Theorem \ref{thm: glc sings are Du Bois} was proven in Subsection \ref{subsec: du bois}.
    \item Theorem \ref{thm: cone theorem gfq} was proven in Subsection \ref{subsec: proof of cone}.
    \item Theorems \ref{thm: lc adjunction foliation nonnqc}, \ref{thm: dcc adjunction is dcc}, and \ref{thm:  ACSS model} were proven in Subsection \ref{subsec: proof of adj}.
    \item Theorem \ref{thm: acc lct alg int gfq} was proven in Subsection \ref{subsec: acc}.
    \item Theorem \ref{thm: global acc alg int gfq} was proven in Subsection \ref{subsec: global acc}.
\end{enumerate}

\begin{thm}[{cf. \cite[Conjecture 4.2(1)]{CS23a}}]\label{thm: fdlt is acss}
    Let $(X,\Ff,B,\Mm)/U$ be a $\Qq$-factorial F-dlt gfq. Then $(X,\Ff,B,\Mm)$ is ACSS.
\end{thm}
\begin{proof}
    Let $f: Y\rightarrow X$ be a foliated log resolution of $(X,\Ff,B,\Mm)$ such that $a(D,\Ff,B,\Mm)>-\epsilon_{\Ff}(D)$ for any prime $f$-exceptional divisor $D$. Let $B_Y:=f^{-1}_*B+(\Supp\Exc(f))^{\Ff_Y}$, then $(Y,\Ff_Y,B_Y,\Mm)$ is $\Qq$-factorial ACSS and $K_{\Ff_Y}+B_Y+\Mm_Y\sim_{\mathbb R,X}E\geq 0$
    for some $f$-exceptional prime divisor $E$ such that $\Supp E=\Supp\Exc(f)$. By Theorem \ref{thm: mmp very exceptional alg int fol}, we may run a $(K_{\Ff_Y}+B_Y+\Mm_Y)$-MMP$/X$ with scaling of an ample$/U$ divisor $A$ which terminates with a good minimal model $(X',\Ff',B',\Mm)/X$ of $(Y,\Ff_Y,B_Y,\Mm)$ such that $E$ is contracted by this MMP. Thus the induced birational morphism $X'\rightarrow X$ is small. Since $X$ is $\Qq$-factorial, $X'\rightarrow X$ is the identity morphism. Notice that the MMP is over $Z$ and hence the theorem follows.
\end{proof}

\begin{proof}[Proof of Theorem \ref{thm: cone theorem nonnqc gpair}]
It follows from Theorems \ref{thm: cone theorem gfq} and \ref{thm: cont thm gpair}.
\end{proof}

%\begin{proof}[Proof of Theorem \ref{thm: qfact nonnqc mmp can run intro}]It follows from Theorems \ref{thm: cone theorem nonnqc gpair} and \ref{thm: eof nonnqc}. \end{proof}

\begin{proof}[Proof of Theorem \ref{thm: kod vanishing gpair intro}]
It immediately follows from Theorem \ref{thm: kod vanishing with lc strata}(2) by letting $U=\{pt\}$.
\end{proof}

\begin{proof}[Proof of Theorem \ref{thm: kv vanishing gpair intro}]
It immediately follows from Theorem \ref{thm: kod vanishing with lc strata}(2).
\end{proof}

\begin{proof}[Proof of Theorem \ref{thm: semi-ampleness intro}]
We write $D=\sum_{i=1}^c r_iD_i$ where $r_1,\dots,r_c$ are linearly independent over $\mathbb Q$ and each $D_i$ is a $\Qq$-divisor. We define $D(\bm{v}):=\sum_{i=1}^cv_iD_i$ for any $\bm{v}=(v_1,\dots,v_c)\in\mathbb R^c$, and let $\bm{r}:=(r_1,\dots,r_c)$. By \cite[Lemma 5.3]{HLS19}, each $D_i$ is $\Qq$-Cartier, so $D(\bm{v})$ is $\Qq$-Cartier for any $\bm{v}\in\mathbb R^c$. 

Let $L:=D-(K_X+B+\Mm_X)$. Since being ample over $U$ is an open condition, there exists an open set $V\ni\bm{r}$ in $\mathbb R^c$ such that 
$\frac{1}{2}L+D(\bm{v})-D$ is ample over $U$ for any $\bm{v}\in V$.

By Theorem \ref{thm: cone theorem gfq}, there exist finitely many $(K_X+B+\Mm_X+\frac{1}{2}L)$-negative extremal rays over $U$ $R_1,\dots,R_l$, and each $R_j=\mathbb R_+[C_j]$ for some rational curve $C_j$ such that
$$-2\dim X\leq (K_X+B+\Mm_X+\frac{1}{2}L)\cdot C_j<0.$$
Since $D$ is nef, $D\cdot C_j\geq 0$ for each $j$. Thus possibly shrinking $V$, we may assume that for any $\bm{v}\in V$, we have that $D(\bm{v})\cdot C_j>0$ for any $j$ such that $D\cdot C_j>0$. Since $r_1,\dots,r_c$ are linearly independent over $\mathbb Q$, for any $j$ such that $D\cdot C_j=0$, we have $D(\bm{v})\cdot C_j=0$ for any $\bm{v}\in\mathbb R^c$. Therefore, $D(\bm{v})\cdot C_j\geq 0$ for any $j$ and any $\bm{v}\in V$.

For any extremal ray $R$ in $\overline{NE}(X/U)$ and any $\bm{v}\in V$, if $R=R_j$ for some $j$, then $D(\bm{v})\cdot R_j\geq 0$. If $R\not=R_j$ for any $j$, then
$$D(\bm{v})\cdot R=(K_X+B+\Mm_X+\frac{1}{2}L)\cdot R+(\frac{1}{2}L+D(\bm{v})-D)\cdot R>0.$$
Therefore, $D(\bm{v})$ is nef over $U$ for any $\bm{v}\in V$. Moreover,
$$D(\bm{v})-(K_X+B+\Mm_X)=\frac{1}{2}L+\frac{1}{2}L+D(\bm{v})-D$$
is ample over $U$.

We let $\bm{v}_1,\dots, \bm{v}_{c+1}\in V\cap\mathbb Q^c$ be rational points such that $\bm{r}$ is in the interior of the convex hull of $\bm{v}_1,\dots,\bm{v}_{c+1}$. Then there exist positive real numbers $a_1,\dots,a_{c+1}$ such that $\sum_{i=1}^{c+1}a_i=1$ and $\sum_{i=1}^{c+1}a_i\bm{v}_i=\bm{r}$. Since $D(\bm{v}_i)$ is a nef $U$-$\Qq$-divisor and $D(\bm{v}_i)-(K_X+B+\Mm_X)$ is ample over $U$, by Theorem \ref{thm:base-point-freeness intro}, $D(\bm{v}_i)$ is semi-ample over $U$ for any $i$. Therefore, $D=\sum a_iD(\bm{v}_i)$ is semi-ample over $U$.
\end{proof}

\begin{proof}[Proof of Theorem \ref{thm: subadj intro}] It is a special case of Definition-Lemma \ref{deflem: subadj minimal lc center}.
\end{proof}

\begin{proof}[Proof of Theorem \ref{thm: mmp fdlt}]
    It follows from Theorem \ref{thm: fdlt is acss} and Proposition \ref{prop: run mmp with scaling gfq}.
\end{proof}

\begin{proof}[Proof of Theorem \ref{thm: +a gmm fdlt}]
    It follows from Theorems \ref{thm: fdlt is acss}, \ref{thm: gmm polarized gfq}, and Proposition \ref{prop: run mmp get mfs}.
\end{proof}

\begin{proof}[Proof of Theorem \ref{thm: +a abundance fdlt}]
    We may assume that $K_{\Ff}+B+A$ is pseudo-effective$/U$. The theorem follows from Theorems \ref{thm: fdlt is acss} and \ref{thm: gmm polarized gfq}.
\end{proof}

\begin{proof}[Proof of Theorem \ref{thm: bpf fdlt}]
   It follows from Theorems \ref{thm: fdlt is acss} and \ref{thm: bpf induced gfq}.
\end{proof}

\begin{proof}[Proof of Theorem \ref{thm: eomfs}]
    It is a special case of Theorem \ref{thm: existence mfs}.
\end{proof}

\begin{proof}[Proof of Theorem \ref{thm: gmm ai num0}] 
First we prove (3). By Theorem \ref{thm: fdlt is acss} and Proposition \ref{prop: projective num 0 mmp}, we may run a $(K_{\Ff}+B)$-MMP with scaling of an ample $\Rr$-divisor, and any such MMP terminates with a log minimal model $(X',\Ff',B')$ of $(X,\Ff,B)$ such that $K_{\Ff'}+B'\equiv 0$. According to \cite[Theorem 1.4]{DLM23}, $K_{\Ff'}+B'\sim_\Rr 0$.

(2) follows from (3) and Theorem \ref{thm:  ACSS model}. (1) follows from (2).
\end{proof}

\begin{proof}[Proof of Theorem \ref{thm: abundance num 0 no restriction to f}]
Let $(Y,\Ff_Y,\bar B_Y;G)/Z$ be a proper ACSS model of $(X,\Ff,B)$ with induced birational morphism $f: Y\rightarrow X$ whose existence is guaranteed by Theorem \ref{thm: property * induction}. Let $K_{\Ff_Y}+B_Y:=f^*(K_{\Ff}+B)$ and $K_Y+B'_Y:=f^*(K_X+B)$. Since $(X,B)$ is lc, the coefficient of any component of $B'_Y$ is $\leq 1$. In particular, any coefficient of $B$ is $\leq 1$. We let
$$\bar B_Y:=f^{-1}_*B+(\Supp\Exc(f))^{\Ff_Y}$$
and $E:=B_Y-\bar B_Y$. Then $E\geq 0$ and $E$ is exceptional over $X$.

Suppose that $K_{\Ff_Y}+\bar B_Y$ is not pseudo-effective. We let
$$F:=\sum_{D\mid D\text{ is an }f\text{-exceptional prime divisor}}D.$$
Then $G\geq F$. Since $(Y,\bar B_Y+G)$ is lc and $G\geq F$, $(Y, \bar B_Y+F)$ is lc. By \cite[Theorem 5.3]{ACSS21}, $K_Y+\bar B_Y+F$ is not pseudo-effective. For any prime $f$-exceptional divisor $D$ such that $D$ is not $\Ff_Y$-invariant, we have $\mult_D\bar B_Y=1$. Therefore, 
$$\bar B_Y+F=f^{-1}_*B+\Supp\Exc(f).$$
Since the coefficient of any component of $B'_Y$ is $\leq 1$, we have
$$E':=f^{-1}_*B+\Supp\Exc(f)-B'_Y\geq 0$$
and $E'$ is exceptional over $X$. Therefore, 
\begin{align*}
    -\infty&=\kappa_{\sigma}(K_Y+\bar B_Y+F)=\kappa_{\sigma}(K_Y+f^{-1}_*B+\Supp\Exc(f))\\
    &=\kappa_{\sigma}(f^*(K_X+B)+E')=\kappa_{\sigma}(K_X+B)\geq 0,
\end{align*}
a contradiction. Thus $K_{\Ff_Y}+\bar B_Y$ is pseudo-effective. Since
$$0\leq\kappa_{\sigma}(K_{\Ff_Y}+\bar B_Y)\leq\kappa_{\sigma}(K_{\Ff_Y}+B_Y)=\kappa_{\sigma}(K_{\Ff}+B)=0,$$
we have $\kappa_{\sigma}(K_{\Ff_Y}+\bar B_Y)=0$. By Theorem \ref{thm: gmm ai num0}(1), $\kappa_{\iota}(K_{\Ff_Y}+\bar B_Y)=0$, so 
$$0=\kappa_{\iota}(K_{\Ff_Y}+\bar B_Y)\leq \kappa_{\iota}(K_{\Ff_Y}+B_Y)= \kappa_{\iota}(K_{\Ff}+B)\leq \kappa_{\sigma}(K_{\Ff}+B)=0.$$
Thus $\kappa_{\iota}(K_{\Ff}+B)=0$ and we are done.
\end{proof}

\begin{proof}[Proof of Theorem \ref{thm: cs23 4.2(1)}]
    It immediately follows from Theorem \ref{thm: fdlt is acss}.
\end{proof}

\begin{proof}[Proof of Theorem \ref{thm: cbf gfq}]
(1-4) follow from Definition-Theorem \ref{defthm: cbf lctrivial morphism}. (5) follows from Definition-Theorem \ref{defthm: cbf lctrivial morphism} and Proposition \ref{prop: gfq cbf preserve sing}. (6) follows from Proposition \ref{prop: gfq cbf preserve sing}. (7) follows from Lemma \ref{lem: td=bd}. (8) follows from Lemma \ref{lem: m preserved under crepant gfq} and Definition-Lemma \ref{deflem: cbf finite}. (9) follows from Definition-Lemma \ref{deflem: cbf gfq}(3) and Definition-Lemma \ref{deflem: cbf finite}(6).
\end{proof}

\begin{proof}[Proof of Theorem \ref{thm: fol adj intro}]
It is a special case of Theorem \ref{thm: lc adjunction foliation nonnqc}.
\end{proof}

\begin{proof}[Proof of Corollary \ref{cor: global acc rank 1 gfq}]  
If $K_{\Ff}$ is pseudo-effective, then $K_{\Ff}\equiv B\equiv\Mm_X\equiv 0$. By Lemma \ref{lem: trivial trace nef imply trivial}, $\Mm\equiv\bm{0}$ so there is nothing left to prove. So we may assume that $K_{\Ff}$ is not pseudo-effective. Since $\rk\Ff=1$, by \cite[Theorem 3.1]{LLM23} (see also \cite[Theorem 1.1]{CP19}), $\Ff$ is algebraically integrable. The corollary follows from Theorem \ref{thm: global acc alg int gfq}.
\end{proof}

\begin{proof}[Proof of Theorem \ref{thm: uniform rational polytope gfq}]
It immediately follows from Theorem \ref{thm: uniform rational polytope}.
\end{proof}

\begin{proof}[Proof of Theorem \ref{thm: gfq mmp very exceptional intro}]
It is a special case of Theorem \ref{thm: mmp very exceptional alg int fol}.
\end{proof}

%\begin{proof}[Proof of Theorem \ref{thm: eof nonnqc}]It is a special case of Theorem \ref{thm: existence of q-factorial glc flips}. \end{proof}

\begin{proof}[Proof of Theorem \ref{thm: ps intro}]
It follows from Theorem \ref{thm: a special b-semiampleness}.
\end{proof}

\end{document}